\documentclass[a4paper,12pt]{amsbook}

\usepackage{amssymb}
\usepackage{amsfonts}
\usepackage{amsmath}
\usepackage{graphicx}
\usepackage{moresize}

\usepackage{tikz-cd}%commutative diagrams

\usepackage{nicematrix}   % for chapters on Hilbert spaces and distributions

\usepackage{hyperref}
\hypersetup{
    colorlinks=true,
    linkcolor=blue,
    filecolor=magenta,      
    urlcolor=cyan,
    pdftitle={Overleaf Example},
    pdfpagemode=FullScreen,
    }

\providecommand{\U}[1]{\protect\rule{.1in}{.1in}}

\newtheorem{theorem}{Theorem}[chapter]
\theoremstyle{plain}

\newtheorem{corollary}{Corollary}[chapter]

\newtheorem{example}{Example}[chapter]
\newtheorem{exercise}{Exercise}[chapter]
\newtheorem{lemma}{Lemma}[chapter]

\newtheorem{proposition}{Proposition}[chapter]
\newtheorem{remark}{Remark}[chapter]

\numberwithin{equation}{section}

\numberwithin{section}{chapter}

\usepackage{imakeidx}
\makeindex

\usepackage{extsizes}

\begin{document}

\title{Functional Analysis and Operator Theory}
\author{Nicola Arcozzi}
%\date{September 2024}

\maketitle

\frontmatter

%\listofchanges   % to check the revisions

%\newpage

%Typeset in \LaTeX with \includegraphics[scale=0.1]{overleaf-logo-primary.png}

\vskip1cm

Copyright ©2025 Nicola Arcozzi

\vskip1cm

\includegraphics[scale=0.25]{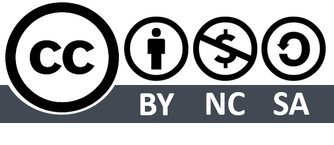}

\vskip1cm

\noindent{\bf License:}

\noindent This work is licensed under the Creative Commons Attribution-Noncommercial Share Alike 4.0 International License. To view a copy of this license, visit 

\noindent \href{https://creativecommons.org/licenses/by-nc-sa/4.0/}{https://creativecommons.org/licenses/by-nc-sa/4.0/}   or send a letter to Creative Commons PO Box 1866,
Mountain View, CA 94042, USA.
You can use, print, duplicate, and share this book as much as you want. You can
base your own notes on it and reuse parts if you keep the license the same. The license is CC-BY-NC-SA. Your derivative work must use  the same licenses.
Derivative works must be prominently marked as such.

\newpage

\section*{Preface}\label{SectPrefare}
These lecture notes merge the contents of two versions of a one-semester course in Functional Analysis and Operator Theory I taught in 2024, then in 2025, for the Master in Mathematics at the University of Bologna. In 2024 I did the basics of Gelfand's theory and spectral properties of compact operators, but the spectral theory of self-adjoint operators for the bounded case only. In 2025 I skipped Gelfand's theory and compact operators, and I did instead the spectral theory of unitary and of possibly unbounded self-adjoint operators. All students were familiar with the basis of Hilbert space theory, and the chapter on Hilbert spaces appears here just to keep the notes reasonably self-contained. The structure of the course, on the other hand, stemmed from a \href{https://www.dmi.unipg.it/didattica/scuola-matematica-interuniversitaria/edizioni-precedenti-corsi-estivi-smi/451-corsi-estivi-smi-2022}{course} in Functional Analysis I taught in 2022 for the "SMI" Summer School in Mathematics of Perugia.

The goal of the course, as it happens with most courses we teach, is to give an introduction to the subject which might be useful (i) as general mathematical culture; (ii) for those who will use these tools or some of their many {\it avatars} in more specialized topics; (iii) as an invitation to study Functional Analysis in greater depth (there are many beautiful treatises, old and new, the study of which might even be pleasant after an introductory course has cleared some ground). Having (i) in mind, I have preferred narrowing the front and moving in depth in the direction of spectral theory, which is one of great achievements of the first thirty or so years of the XX century. Concerning (ii), I unfortunately had to sacrifice the "quick introduction to Calculus of Variations from the viewpoint of functional analysis" I always have in the back of my mind. Time is short and choices have to be made.

The choice of much background material is based on the concrete groups of students attending the class. What some of them had not seen before, I stated and proved (the exception being some notions on Borel measures, which I just enumerated without proofs). All of them were familiar with Fourier transforms, for instance, but many had just a vague idea of what Fourier series are about.

All this accounts for the Frankensteinish structure of the lecture notes, where a number of repetitions and inversions of topics occur. For instance, in the chapter on preliminaries I use Banach-Alaoglu theorem, which appears two chapters downstream. In cases like this, I have preferred to keep most contents of the same kind in the same chapter, where the student can find them better contextualized. The chapters on Hilbert and Banach spaces are copied from notes for other courses and I did not make much of an effort to harmonize them with the rest. 

Some chapters have introductions explaining why we study this or that and which sources I am using.
Let me just mention the main references:
\begin{enumerate}
    \item Hilbert and Banach spaces: Michael Reed, Barry Simon,  \href{https://books.google.it/books/about/Functional_Analysis.html?id=16jFzQEACAAJ&redir_esc=y}{Methods of Modern Mathematical Physics: Functional analysis}, Academic Press, 1972;
    \item compact operators and their spectra: Reed and Simon;
    \item Banach algebras: Functional analysis, by Peter D. Lax, Wiley-Interscience, 2002, see this \href{https://www.ams.org/journals/bull/2006-43-01/S0273-0979-05-01073-6/S0273-0979-05-01073-6.pdf}{review by Meijun Zhu};
    \item spectral theory: for the unitary operators, Michael Taylor's lecture notes \href{https://mtaylor.web.unc.edu/wp-content/uploads/sites/16915/2018/04/specthm.pdf}{The Spectral Theorem for Self-Adjoint and Unitary Operators (2018)}; for the unbounded, self-adjoint operators, \href{https://terrytao.wordpress.com/2011/12/20/the-spectral-theorem-and-its-converses-for-unbounded-symmetric-operators/}{The spectral theorem and its converses for unbounded symmetric operators (2011)} by Terence Tao.
\end{enumerate}

\vskip0.5cm

Functional analysis, unlike Athena, was not born fully grown and armored from the head of Zeus. Its origins lie in rather concrete and specific problems, especially those concerning integral and differential equations, and calculus of variations. Then quantum mechanics was invented, and it needed what functional analysis was available, and much more. The presence, and the usefulness, of abstract theories emerged gradually, but not so slowly. There are several places where the story is told. One of them is \href{https://home.agh.edu.pl/~rudol/History_of_F_A_beginings.pdf}{The Establishment of Functional Analysis} (Historia Mathematica 11(3) 1984, p. 258-321) by Garrett Birkhoff and Erwin Kreyszig.

\vskip0.5cm

It is a pleasure to thank Nikolaos Chalmoukis, Giovanna Citti, Giovanni Dore, and Davide Guidetti, with whom I had a number of conversations on the material while teaching the class. Some material was included to answer questions from attendees of the Spring 2025 class; often ones I would not have ever thought of (and point in fact to important chapters of the theory). 

\vskip0.5cm 

\noindent The \LaTeX source is dowloadable, in case you wish to rearrange and modify the material according to your taste and needs. For comments and pointing out mistakes, you can write me at \href{mailto:nicola.arcozzi@unibo.it}{nicola.arcozzi@unibo.it}.

\vskip0.5cm

Bologna 2025

Nicola Arcozzi

\tableofcontents
%\listoffigures
%\listoftables

\mainmatter

\chapter{Some tools from real and complex analysis}\label{chaptprelim}

This section recalls some results from real and complex analysis which will be used in the sequel. Many of them are probably well known by most students with a bachelor in maths or physics, but perhaps not all of them. When in doubt, I provide proofs. I recommend browsing through the material and, only when needed, reading more carefully the parts which are necessary to prove this or that result in spectral theory. 
The prerequisites are basic notions from advanced calculus, and real and complex analysis. From advanced calculus: Green's formula, closed and exact fields (or $1$-forms); from real analysis: measure theory and Borel measures (including F. Riesz representation theorem for measures, which is the only non-elementary prerequisite), the basic density theorems, Fourier transforms; from complex analysis: the Cauchy-Riemann equations, Cauchy theorem and Cauchy formula, the open mapping theorem, expansion of a holomorphic function as a power series. In one proof we use the Banach-Alaoglu theorem, which is proved in the chapter on Banach spaces. Some basics on Fourier series are also in the chapter on Hilbert spaces, with slightly different proofs.

\vskip0.5cm 

\noindent{\bf Notation.} If $X$ is a topological space, we denote by $\overline{A}$ the closure of a subset $A$ of $X$. We make an exception for the complex plane $\mathbb{C}$, where the closure of $A\subseteq\mathbb{C}$ is $\text{Cl}A$, while $\overline{A}=\{\overline{z}:z\in A\}$.

\section{The definition of Banach space}\label{SectCompleMS}

\subsection{Normed linear spaces}\label{SSectNormedLinearSpaces}
An especially important family of metric spaces is that of the normed linear spaces\index{normed linear spaces}. 
A {\it{normed linear space}} is a vector space $X$ over $\mathbb{C}$ (or
over $\mathbb{R}$) endowed with a {\it{norm}}: a function
\[ \| \cdot \| \colon X \rightarrow [0, \infty), \]
satisfying the properties:
\begin{enumerate}
  \item[(i)] $\| x \| = 0$ if and only if $x = 0$;
  \item[(ii)] $\| \lambda x \| = | \lambda |  \| x \|\ \text{if}\ x \in X$ and
  $\lambda \in \mathbb{C}$;
  \item[(iii)] $\| x + y \| \leq \| x \| + \| y \|$.
\end{enumerate}
In practice, when we introduce a perspective norm on some concrete vector
space $X$, we have to verify that it has a finite value for each $x$ in $X$,
and we have to verify the conditions above, (iii) being sometimes subtle.

A normed linear space becomes a metric space when endowed with the distance
\[
d(x,y):=\|x-y\|.
\]
In addition to (i-iii), this distance satisfies properties linking it with the algebraic structure.
\begin{enumerate}
    \item[(iv)] $d(\lambda x,\lambda y)=|\lambda| d(x,y)$ if $x,y\in X$ and $\lambda\in{\mathbb C}$; 
    \item[(v)] $d(x+a,y+a)=d(x,y)$ if $x,y,a\in X$.
\end{enumerate}
The map $x\mapsto\lambda x$ is a {\it (complex) homothety}, and $x\mapsto x+a$ is a {\it translation}.
\begin{footnotesize}\begin{exercise}
Show that if a metric $d$ on a vector space $X$ satisfies (iv-v), then $\|x\|:=d(x,0)$ defines a norm, and that $d(x,y)=\|x-y\|$.
\end{exercise}\end{footnotesize}
For a continuous function $f\colon[a,b]\to{\mathbb C}$ we define:
\begin{align*}
\label{eqPnormCont}
    \|f\|_{L^p}&=\left(\int_a^b|f(x)|^p\right)^{1/p} \text{ for }1\le p<\infty, \text{ the $L^p$ {\it norm}},\crcr
    \|f\|_{L^\infty_0}&=\max_{x\in[a,b]}|f(x)|=\sup_{x\in[a,b]}|f(x)|,\text{ the {\it uniform norm}}.
\end{align*}
They make $(C[a,b],\|\cdot\|_{L^p})$ and $(C[a,b],\|\cdot\|_{u})$ into normed linear spaces.

\subsection{Complete metric spaces and Banach spaces}\label{SSectBanachComplete}
A basic problem in applications is the convergence of a sequence of objects to an object. Consider for instance the convergence of an algorithm. This notion can be formalized as the convergence of a sequence of points to a point in a metric space. Often, however, the nature of the objects in the sequence is clear, but so is not that of the limiting object. Think of the approximation of $\sqrt{2}$ by means of decimal (or binary) numbers having finitely many digits. The notion of "familiar" objects converging to a "nonfamiliar", "ghost" one is encoded in the notion of Cauchy sequence. The conceptual tool to make "ghosts" into "real" objects is the completion of a metric space. 

A sequence $\{ x_n \}_{n = 1}^{\infty}$ in $(X, d)$ {\it converges} to the {\it limit} $a
\in X$, $\lim_{n \rightarrow \infty} x_n = a$, if for all $\epsilon > 0$ there
is $n (\epsilon) > 0$ such that if $n > n (\epsilon)$, then $d (x_n, a)
\leqslant \epsilon$. The sequence $\{ x_n \}_{n = 1}^{\infty}$ is
{\it{Cauchy}} if for all $\epsilon > 0$ there is $n (\epsilon) > 0$ such
that if $n > n (\epsilon)$ and $j > 0$ one has $d (x_n, x_{n + j}) \leqslant
\epsilon$. 
\begin{footnotesize}\begin{exercise}\label{ExeDistCont}
 Show that, if $\lim_{n\to\infty}x_n=x$ and    $\lim_{n\to\infty}y_n=x$ in $(X,d)$, then $\lim_{n\to\infty}d(x_n,y_n)=d(x,y)$ exists in $\mathbb{R}$ in the usual sense.
\end{exercise}\end{footnotesize}
All convergent sequences are Cauchy, but the opposite implication
generally fails. The space $(X, d)$ is {\it{complete}} if all Cauchy
sequences in it converge. 

For instance, if ${\mathbb Q}\subset{\mathbb R}$ is the set of the rational numbers and $\{ x_n \}_{n = 1}^{\infty}$ is a sequence in ${\mathbb Q}$ converging to $\sqrt{2}$, then (with respect to the usual distance function $d(x,y)=|x-y|$) the sequence $\{ x_n \}_{n = 1}^{\infty}$ 
is Cauchy, but not convergent, in ${\mathbb Q}$. Familiar examples of complete metric spaces are ${\mathbb R}$ and ${\mathbb C}$, with respect to the distance $d(x,y)=|y-x|$.

Any metric space can be canonically imbedded in a complete one.
\begin{theorem}\label{TheoCompletenessMetric}
  Let $(X, d)$ be a metric space. Then there exists a complete metric space
  $(\tilde{X}, \tilde{d})$, the {\bf completion} of $(X, d)$, and an
  injective map $i \colon X \rightarrow \tilde{X}$ such that:
  \begin{enumerate}
    \item[(i)] $\tilde{d} (i (x), i (y)) = d (x, y)$;
    \item[(ii)] $i (X)$ is dense in $\tilde{X}$.
  \end{enumerate}
  The completion is unique in the following sense. For any other complete
  metric space $(Z, \delta)$ endowed with an injection $j \colon X \rightarrow Z$
  satisfying properties (i) $\delta(j(x),j(y))=d(x,y)$, and (ii) $j(X)$ is dense in $Z$, there is a unique map $F \colon \tilde{X}
  \rightarrow Z$ which is a surjective isometry, $\delta (F (\tilde{x}), F
  (\tilde{y})) = \tilde{d} (\tilde{x}, \tilde{y})$.
\end{theorem}

\begin{proof} We start with the construction of $(\tilde{X},\tilde{d})$.
  Let $C$ be the set of all Cauchy sequences in $X$ and for $\{ x_n \}, \{ y_n
  \} \in C$ set
  \[ \{ x_n \} \sim \{ y_n \} \Leftrightarrow \lim_{n \rightarrow \infty} d
     (x_n, y_n) = 0. \]
  The relation $\sim$ is an equivalence relation (check it!). Write $[\{ x_n \}]$ for the
  equivalence class of the Cauchy sequence $\{ x_n \}$, and let $\tilde{X} =
  C / \sim$, \ be the corresponding quotient space. Define
  \[ \tilde{d} ([\{ x_n \}], [\{ y_n \}]) = \lim_{n \rightarrow \infty} d
     (x_n, y_n).\] 
     It is easy to see that $\tilde{d}$ is well defined, and that it defines a distance on $\tilde{X}$.
  (check that it is independent of the representatives). Finally, for $x \in X$ set $i
  (x) = [\{ x_n = x \}]$, the class of the corresponding constant function.
  Properties (i) and (ii) are easily verified (exercise). Some work is needed to prove completeness, using a
  diagonal trick.
  
  Let $\{ \tilde{x}_n \}$, $\tilde{x}_n = [\{ x^n_m \}_{m = 1}^{\infty}]$, be
  a Cauchy sequence in $\tilde{X}$. For $k \geqslant 1$ select integer $n (k)
  > n (k - 1)$ (if $n (k - 1)$ was already selected) such that, for $n \geq n
  (k)$ and $j \geqslant 1$,
  \[ \frac{1}{2^k} \geqslant \tilde{d} (\tilde{x}_n, \tilde{x}_{n + j}) : =
     \lim_{m \rightarrow \infty} d (x^n_m, x^{n + j}_m) . \]
  Select then $m (k) > m (k - 1)$ such that, for $i = 0, 1$, $m \geqslant m
  (k)$, and $j \geqslant 1$,
  \[ (\ast) \quad d (x^{n (k + i)}_m, x^{n (k + i)}_{m + j}) \leqslant
     \frac{1}{2^k}, \]
  which we can ask because each $\{ x^n_m \}_{m = 1}^{\infty}$ is Cauchy, and
  \[ (\ast \ast) \quad d (x^{n (k)}_{m (k)}, x^{n (k + 1)}_{m (k)}) \leqslant
     \frac{2}{2^k}, \]
  which we can ask because $\frac{1}{2^k} \geqslant \tilde{d} (\tilde{x}_n,
  \tilde{x}_{n + j})$. Set then $a_k = x^{n (k)}_{m (k)}$.
  
  The sequence $\{ a_k \}_{k = 1}^{\infty}$ is Cauchy in $X$:
  \begin{eqnarray*}
    d (a_k, a_{k + 1}) & = & d (x^{n (k)}_{m (k)}, x^{n (k + 1)}_{m (k +
    1)})\\
    & \leqslant & d (x^{n (k)}_{m (k)}, x^{n (k + 1)}_{m (k)}) + d (x^{n (k +
    1)}_{m (k)}, x^{n (k + 1)}_{m (k + 1)})\\
    & \leqslant & \frac{2}{2^k} + \frac{1}{2^k}
  \end{eqnarray*}
  by $(\ast)$ and $(\ast \ast)$, and, by geometric sums, $d (a_k, a_{k + j})
  \leqslant \frac{6}{2^k}$. Let $a = [\{ a_k \}]$.
  
  We want to prove that $\tilde{d} (\tilde{x}_{n (k)}, a) \rightarrow 0$.
  Using the fact that $\lim_{n \rightarrow \infty} d (x_n, y_n) = \lim_{m
  \rightarrow \infty} d (x_{j_m}, y_{k_m})$ for all subsequences $\{ j_m \}$
  and $\{ k_m \}$ of the positive integers, provided the initial sequences $\{
  x_n \}$ and $\{ y_n \}$ are Cauchy (check it), we have
  \begin{eqnarray*}
    \tilde{d} (\tilde{x}_{n (k)}, a) & = & \lim_{l \rightarrow \infty} d (x^{n
    (k)}_{m (l)}, x^{n (l)}_{m (l)})\\
    & \leqslant & \limsup_{l \rightarrow \infty} d (x^{n (k)}_{m (l)}, x^{n
    (k)}_{m (k)}) + \lim_{l \rightarrow \infty} d (x^{n (k)}_{m (k)}, x^{n
    (l)}_{m (l)})\\
    & \leqslant & \frac{1}{2^k} + \frac{6}{2^k}.
  \end{eqnarray*}
  by $(\ast)$ and the fact that $\{ a_k = x^{n (k)}_{m (k)} \}$ is Cauchy (the precise value $6$ is indeed unimportant). 

  We now come to uniqueness. Let $(Z,\delta)$, $j$ as in the hypothesis, and let $\Tilde{x}=[\{x_n:\ n\ge1\}]$ be an element in $\Tilde{X}$, where $\{x_n\}$ is a Cauchy sequence in $X$. Then, $\{i(x_n):\ n\ge1\}$ is Cauchy in $\Tilde{X}$ and $\{j(x_n):\ n\ge1\}$ is Cauchy in $Z$. Define
  $F(\Tilde{x})=\lim_{n\to\infty}j(x_n)$, which exists because $Z$ is complete. The definition is well posed, since for two equivalent Cauchy sequences $\{x_n\},\ \{y_n\}$ in $Z$ we have
\[
\delta(j(x_n),j(y_n))=d(x_n,y_n)\to0 \text{ as }n\to\infty,
\]
hence, they have the same limit in $(Z,\delta)$. Moreover, for $\Tilde{x}=[\{x_n\}]$ and $\Tilde{y}=[\{y_n\}]$ in $\Tilde{X}$,
\begin{eqnarray*}
\delta(F(\Tilde{x},\Tilde{y}))&=&\delta(F([\{x_n\}]),F([\{y_n\}]))\crcr 
&=&\delta(\lim_{n\to\infty}j(x_n),\lim_{n\to\infty}j(y_n))\crcr 
&=&\lim_{n\to\infty}\delta(j(x_n),j(y_n))=\lim_{n\to\infty}d(x_n,y_n)=\lim_{n\to\infty}\Tilde{d}(i(x_n),i(y_n))\crcr 
&=&\Tilde{d}(\Tilde{x},\Tilde{y}).
\end{eqnarray*}
The equalities from second to third, and third to fourth line, follow from Exercise \ref{ExeDistCont} applied to $\delta$ and $\Tilde{d}$, respectively.

The isometry $F\colon\Tilde{X}\to Z$ is surjective. If $z\in Z$, by hypothesis there is a sequence $\{x_n\}$ in $X$ such that $\delta(j(x_n), z)\to 0$, so that $\{x_n\}$ is Cauchy in $X$.
Hence, $\{i(x_n)\}$ is Cauchy in $\Tilde{X}$ and $\lim_{n\to\infty}F(i(x_n))=\lim_{n\to\infty}j(x_n)=z$.
\end{proof}

A normed linear space is a {\it \index{Banach space} Banach space} if it is complete with respect to the distance induced by the norm.

Let $(X,\|\cdot\|)$ be a normed linear space, and let $(\tilde{X},\tilde{d})$ be its completion with respect to the distance $d(x,y)=\|x-y\|$ on $X$. We want to introduce on $(\tilde{X},\tilde{d})$ a structure of normed, linear space, in such a way that $\tilde{X}$ becomes a Banach space. To this aim, we define algebraic operations between equivalence classes of Cauchy sequences,
\[
[\{x_n\}]+[\{y_n\}]:=[\{x_n+y_n\}],\ \lambda [\{x_n\}]:=[\{\lambda x_n\}];
\]
and of the norm:
\[
\|[\{x_n\}]\|:=\lim_{n\to\infty}\|x_n\|.
\]
\begin{proposition}\label{PropLinCompl}
Sum, multiplication times scalar, and norm are well defined, and make $\tilde{X}$ into a normed linear space. Moreover, if $z,w$ belong to $\tilde{X}$, then
\[
\tilde{d}(z,w)=\|z-w\|,
\]
so that, in particular, $(\tilde{X},\|\cdot\|)$ is a Banach space. Also, the imbedding map $i\colon X\to\tilde{X}$ is linear.
\end{proposition}
\begin{proof}
The statement can be split in a number of statements, whose proof is left to the reader.
\begin{enumerate}
    \item[(a)] About the sum, we have to show that (i) $\{x_n+y_n\}$ is Cauchy if $\{x_n\}$ and $\{y_n\}$ are; (ii) that $[\{x_n+y_n\}]=[\{z_n+w_n\}]$ if $[\{x_n\}]=[\{z_n\}]$ and $[\{y_n\}]=[\{w_n\}]$.
    \item[(b)] Similar statements must be verified for the product of a vector with a scalar.
    \item[(c)] We have to verify that for $[\{x_n\}]$ in $\tilde{X}$, $\lim_{n\to\infty}\|x_n\|$ exists and does not depend on the particular representative.
    \item[(d)] For $z,w$ in $\tilde{X}$, $\tilde{d}(z,w)=\|z-w\|$.
    \item[(e)] The imbedding $i$ is linear.
\end{enumerate}
\end{proof}
\begin{footnotesize}\begin{exercise}
Prove the assertions above.
\end{exercise}\end{footnotesize}

\section{Borel measures}\label{SectBorelMeas} 
\subsection{Borel measures on the real line}\label{SSectBorelR}
We quickly review some facts from measure theory. A good reference is Gerald B. Folland, \href{https://www.wiley.com/en-us/Real+Analysis%3A+Modern+Techniques+and+Their+Applications%2C+2nd+Edition-p-9780471317166}{Real Analysis: Modern Techniques and Their Applications}, Wiley 2nd Edition (1999).

Let $X$ be a topological space. The {\it \index{Borel!$\sigma$-algebra}Borel $\sigma$-algebra} $\mathcal{B}(X)$ is the smallest one for which the open sets are measurable. A function $X\xrightarrow{f}Y$ from $X$ to another topological space $Y$ is {\it \index{Borel!measurable function} Borel measurable} if the preimage $f^{-1}(A)$ of any open subset $A$ of $Y$ is Borel measurable. In particular, a continuous function is measurable. The composition of two Borel measurable functions is measurable (it is not true, in general, that $g\circ f$ is measurable if $\mathbb{R}\xrightarrow{f}\mathbb{R}\xrightarrow{g}\mathbb{R}$, and $f$ and $g$ are Lebesgue measurable: Borel measurability is more stable with respect to composition). A {\it\index{Borel!measure} Borel measure} is a measure on the Borel $\sigma$-algebra. (Unless otherwise stated, a measure is positive).

A topological space $X$ is {\it locally compact} if any point $x$ in it has a compact neighborhood. A continuous function $X\xrightarrow{\varphi}\mathbb{C}$ {\it vanishes at infinity}, $\varphi\in C_0(X)$, if for each $\epsilon>0$ there exists $K$ compact such that $|\varphi(x)|<\epsilon$ if $x\in X\setminus K$. The space $C_0(X)$ is a Banach space with respect to the {\it uniform norm},
\[
\|\varphi\|_{L^\infty_0}=\sup_{x\in X}|\varphi(x)|,
\]
and the space $C_c(X)$ of the functions having compact support is dense in it.

A Borel measure on a locally compact space is a {\it\index{Radon measure (finite): finite Borel measure} Radon measure} if
\begin{enumerate}
    \item[(i)] $\mu(K)<\infty$ if $K$ is compact;
    \item[(ii)] $\mu(E)=\inf\{\mu(V):V\supseteq E,\ V\text{ open in }X\}$ for all measurable $E$;
    \item[(iii)] $\mu(V)=\sup\{\mu(K):K\subseteq E,\ K\text{ compact in }X\}$.
\end{enumerate}
Let $\mu$ be a finite, positive Borel measure on $X$ which is bounded on compact sets, and define $\Lambda:C_0(X)\to\mathbb{C}$,
\begin{equation}\label{eqRMKRTM}
    \Lambda(\varphi)=\int_X\varphi d\mu.
\end{equation}
Then
\begin{enumerate}
    \item[(i)] $\Lambda$ is linear and {\it positive}: if $\varphi\ge0$, then $\Lambda(\varphi)\ge0$;
    \item[(ii)] $\Lambda$ is {\it a bounded functional},
    \[
    \|\Lambda\|_{C_0(X)^\ast}:=\sup\{|\Lambda(\varphi)|:\varphi\in C_0(X)\}\le C\|\varphi\|_{L^\infty_0}.
    \]
\end{enumerate}
The following is a version of the {\it F.Riesz-Markov-Kakutani representation theorem for measures} provides a converse to the fact just stated.
\begin{theorem}\label{theoFRrieszMK}[F.Riesz-Markov-Kakutani theorem\index{theorem!Riesz-Markov-Kakutani for positive measures} for positive measures]
    Let $\Lambda:C_0(X)\to\mathbb{C}$ be a positive, linear functional. Then, there exists a unique Radon measure $\mu$ such that \eqref{eqRMKRTM} holds. Moreover, $\mu$ is finite and $\|\Lambda\|_{C_0(X)^\ast}=\mu(X)$.
\end{theorem}
Two positive measures $\mu,\nu$ on the same $\sigma$-algebra on $X$ are {\it mutually singular}, $\mu\perp\nu$, if $X=A\cup B$ with $A,B$ measurable, $A\cap B=\emptyset$, $\mu(B)=0=\nu(A)$.

Radon measures on the real line, or on the unit circle, play an important role in spectral theory, and it can help having an insight on the vastity of their catalogue. {\it All positive Borel measures on $\mathbb{R}$ (in fact, on $\mathbb{C}$, and on a large class of locally compact spaces) which are finite on the compacts are Radon measures, and viceversa. We will use the terms "finite Borel" and "finite Radon" as synonimous.} If $\mu$ is one such measure, and it is finite, then it can be uniquely decomposed as the sum of three positive, mutually singular, finite Borel measures,
\begin{equation}\label{eqMeasDec}
    \mu=\mu_{ac}+\mu_{sc}+\mu_d,
\end{equation}
where
\begin{enumerate}
    \item[(i)] $d\mu_{ac}=fdm$ is the {\it absolutely continuous part of $\mu$}, $m$ is the Lebesgue measure and $f\ge0$, $f\in L^1(m)$; 
    \item[(ii)] $\mu_d=\sum_{n}\lambda_n\delta_{x_n}$, the {\it discrete part of $\mu$}, is the sum of positive multiples ($\lambda_n>0$) of point masses $\delta_{x_n}$ at points $x_n\in\mathbb{R}$, with $\sum_n\lambda_n<\infty$;
    \item[(iii)] $\mu_{sc}$, the {\it singular continuous part of $\mu$}, is singular with respect to the Lebesgue measure and it has no point mass, $\mu_{sc}(\{x\})=0$ for all $x$.
\end{enumerate}
The natural Hausdorff measure supported on the Cantor set is a nontrivial example of a continuous, singular measure. If $\mu$ is a finite Radon measure on $\mathbb{R}$ and $\alpha(t)=\mu(-\infty,t]$ is the (increasing, right continuous) distribution function of $\mu$, then the decomposition in \eqref{eqMeasDec} corresponds to the decomposition $\alpha=\alpha_{ac}+\alpha_{sc}+\alpha_d$, where the functions on the right are increasing, $\alpha_{ac}$ is {\it absolutely continuous}, $\alpha_{d}$ is {\it pure jump}, and $\alpha_{sc}$ is continuous and it has derivative vanishing $a.e.$. with respect to Lebesgue's measure.

A {\it finite complex measure} on a measurable space $(X,\mathcal{F})$ (where $\mathcal{F}$ is a $\sigma$-algebra) is a map $\mu:\mathcal{F}\to\mathbb{C}$ such that, if $\{E_n\}_{n=1}^\infty$ is a sequence of disjoint measurable sets, then
\[
\mu\left(\bigcup_{n=1}^\infty E_n\right)=\sum_{n=1}^\infty\mu(E_n),
\]
and the series converges absolutely. A {\it finite, complex Borel measure\index{Borel!complex measure}} is one defined on a Borel $\sigma$-algebra.

Any finite complex measure can be uniquely be decomposed as a linear combination ({\it Jordan decomposition}) of four positive finite measures,
\[
\mu=\mu_{R+}-\mu_{R-}+i(\mu_{I+}-\mu_{I-}),
\]
with the requirement that $\mu_{R+}\perp\mu_{R-}$ and $\mu_{I+}\perp\mu_{I-}$. The four measures are Radon if $\mu$ is Radon. The integral with respect to $\mu$ is defined as the linear combination of the integrals with respect to the measures in the decomposition.

The {\it total variation} $|\mu|$ of a real valued measure $\mu=\mu_{R+}-\mu_{R-}$ is the positive measure $|\mu|=\mu_{R+}+\mu_{R-}$; that of a complex measure $\mu$, more generally, is defined by
\[
|\mu|(E)=\sup\left\{\sum_{n=1}^\infty |\mu(E_j)|:E=\bigcup_{n=1}^\infty E_n\text{ is a disjoint decomposition of }E\right\},
\]
and clearly $E$ and all $E_n$'s in the definition are measurable.
Although it is not obvious, the two definitions agree for a real valued measure.

Complex measures intervene in the following version of F. Riesz representation theorem.
\begin{theorem}\label{theoRMKcompMeas}[F.Riesz-Markov-Kakutani theorem\index{theorem!Riesz-Markov-Kakutani for complex measures} for complex measures]
    Let $X$ be a locally compact space, and $\Lambda:C_0(X)\to\mathbb{C}$ be a bounded, linear functional. Then, there exists a unique finite, complex Radon measure $\mu$ such that \eqref{eqRMKRTM} holds. Moreover,
    \[
    \|\Lambda\|_{C_0(X)^\ast}=|\mu|(X).
    \]
    Conversely, \eqref{eqRMKRTM} defines a bounded, linear functional on $C_0(X)$ if $\mu$ is finite, complex Radon measure.
\end{theorem}
\subsection{Borel measurable functions and the Lebesgue-Hausdorff theorem}\label{SSectBorelHierarchy}
The content of this section will only be used to prove some uniqueness statements in spectral theory, and to provide one of two proofs of Stone's theorem. I suggest you just read the statement of theorem \ref{theoLebesgueHausdorff} and of its corollary. There are properties which we will first prove for continuous functions on $\mathbb{R}$, and that are preserved under pointwise limits. The punchline is that such properties hold for all Borel measurable functions.
My source for this material is \cite{Sr1998}, S.M. Srivastava \href{https://books.google.it/books?id=FhYGYJtMwcUC}{A Course on Borel Sets}, Springer 1998.

\begin{theorem}\label{theoLebesgueHausdorff}[Lebesgue-Hausdorff theorem\index{theorem!Lebesgue-Hausdorff}] Let $\mathcal{F}$ be a class of functions $\mathbb{R}\xrightarrow{f}\mathbb{C}$ which (i) contains the continuous functions; (ii) is closed under pointwise limits: if $f_n\in\mathcal{F}$ and $f_n(x)\to f(x)$ for all $x\in\mathbb{R}$, then $f\in\mathcal{F}$.

Then, $\mathcal{F}$ contains all Borel measurable functions.
\end{theorem}
Clearly, there is no loss of generality in assuming that we consider functions $\mathbb{R}\xrightarrow{f}\mathbb{R}$, and in fact the theorem holds in a much broader context. A slight modification of the proof gives the following variant.
\begin{corollary}\label{corLebesgueHausdorff}
    Let $\mathcal{F}$ be a class of functions $\mathbb{R}\xrightarrow{f}\mathbb{C}$ which (i) contains the bounded continuous functions; (ii) is closed under pointwise limits of uniformly bounded functions: if $f_n\in\mathcal{F}$, $\|f_n\|_{L^\infty_0}\le C$, and $f_n(x)\to f(x)$ for all $x\in\mathbb{R}$, then $f\in\mathcal{F}$.

    Then, $\mathcal{F}$ contains $L^\infty_0(\mathbb{R})$.
\end{corollary}
The smallest class $\mathbf{B}=\mathbf{B}(\mathbb{R},\mathbb{R})$ which contains $C(\mathbb{R})$ and is closed under pointwise limits is called the {\it Baire class} of $\mathbb{R}$. The symbol $\mathbf{B}^\infty$ is the smallest one containing bounded continuous functions, and closed under pointwise limits of uniformly bounded functions. 

In order to prove theorem \ref{theoLebesgueHausdorff} we need to define a hierarchy into the Baire class. If $\mathcal{A}$ is a class of functions $f:\mathbb{R}\to\mathbb{R}$, let $\lim\mathcal{A}$ be the set of the pointwise limits of sequences of functions in $\mathcal{A}$. If $\mathcal{A}$ is a class of bounded function, $\lim^\infty\mathcal{A}$ is defined similarly, but requiring that the functions' sequence is uniformly bounded.

Let $\omega_1$ be the smallest uncountable ordinal.
\begin{itemize}
    \item[(i)] Set $\mathbf{B}_0=C(\mathbb{R})$.
    \item[(ii)] If $\alpha$ is a countable ordinal,
    \[
    \mathbf{B}_\alpha=\lim\left(\bigcup_{\beta<\alpha}\mathbf{B}_\beta\right).
    \]
    \item[(iii)] The classes $\mathbf{B}^\infty_\alpha$ are defined similarly, starting with $\mathbf{B}^\infty_0=C_b(\mathbf{R})$, the bounded and continuous functions, and using $\lim^\infty$ instead of $\lim$.
\end{itemize}
With this definition, $\mathbf{B}_1\supseteq\mathbf{B}_0$ contains pointwise limits of continuous functions, $\mathbf{B}_2\supseteq\mathbf{B}_1$ contains pointwise limits of functions which are of the form
\[
f(x)=\lim_{n\to\infty}\lim_{m\to\infty}f_{m,n}(x),
\]
where each $f_{m,n}$ is continuous, and so on. The definition looks more interesting at {\it limit ordinals}\footnote{Recall that an ordinal $\alpha$ is a limit ordinal if it is not the successor of any ordinal, $\alpha\ne\beta+1$ for all $\beta<\alpha$.} such as
\[
\omega_0=\sup\{n:\text{ $n$ is a finite ordinal}\},
\]
the first infinite ordinal. If $\alpha$ is a limit ordinal, then $f\in \mathbf{B}_\alpha$ if and only if there is a sequence $f_n\in\mathbf{B}_{\alpha_n}$ with $\alpha_n<\alpha$, such that $f$ is the pointwise limit of the $f_n$'s.

Let $\omega_1=\sup\{\alpha:\text{ $\alpha$ is countable ordinal}\}$ be the first uncountable ordinal.
\begin{lemma}\label{lemmaBaireDecomposed}
We have the decomposition:
    \[
    \mathbf{B}=\bigcup_{\alpha<\omega_1}\mathbf{B}_\alpha, \text{ and }\mathbf{B}^\infty=\bigcup_{\alpha<\omega_1}\mathbf{B}^\infty_\alpha.
    \]
\end{lemma}
\begin{proof}
We first prove the trivial $\supseteq$. Surely $\mathbf{B}_0=C(\mathbb{R})\subset\mathbf{B}$. We have to show that the definition of $\mathbf{B}_\alpha$ from the  form the preceding $\mathbf{B}_\beta$'s, $\beta<\alpha$, preserves membership in $\mathbf{B}$ if $\mathbf{B}_\beta\subseteq\mathbf{B}_\beta$ for all such $\beta$'s. Let $f_n\in\mathbf{B}_{\alpha_n}$ with $\alpha_n<\alpha$, and $f=\lim f_n$. Then $f$ is pointwise limit of functions in $\mathbf{B}$, hence it belongs to $\mathbf{B}$. 

To prove $\subseteq$, it suffices to show that $\bigcup_{\alpha<\omega_1}\mathbf{B}_\alpha$ is closed under pointwise limits. 
Let $f_n\in\mathbf{B}_{\alpha_n}$ with $\alpha_n<\omega_1$, and $f=\lim f_n$. The supremum of a countable family of countable ordinals is a countable ordinal\footnote{We can identify each ordinal $\alpha$ with the set $[0,\alpha)$ of the ordinals less that $\alpha$: $\alpha$ is countable if and only if $[0,\alpha)$ is countable of a set. The ordinal $\sup_n\alpha_n$ is the one associated with $\bigcup_{n\ge1}[0,\alpha_n)$, which is a countable set. Hence, $\sup_{n\ge1}\alpha_n$ is countable.}, $\alpha:=\sup\{\alpha_n:n\ge1\}<\omega_1$. Hence, $f_n\in\mathbf{B}_\alpha$ for all $n$, which implies that $f\in\mathbf{B}_{\alpha+1}$, and $\alpha+1<\omega_1$.

The proof for the class $\mathbf{B}^\infty$ is the same.
\end{proof}
\begin{corollary}\label{corollaryBaireIs} The Baire class $\mathbf{B}$ is closed under the following operations.
    \begin{enumerate}
        \item[(i)] If $f,g\in\mathbf{B}$, then $\max(f,g),\min(f,g)\in\mathbf{B}$.
        \item[(ii)] If $f,g\in\mathbf{B}$ and $a,b\in\mathbb{R}$, then $af+bg\in\mathbf{B}$.
    \end{enumerate}
    Similar statements hold for $\mathbf{B}^\infty$.
\end{corollary}
\begin{proof}
    The two statements hold for $\mathbf{B}_0=C(\mathbb{R})$. Suppose they hold in $\mathbf{B}_\beta$ for all $\beta\in[0,\alpha)$, and that 
    $f,g\in\mathbf{B}_\alpha$. Then, $f=\lim f_n$ with $f_n\in\mathbf{B}_{\alpha_n}$ ($\alpha_n<\alpha$), $g=\lim g_n$ with 
    $g_n\in\mathbf{B}_{\beta_n}$ ($\beta_n<\alpha$), so that $af_n+bg_n\in\mathbf{B}_{\max(\alpha_n,\beta_n)}$ and 
    $\lim(af_n+bg_n)=af+bg\in\mathbf{B}_\alpha$, by definition of $\mathbf{B}_\alpha$. The same reasoning can be done with $\max(f,g),\min(f,g)$. 

    By transfinite induction, they hold for all $\mathbf{B}_\alpha$ with $\alpha<\omega_1$. By lemma \ref{lemmaBaireDecomposed}, they hold the Baire class $\mathbf{B}$.
\end{proof}
Next, we have to formulate a different definition of the Borel class.
\begin{proposition}\label{propBorelRenew} Recall that the Borel class $\mathcal{B}=\mathcal{B}(\mathbb{R})$ on $\mathbb{R}$ is the smallest $\sigma$-algebra containing the open sets. Then, $\mathcal{B}=\mathcal{C}$, where $\mathcal{C}$ is the smallest family of subsets of $\mathbb{R}$ which contains all open sets and it is closed under countable intersections and countable disjoint unions.
\end{proposition}
\begin{proof}
   It is clear that $\mathcal{B}$ has the two stated properties, hence $\mathcal{B}\supseteq\mathcal{C}$. To show the opposite inclusion, it suffices to show that $\mathcal{C}$ is closed under complements, if $A\in\mathcal{C}$, then $\mathbb{R}\setminus A\in\mathcal{C}$. If this holds, $\mathcal{C}$ is a $\sigma$-algebra containing the open sets, hence it contains $\mathbf{B}$. Let
   \[
   \mathcal{D}=\{A\in\mathcal{C}:\mathbb{R}\setminus A\in\mathcal{C}\}.
   \]
   Since all closed sets in $\mathbb{R}$ are countable intersections of open sets (the universal notation for the class of such intersections is $G_\delta$), they belong to $\mathcal{C}$ and they complements do as well, hence, $\mathcal{D}$ contains open and closed sets.

   We next show that $\mathcal{D}$ is closed under countable intersection. Suppose $A_n\in\mathcal{D}$ for $n\ge0$, so that $A_n,\mathbb{R}\setminus A_n\in\mathcal{D}$, and $\cap_{n\ge0}A_n\in\mathcal{C}$. The sets $B_0=\mathbb{R}\setminus A_0$, and $B_n=(\mathbb{R}\setminus A_n)\cap A_1\cap\dots\cap A_{n-1}$ ($n\ge1$) are pairwise disjoint and belong to $\mathcal{D}$. Then, since $\mathcal{C}$ is closed under disjoint unions,
   \[
   \mathcal{C}\ni\bigcup_{n\ge0}B_n=\mathbb{R}\setminus\left(\bigcap_{n\ge0}A_n\right).
   \]
   Thus, $\bigcap_{n\ge0}A_n\in\mathcal{D}$.

   It is now clear that $\mathcal{D}$ is closed under countable, disjoint unions. If $B_n\in\mathcal{D}$ are pairwise disjoint, $n\ge0$, then $\bigcup_{n\ge0}B_n\in\mathcal{C}$ by definition, and $\mathbb{R}\setminus\left(\bigcup_{n\ge0}B_n\right)=\bigcap_{n\ge0}(\mathbb{R}\setminus B_n)\cap\mathbb{C}$ by the previous step. Hence,  $\bigcup_{n\ge0}B_n\in\mathcal{D}$.

   Since $\mathcal{D}$ contains the open sets and it is closed under countable intersections and countable unions of disjoint sets, $\mathcal{D}=\mathcal{C}$. Hence, $\mathcal{C}$ is also closed under complements, hence $\mathcal{C}=\mathcal{B}$.
\end{proof}
We are now ready to prove the Lebesgue-Hausdorff theorem for characteristic functions of sets.
\begin{proposition}\label{propLebHausSets}
    If $B$ is a Borel set in $\mathbb{R}$, then $\chi_B\in \mathbf{B}^\infty\subset\mathbf{B}$.
\end{proposition}
\begin{proof}
    Let 
    \[
    \mathcal{U}=\{B\subset \mathbb{R}:\chi_B\in \mathbf{B}^\infty\}.
    \]
    If $U$ is open, there are closed sets $F_n$ such that $F_n\nearrow U$, and for each of them, by Urysohn lemma, there is $f_n:\mathbb{R}\to[0,1]$ continuous, with $f_n|_{F_n}=1$ and $f_n|_{\mathbb{R}\setminus U}=0$. Since $\{f_n\}$ is uniformly bounded and $f_n\to\chi_U$, we have $U\in\mathcal{U}$.

    Suppose $\{B_n:n\ge0\}$ is a family of disjoint subsets in $\mathcal{U}$, having union $B$, and set
    \[
    f_n=\sum_{j=0}^n\chi_{B_n}.
    \]
    Then, $0\le f_n\le 1$, and $f_n\to\chi_{B}$ by corollary \ref{corollaryBaireIs}, hence $B\in \mathcal{U}$.

    Let $\{B_n:n\ge0\}$ is a family of subsets in $\mathcal{U}$, having intersection $B$, and set
    \[
    f_n=\min_{j\le n}^n\chi_{B_n}.
    \]
    Then, $0\le f_n\le 1$, and $f_n\to\chi_{B}$ by corollary \ref{corollaryBaireIs}, hence $B\in \mathcal{U}$.

    By proposition \ref{propBorelRenew}, we have that $\mathcal{U}=\mathbf{B}^\infty$.
\end{proof}

\begin{proof}[Proof of the Lebesgue-Hausdorff theorem] The characteristic functions of a Borel set belong to $\mathbf{B}^\infty$ by proposition \ref{propLebHausSets}. Any (Borel) measurable function $f$ is the pointwise limit of finite linear combinations $f_n$ of characteristic functions of (Borel) measurable sets, hence $f$ lies in $\mathbf{B}$. If $f$ is also bounded, we can take the $f_n$ to be uniformly bounded, hence $f\in\mathbf{B}^\infty$. 
\end{proof}
\section{Fourier transforms and series}\label{SsectFourier}
In this subsection we recall without proofs some basic facts in Fourier theory. These are well known by virtually all students, but there are various equally popular normalizations for the Fourier transform and the Fourier series, and it might be helpful to have the formulas for the choice we make here. 
\subsection{The Fourier transform}\label{SSSectFT}
The {\it Fourier transform\index{Fourier transform}}, defined in the first line, enjoys some properties. Below, $\mathbf{D}=i^{-1}\frac{\partial\ }{\partial x}$ is the symmetric version of the derivative (the {\it \index{momentum operator}momentum operator} of quantum mechanics).
\begin{align}\label{alignBFT}
    \widehat{\varphi}(\omega)&=(\mathcal{F}\varphi)(\omega)=\int_{\mathbb{R}}\varphi(x)e^{-i\omega x}dx;\\ 
    \varphi(x)&=\frac{1}{2\pi}\int_{\mathbb{R}}\widehat{\varphi}(\omega)e^{i\omega x}d\omega;\\ 
    \int_{\mathbb{R}}\overline{\varphi(x)}\psi(x)dx&=\frac{1}{2\pi}\int_{\mathbb{R}}\overline{\widehat{\varphi}(\omega)}\widehat{\psi}(\omega)d\omega;\\ 
    \widehat{\mathbf{D}\varphi}(\omega)&=\omega\widehat{\varphi}(\omega)=M_\omega\widehat{\varphi}(\omega).
\end{align}
The formula in the third line is the {\it Plancherel identity\index{Plancherel formula!Fourier transform}}, which is especially interesting when $\varphi=\psi$.
The definition makes sense for $f\in L^1(\mathbb{R})$, but it can be extended to $L^2(\mathbb{R})$ in the appropriate sense.
The two central identities hold, in the appropriate sense, if $\varphi\in L^2(\mathbb{R})$, while the fourth holds if, e.g., $\varphi\in C_c^1(\mathbb{R})$. For functions $f$ in $L^1(\mathbb{R})$ we the following hold:
\begin{align}\label{alignRL}
    &\widehat{f}\in C_0(\mathbb{R})\ (\text{i.e. $\widehat{f}\in C(\mathbb{R})$ and }\lim_{\omega\to\pm\infty}\widehat{f}(\omega)=0);\\ 
    &\widehat{f}=0\ \implies\ f=0.
\end{align}
The {\it convolution} $f\ast g$ of $f,g:\mathbb{R}\to\mathbb{C}$ is defined as
\[
(f\ast g)(x)=\int_{\mathbb{R}}f(x-y)g(y)dy.
\]
For $1\le p\le\infty$ {\it Young's inquality} holds,
\begin{equation}\label{eqBabyYoung}
    \|f\ast g\|_{L^p}\le\|f\|_{L^1}\|g\|_{L^p}.
\end{equation}
We have the fundamental:
\begin{equation}\label{eqConvFT}
 \widehat{f\ast g}(\omega)=\widehat{f}(\omega)\widehat{g}(\omega).   
\end{equation}
This holds if, e.g., $f,g\in L^2\cup L^1$. 

The convolution between a function in $C_0(\mathbb{R})$ and a finite Borel measure is defined in the natural way,
    \begin{equation}\label{eqConvMeasFunct}
    (\varphi\ast\mu)(x)=\int_{\mathbb{R}}\varphi(x-u)d\mu(u).
    \end{equation}
The {\it Fourier transform of the measure\index{Fourier transform!of a measure}} $\mu$ is obviously
    \[
    \widehat{\mu}(\omega)=\int_{\mathbb{R}}e^{-ix\omega}d\mu(x),
    \]
the convolution-to-product relation still holds, with the same proof, if e.g. $\varphi\in L^1(\mathbb{R})\cap C_0(\mathbb{R})$.
    \begin{equation}\label{eqConvProdMeas}
    \widehat{\varphi\ast\mu}=\widehat{\varphi}\widehat{\mu}.
    \end{equation}
We might define the convolution between a finite Borel measure and a function in $L^1(\mathbb{R})$, or even between finite Borel measures. The definition in \eqref{eqConvMeasFunct} can not be generalized naively, since a function in $L^1(\mathbb{R})$ is defined a.e with respect to the Lebesgue measure, while the measure $\mu$ might be orthogonal to the Lebesgue measure.

We end with some elementary, but useful formulas. Let $(J\varphi)(x)=\varphi(-x)$ the operator "reflection with respect to the origin" in $\mathbb{R}$, and $(C\varphi)(x)=\overline{\varphi(x)}$. Then,
\begin{equation}\label{eqEasyFT}
    CJ=JC,\ C\mathcal{F}=J\mathcal{F}C,\ J\mathcal{F}=\mathcal{F}J.
\end{equation}

\subsection{Fourier series}\label{SSSectFourierseries}
We briefly review some useful formulas concerning the Fourier series of $\varphi\in L^2[-\pi,\pi]$, where integrals are taken with respect to the Lebesgue measure. We might identify $(-\pi,\pi]$ with the {\it torus} $\mathbb{T}=\{e^{it}:t\in(-\pi,\pi]\}$ (the unit circle in the complex plane), and we can think of functions on $(-\pi,\pi]$ as of functions on $\mathbb{T}$ ($f(t)\equiv f(e^{it})$), or as $2\pi$-periodic functions on $\mathbb{R}$.

First, we state an important result.
\begin{theorem}\label{theoONB}
    Let $e_n(t)=e^{int}$, $e_n:[-\pi,\pi]$, $n\in \mathbb{Z}$. Then, $\mathcal{T}=\left\{(2\pi)^{-1/2}e_n:n\in\mathbb{Z}\right\}$ is a orthonormal basis for $L^2[-\pi,\pi]$.
\end{theorem}
If you have never seen any of the numerous proofs of this result, I provide one after we have introduced the Poisson kernel on $[-\pi,\pi]$. It is readily verified that $\mathcal{T}$ is a orthonormal system,
\begin{equation}\label{eqONStrig}
    \langle e_m,e_n\rangle_{L^2}=\int_{-\pi}^{\pi}e^{-imt}e^{int}dt=\begin{cases}
       2\pi\text{ if }m=n,\\
       \left[\frac{e^{i(n-m)t}}{i(n-m)}\right]_{-\pi}^\pi=0\text{ if }n\ne m.
    \end{cases}
\end{equation}
By general Hilbert space theory, the {\it trigonometric system\index{trigonometric system}} $\mathcal{T}$ is a orthonormal basis for a closed subspace $H(\mathcal{T})$ of $L^2[-\pi,\pi]$. Theorem \ref{theoONB} says that 
\[
H(\mathcal{T})=L^2[-\pi,\pi].
\]
For a function $\varphi\in L^1[-\pi,\pi]$, we define its {\it Fourier coefficients\index{Fourier!coefficients}} $(\mathcal{F}\varphi)(n)=\widehat{\varphi}(n)$, and list some of their basis properties.
\begin{align}\label{eqFTproperties}
    \widehat{\varphi}(n)=\int_{-\pi}^{\pi}\varphi(x)e^{-inx}dx \text{ for }n\in\mathbb{Z};\\ 
    \varphi(x)=\frac{1}{2\pi}\sum_{n\in\mathbb{Z}}\widehat{\varphi}(n)e^{in x};\\ 
    \langle\varphi,\psi\rangle_{L^2[-\pi,\pi]}=\frac{1}{2\pi}\langle\widehat{\varphi},\widehat{\psi}\rangle_{\ell^2(\mathbb{Z})};\\ 
    \widehat{\mathbf{D}\varphi}(n)=n\widehat{\varphi}(n)=(M_n\widehat{\varphi})(n).
\end{align}
The first equality is the definition of the {\it Fourier coefficients} of $\varphi$. The second holds for $\varphi\in H(\mathcal{T})$, and the equality has to be interpreted in the $L^2$ norm. The third relation is equivalent to the second, and it holds if either $\varphi$ or $\psi$ belong to $H(\mathcal{T})$. while the last relation holds if $\varphi\in C^1_{per}$, the space of the $C^1$, $2\pi$-periodic functions $\varphi:\mathbb{R}\to\mathbb{C}$. The second and third relations say that
\[
\mathcal{F}:H(\mathcal{T})\to\ell^2(\mathbb{Z})
\]
is a bijective isometry, but for a factor,
\[
\|\mathcal{F}\varphi\|_{\ell^2}^2=2\pi\|\varphi\|_{L^2[-\pi,\pi]}^2.
\]
This is in perfect agreement with the Fourier transform, and it can be extended to the context of {\it locally compact Abelian groups}, and to a more general non-commutative universe, with some highly nontrivial adjustments. For the commutative theory, see e.g. \href{https://onlinelibrary.wiley.com/doi/book/10.1002/9781118165621}{Fourier Analysis on Groups (1962)}, Wiley‐Interscience, by Walter Rudin.

The fourth relation in \eqref{eqFTproperties} can be used to give the Fourier series expansion a more concrete meaning. We give here a statement having an elementary proof.
\begin{lemma}\label{lemmaFTconverges}
    Let $\varphi\in H(\mathcal{T})\cap C^2_{per}[-\pi,\pi]$. Then,
    \[
    \varphi(t)=\frac{1}{2\pi}\sum_{n=-\infty}^\infty\widehat{\varphi}(n)e^{int},
    \]
    and the series converges uniformly.
\end{lemma}
Here $C^2_{per}([-\pi,\pi])$ is the space of $2\pi$-periodic functions in $C^2$.
\begin{proof}
    The series converges in $L^2$, and it converges to $\varphi$ because $\varphi\in H(\mathcal{\tau})$. Hence a subsequence converges a.e.,
    \begin{equation}\label{eqSubSeqConvFT}
    \varphi(t)=\lim_{j\to+\infty}\frac{1}{2\pi}\sum_{n=-N_j}^{N_j}\widehat{\varphi}(n)e^{int}\ a.e.\ t
    \end{equation}
    for some sequence $N_j\nearrow \infty$, and $|\widehat{\varphi''}(n)|\le\|\varphi''\|_{L^1[-\pi,\pi]}$, hence the sequence $\{\widehat{\varphi''}(n):n\in\mathbb{Z}\}$ is bounded. On the other hand,
    \[
    \widehat{-\varphi''}(n)=\widehat{\textbf{D}^2\varphi}(n)=n^2\widehat{\varphi}(n),
    \]
    hence, the series 
    \[
    S_\infty(\varphi)(t)=\frac{1}{2\pi}\sum_{n=-\infty}^{\infty}\widehat{\varphi}(n)e^{int}
    \]
    converges uniformly. By \eqref{eqSubSeqConvFT}, $S_\infty(\varphi)$ converges uniformly to $\varphi$.
\end{proof}
We can define the convolution of functions on $[-\pi,\pi]$ by thinking of them as $2\pi$-periodic functions (or functions on $\mathbb{T}$):
\[
(f\ast g)(t)=\int_{-\pi}^\pi f(t-s)g(s)ds.
\]
\begin{lemma}\label{lemmaConvFT}
    If $f\in H(\mathcal{T})$ and its Fourier series converges absolutely (e.g. if $f\in C^2_{per}[-\pi,\pi]$), and $g\in L^1[-\pi,\pi]$, then $f\ast g\in H(\mathcal{T})$ and $\widehat{f\ast g}(n)=\widehat{f}(n)\widehat{g}(n)$.
\end{lemma}
\begin{proof}
\begin{eqnarray*}
    (f\ast g)(t)&=&\int_{-\pi}^\pi \frac{1}{2\pi}\sum_{n=-\infty}^\infty \widehat{f}(n)e^{in(s-t)}g(s)ds\\ 
    &=&\frac{1}{2\pi}\sum_{n=-\infty}^\infty \widehat{f}(n) \int_{-\pi}^\pi e^{-ins} g(s)dse^{int}\\ 
    &=&\frac{1}{2\pi}\sum_{n=-\infty}^\infty \widehat{f}(n)\widehat{g}(n)e^{int},
\end{eqnarray*}
which converges absolutely to a function in $H(\mathcal{T})$, having Fourier coefficients $\widehat{f}(n)\widehat{g}(n)$.
\end{proof}
\section{Basic theory of harmonic functions}\label{SSectHarmFunct}
A function $U:\Omega\to\mathbb{C}$ is {\it harmonic\index{harmonic function}} on the open subset $\Omega\subseteq\mathbb{C}$ if $U\in C^2(\Omega)$ and 
\begin{equation}\label{eqLaplaceEqPlane}
\Delta U(z):=\partial_{xx}U(z)+\partial_{yy}U(z)=0
\end{equation}
for all $z\in \Omega$. In this two-dimensional context, it is well known that, if $\Omega$ is also simply connected, then {\it a real valued $U$ is harmonic if and only if it is the real part of a function $F$ holomorphic in $\Omega$} (the "if" holds with no special assumption on $\Omega$). The proof of this fact easily reduces to the {\it Cauchy-Riemann equations} for a holomorphic $F=U+iV$, 
$\begin{cases}
&\partial_xU=\partial_yV\\ 
&\partial_yU=-\partial_xV.
\end{cases}$
Equation \eqref{eqLaplaceEqPlane} implies that $-\partial_yUdx+\partial_xUdy$ is a closed $1$-form, which is then exact. Let $V$ be one of its potentials. Then, $\partial_xV=-\partial_yU$ and $\partial_yV=\partial_xU$, which are the Cauchy-Riemann equations.
Viceversa, that the Cauchy-Riemann equations imply \eqref{eqLaplaceEqPlane} for both $U$ and $V$ is easily verified.

We list here the few facts that we will need.
\begin{enumerate}
    \item[(i)] [{\bf Series expansion}] {\it If $U$ is harmonic in $D(a,R)=\{z:|z-a|<R\}$, with $0<R\le\infty$, then $U$ can be written in the form
    \[
    U(a+re^{it})=\sum_{n=-\infty}^{+\infty}a_n r^{|n|}e^{int},
    \] 
    for $0\le r<R$. Moreover, if $0<r<R$, then
    \begin{equation}\label{eqCoeff}
    a_n=\frac{1}{2\pi r^{|n|}}\int_0^{2\pi}U(a+re^{it})e^{-int}dt.
    \end{equation}}
    \begin{proof}
    We can assume $a=0$ and that $U=\text{Re}F$ is the real part of a holomorphic $F$. If $  F(z)=\sum_{n=0}^\infty b_nz^n$ is the Taylor expansion of $F$, then
    \[
    U(re^{it})=\frac{F(e^{it})+\overline{F(e^{it})}}{2}=\text{Re}(b_0)+\sum_{n=-\infty}^{-1}\overline{b_{-n}}r^{|n|}e^{int}+\sum_{n=1}^\infty b_n r^n e^{int},
    \]
    The formula for the coefficients follows by integrating term by term.
    \end{proof}
    \item[(ii)] [{\bf Liouville theorem}] {\it A bounded harmonic function on $\mathbb{C}$ is constant.} 
    \begin{proof}
    For $n\ne0$, by \eqref{eqCoeff},
    \[
    |a_n|\le\frac{\|U\|_{L^\infty} }{r^{|n|}}\to0
    \]
    as $r\to\infty$.
    \end{proof}
    \item[(iii)] [{\bf Identity principle}] {\it If $U$ is harmonic on $\Omega$, connected, and if $U$ vanishes on an open set, then it vanishes identically. More generally, if $U$ is constant on an open subset of $\Omega$, then it is globally constant.}
    \begin{proof}
       Let 
       \[
       Z=\{a\in\Omega:\exists D(a,r)\subseteq\Omega \text{ such that $U$ vanishes on }D(a,r)\}.
       \]
       By \eqref{eqCoeff}, if $a\in Z$, then all series coefficients of $U$ vanish, $Z$ contains the largest disc centered at $a$ which is contained in $\Omega$. The set $Z$ is nonempty and open. Let $Z\ni a_n\to b\in\Omega$, and let $D(b,r)\subset\Omega$. Pick $a_n$ such that $|b-a_n|\le r/4$. Then $b\in D(a_n,r/2)\subseteq Z$. Being $\Omega$ connected, $Z=\Omega$.
    \end{proof}

    \begin{proof}[A proof using the open mapping theorem] If $\Omega$ is simply connected, then $U=\text{Re}F$ with $F$ holomorphic in $\Omega$. Since $F(D(a,r))$ fails to be open, then $F$, hence $U$, is constant.
        Suppose $\Omega$ is just connected, and let $b\in\Omega$. Consider a simple, broken line $l$ in $\Omega$ which joins $a$ and $b$, and let $A=\{z\in\Omega:\text{dist}(z,l)<\epsilon\}$. If $\epsilon$ is small, then $A$ is simply connected (exercise) and contained in $\Omega$. The previous case applies, then $U(b)=U(a)=0$. 
    \end{proof}
    \item[(iv)] [{\bf Mean value property.}] {\it If $U$ is harmonic in $\Omega$ and $\text{Cl}D(a,r)\subset\Omega$, then
    \[
    U(a)=\frac{1}{2\pi}\int_{0}^{2\pi}U(a+re^{it})dt.
    \]}
    \begin{proof}
        It follows from (i).
    \end{proof}
    \item[(v)] [{\bf Maximum principle}] {\it Let $U$ be real valued and continuous on $\Omega$, open and connected. If $U$ satisfies the mean value property and it achieves its maximum, then it is constant.} 
    \begin{proof}
    Let $M$ be the maximum of $U$, and $Z=\{a\in\Omega:\ U(z)=M\}$. By the mean value property and intuition, if $a\in Z$ and $D(a,r)\subseteq\Omega$, then $D(a,r)\subseteq Z$. This shows that $Z$ is open, but $Z$ is also closed in the relative topology of $\Omega$, then $Z=\Omega$.
    \end{proof}

    \item[(vi)][{\bf Local maximum principle}] {\it Let $U$ be real valued and harmonic on $\Omega$, open and connected. If $U$  has a local maximum in $\Omega$, then it is constant.} 
    \begin{proof} Suppose $a$ is a maximum in $D(a,r)\subset\Omega$, and let $F=U+iV$ be holomorphic in $D(a,r)$. Then, $F(D(a,r))\subseteq C:=\{w=u+iv:\ u\le U(a)\}$ and $F(a)\in \partial C$, which contradicts the open mapping theorem. Thus $F$ is constant, hence $U$ is constant on $D(a,r)$. By the identity principle, $U$ is constant.
    \end{proof}
    
    \item[(vii)] [{\bf Poisson representation}] {\it If $U$ is harmonic in $\Omega$ and $\text{Cl}D(a,r)\subset\Omega$, then for $z\in D(a,r)$:
    \begin{equation}\label{eqPoissonDisc}
        U(a+z)=\frac{1}{2\pi}\int_0^{2\pi}\frac{r^2-|z|^2}{|re^{it}-z|^2}U(a+re^{it})dt.
    \end{equation}}
    The expression
    \[
    P_z^r(e^{it})=\frac{1}{2\pi}\frac{r^2-|z|^2}{|re^{it}-z|^2}dt
    \]
    is the {\it Poisson kernel at $z$}.
    \begin{proof}
     It suffices to show it for $U=F$ holomorphic, and $a=0$. We start with Cauchy formula,
     \[
     F(z)=\frac{1}{2\pi i}\int_{|w|=r}\frac{F(w)}{w-z}dw=\frac{1}{2\pi}\int_0^{2\pi}\frac{re^{it}}{re^{it}-z}F(re^{it})dt.
     \]
     The difference between the kernel in the last expression and the Poisson kernel is 
     \begin{eqnarray*}
     \frac{re^{it}}{re^{it}-z}-\frac{r^2-|z|^2}{|re^{it}-z|^2}&=&\frac{re^{it}(re^{-it}-\overline{z})-(r^2-|z|^2)}{(re^{it}-z)(re^{-it}-\overline{z})}\\ 
     &=&\frac{|z|^2-\overline{z}re^{it}}{(re^{it}-z)(re^{-it}-\overline{z})}\\ 
     &=&-\frac{\overline{z}}{re^{-it}-\overline{z}}\\ 
     &=&-\sum_{n=0}^\infty\left(\frac{\overline{z}}{r}\right)^{n+1}e^{i(n+1)t}.
     \end{eqnarray*}
     Integrating term by term, we see that
     \[
     \frac{1}{2\pi}\int_0^{2\pi}\left[\frac{re^{it}}{re^{it}-z}-\frac{r^2-|z|^2}{|re^{it}-z|^2}\right]F(re^{it})dt=0, 
     \]
     hence that 
     \[
     F(z)=\frac{1}{2\pi}\int_0^{2\pi}\frac{re^{it}}{re^{it}-z}F(re^{it})dt=\frac{1}{2\pi}\int_0^{2\pi}\frac{r^2-|z|^2}{|re^{it}-z|^2}F(re^{it})dt,
     \]
     as wished. Since the Poisson kernel is positive (thus real), the relation holds for the real and the imaginary part of $F$ separately.
    \end{proof}
    There are indeed other ways to guess the Poisson kernel and to prove the Poisson representation. Here is one using Fourier series, starting from (i). It does not really make use of the completeness of the trigonometric system, since in the Plancherel formula below one function (there denoted by $Q$) belongs to the span $H(\mathcal{T})$ of the trigonometric system.
    \begin{proof}
     Let $a=0$ and $z=\rho e^{is}$, with $0\le \rho<r$, and $\text{Cl}D(0,r)\subset\mathbb{C}$. By Plancherel's formula,
     \begin{eqnarray*}
     U(\rho e^{is})&=&\sum_{n=-\infty}^{+\infty}a_n\rho^{|n|}e^{ins}=\sum_{n=-\infty}^{+\infty}a_nr^{|n|}\overline{(\rho/r)^{|n|}e^{-is}}\\ &=&2\pi\langle Q,U_r  \rangle_{L^2[0,2\pi]},
     \end{eqnarray*}
     where 
     \begin{eqnarray*}
     U_r(e^{it})&=&\frac{1}{2\pi}\sum_{n=-\infty}^{+\infty}a_nr^{|n|}e^{int}\\ 
     &=&\frac{1}{2\pi}U(re^{it}),
     \end{eqnarray*}
     and (by computing geometric series)
     \begin{eqnarray}\label{eqPoissonSeries}
         Q(e^{it})&=&\frac{1}{2\pi}\sum_{n=-\infty}^{+\infty}(\rho/r)^{|n|}e^{-ins}e^{int}\\
         &=&\frac{1}{2\pi}\left[\frac{1}{1-\frac{\rho}{r}e^{i(t-s)}}+\frac{1}{1-\overline{\frac{\rho}{r}e^{i(t-s)}}}-1\right]\\  
         &=&\frac{1}{2\pi}\frac{r^2-\rho^2}{|re^{it}-\rho e^{is}|^2}.
     \end{eqnarray}
     Summarizing,
     \[
     U(\rho e^{is})=2\pi\langle Q,U_r  \rangle_{L^2[0,2\pi]}=\frac{1}{2\pi}\int_0^{2\pi}\frac{r^2-\rho e^{is}}{|re^{it}-\rho e^{is}|^2}F(re^{it})dt,
     \]
     as promised.
    \end{proof}
    \item[(viii)] [{\bf The Poisson extension of a continuous function on the circle}] {\it Let $\varphi=\varphi(t):\partial D(a,r)$ be continuous, and define
    \begin{equation}\label{eqPoissonKernelTrue}
    U(a+\rho e^{is})=\begin{cases}
        \frac{1}{2\pi}\int_0^{2\pi}\frac{r^2-\rho^2}{|re^{it}-\rho e^{is}|^2}\varphi(a+re^{it})dt\text{ if }0\le\rho<r,\\
        \varphi(r e^{is})\text{ if }\rho=r.
    \end{cases}
    \end{equation}
    Then, $U$ is harmonic in $D(a,r)$ and continuous on $\text{Cl}D(a,r)$.}
    
    The proof is similar to that of proposition \ref{propDirProbLap} below.
    \begin{proof}
        We have the following three properties.
        \begin{enumerate}
            \item[(a)] The function $z\mapsto P_z^r(e^{is})$ is harmonic and positive in $D(0,r)$.
            \item[(b)] Consider the function  $s\mapsto P_\rho^r(e^{is})$ on $[-\pi,\pi]$. Then, 
            \[
            \int_{-\pi}^\pi P^r_\rho(e^{is})ds=1.
            \]
            \item[(c)] For each $\delta>0$,
            \[
            \lim_{\rho\to r}\sup_{\delta\le|s|\le\pi}P_\rho^r(e^{is})=0.
            \]
        \end{enumerate}
        (a) Positivity is clear, and harmonicity holds because, by \eqref{eqPoissonSeries},
        \begin{equation}\label{eqSeriesPoiss}
        P_z^r(e^{is})=\frac{1}{2\pi}\frac{r^2-|z|^2}{|re^{is}-z|^2}=\sum_{n=0}^\infty r^{-|n|}e^{-ins} z^n+\sum_{n=1}^{\infty} r^{-|n|}e^{ins} \overline{z}^n,
        \end{equation}
        and the functions $z^n$, $\overline{z}^n$ are harmonic.

        (b) Integrate the series \eqref{eqSeriesPoiss} term by term.

        (c) We have 
        \begin{eqnarray*}
        \frac{r^2-\rho^2}{|re^{is}-\rho|^2}&=&\frac{r^2-\rho^2}{r^2-2r\rho\cos(s)+\rho^2}=\frac{r^2-\rho^2}{(r-\rho)^2+2r\rho[1-\cos(s)]}\\ 
        &\le&\frac{r^2-\rho^2}{(r-\rho)^2+2r\rho[1-\cos(\delta)]}\\ 
        &\to&0\text{ as }\rho\to r.
        \end{eqnarray*}

        Differentiating with respect to $z=\rho e^{is}$ under the integral on the right side of \eqref{eqPoissonKernelTrue}, and using (a), we see that $U$ is harmonic. We only have to show continuity at points $a+re^{it}$. We can assume $a=0$, and verify continuity at $r$. Fix $\epsilon>0$. We estimate separately the summands on the right of:
        \[
        |U(\rho e^{is})-\varphi(r)|\le|U(\rho e^{is})-\varphi(re^{is})|+|\varphi(re^{is})-\varphi(r)|.
        \]
        The second summand can be made less than $\epsilon$ by choosing $|s|\le\eta(\epsilon)$, by the uniform continuity of $\varphi$.
        For the first, using (b),
        \begin{eqnarray*}
            |U(\rho e^{is})-\varphi(re^{is})|&=&\frac{1}{2\pi}\left|\int_{-\pi}^{\pi}\frac{r^2-\rho^2}{|re^{it}-\rho e^{is}|^2}[\varphi(r^{it})-\varphi(re^{is})]dt\right|\\ 
            &\le&\frac{\epsilon}{2\pi}\int_{|s-t|\le\eta(\epsilon)}\frac{r^2-\rho^2}{|re^{it}-\rho e^{is}|^2}dt\\
            &\ &+\frac{2\|\varphi\|_{L^\infty}}{2\pi}\int_{|s-t|\ge\eta(\epsilon)}\frac{r^2-\rho^2}{|r-\rho e^{i(s-t)}|^2}dt\\ 
            &\le&\epsilon+\frac{\|\varphi\|_{L^\infty}}{\pi}\epsilon,
        \end{eqnarray*}
        if we choose $\rho(\epsilon)\le\rho<r$, by (c).
    \end{proof}
     \item[(ix)] [{\bf The completeness of the trigonometric system}\index{trigonometric system!completeness}] We prove here theorem \ref{theoONB}. Let $f\in L^2[-\pi,\pi]$ and fix $\epsilon>0$. We want to show that there is a trigonometric polynomial $p=p(e^{it})$ such that $\|f-p\|_{L^2}\le C\epsilon$. By density of $C_{per}[-\pi,\pi]$ in $L^2[-\pi.\pi]$, there is a continuous function $\varphi$ such that $\|f-\varphi\|_{L^2}\le\epsilon$. Consider now the Poisson extension of $\varphi$ in theorem \ref{theoONB} (with $a=0$ and $r=1$), and set $U_\rho(e^{it})=U(\rho e^{it})$. There is $\rho=\rho(\epsilon)$ such that
     \[
     \|U_\rho-\varphi\|_{L^2}\le\sqrt{2\pi}\|U_\rho-\varphi\|_{L^\infty}\le \epsilon,
     \]
     since by our estimates the convergence of $U_\rho$ to $\varphi$ is uniform.

     Last, we have that
     \[
     U(\rho e^{is})=\frac{1}{2\pi}\sum_{n=-\infty}^\infty \widehat{\varphi}(n)\rho^{|n|}e^{ins}.
     \]
     Let $p(e^{is})=\frac{1}{2\pi}\sum_{n=-N}^N \widehat{\varphi}(n)\rho^{|n|}e^{ins}$, which is a trigonometric polynomial. Then
     \[
     \|U_\rho-p\|_{L^\infty}\le\frac{\|\varphi\|_{L^1}}{2\pi}\sum_{|n|>N}\rho^{|n|}\le\epsilon,
     \]
     provided we choose $N$ large enough.
    \item[(x)] [{\bf The inverse mean value property}] {\it Let $V$ be continuous on $\Omega$, open, and suppose it satisfies the mean value property. Then, $V$ is harmonic in $\Omega$.}
    \begin{proof}
        Consider $\text{Cl}D(a,r)\subset\Omega$, and consider the Poisson extension $U$ of the restriction of $V$ to $\partial D(a,r)$, which is harmonic in $D(a,r)$, and continuous in $\text{Cl}D(a,r)$ by (viii). The function $V-U$ satisfies the mean value property in $D(a,r)$ by (iv), hence its maximum is achieved at the boundary by (v). Since the boundary values vanish, $V\le U$. The same reasoning applies to $U-V$, then $V=U$ is harmonic in $D(a,r)$. 
    \end{proof}
    \item[(xi)] [{\bf The reflection principle}] {\it Let $\Omega\subseteq\mathbb{C}_+$ be a domain such that $\partial\Omega\cap\mathbb{R}$ is a finite union of the closures of the open intervals $I_1,\dots,I_n$, and let $U:\Omega\to\mathbb{R}$ be harmonic and such that $\lim_{z\to x}U(z)=0$ for all $x\in\partial\Omega\cap\mathbb{R}$. Let $\overline{\Omega}=\{\overline{z}:z\in\Omega\}$ be the reflection of $\Omega$ in the real axis, and define $U(\overline{z})=-U(z)$ for $\overline{z}\in\overline{\Omega}$, and $U=0$ on $I_1\cup\dots\cup I_n$. Then, $U$ is harmonic on $\Omega\cup\overline{\Omega}\cup I_1\cup\dots\cup I_n$.} 
    \begin{proof} We apply (x).
        The extended function is harmonic in $\Omega\cup\overline{\Omega}$, hence it satisfies the mean value property there. We only have to verify the mean value property at the points $x\in I_1\cup\dots\cup I_n$, where the function vanishes: the property there holds by symmetry.
    \end{proof}
\end{enumerate}

\section{The Dirichlet problem for the Laplace equation in the upper half plane}\label{SSectDirLap}
Let $f:\mathbb{R}\to\mathbb{C}$ continuous and bounded. A solution of the {\it Dirichlet problem} for the {\it Laplace equation} in the upper half plane $\mathbb{C}_+=\{z=x+iy:\ y>0\}$ with {\it boundary conditions} $f$ is a continuous function $U:\mathbb{C}_+\cup\mathbb{R}\to\mathbb{C}$ such that
\begin{align}\label{eqDirLap}
    &\Delta U(z)=(\partial_{xx}+\partial_{yy})U(x+iy)=0\text{ for }z\in\mathbb{C}_+\\
    &U(x)=f(x) \text{ for }x\in\mathbb{R}\\ 
    &U\text{ is bounded}.
\end{align}
For $z=x+iy\in\mathbb{C}_+$, let
\begin{equation}\label{eqPoissonKernel}
    P_y(x)=\frac{1}{\pi}\frac{y}{x^2+y^2}.
\end{equation}
be the {\it Poisson kernel}.
\begin{proposition}\label{propDirProbLap}
The unique solution of \ref{eqDirLap} is given by 
\begin{equation}\label{eqPossonIntegral}
U(x+iy)=(P_y\ast f)(x)=\frac{1}{\pi}\int_{\mathbb{R}}f(u)\frac{y}{(x-u)^2+y^2}du.
\end{equation}
\end{proposition}
There are various ways to guess the formula \eqref{eqPossonIntegral} and the form of the Poisson kernel, which are interesting in that they point to extensions and generalizations. We provide one of them here below, and two more after the proof of proposition \ref{propDirProbLap}. 

\vskip0.5cm

\noindent{\bf The Poisson kernel via Fourier transforms.} By Fourier inversion formula, also assuming for the moment that $f\in L^1(\mathbb{R})$,
\[
f(x)=\int_{\mathbb{R}}\widehat{f}(\omega)e^{i\omega x}d\omega
\]
is the superposition of the bounded functions $x\mapsto e^{i\omega x}$ (with $\omega\in\mathbb{R}$ fixed). If we extend each of them to a bounded harmonic function on the upper half-plane, by superposition we can extend $f$ to the same. A natural guess is $z\mapsto e^{i\omega z}=e^{i\omega x}e^{-\omega y}$, which is holomorphic, hence harmonic, for $z\in\mathbb{C}$. We have boundedness on $\mathbb{C}_+$, however, only when $\omega\ge0$. When $\omega<0$, the function $z\mapsto e^{i\omega\overline{z}}$ is bounded on $\mathbb{C}_+$ and anti-holomorphic, hence harmonic. We are then led to look for a solution to \eqref{eqDirLap} having the form
\begin{equation}\label{eqDorFourierSolved}
    U(x+iy)=\frac{1}{2\pi}\int_{\mathbb{R}}\widehat{f}(\omega)e^{-y|\omega|}e^{i\omega x}d\omega=(P_y\ast f)(x),
\end{equation}
where
\begin{eqnarray}\label{eqPoissonViaFourier}
    P_y(x)&=&\frac{1}{2\pi}\int_{\mathbb{R}}e^{-y|\omega|}e^{i\omega x}d\omega\\ 
    &=&\frac{1}{2\pi}\int_0^\infty e^{iz\omega}d\omega+\frac{1}{2\pi}\int_{-\infty}^0e^{i\overline{z}\omega}d\omega\\ 
    &=&\frac{1}{2\pi}\left(\left[\frac{e^{iz\omega}}{iz}\right]_{\omega=0}^{+\infty}+\left[\frac{e^{i\overline{z}\omega}}{i\overline{z}}\right]_{\omega=-\infty}^0\right)\\ 
    &=&\frac{1}{2\pi i}\left(-\frac{1}{z}+\frac{1}{\overline{z}}\right)=\frac{z-\overline{z}}{2\pi i|z|^2}\\ 
    &=&\frac{1}{\pi}\frac{y}{x^2+y^2},
    \end{eqnarray}
which is the formula we wanted.

\vskip0.5cm

    The {\it Poisson kernel }  satisfies the following crucial properties:
\begin{enumerate}
    \item[(i)] $P_y>0$ and $\int_{\mathbb{R}}P_y(x)dx=1$;
    \item[(ii)] for each $\delta>0$, $P_y(x)\to0$ as $y\to0$ uniformly in $|x|\ge\delta$;
    \item[(iii)] $z=x+iy\mapsto P_y(x)$ is {\it harmonic} on $\mathbb{C}_+$, and the same holds for $z\mapsto P_y(x-u)$ for each fixed real $u$. This means that $[\partial_{xx}+\partial_{yy}](P_y(x))=0$ for all $z=x+iy\in\mathbb{C}_+$.
\end{enumerate}
The first two properties say that $\{P_y\}_{y>0}$ is an {\it approximation of the identity} on $\mathbb{R}$, whose special nature is encoded in the third.
Property (i) follows from the expression for $P_y$ in the "time" $x$ variable and from that in the "frequency" $\omega$ one. Property (ii) is easily verified, and property (iii) follows by differentiating \eqref{eqPoissonViaFourier} under the integral sign.

\begin{proof}[Proof of proposition \ref{propDirProbLap}] We start by showing that the expression in \eqref{eqPossonIntegral} is a solution to the Dirichlet problem. By differentiating under the integral and property (iii) of the Poisson kernel, $U$ is harmonic on $\mathbb{C}_+$. It is bounded by property (i) and an elementary estimate. We have to show that for all $x$ in $\mathbb{R}$,
\begin{equation}\label{eqPoissonLimit}
    \lim_{\mathbb{C}_+\ni z\to x}U(z)=f(x).
\end{equation}
We estimate each of the two summands on the right of
\[
|U(x+iy)-f(x_0)|\le|U(x+iy)-U(x_0+iy)|+|U(x_0+iy)-f(x_0)|.
\]
Fix $\epsilon>0$. For the first term, fix $\delta\le1$ such that $\sup_{|u|\le1}|f(x-u)-f(x_0-u)|\le\epsilon$ if $|x-x_0|\le\delta$, so that (by uniform continuity of $f$ on compact intervals)
\begin{eqnarray*}
    |U(x+iy)-U(x_0+iy)|&=&\frac{1}{\pi}\left|\int_{\mathbb{R}}[f(x-u)-f(x_0-u)]\frac{y}{u^2+y^2}\right|\\ 
    &\le&\frac{\sup_{|u|\le1}|f(x-u)-f(x_0-u)|}{\pi}\\ 
    &\ &+\frac{2\|f\|_{L^\infty_0}}{\pi}\int_{|u|\ge1}\frac{y}{u^2+y^2}du\\ 
    &\le&\frac{\epsilon}{\pi}+\frac{2\|f\|_{L^\infty_0}}{\pi}\epsilon,
\end{eqnarray*}
provided $0<y<y_{\epsilon}$.

For the second term,
\begin{eqnarray*}
    |U(x_0+iy)-f(x_0)|&=&\frac{1}{\pi}\left|\int_{\mathbb{R}}[f(x_0-u)-f(x_0)]\frac{y}{u^2+y^2}du\right|\\ 
    &\le&\frac{\epsilon}{\pi}+\frac{2\|f\|_{L^\infty_0}}{\pi}\int_{|u|\ge\delta}\frac{y}{u^2+y^2}du\\ 
    &\le&\frac{\epsilon}{\pi}+\frac{2\|f\|_{L^\infty_0}}{\pi}\epsilon
\end{eqnarray*}
provided $0<y<y_\epsilon'$. Overall, if $0<y<\max(y_\epsilon,y_\epsilon')$ and $|x-x_0|\le \delta$, then $|U(x+iy)-f(x_0)|\le C\epsilon$.

To show uniqueness it suffices to show that a bounded function $U$ which is harmonic in $\mathbb{C}_+$, continuous on $\text{Cl}\mathbb{C}_+$, and that it vanishes on $\mathbb{R}$, has to be constant. Extend such function to one which is defined on the whole of $\mathbb{C}$ by setting $U(\overline{z})=-U(z)$. By the reflection principle, it is harmonic and bounded on $\mathbb{C}$, hence is has to vanish identically by Liouville's theorem.
\end{proof}

\vskip0.5cm

\noindent{\bf The Poisson kernel via (at this point, formal) spectral theory.} The operator $\Delta=\partial_{xx}$ is symmetric on the domain $C_c^2(\mathbb{R})$ which is dense in $L^2(\mathbb{R})$, and $-\Delta$ is positive,
\[
\langle -\Delta\varphi,\varphi\rangle_{L^2}=\int_\mathbb{R}|\varphi'(x)|^2dx\ge0.
\]
We are looking for a function $U:[0,\infty)\to L^2(\mathbb{R})$ such that
\begin{align}\label{eqDirInOT}
    &U''(y)+\Delta U(y)=0,\\ 
    &U(0)=f,\\ 
    &U(+\infty)=0.
\end{align}
Interpreting $\-\Delta$ as a "negative constant", the first equation in \ref{eqDirInOT} has formal solution
\[
U(y)=e^{-\sqrt{-\Delta}y}A+e^{\sqrt{-\Delta}y}B,
\]
with $A,B$ "constants" in $L^2(\mathbb{R})$. Imposing the boundary conditions we have $B=0$ (at least on an intuitive level), and $A=f$. Thus,
\begin{equation}\label{eqDirAbstract}
    U(y)=e^{-\sqrt{-\Delta}y}f,
\end{equation}
or, thinking of $U$ as a function of $(x,y)$,
\[
U(x,y)=\left[e^{-\sqrt{-\Delta}y}f\right](x).
\]
This simple, so far formal, calculation has the advantage that it works for any positive operator $L$ instead of $-\Delta$. For instance, $L$ might be minus the Laplace-Beltrami operator on a geodesically complete Riemannian manifold, or the generator of a random walk on a graph, or minus a fractional Laplacian, etcetera.

In our case, the Fourier transform provides a "spectral representation" for the operator $-\Delta$; 
\[
[\mathcal{F}(-\Delta f)](\omega)=\omega^2[\mathcal{F}f](\omega).
\]
which is another way to say that $-\Delta$ acts like multiplication times $\omega^2$ on the frequency side. The operator $e^{-y\sqrt{-\Delta}}$, then, acts as multiplication times $e^{-y|\omega|}$, and we are back to the previous deduction of the the Poisson kernel's expression.

\vskip0.5cm

\noindent{\bf The Poisson kernel via the fundamental solution for the Laplace operator.} A third way to construct Poisson kernels is by means of Green's formula and the fundamental solution for the Laplace operator. Let $G(z,w)=\frac{1}{2\pi}\log|z-w|^{-1}$, $z,w\in\mathbb{C}$, which satisfies
\begin{align}\label{eqFundSol}
    &\Delta_zG(z,w)=\delta_z(w),\\
    &G(z,w)=O_{z\to\infty}(\log|z|).
\end{align}
The first equation means that
\[
\int_\mathbb{C}G(z,w)\Delta\varphi(w)dudv=\varphi(z)
\]
for all $\varphi\in C^2_c(\mathbb{C})$.

Let $\Omega$ be a bounded domain in $\mathbb{C}$ with smooth boundary $\partial\Omega$, and denote by $\nu$ the external unit normal. Let $D(z,\epsilon)$ be the disc with center $z$ and radius $\epsilon$. Suppose $H(z,w)$ is defined for $z,w\in \text{Cl}\Omega$, $z\ne w$, and has the properties:
\begin{align}\label{eqProGreenOmega}
    &\Delta_wH(z,w)=0\text{ if }z\ne0;\\ 
    &H(z,w)=O_{w\to z}(\log|z-w|);\\ 
    &\partial_r H(z,z+r e^{it})=\frac{1}{2\pi r}(1+o_{r\to 0}(1));\\ 
    &H(z,\xi)=0\text{ if $\xi\in\partial\Omega$ and $z\in\Omega$}.
\end{align}
The first three properties also hold for $G$, while the fourth is linked with the domain $\Omega$. If a function $U:\text{Cl}\Omega\to\mathbb{C}$ is $C^2$ and $\Delta U=0$ in $\Omega$, then, by Green's formula,
\begin{eqnarray*}
    0&=&\int_{\Omega\setminus D(z,\epsilon)}[H(z,w)\Delta U(w)-\Delta_wH(z,w)U(w)]dudv\\ 
    &=&\int_{\partial\Omega-\partial D(z,\epsilon)}\left[H(z,\zeta)\partial_\nu U(\zeta)-\partial_\nu H(z,\zeta) U(\zeta)\right]d\sigma(\zeta)\\ 
    &=&\int_{\partial D(z,\epsilon)}\partial_\nu H(z,\zeta) U(\zeta)d\sigma(\zeta)-\int_{\partial\Omega}\partial_\nu H(z,\zeta) U(\zeta)d\sigma(\zeta)\\ 
    &\to&U(z)-\int_{\partial\Omega}\partial_\nu H(z,\zeta)U(\zeta)d\sigma(\zeta)
\end{eqnarray*}
as $\epsilon\to0$. (We did not use the second equation in \eqref{eqProGreenOmega}, which becomes relevant when solving the Poisson equation: the nonhomogeneous version of the Laplace equation). At the end of the day, we have the inspiring formula:
\begin{equation}\label{eqPoissonOmega}
    U(z)=\int_{\partial\Omega}\partial_\nu H(z,\zeta)U(\zeta)d\sigma(\zeta),
\end{equation}
which reconstructs the interior values of a function $U$ harmonic in $\Omega$ by its boundary values. We set $P^\Omega_z(\zeta):=\partial_\nu H(z,\zeta)$ the {\it Poisson kernel for $\Omega$ at $z$}. 

This formula can not be immediately applied to the upper-half plane, which is not bounded. We might however borrow the general idea. The function $H:\text{Cl}\mathbb{C}_+\times\text{Cl}\mathbb{C}_+\setminus\{(z,z):z\in\text{Cl}\mathbb{C}_+\}\to\mathbb{R}$,
\begin{equation}\label{eqGreenHalfPlane}
    H(z,w)=G(z,\overline{w})-G(z,w)=\frac{1}{2\pi}\left(\log|z-\overline{w}|^{-1}-\log|z-w|^{-1}\right)
\end{equation}
satisfies all conditions in \eqref{eqProGreenOmega} with $\Omega=\mathbb{C}_+$.
The normal derivative at the boundary is $-\partial_v$,
\begin{eqnarray*}
    -\partial_v|_{v=0}H(z,u+iv)&=&\frac{1}{2\pi}\left(\frac{y}{|z-u|^2}-\frac{-y}{|z-u|^2}\right)\\ 
    &=&\frac{1}{\pi}\frac{y}{(x-u)^2+y^2}.
\end{eqnarray*}

\section{The Fourier transforms of Borel measures}\label{SSectBorelMeas}
We start with an approximation lemma for finite Borel measures. I include the approximation by point masses because it has many applications and generalizations, although we will only make use of the approximation by $L^1$ functions.
\begin{lemma}\label{lemmaApproxMeas}
Let $\mu$ be a finite Borel measure on $\mathbb{R}$, not necessarily finite. We can find a sequence $\{\mu_n\}$ of measures made of a finite number of atoms such that $\mu_n\xrightarrow[w^\ast]{n\to\infty}\mu$. 

Another approximation of $\mu$ is provided by the Poisson integral.
\begin{equation}\label{eqAppPoissonMeasure}
    P_y[\mu]\xrightarrow[w^\ast]{y\to0}\mu.
\end{equation}
We also have that $\|P_y[\mu]\|_{L^1}\le|\mu|(\mathbb{R})$, hence $\mu$ can be $w^\ast$-approximated by $L^1$ functions.
\end{lemma}
\begin{proof} It suffices to show the statement for a positive $\mu$, since any finite complex Borel measure is a linear combination of four positive ones.
 Let $\alpha_\mu(t)=\alpha(t)=\mu(-\infty,t]$ be the distribution function of $\mu$, which is increasing, right continuous, $\alpha(t)\to0$ as $t\to-\infty$, and $\alpha(t)\to\mu(\mathbb{R})$ as $t\to+\infty$. For $\varphi\in C_0(\mathbb{R})$, the following limit is easily proved using its uniform continuity:
\begin{align}\label{eqApproxMeas}
    &\int_{\mathbb{R}}\varphi d\mu=\lim_{n\to\infty}\sum_{j:|j2^{-n}|\le n}\varphi\left(j2^{-n}\right)[\alpha(j2^{-n})-\alpha((j-1)2^{-n})],\\ 
\end{align}
Equation \eqref{eqApproxMeas} provides the desired limit with
\begin{align*}
&\mu_n=\sum_{j:|j2^{-n}|\le n}[\alpha(j2^{-n})-\alpha((j-1)2^{-n})]\delta_{j2^{-n}}\\ 
\end{align*}

About the second statement, let $\varphi\in C_0(\mathbb{R})$. Then,
\begin{eqnarray*}
    \int_{\mathbb{R}}\varphi(x)P_y[\mu](x)dx&=&\int_{\mathbb{R}}\varphi(x)\frac{1}{2\pi}\int_{\mathbb{R}}\frac{y}{(x-u)^2+y^2}d\mu(u)dx\\ 
    &=&\int_{\mathbb{R}}\frac{1}{2\pi}\int_{\mathbb{R}}\frac{y}{(x-u)^2+y^2}\varphi(x)dxd\mu(u)\\ 
    &=&\int_{\mathbb{R}}P_y[\varphi](u)d\mu(u).
\end{eqnarray*}
From the proof of proposition \ref{propDirProbLap} it is clear that $P_y[\varphi]$ converges uniformly to $\varphi$ as $y\to0$ if $\varphi$ is bounded and uniformly continuous, which is the present case since $\varphi\in C_0$. Then, the last term in the chain of equalities tends to $\int_\mathbb{R}\varphi(x)dx$, as wished. The estimate $\|P_y[\mu]\|_{L^1}\le|\mu(\mathbb{R})|$ holds because the Poisson kernel has integral one.
\end{proof}
\begin{proposition}\label{propFTmeasures}
    Let $\mu$ be a finite Borel measure on $\mathbb{R}$. Then, $\widehat{\mu}:\mathbb{R}\to\mathbb{C}$ is bounded and uniformly continuous.
\end{proposition}
\begin{proof}
   We have $|\widehat{\mu}(\omega)|\le|\mu(\mathbb{R})|$, hence $\widehat{\mu}$ is bounded. Uniform continuity is proved as in the case of $L^1$ functions,
   \begin{eqnarray*}
       |\widehat{\mu}(\omega)-\widehat{\mu}(\eta)|&\le&\int_{\mathbb{R}}\left|e^{i\omega x}-e^{i\eta x}\right|d|\mu|(x)\\ 
       &\le&2\int_{x:|x|> R}d|\mu|(x)+\int_{x:|x|\le R}\left|e^{i(\omega-\eta)x}\right|d|\mu|(x).
   \end{eqnarray*}
   For fixed $\epsilon>0$, choose first $R$ such that the integral in first summand is bounded by $\epsilon$ (dominated convergence), then choose $\delta$ such that $\left|e^{i(\omega-\eta)x}\right|\le\epsilon$ if $|\omega-\eta|\le \delta$.
\end{proof}
\begin{theorem}\label{theoUniqueFTmeas}
    Let $\mu$ be a finite Borel measure on $\mathbb{R}$. If $\widehat{\mu}=0$, then $\mu=0$.
\end{theorem}
\begin{proof}
    Let $P_y[\mu]=P_y\ast\mu$ be the Poisson approximation of $\mu$. For all $\varphi\in C_0(\mathbb{R})$:
    \begin{eqnarray*}
        \int_{\mathbb{R}}P_y[\mu](x)\overline{\varphi(x)}dx&=&\frac{1}{2\pi}\int_{\mathbb{R}}\widehat{P_y[\mu]}(\omega)\overline{\widehat{\varphi}}(\omega)d\omega\\ 
        &=&\frac{1}{2\pi}\int_{\mathbb{R}}\widehat{\mu}(\omega)e^{-y|\omega|}\overline{\widehat{\varphi}}(\omega)d\omega\\ 
        &=&0.
    \end{eqnarray*}
    As $y\to0$, by lemma \ref{lemmaApproxMeas} we have that
    \[
    \int_{\mathbb{R}}\overline{\varphi(x)}d\mu(x)=\lim_{y\to0}\int_{\mathbb{R}}P_y[\mu](x)\overline{\varphi(x)}dx=0.
    \]
    Since this holds for all $\varphi\in C_0$, $\mu=0$ as a measure.
\end{proof}
We close with a useful variation on Plancherel formula, which holds if $\mu$ is a finite measure and $\widehat{\varphi}\in L^1(\mathbb{R})$:
\begin{equation}\label{eqPlachMeas}
   \int_{\mathbb{R}}\varphi(-x)d\mu(x)=\frac{1}{2\pi}\int_{\mathbb{R}}\widehat{\varphi}(\omega)\widehat{\mu}(\omega)d\omega..
\end{equation}
In fact,
\begin{eqnarray*}
    \int_{\mathbb{R}}\varphi(-x)d\mu(x)&=&\int_{\mathbb{R}}\frac{1}{2\pi}\int_{\mathbb{R}}\widehat{\varphi}(\omega)e^{-i\omega x}d\omega d\mu(x)\\ 
    &=&\frac{1}{2\pi}\int_{\mathbb{R}}\widehat{\varphi}(\omega)\widehat{\mu}(\omega)d\omega.
\end{eqnarray*}

\section{Herglotz representation}\label{SsectHerglotz} The theorem of Herglotz we present here will be the key to determine the spectral measures for a self-adjoint operator, which will be the ones associated with holomorphic functions indexed by vectors in a Hilbert space.

Let $\mu$ be a finite Borel measure on $\mathbb{R}$. Its {\it Poisson integral} on $\mathbb{C}_+$ is defined in the natural way.
\begin{equation}\label{eqPoiIntMeas}
    P[\mu](z):=(P_y\ast\mu)(x)=\int_{\mathbb{R}}P_y(x-u)d\mu(u)=\frac{1}{\pi}\int_{\mathbb{R}}\frac{y}{(x-u)^2+y^2}d\mu(u).
\end{equation}
If we want to fix the height $y$, we write
\[
P_y[\mu](x):=P[\mu](x+iy),\ P_y[\mu]:\mathbb{R}\to\mathbb{R}.
\]
The function $P[\mu]$ is (i) harmonic in $\mathbb{C}_+$ (differentiate under the integral);  (ii) positive, if the measure $\mu$ is positive.
\begin{remark}\label{remPoiKerPosMeas}
    Let $\mu\ge0$ be a finite Borel measure on $\mathbb{R}$. Then $P[\mu]:\mathbb{C}_+\to\mathbb{R}$ is a positive harmonic function, and
    \begin{equation}\label{eqPoiIntMeasDecay}
        P[\mu](x+iy)\le\frac{\mu(\mathbb{R})}{\pi y},
    \end{equation}
    since $P_y(x-u)\le\frac{1}{\pi y}$.
\end{remark}
The remark has a surprising converse.
\begin{theorem}\label{theoHerglotzHarm}[Herglotz theorem for positive harmonic functions\index{theorem!of Herglotz for harmonic functions}]
    Let $U:\mathbb{C}_+\to[0,\infty)$ be a positive, harmonic function, and suppose that $U(x+iy)\le\frac{C}{y}$ for some constant independent of $x,y$. Then, there exists a positive, finite measure $\mu$ on $\mathbb{R}$ such that $U=P[\mu]$. Moreover, such measure is unique.
\end{theorem}
    \begin{proof}
        We do some bootstrapping followed by an application of the Hahn-Banach theorem. For fixed $\epsilon>0$, let $U_\epsilon(x)=U(x+i\epsilon)$, which is positive and harmonic for $z$ in $\mathbb{C}_+$, continuous up to the real line, and bounded, $U_\epsilon(x)\le\frac{C}{\epsilon}$. Since the Poisson kernel is integrable, we can form
        \[
        P[U_\epsilon](z)=P_y[U_\epsilon](x)=\frac{1}{\pi}\int_{\mathbb{R}}\frac{y}{(x-u)^2+y^2}U_\epsilon(u)du.
        \]
        \noindent{\bf Step I.} We will show the semigroup law $U_{\epsilon+y}=P_y[U_\epsilon]$ for all $y>0$. Let $G_\epsilon(x+iy)=P_y[U_\epsilon](x)$ and $H_\epsilon(x+iy)=U_{\epsilon+y}(x)$. Both functions are positive and harmonic on $\mathbb{C}_+$; they continuous on $\mathbb{C}_+\cup\mathbb{R}$, the closure of the upper-half plane; $G_\epsilon(x)=U(x+i\epsilon)=H_\epsilon(x)$ for $x\in\mathbb{R}$; they are both bounded. Their difference $K_\epsilon=H_\epsilon-G_\epsilon$ is then harmonic on the upper-half plane, bounded, continuous on its closure, and vanishing on $\mathbb{R}$. We can use the reflection principle to extend $K_\epsilon$ to a harmonic, bounded ($K_\epsilon(\overline{z}):=-K_\epsilon(z)$) function defined on $\mathbb{C}$, which vanishes on the real line.
        By Liouville theorem, a bounded harmonic function on $\mathbb{C}$ must be constant, hence $K_\epsilon(z)=K_\epsilon(0)=0$ for all complex $z$. This proves our claim.

        \noindent{\bf Step II.} We show that $1/\pi\int_{\mathbb{R}}U_\epsilon(u)du\le C$, where $C$ is the constant in the hypothesis, which is independent of $\epsilon>0$. In fact, by the assumption on the decay of $U(x+iy)$ with respect to to $y$, and using Step I,
        \begin{eqnarray*}
            \frac{CR}{\epsilon+R}&\ge&R U_{\epsilon+R}(x)=R P_R[U_\epsilon](x)\\ 
            &=&\frac{1}{\pi}\int_{\mathbb{R}}\frac{R^2}{(x-u)^2+R^2}U_\epsilon(u)du\\
            &=&\frac{1}{\pi}\int_{\mathbb{R}}\left(1-\frac{(x-u)^2}{(x-u)^2+R^2}\right)U_\epsilon(u)du\\ 
            &\to&\frac{1}{\pi}\int_{\mathbb{R}}U_\epsilon(u)du\text{ as }R\to\infty,
        \end{eqnarray*}
        by monotone convergence. We proved the claim, since the limit as $R\to\infty$ of the first term in the chain is $C$.

        \noindent{\bf Step III.} Consider the family of finite positive Borel measures $d\mu_\epsilon(u)=U_\epsilon(u)du$. Since their norms $\mu_\epsilon(\mathbb{R})\le C$ are uniformly bounded, by the Banach-Alaoglu theorem there exist a sequence $\epsilon_n\to0$ and a positive Borel measure $\mu$ with $\mu(\mathbb{R})\le C$  such that $\mu_{\epsilon_n}\xrightarrow[w^\ast]{}\mu$. Recall that (by the Riesz representation of measures) the space of the finite Borel measures on $\mathbb{R}$ is the dual of $C_0(\mathbb{R})$, the space of the continuous functions tending to zero as $x\to\pm\infty$. In particular, for all $x+iy$ in $\mathbb{C}_+$,
        \begin{eqnarray*}
            U(x+iy)&=&\lim_{n\to\infty}U_{\epsilon_n}(x+iy)=\lim_{n\to\infty}\frac{1}{\pi}\int_{\mathbb{R}}\frac{y}{(x-u)^2+y^2}U_{\epsilon_n}(u)du\\ 
            &=&\frac{1}{\pi}\int_{\mathbb{R}}\frac{y}{(x-u)^2+y^2}d\mu(u),
        \end{eqnarray*}
        as wished.

        \noindent{\bf Step IV.} We prove uniqueness. The Herglotz representation formula can be written as $U_y=P_y\ast\mu$. 
        It is well known (and easy to verify) that 
        \begin{equation}\label{eqPoissonFourier}
            \widehat{P_y}(\omega)=e^{-y|\omega|}.
        \end{equation}
        If $U_y=P_y\ast\mu$, then $\widehat{U_y}(\omega)=e^{-y|\omega|}\widehat{\mu}(\omega)$. This shows that $\widehat{\mu}$ is uniquely determined by $U$, and it is known that $\mu$ is determined by $\widehat{\mu}$.
    \end{proof}
    \begin{corollary}\label{corHerglotzComplex}[Herglotz theorem for holomorphic functions $f:\mathbb{C}_+\to\mathbb{C}_+\cup\mathbb{R}$\index{theorem!of Herglotz for holomorphic functions}] Let $f:\mathbb{C}_+\to\mathbb{C}_+\cup\mathbb{R}$ be holomorphic and such that
    \begin{equation}\label{eqHeglotzDecayC}
        \text{Im}f(x+iy)\le\frac{C}{y}.
    \end{equation}
        Then, there exist a unique finite, positive Borel measure $\mu$ and a real number $a$ such that 
        \begin{equation}\label{eqHerglotzThesisC}
            f(z)=\frac{1}{\pi}\int_{\mathbb{R}}\frac{d\mu(u)}{z-u}+a.
        \end{equation}
    \end{corollary}
\begin{proof}
    We have that
    \[
    \frac{1}{\pi}\frac{1}{u-(x+iy)}=\frac{1}{\pi}\frac{(u-x)+iy}{(u-x)^2+y^2}=Q_y(x-u)+iP_y(x-y)
    \]
    is holomorphic with respect to $z=x+iy$. If $f=\varphi+i\psi$ is as in the hypothesis, then, $\psi$ satifies the hypothesis of Herglotz theorem, then there exists $\mu\ge0$ (Borel, finite) such that $\psi=P[\mu]$. The function
    \[
    Q[\mu](x+iy)=\frac{1}{\pi}\int_{\mathbb{R}}\frac{u-x}{(u-x)^2+y^2}d\mu(u)
    \]
    is harmonic in $\mathbb{C}$ and $Q[\mu]+iP[\mu]=Q[\mu]+i\psi$ is holomorphic. Hence, $Q[\mu]-\varphi$ is a constant.

    Uniqueness is easily checked.
\end{proof}
\section{Bochner's theorem}\label{SSectBochner}
This is here because it's the key for one of the proofs of Stone's theorem on groups of unitary operators. The result itself is one of the cornerstones of Fourier analysis and has applications to analysis, probability, and statistics.

A function $f:\mathbb{R}\to\mathbb{C}$ is {\it positive definite\index{positive definite function}} if
\[
\sum_{j,k=1}^n\varphi(x_k-x_j)c_j\overline{c_k}\ge0
\]
for all choices of $n$ points $x_1,\dots,x_n$ in $\mathbb{R}$ and $n$ complex scalars $c_1,\dots,c_n$. This is the same as saying that for all choices of $x_1,\dots,x-n$, the matrix $[\varphi(x_k-x_j)]_{k,j=1}^n$ is definite positive. 

Consider the points $0,x$. The matrix
\[
\begin{pmatrix}
    \varphi(0)&\varphi(x)\\ 
    \varphi(-x)&\varphi(0)
\end{pmatrix}
\]
is positive definite (hence, Hermitian, since the scalar field is $\mathbb{C}$), thus
\begin{align}\label{eqPDFbasics}
&\varphi(-x)=\overline{\varphi(x)},\\
&\varphi(0)\ge|\varphi(x)|.
\end{align}
\begin{theorem}\label{theoBochner}[Bochner's theorem\index{theorem!Bochner}\index{Borel!measure}]
    For a function $f:\mathbb{R}\to\mathbb{C}$ the following are equivalent:
    \begin{enumerate}
        \item[(i)] $f$ is continuous and positive definite;
        \item[(ii)] there exists a finite, positive Borel measure $\mu$ on the real line such that $\widehat{\mu}=f$.  
    \end{enumerate}
\end{theorem}
The idea of the proof can be best appreciated in a finite environment. Let $\mathbb{Z}_N=\{0,1,\dots,N-1\}$ be the group of the integers modulo $N$, and for $f:\mathbb{Z}_N\to\mathbb{C}$, let its {\it discrete Fourier transform} $\widehat{f}:\mathbb{Z}_N\to\mathbb{C}$ be
\[
\widehat{f}(\omega)=\sum_{n=0}^{N-1}e^{-2\pi i n\omega/N}f(n).
\]
The usual formulas hold:
\begin{align*}
  &f(n)=\frac{1}{N}\sum_{\omega=0}^{N-1}e^{2\pi i n\omega/N}\widehat{f}(\omega);\\
  &\sum_{n=0}^{N-1}\overline{\nu(n)}f(n)=\frac{1}{N}\sum_{\omega=0}^{N-1}\overline{\widehat{\nu}(\omega)}\widehat{f}(\omega);\\
  &\widehat{f\ast\nu}(\omega)=\widehat{f}(\omega)\widehat{\nu}(\omega),
\end{align*}
where 
\[
(f\ast\nu)(n)=\sum_{j=0}^{N-1}f(n-j)\nu(j).
\]
Using all this, if $f,c:\mathbb{Z}_N\to\mathbb{C}$,
\begin{eqnarray}\label{eqBochnerProofD}
    \frac{1}{N}\sum_{\omega=0}^{N-1}\widehat{f}(\omega)|\widehat{c}(\omega)|^2&=&\frac{1}{N}\sum_{\omega=0}^{N-1}\overline{\widehat{c}(\omega)}[\widehat{f}\widehat{c}](\omega)\\ 
    &=&\sum_{n=0}^{N-1}\overline{c(n)}[f\ast c](n)\\ 
    &=&\sum_{m,n=0}^{N-1}f(n-m)\overline{c(n)}c(m).
\end{eqnarray}
Comparing the first and the last terms in the chain of equalities, it is clear that if $\widehat{f}$ is positive, then $f$ is positive definite. Conversely, if $f$ is positive definite, then the first term is positive for each choice of $c$. For any fixed $\omega_0$ in $\mathbb{Z}_N$, we can choose $c$ such that $\widehat{c}(\omega)=\delta_{\omega_0,\omega}$, obtaining $\widehat{f}(\omega_0)\ge0$.

In the general case, the idea is the same, but there are complications. We start with a useful result in linear algebra.

If $A,B$ are $N\times N$ matrices with complex entries, their {\it Schur (or Hadamard) product} is the $N\times N$ matrix $A\circ B$ having entries $[A\circ B]_{ij}=A_{ij}B_{ij}$. The matrix $A$ is {\it positive definite} if 
\[
0\le \sum_{ij}\overline{c_i}A_{ij}c_j=\langle c,Ac\rangle_{\mathbb{C}^N}
\]
for all $c\in\mathbb{C}^N$. Sometimes this property is called semi-positive definiteness. We say that $A$ is {\it strictly positive definite} if $\langle c,Ac\rangle_{\mathbb{C}^N}>0$ whenever $c\ne0$.
\begin{theorem}[Schur product theorem]\label{theoProdOcPDF}
    If $A,B$ are positive definite $N\times N$ matrices, then $A\circ B$ is positive definite.
\end{theorem}
For the proof, we lift the problem to a more general one living on the {\it tensor product} $\mathbb{C}^N\otimes\mathbb{C}^N$, the space having as basis $\{e_i\otimes e_j:1\le i,j\le N\}$. The {\it tensor product} $\otimes$ is bilinear by definition,
\begin{align*}
    (ax+by)\otimes z=a(x\otimes z)+b(y\otimes z),\\
    z \otimes (ax+by)=a(z\otimes x)+b(z\otimes y),
\end{align*}
for all $x,y,z\in\mathbb{C}^N$ and $a,b\in\mathbb{C}$; and no other relation is imposed.

Here is a concrete model of the tensor product. The elements of $\mathbb{C}^N\otimes\mathbb{C}^N$ might be seen as linear combinations of vectors of the form
\[
h\otimes k=h k^t\in \mathbb{C}^{N\times N},
\]
where $k^t$ is the transpose of the column vector $k$. This way,  we have a linear isomorphism $\mathbb{C}^N\otimes\mathbb{C}^N\to\mathbb{C}^{N\times N}$, which identifies elements in the tensor product with matrices.

A linear map $A:\mathbb{C}^N\to\mathbb{C}^N$ induces two linear maps 
\[
\mathbb{C}^N\otimes\mathbb{C}^N\xrightarrow{A\otimes I,I\otimes A}\mathbb{C}^N\otimes\mathbb{C}^N,
\]
where
\begin{align*}
    (A\otimes I)(x\otimes y)=(Ax)\otimes y,\\ 
    (I\otimes A)(x\otimes y)=x\otimes(Ay).
\end{align*}
More generally,
\[
(A\otimes B)(x\otimes y)=(A\otimes I)(I\otimes B)(x\otimes y)=(Ax\otimes By).
\]
In the column$\times$row model, you might think of $A\otimes B$ as the linear map $xy^t\mapsto (Ax)(By)^t=Axy^tB^t$.

The vector space $\mathbb{C}^N\otimes\mathbb{C}^N$ is endowed with the inner product\index{inner product} making tensor products $e_i\otimes e_j$ of standard basis elements of $\mathbb{C}^N$ into a orthonormal basis. This gives
\begin{equation}\label{eqOinnerproduct}
\langle x\otimes u,y\otimes v\rangle=\langle x,y \rangle\cdot \langle u,v\rangle.
\end{equation}
In fact,
\begin{eqnarray*}
    \langle x\otimes u,y\otimes v\rangle&=&\sum_{ijlm}\overline{x_i}u_j\overline{y_l}v_m\langle e_i\otimes e_j,e_l\otimes e_m\rangle\\ 
    &=&\sum_{ij}\overline{x_i}u_j\overline{y_i}v_j\\
    &=&\langle x,y \rangle\cdot \langle u,v\rangle.
\end{eqnarray*}
In the column$\times$row model,
\[
\langle xu^t,yv^t \rangle=y^t\overline{x}v^t\overline{u}.
\]
A consequence of \eqref{eqOinnerproduct} is that we can replace the standard basis of $\mathbb{C}^N$ by any other orthonormal basis in the calculations. In fact, \eqref{eqOinnerproduct} might be takes as the definition of the inner product on $\mathbb{C}^N\times\mathbb{C}^N$.
\begin{lemma}[Tensor product preserves positivity]\label{lemmaPosTensProd}
    If $A$ and $B$ are positive, then $A\otimes B$ is positive.
\end{lemma}
\begin{proof}[Proof of lemma \ref{lemmaPosTensProd}] Choose a orthonormal basis $\{\epsilon_i\}$ for $\mathbb{C}^N$ with respect to which $B$ is diagonal,
$B_{ij}=\delta_{ij}\lambda i$ with $\lambda_1,\dots,\lambda_N\ge0$.
For a vector $V=\sum_{ij}v_{ij}\epsilon_i\otimes \epsilon_j\in \mathbb{C}^N\otimes\mathbb{C}^N$, we have
\begin{eqnarray*}
    \langle V,(A\otimes B) V\rangle&=&\sum_{ijlm}\overline{v_{ij}}v_{lm}\langle \epsilon_i\otimes \epsilon_j,(A\otimes B) (\epsilon_l\otimes \epsilon_m)\rangle\\ 
    &=&\sum_{ijlm}\overline{v_{ij}}v_{lm}\langle \epsilon_i,A\epsilon_l\rangle \langle \epsilon_j,B\epsilon_m\rangle\\ 
    &=&\sum_{ijlm}\overline{v_{ij}}v_{lm}A_{il}B_{jm}\\
    &=&\sum_{jm}B_{jm}\left(\sum_{il}\overline{v_{ij}}v_{lm}A_{il}\right)\\ 
    &=&\sum_j\lambda_j\left(\sum_{il}\overline{v_{ij}}v_{lj}A_{il}\right)\\
    &\ge&0,
\end{eqnarray*}
since the sums in the parenthesis are positive and so are the $\lambda_j$'s.
\end{proof}
\begin{proof}[Proof of Schur's theorem given lemma \ref{lemmaPosTensProd}]
    By the lemma above,
    \[
    0\le\langle V,(A\otimes B) V\rangle=\sum_{ijlm}\overline{v_{ij}}v_{lm}A_{il}B_{jm}.
    \]
    Choose in particular $v_{ij}=\delta_{ij}c_i$, to obtain
    \[
    0\le \sum_{il}\overline{c_i}c_lA_{il}B_{il},
    \]
    as wished.
\end{proof}
The topic we have just sketched can be interpreted in terms of, and applied to, the theory or {\it Hilbert spaces with reproducing kernel}, which is a very active area in pure mathematics and in its applications to a number of areas, including statistics and machine learning. See e.g. \href{https://www.cambridge.org/core/books/an-introduction-to-the-theory-of-reproducing-kernel-hilbert-spaces/C3FD9DF5F5C21693DD4ED812B531269A}{An Introduction to the Theory of Reproducing Kernel Hilbert Spaces (2016)} by 
Vern I. Paulsen, Mrinal Raghupathi, Cambridge University Press.
\begin{proof}[Proof of Bochner's theorem] If $\mu\ge0$ is a finite Borel measure, then $\widehat{\mu}$ is continuous and bounded. About positive definiteness,
\begin{eqnarray*}
    \sum_{lj}\widehat{\mu}(\omega_l-\omega_j)c_l\overline{c_j}&=&\int_{\mathbb{R}}\sum_{lj}e^{-i(\omega_l-\omega_j)x}c_l\overline{c_j}d\mu(x)\\ 
    &=&\int_\mathbb{R}\left|\sum_l e^{-i\omega_l x}\right|^2d\mu(x)\\ 
    &\ge0&.
\end{eqnarray*}
    In the opposite direction, suppose that $f:\mathbb{R}\to\mathbb{C}$ is positive definite, hence bounded, and continuous. For $y>0$, let $f_y(\omega)=f(\omega)e^{-y|\omega|}$, which is still positive definite by Schur's theorem, and it lies in $L^1(\mathbb{R})$ (hence in $L^2(\mathbb{R})$). Set $m_y(x)=\frac{1}{2\pi}\int_{\mathbb{R}}f_y(\omega)e^{i\omega x}d\omega=[\mathcal{F}^{-1}f_y](x)$, so that $m_y\in L^2(\mathbb{R})\cap C_0(\mathbb{R})$. For all functions $\varphi\in\mathcal{S}(\mathbb{R})$ (the Schwartz class, but if you are not familiar with that object, you might even consider $\varphi\in C_c^\infty(\mathbb{R})$), we have:
    \begin{eqnarray*}
        \int_{\mathbb{R}}m_y(x)|\varphi(x)|^2dx&=&\int_{\mathbb{R}}[m_y(x)\varphi(x)]\overline{\varphi(x)}dx\\ 
        &=&\frac{1}{2\pi}\int_{\mathbb{R}}\widehat{m_y\varphi}(\omega)\overline{\widehat{\varphi}(\omega)}d\omega\\ 
        &=&\frac{1}{2\pi}\int_{\mathbb{R}}[\widehat{m_y}\ast\widehat{\varphi}](\omega)\overline{\widehat{\varphi}(\omega)}d\omega\\ 
        &=&\frac{1}{2\pi}\int_{\mathbb{R}}\int_{\mathbb{R}}\widehat{m_y}(\omega-\eta)\widehat{\varphi}(\omega)\overline{\widehat{\varphi}(\omega)}d\eta d\omega\\ 
        &=&\frac{1}{2\pi}\int_{\mathbb{R}}\int_{\mathbb{R}}f_y(\omega-\eta)\widehat{\varphi}(\omega)\overline{\widehat{\varphi}(\omega)}d\eta d\omega\\ 
        &=&\frac{1}{2\pi}\lim_{N\to\infty}\sum_{l,m=-N^2}^{N^2}f_y\left(l/N-m/N\right)\widehat{\varphi}(l/N)\overline{\widehat{\varphi}(m/N)}\frac{1}{N^2}\\
        &\ge&0
    \end{eqnarray*}
    because $f_y$ is positive definite. Hence, $m_y\ge0$ is a positive function with $\|m_y\|_{L^1}=\widehat{m_y}(0)=f(0)$.

    For $\varphi\in \mathcal{S}(\mathbb{R})$, let
    \begin{eqnarray*}
    \Lambda_y(\varphi)&:=&\int_{\mathbb{R}}m_y(x)\varphi(x)dx=\frac{1}{2\pi}\int_{\mathbb{R}}f(\omega)e^{-y|\omega|}\widehat{\varphi}(-\omega)d\omega\\
    &\to&\frac{1}{2\pi}\int_{\mathbb{R}}f(\omega)\widehat{\varphi}(-\omega)d\omega.
    \end{eqnarray*}
    as $y\to0$.

    On the other hand, by Banach-Alaoglu there exist a sequence $y_n\to0$ and a finite (positive) Borel measure $\mu$ such that $\mu(\mathbb{R})\le f(0)$ and $m_{y_n}dm\xrightarrow[w^\ast]{}\mu$, hence
    \[
    \lim_{n\to\infty}\int_{\mathbb{R}}m_{y_n}(x)\varphi(x)dx=\int_{\mathbb{R}}\varphi(x)d\mu(x).
    \]
    We have shown that
    \[
    \frac{1}{2\pi}\int_{\mathbb{R}}f(\omega)\widehat{\varphi}(-\omega)d\omega=\int_{\mathbb{R}}\varphi(x)d\mu(x)=\frac{1}{2\pi}\int_{\mathbb{R}}\widehat{\mu}(\omega)\widehat{\varphi}(-\omega)d\omega
    \]
    holds for all $\varphi\in \mathcal{S}$ (use \eqref{eqPlachMeas} for the second equality), hence $f=\widehat{\mu}$.
\end{proof}

\chapter{Hilbert spaces}\label{chaptHilbert}
\label{chap:Hilbert_spaces}

Among function spaces, Hilbert spaces are probably the most ubiquitous in mathematics and its applications. At its inception, "Hilbert theory" dealt with concrete $L^2$ spaces, but soon an abstract theory was developed. The advantage of the abstract theory is that it makes it easier to recognize when a Hilbert space structure underlies a cluster of mathematical objects and phenomena, and it helps in generating Hilbert spaces in which better state, and solve, theoretical, as well as practical problems. Although all Hilbert spaces are isomorphic to some $L^2$ space, the point of view that knowledge of $L^2$ is all that's needed is far too simplistic, and basically incorrect. Most Hilbert spaces are spaces of functions defined on some set of points, and such "points" constitute further structure, which often is at the core of the problem we have at hands. Is this chapter, however, we are mostly interested in the abstract theory.

Another useful way of thinking is that of viewing at Hilbert spaces as generalizations (often infinite-dimensional, with complex rather than real scalars) of the Euclidean space. Even if in applications each vector in the Hilbert space represents a functions, we can think of each of them as a "point" in a linear space where notions like orthogonality make perfect sense. This intuition is helpful in translating complex phenomena in simple pictures, and pictures into statements (which are sometimes true, sometimes false).

\section{Basic geometry of Hilbert spaces and Riesz Lemma}\label{SectBasicHilbert}

\subsection{Definition and basic properties}\label{SSectBasicHilbert}
An {\it inner product\index{inner product}} on a vector space $V$ over $\mathbb{C}$ (or $\mathbb{R}$ ) is a map $\langle\cdot , \cdot\rangle: V \times V \rightarrow \mathbb{C}$ such that:
\begin{enumerate}
\item[(i)] $\langle x , x\rangle \geqslant 0$ and $\langle x , x\rangle=0$ if and only if $x=0$;

\item[(ii)] $\langle z , \alpha x+\beta y\rangle=\alpha\langle z , x\rangle+\beta\langle z , y\rangle$ for $x, y, z \in V$ and $\alpha, \beta \in \mathbb{C}$;

\item[(iii)] $\overline{\langle x , y\rangle}=\langle y , x\rangle$.
\end{enumerate}

Two vectors $x, y$ are orthogonal, $x \perp y$, if $\langle x , y\rangle=0$. We define $\|x\|:=\langle x , x\rangle^{1 / 2}$ to be the norm of $x \in V$. A simple calculation gives:
\begin{lemma}\label{LemmaPyth}
[Pythagorean relation] Let $x, y \in V, x \neq 0$. Then,
$$
\|y\|^{2}=\left\|y-\frac{x}{\|x\|}\left\langle  \frac{x}{\|x\|}, y\right\rangle\right\|^{2}+\left|\left\langle  \frac{x}{\|x\|}, y\right\rangle\right|^{2} .
$$    
\end{lemma}
\begin{proof}
    Expand the right hand side using the definition of norm and the properties of the inner product.
\end{proof}
\begin{footnotesize}\begin{exercise}\label{ExeHOne}
    Explain why the lemma is called as it is with a picture, or observing that the vectors $\frac{x}{\|x\|}\left\langle \frac{x}{\|x\|}, y\right\rangle$ and $y-\frac{x}{\|x\|}\left\langle  \frac{x}{\|x\|}, y\right\rangle$ are orthogonal and their sum is $y$.
\end{exercise}\end{footnotesize}
 
\begin{corollary}\label{CorCS}
[Cauchy-Schwarz inequality]For $x, y \in V$ we have
$$
|\langle x , y\rangle| \leqslant\|x\|\|y\| \text {, }
$$
with equality if and only if $x, y$ are linearly dependent.    
\end{corollary}
\begin{proof}
If $x=0$, there is nothing to prove.  Otherwise, the inequality follows from dropping the first summand from the right of the Pythagorean relation. In case of equality, either $x=0$ or the first term in the right of the Pythagorean relation vanishes, in which case $x, y$ are linearly dependent. Conversely, if they are linearly dependent it is easy to see that equality holds in Cauchy-Schwarz. 
\end{proof}
\begin{proposition}\label{PropHilbertNorm}
    The function $x\mapsto\|x\|$ defines a norm on $V$.
\end{proposition}
\begin{proof} It is an immediate consequence of Cauchy-Schwarz:
$$
\begin{aligned}
\|x+y\|^{2} &=\langle x+y , x+y\rangle \\
&=\|x\|^{2}+\langle x , y\rangle+\langle y , x\rangle+\|y\|^{2} \\
&=\|x\|^{2}+2 \operatorname{Re}(\langle x , y\rangle)+\|y\|^{2} \\
& \leqslant\|x\|^{2}+2\|x\|\|y\|+\|y\|^{2} \\
&=(\|x\|+\|y\|)^{2} .
\end{aligned}
$$    
\end{proof}
An inner product space $(H,\langle\cdot , \cdot\rangle)$ is a {\it Hilbert space\index{Hilbert space}} if it is complete with respect to the norm induced by the inner product.
\begin{footnotesize}\begin{exercise}\label{ExeHThree}
Show that the inner product can be written in terms of norms by means of the {\bf polarization identity\index{polarization identity}}:
\begin{equation}\label{eqPolarization}
\langle x , y\rangle=\frac{1}{4}\left[\left(\|x+y\|^{2}-\|x-y\|^{2}\right)-i\left(\|x+i y\|^{2}+\|x-i y\|^{2}\right)\right] .
\end{equation}
\end{exercise}\end{footnotesize}
\begin{lemma}\label{LemmaParallelogram}
    [Parallelogram law] Let $x, y \in V$. Then,
$$
\|x+y\|^{2}+\|x-y\|^{2}=2\left(\|x\|^{2}+\|y\|^{2}\right) .
$$
\end{lemma}
\begin{proof}
    Expand the expression on the left using the definition of $\|\cdot\|$.
\end{proof}
\begin{footnotesize}\begin{exercise}\label{ExeHFour}
Show that, if $\|\cdot\|$ is a norm on a vector space which satisfies the parallelogram law, then it comes from an inner product as defined by the polarization identity. 
\end{exercise}\end{footnotesize}
\begin{footnotesize}\begin{exercise}\label{ExeHFive}
Show that the (two-dimensional) Banach spaces $\ell^{\infty}(\{0,1\}), \ell^{1}(\{0,1\})$ are not inner product spaces (i.e. that their norm does not come from an inner product). A {\bf Hint.} The parallelogram law fails.
\end{exercise}\end{footnotesize}
\subsection{$L^2$ as a Hilbert space} Let $(X,\mathcal{F},\mu)$ be a measure space. Observe that, if $f,g\in L^2(\mu)$, then $\overline{f}g\in L^1(\mu)$:
\[
\int_x|\overline{f}g|d\mu\le\|f\|_{L^2}\|g\|_{L^2}
\]
by H\"older's inequality with $p=p'=2$. We can then define
\[
\langle f,g\rangle_{l^2}=\int_X\overline{f}gd\mu,
\]
which satisfies the properties of a inner product, with associated norm $\langle f,f\rangle_{l^2}=\|f\|_{L^2}^2$. We saw that $L^2(\mu)$ is complete with respect to this norm.
\subsection{Projections onto subspaces}\label{SSectHilbertProjectionBasic}
Let $M \subseteq H$ be a closed, linear subspace of $H$. Its {\it orthogonal complement}, $M^{\perp}=H \ominus M$, is:
\[
M^{\perp}=\{x \in H: x \perp y \text { for all } y \in M\}.
\]
\begin{lemma}[Projection Lemma\index{projection!lemma}]\label{LemmaProjection}
 Let $M \subseteq H$ be a closed, linear subspace of $H$. Then, for each $x \in H$ there exists a unique $u \in H$ such that
$$
\|x-u\|=\inf \{\|x-y\|: y \in M\} .
$$    
\end{lemma}
\begin{proof}
Let $\left\{y_{n}\right\}_{n \geqslant 1}$ be a sequence in $M$ such that $\left\|y_{n}-x\right\| \rightarrow \delta:=\inf \{\|x-y\|: y \in M\}$. We show that it is a Cauchy sequence- By the parallelogram law,
$$
\begin{aligned}
\left\|y_{m}-y_{n}\right\|^{2}=&\left\|\left(y_{m}-x\right)-\left(y_{n}-x\right)\right\|^{2} \\
=& 2\left(\left\|y_{m}-x\right\|^{2}+\left\|y_{n}-x\right\|^{2}\right)-\left\|y_{m}+y_{n}-2 x\right\|^{2} \\
=& 2\left(\left\|y_{m}-x\right\|^{2}+\left\|y_{n}-x\right\|^{2}\right)-4\left\|\frac{y_{m}+y_{n}}{2}-x\right\|^{2} \\
\leqslant & 2\left(\left\|y_{m}-x\right\|^{2}+\left\|y_{n}-x\right\|^{2}\right)-4 \delta^{2} \\
& \text { because } \frac{y_{m}+y_{n}}{2} \in M \\
\rightarrow & 4 \delta^{2}-4 \delta^{2}=0 .
\end{aligned}
$$
We now let $u=\lim y_{n} \in M$ (because $M$ is closed).

If $u^{\prime}$ is another point with the same property, then $\frac{u+u^{\prime}}{2} \in M$ and
$$
\left\|u-u^{\prime}\right\|^{2}=2\left(\|u-x\|^{2}+\left\|u^{\prime}-x\right\|^{2}\right)-4\left\|\frac{u+u^{\prime}}{2}-x\right\| \leqslant 4 \delta^{2}-4 \delta^{2}=0,
$$
hence, $u=u^{\prime}$.
\end{proof}
We call $u=\pi_{M}(x)$ the {\it orthogonal projection\index{projection}} of $x$ onto $M$.
\begin{proposition}\label{PropHSix}
Let $E \subseteq H$ and define $E^{\perp}=\{x \in H: x \perp y$ for all $y \in E\}$. Then,
\begin{enumerate}
    \item[(i)] $E^{\perp}$ is a closed, linear subspace of $H$, and $E\cap E^\perp=\{0\}$;
    \item[(ii)] if $M$ is a linear subspace of $H$, then $M^\perp=\overline{M}^\perp$;
    \item[(iii)] if $M$ is a linear subspace of $H$, then $\left(M^{\perp}\right)^{\perp}=\bar{M}$ is the closure of $M$ in $H$.
\end{enumerate}
\end{proposition}
\begin{proof}
(i) Let $x,y\in E^\perp$, $\alpha,\beta\in\mathbb{C}$, and let $z\in E$. Then,
\[
\langle z,\alpha x+\beta y\rangle=\alpha\langle z,x\rangle+\beta\langle z,y\rangle=0,
\]
hence, $\alpha x+\beta y\in E^\perp.$ 

Let now $E^\perp\ni x_n$ and $\lim_{n\to\infty}\|x_n-x\|=0$. Then, for $z\in E$,
\[
\left|\langle z,x\rangle\right|=\left|\langle z,x\rangle-\langle z,x_n\rangle\right|\le \|z\|\cdot\|x-x_n\|\to0,
\]
as $n\to\infty$. So, $\left|\langle z,x\rangle\right|$ is smaller than any positive $\epsilon$, hence $\langle z,x\rangle=0$, i.e. $x\in E^\perp$.

Finally, if $x\in E\cap E^\perp$, then $0=\langle v,v\rangle=\,v\,^2$, hence, $v=0$.

(iii) $\left(M^{\perp}\right)^{\perp}$ is linear and closed by (i), and $\left(M^{\perp}\right)^{\perp}\supseteq M$: if $x\in M$, then, for $z\in M^\perp$ we have that $\langle z,x\rangle=0$, then $x\in \left(M^{\perp}\right)^{\perp}$. We have thus proved that $\overline{M}$, the smallest closed linear subspace containing $M$, is a subset of $\left(M^{\perp}\right)^{\perp}$.

Consider $x\in \left(M^{\perp}\right)^{\perp}\setminus\overline{M}$. By the Projection Lemma, $x=u+v$ with $u\in \overline{M}$, and $v\in\overline{M}^\perp=M^\perp$. So, $v\in M^\perp\cap \left(M^{\perp}\right)^{\perp}$, hence, $v=0$ and $x\in\overline{M}$, contradicting the assumption.
\end{proof}
\begin{footnotesize}\begin{exercise}\label{exeFromMtoE}
    Let $E$ be a subset of $H$. Show that $(E^\perp)^\perp=\overline{\text{span}(E)}$ is the closure of the linear span of $E$ in $H$.
\end{exercise}\end{footnotesize}
\begin{theorem}\label{TheoOrthDec}
    [Orthogonal Decomposition\index{theorem!orthogonal decomposition}] Let $M \subseteq H$ be a closed, linear subspace of $H$. Then, for each $x \in H$ there is a unique decomposition $x=u+v$ with $u \in M$ and $v \in M^{\perp}$.
\end{theorem}
\begin{proof}
Let $u\in M$ be as in the Projection Lemma, and for $t \in \mathbb{R}$ and $w\in M$, so that $u+tw\in M$, define $f(t)=\,u+t w-x\|^{2}$. By minimality, $f^{\prime}(0)=0$, and
$$
\begin{aligned}
f^{\prime}(0) &=\left.\frac{d}{d t}\right|_{t=0}\left(\|u-x+tw\|^{2}\right) \\
&=\left.\frac{d}{d t}\right|_{t=0}\left(\|u-x\|^{2}+2 t \operatorname{Re}(\langle u-x , w\rangle)+t^{2}\|w\|^{2}\right) \\
&=\left[2 \operatorname{Re}(\langle u-x , w\rangle)+2 t\|w\|^{2}\right]_{t=0} \\
&=2 \operatorname{Re}(\langle u-x , w\rangle) .
\end{aligned}
$$
Hence, $\operatorname{Re}(\langle u-x , w\rangle)=0$. The same reasoning with $g(t)=\|u+i t w-x\|^{2}$ gives $\operatorname{Im}(\langle u-x , w\rangle)=0$. Hence, $\langle u-x , w\rangle$ for all $w\in M$, i.e. $v:=x-u\in M^\perp$.

If I have two decompositions $x=u+v=u'+v'$ with $u,u'\in M$ and $v,v'\in M^\perp$, then $u-u'=v'-v\in M\cap M^\perp$, hence 
$u-u'=v'-v=0$.
\end{proof}
In the theorem on the orthogonal decomposition of a vector, we saw that, if $M$ is a closed subspace of a Hilbert space $H$, and $x\in H$, then there exist unique $u\in M$ and $v\in M^\perp$ such that $x=u+v$. Define $\pi_M:H\to H$ by $\pi_M(x)=u$, the {\it (orthogonal) projection} of $x$ onto $M$.
\begin{footnotesize}\begin{exercise}\label{ExeHSeven}
    Let $M$ be a closed, linear subspace of a Hilbert space $H$.
    \begin{enumerate}
        \item[(i)] Show that $\pi_{M}: H \rightarrow H$ is linear and that $\left\|\pi_{M}x\right\|\le\|x\|$ for all $x$ in $H$.
        \item[(ii)] Show that $\pi_{M}^{2}=\pi_{M}$ and that $\pi_{M}$ is {\bf self-adjoint}, i.e. that
$$
\left\langle\pi_{M}x , y\right\rangle=\left\langle x , \pi_{M}y\right\rangle
$$
for $x, y \in H$.
       \item[(iii)] Verify that $\pi_M+\pi_{M^\perp}=I$ is the identity operator ($I(x)=x$), and that $\pi_M\circ\pi_{M^\perp}=0$. 
    \end{enumerate}
\end{exercise}\end{footnotesize}
It is an interesting fact that the Exercise \ref{ExeHSeven} has a converse: all operators sharing the properties of projections are in fact projections on closed subspaces. We can, that is, indifferently work with the {\it lattice} of the closed subspaces of $H$, or with the lattice of the projection operators. Since projections operators are spacial instances of linear operators on $H$, this identification gives much more flexibility in theory and calculations. (This is similar to the case of $\sigma$-algebras, where characteristic functions of measurable sets are special instances of measurable functions, and measurable functions have a rich structure).
\begin{proposition}\label{PropHSixInverse} Let $\pi:H\to H$ be a linear operator satisfying:
\begin{enumerate}
\item[(i)] $\pi^2=\pi$;  
\item[(ii)] $\pi^\ast=\pi$;  
\item[(iii)] $\|\pi(x)\|\le C\|x\|$ for some $C>0$.  
\end{enumerate}
    Then, there exists a closed, linear subspace $M$ of $H$ such that $\pi=\pi_M$.
\end{proposition}
\begin{proof}
   Let $M=\{\pi(x):\ x\in H\}$ be the range of $\pi$. For $\pi(x),\pi(y)\in M$ and $\alpha,\beta\in\mathbb{C}$ we have $\alpha\pi(x)+\beta\pi(y)=\pi(\alpha x+\beta y)\in M$, hence $M$ is a linear space. Let $M\ni\pi(x_n)\to y$ in $H$. Then, $\pi(x_n)=\pi^2(x_n)\to\pi(y)$ because
   \[
   \|\pi^2(x_n)-\pi(y)\|\le C\|\pi(x_n)-y\|\to0\text{ as }n\to\infty.
   \]
   Thus, $y=\pi(y)\in M$, showing that $M$ is closed.

   Finally, let $x$ be an element of $H$, and write $x=u+v=\pi_M(x)+\pi_{M^\perp}(x)$, as in Theorem \ref{TheoOrthDec}. By definition of $M$, $u=\pi(w)$ for some $w$. For $y\in H$ we have:
   \[
   \langle \pi(v),y\rangle=\langle v,\pi(y)\rangle=0, 
   \]
   since $\pi(y)\in M$ and $v\in M^\perp$. Hence, $\pi(v)=0$, and so
   \[
   \pi(x)=\pi(u)=\pi^2(w)=\pi(w)=u=\pi_M(x),
   \]
   showing $\pi=\pi_M$.
\end{proof}
\subsection{F. Riesz representation in Hilbert spaces}\label{SSectRieszRepresentationHilbert}
F. Riesz representation theorem in Hilbert spaces, like the one we saw concerning measures, shows that any element from an abstract collection of objects can be represented as an object from a very special and concrete subclass. The first theorem we saw represented positive functionals on $C_c(X)$ in terms of measures; the present one represents bounded, linear functionals on $H$ in terms of inner products.
\begin{footnotesize}\begin{exercise}\label{ExeHEight}
    Let $y \in H$ and define $T_{y}(x)=\langle y , x\rangle$. Then, $T_{y}: H \rightarrow \mathbb{C}$ is a linear operator (a linear functional) and $\left\|T_{y}\right\|=\|y\|$. 
\end{exercise}\end{footnotesize}

\begin{theorem}
    [F. Riesz Lemma\index{theorem!Riez representation in Hilbert spaces}] Let $T: H \rightarrow \mathbb{C}$ be a bounded, linear functional. Then, there is a unique $y_{T}$ in $H$ such that $T(x)=\left\langle y_T , x\right\rangle$.
\end{theorem}
\begin{proof}
Let $N=\operatorname{Ker} T$. If $N=H$, set $y_{T}=0$. If not, by the Orthogonal Decomposition Theorem there is some $x_{0} \in N^{\perp} \backslash\{0\}$, and we can assume that $\left\|x_{0}\right\|=1$. We write
$$
x=\left(x-\frac{T(x)}{T\left(x_{0}\right)} x_{0}\right)+\frac{T(x)}{T\left(x_{0}\right)} x_{0}=u+\lambda x_{0}.
$$
The second summand is in $\operatorname{span}\left(\left\{x_{0}\right\}\right) \subseteq N^{\perp}$, with $\lambda=\left\langle x_0 , x\right\rangle$, while $u$ clearly lies in $N$ (incidentally, this shows that $\operatorname{Ker}(N)$ has codimension one in $H)$.
$$
T(x)=\lambda T\left(x_{0}\right)=\left\langle x_0 , x\right\rangle T\left(x_{0}\right)=\left\langle  \overline{T\left(x_{0}\right)} x_{0}, x\right\rangle,
$$
hence, $y_{T}=\overline{T\left(x_{0}\right)} x_{0}$.
\end{proof}

\begin{footnotesize}\begin{exercise}\label{ExeHNine}
Riesz Lemma provides an identification, in fact an isometry, $C: T \mapsto y_{T}$ of the dual space $H^{*}=\{T$ : $H \rightarrow \mathbb{C}$ such that $\|T\|<\infty\}$ with $H$ itself. The map $C$ is conjugate linear,
$$
C(a S+b T)=\bar{a} C(S)+\bar{b} C(T) .
$$
Prove this and prove that, if $A: H \rightarrow K$ is a bounded, linear map between Hilbert spaces, then $C(T A)=A * C(T)$, where the {\it (Hilbert space) adjoint operator $A^\ast: K \rightarrow H$} is defined by $\langle A x , y\rangle=\left\langle x , A^\ast y\right\rangle$.   \index{operator!adjoint} 
\end{exercise}\end{footnotesize}

\section{Orthonormal systems}\label{SectOrthogonal}
\subsection{Orthogonal vectors}\label{SSectOrthigonal}
A family $S=\left\{f_{a}: a \in I\right\}$ of vectors is an {\it orthogonal system\index{orthogonal system}} if any two vectors in it are orthogonal, and an {\it orthonormal system\index{orthonormal system}} if, in addition, each $f_{a}$ has unit norm. An {\it orthonormal basis (o.n.b.)\index{orthonormal basis}} for $H$ is an orthonormal system which is maximal: no other vector can be added to it without breaking the orthonormality condition.

\begin{footnotesize}\begin{exercise}\label{ExeHEleven}. Show that the orthonormal system $\left\{e_{a}: a \in I\right\}$ is a o.n.b. for $H$ is and only if span $\left\{e_{a}: a \in I\right\}$ is dense in $H$.
\end{exercise}\end{footnotesize}

\begin{theorem}\label{eqBessel}[Bessel Inequality\index{theorem!Bessel inequality}] Let $\left\{e_{i}: i=1, \ldots, n\right\}$ be an orthonormal system in the Hilbert space $H$, and let $x\in H$. Then, the vector $x-\sum_{i=1}^{n}\left\langle x , e_{i}\right\rangle e_{i}$ is orthogonal to $\operatorname{span}\left\{e_{1}, \ldots, e_{n}\right\}$, and
$$
\|x\|^{2} \geqslant \sum_{i=1}^{n}\left|\left\langle  e_{i}, x\right\rangle\right|^{2} .
$$
If $\left\{e_{a}: a \in I\right\}$ is an o.n.b., the numbers $\left\langle  e_{a}, x\right\rangle=\hat{x}(a)$ are the Fourier coefficients of $x$ w.r.t. the basis. 
\end{theorem}
\begin{proof}
The first assertion is clear: 
\[
\left\langle x-\sum_{i=1}^{n}\left\langle x , e_{i}\right\rangle e_{i} , e_j \right\rangle=\langle x, e_j\rangle-\sum_{i=1}^n \left\langle x , e_{i}\right\rangle \left\langle e_i , e_j\right\rangle=0,
\]
by orthonormality of the system. As a consequence,
$$
\begin{aligned}
\|x\|^{2} &=\left\|x-\sum_{i=1}^{n}\left\langle  e_{i}, x\right\rangle e_{i}\right\|^{2}+\left\|\sum_{i=1}^{n}\left\langle  e_{i}, x\right\rangle e_{i}\right\|^{2} \\
&=\left\|x-\sum_{i=1}^{n}\left\langle  e_{i}, x\right\rangle e_{i}\right\|^{2}+\sum_{i=1}^{n}\left|\left\langle e_{i}, x\right\rangle\right|^{2}
\end{aligned}
$$
Before proceeding, we need to clarify some concepts about infinite sums. If $\left\{c_{a}: a \in I\right\}$ is a family of positive numbers, then
$$
\sum_{a \in I} c_{a}:=\sup \left\{\sum_{i=1}^{n} c_{a_{i}}:\left\{a_{1}, \ldots, a_{n}\right\} \subseteq I\right\} .
$$
\end{proof}
\begin{corollary}\label{corProjectionFinite}
    Let $M=\text{span}\{e_1,\dots,e_n\}$. Then,
    \[
    \pi_M(x)=\sum_{i=1}^{n}\left\langle x , e_{i}\right\rangle e_{i}.
    \]
\end{corollary}
\begin{proof}
    The first assertion in Bessel's inequality says that $x=u+v$ with $u=\sum_{i=1}^{n}\left\langle x , e_{i}\right\rangle e_{i}\in M$ and $v\in M^\perp$, hence, $u=\pi_M(x)$.
\end{proof}
\begin{footnotesize}\begin{exercise}\label{ExeHTen}
     Show that, if $\sum_{a \in I} c_{a}<\infty$, then $\left\{a \in I: c_{a} \neq 0\right\}$ is at most countable. If $\left\{c_{a}: a \in I\right\}$ is a family of complex numbers, we say that $\sum_{a \in I} c_{a}$ converges absolutely if $\sum_{a \in I}\left|c_{a}\right|<\infty$ converges. In this case we say that $\sum_{a \in I} c_{a}$ converges in $\mathbb{C}$ absolutely, hence irrespective of how $I$ is ordered (prove it if you have never done it before!).
\end{exercise}\end{footnotesize}
\subsection{Fourier analysis and synthesis of a vector}\label{SSectSpectralAnSynth}\index{Fourier analysis!in Hilbert spaces}\index{Fourier synthesis!in Hilbert spaces}
\begin{theorem}\label{theoSpectralAnSynth}
    Let $\{e_a\}_{a\in I}$ be a o.n.b. of $H$. Then,
    \begin{enumerate}
        \item[(i)] For all $x$ in $H$,
        \begin{equation}\label{eqSpectAn}
            x=\sum_{a\in I}\langle e_a, x\rangle e_a \text{ (Analysis of $x$)}
        \end{equation}
        converges in $H$;
        \item[(ii)] we have\index{Plancherel formula!abstract}
        \begin{equation}\label{eqPlancherelAbstract}
            \|x\|^2=\sum_{a\in I}|\langle e_a, x\rangle|^2 \text{ (Plancherel Isometry)};
        \end{equation}
        \item[(iii)] if $\{c_a\}_{a\in I}$ is a sequence in $H$ such that $\sum_{a\in I}|c_a|^2$ converges, then
        \begin{equation}\label{eqSpectSynthesis}
            \sum_{a\in I}c_ae_a\text{ converges in }H. \text{ (Analysis of $x$)}
        \end{equation}
    \end{enumerate}
\end{theorem}     
\begin{proof}
    By Bessel inequality, if $A\subseteq I$ with $\sharp(A)<\infty$, then $\sum_{a\in A}|\langle e_a, x\rangle|^2\le \|x\|^2$. As a consequence, $I(x):=\{a\in I:\ \langle e_a, x\rangle\}=\{a_n\}_{n}$ is countable, hence, \eqref{eqPlancherelAbstract} holds with $\le$. Set now $x_n=\sum_{j=1}^n\langle e_{a_j}, x\rangle e_{a_j}$. We show that is defines a Cauchy sequence in $H$:
    \[
        \|x_{n+l}-x_n\|^2=\sum_{j=1}^l|\langle e_{a_j}, x\rangle|^2\to0
    \]
    as $n\to\infty$, by the convergence of the series on the right of \eqref{eqPlancherelAbstract}. Let $y=\lim_{n\to\infty}x_n$ (the limit is taken in $H$-norm):
    \begin{eqnarray*}
        \langle x-y, e_{a_m}\rangle&=&\left\langle x-\lim_{n\to\infty}\sum_{j=1}^n\langle e_{a_j}, x\rangle e_{a_j}  ,  e_{a_m}\right\rangle\crcr 
        &=&\lim_{n\to\infty}\left\langle x-\sum_{j=1}^n\langle e_{a_j}, x\rangle e_{a_j}  ,  e_{a_m}\right\rangle\crcr 
        &=&\lim_{n\to\infty}\left(\langle e_{a_m}, x\rangle-\langle e_{a_m}, x\rangle\right)=0.
    \end{eqnarray*}
    The same argument shows that $\langle x-y, e_{a}\rangle=0$ if $a\notin I(x)$. Hence, $x-y\perp e_a$ for all elements of the orthonormal basis, hence $x=y$, which shows \eqref{eqSpectAn}.

    We finish the proof of Plancherel's identity:
    \begin{eqnarray*}
        0&=&\lim_{n\to\infty}\left\|x-\sum_{j=1}^n\langle e_{a_j}, x\rangle e_{a_j}\right\|^2\crcr 
        &=&\left(\|x\|^2-\sum_{j=1}^n|\langle e_{a_j}, x\rangle|^2\right)\crcr 
        &=&\|x\|^2-\sum_{n=1}^\infty|\langle e_{a_n}, x\rangle|^2\crcr 
        &=&\|x\|^2-\sum_{a\in I}|\langle e_{a}, x\rangle|^2,
    \end{eqnarray*}
    as wished.
\end{proof}
\subsection{Orthonormal basis in separable Hilbert spaces}\label{SSectGramSchmidt}
\subsubsection{Gram-Schmidt algorithm}\label{SSSectGramSchmidt}\index{Gram-Schmidt algorithm}
A Hilbert space is {\it separable} if it separable as a metric space, i.e. if it has a countable, dense set. Most Hilbert spaces encountered in theory and applications are separable, and in this case the existence of a orthonormal basis is constructive. We will then treat first the separable case, then we will see how things work in the general case, where Zorn's Lemma is required, and we have no way to "see" the basis.
\begin{theorem}[Gram-Schmidt]\label{theoGramSchmidt}
    Let $H$ be a separable Hilbert space. Then, $H$ has a countable basis.
\end{theorem}
\begin{proof} We start with a countable, dense subset $G=\{g_n:\ n\ge1\}$ of $H$. Inductively removing vectors linearly dependent from the previously chosen ones, we obtain a linearly independent subfamily $F=\{f_n:\ n\ge1\}$ such that ${\text span}(F)={\text span}(G)$ is still dense in $H$. We then transform $F$ into an orthonormal system:
\begin{eqnarray*}
e_1&=&\frac{f_1}{\|f_1\|},\crcr
e_2&=&\frac{f_2-\langle  e_1, f_2\rangle e_1}{\|f_2-\langle  e_1, f_2\rangle e_1\|},\crcr
&\dots&\crcr
e_n&=&\frac{f_n-\sum_{j=1}^{n-1}\langle  e_j, f_n\rangle e_j}{\|f_n-\sum_{j=1}^{n-1}\langle  e_j, f_n\rangle e_j\|},\crcr
&\dots&
\end{eqnarray*}
The denominators do not vanish because $F$ is a linearly independent family, and inductively we see that
${\text span}\{e_1,\dots,e_n\}={\text span}\{f_1,\dots,f_n\}$.
Hence, $\{e_n:\ n\ge1\}$ is a orthonormal basis for $H$.
    
\end{proof}

\begin{footnotesize}\begin{exercise}\label{ExeHTwelve}
    Let $H=L^{2}[0,1]$ with the Lebesgue measure, and consider the functions $1, x, x^{2}$. Apply Gram-Schmidt's algorithm to find three orthonormal vectors in $H$.
\end{exercise}\end{footnotesize}

\subsubsection{The classification of separable Hilbert spaces}\label{SSSclassifyHilbertSeparable}
\begin{theorem}\label{theoClassifyHilbertSeparable}
    A Hilbert space $H$ is separable if and only if it has one countable orthonormal basis $\mathcal{B}$. If it does, all o.n.b. of $H$ have the same cardinality. Moreover,
    \begin{enumerate}
        \item[(i)] if $\sharp(\mathcal{B})=d<\infty$, then $H$ is isometrically isomorphic to $\mathbb{C}^d$;
        \item[(ii)]  if $\sharp(\mathcal{B})=\infty$, then $H$ is isometrically isomorphic to $\ell^2(\mathbb{N})$.
    \end{enumerate}
\end{theorem}
\begin{proof}
    If $H$ is separable, then Theorem \ref{theoGramSchmidt} provides a countable o.n.b. $\mathcal{B}$. Suppose viceversa that $H$ has a countable basis $\mathcal{B}$. If $\sharp(\mathcal{B})=d<\infty$, $\mathcal{B}=\{e_1,\dots,e_d\}$, then each $x$ in $H$ can be written as $x=\sum_{j=1}^d\langle e_j, x\rangle e_j$ by the Analysis part of Theorem \ref{theoSpectralAnSynth} (or by a much more elementary argument). The map
    \[
    L_\mathcal{B}:x\mapsto (\langle e_j, x\rangle)_{j=1}^d
    \]
    is an isometric isomorphism of $H$ onto $\mathbb{C}^d$. Since $H$ is finite dimensional, all its basis have the same dimension.
    In the countable case, $\mathcal{B}=\{e_n\}_{n=1}^\infty$, and we define $L_\mathcal{B}:H\to \ell^2(\mathbb{N})$ in the same way, $L_\mathcal{B}:x\mapsto \left\{\langle e_n, x\rangle\right\}_{n=1}^\infty$. Theorem \ref{theoSpectralAnSynth} implies that $L_\mathcal{B}$ is a isometric isomorphism.

    In the opposite direction, suppose $H$ has a countable o.n.b., and consider the countable $Q$ family of the vectors $q=\sum_{j=1}^n q_j e_j$, where $n\ge1$, the $q_j$'s are complex rationals, $q_j\in\mathbb{Q}+i\mathbb{Q}$. If $x=\sum_{j=1}^\infty c_j e_j$ and $\epsilon>0$, then
    \[
        \left\|\sum_{j=1}^\infty c_j e_j-\sum_{j=1}^n q_j e_j\right\|^2=\sum_{j=n+1}^\infty|c_j|^2+\sum_{j=1}^n|c_j-q_j|^2.
    \]
    Choose first $n$ such that the first sum is dominated by $\epsilon$, then $q_1,\dots,q_n$ such that each summand in the second sum is dominated by $\epsilon/n$. This shows separability.

    We are left with proving that, if one o.n.b. $\mathcal{B}$ is infinite countable, then any other o.n.b. $\mathcal{B}_1$ is. Let $\mathcal{B}_1=\{e_a\}_{a\in I}$, and observe that, for any $e_a\ne e_b$ in it, $\|e_a-e_b\|^2=2$. Let $\{q_n\}_{n=1}^\infty$ be any countable, dense set in $H$. Then, for each $a\in I$ there is $n=n(a)$ such that $\|e_a-q_n\|<\sqrt{2}/2$, and if $a\ne b$, then $n(a)\ne n(b)$. This provides an injective map $I\to\mathbb{N}$, hence $I$ is countable.
\end{proof}

\begin{footnotesize}
    
\subsection{Intermezzo: orthonormal basis in general Hilbert spaces}\label{SSectZornOrth} 

\subsubsection{Existence of o.n.b.}\label{SSSectONBveryinfinite}
A widely used version of the Axiom of Choice is Zorn's Lemma, which we are stating below. We work with a partially ordered set $(A,\le)$. A {\it chain} in $A$ is a subset $B\subseteq A$ on which the partial order is in fact a total order: for any $x,y$ in $B$, either $x\le y$ or $y\le x$. An element $m$ in $A$ is {\it maximal} for $C\subseteq A$ if $x\le m$ for all $x$ in $C$ and if $n\le m$ has the same property, then $n=m$.
\begin{theorem}[Zorn Lemma]\label{TheoZorn}
Let $(A,\le)$ be a (nonempty) partially ordered set with the property that any chain $B$ in $A$ has a maximal element in $A$. Then, $A$ has a maximal element.
\end{theorem}
As a consequence, we have that an orthornormal basis exists for any Hilbert space.
\begin{theorem}\label{TheoOnbZorn}
Every Hilbert space $H\ne\{0\}$ has a orthonormal basis.
\end{theorem}
\begin{proof}
Let $A$ be the set of all orthonormal systems of $H$, ordered by inclusion. It is nonempty because $\{x/\|x\|\}\in A$ if $x\ne 0$ is an element of $H$. If $B$ is a chain in $A$, then
\[
{\mathcal G}=\cup_{{\mathcal F}\in B}{\mathcal F}\in A
\]
is a maximal element for $B$: it is an orthonormal system, and all orthonormal system containing all ${\mathcal F}$ in $B$ contain ${\mathcal G}$.

By Zorn's Lemma, $A$ has a maximal element ${\mathcal H}$ and, unravelling definitions, ${\mathcal H}$ is a maximal orthonormal system in $H$, hence a basis.
\end{proof}

\subsubsection{The dimension of a Hilbert space}\label{SSSectHilbertDimension}\index{Hilbert space!dimension}
By Theorem \ref{theoSpectralAnSynth}, if $H$ be a Hilbert space and $\mathcal{B}=\{e_a:a\in I\}$ an orthonormal basis for it, then the map
\[
L_\mathcal{B}:x=\sum_{a\in I}\langle  e_a, x\rangle e_a\mapsto\{\langle  e_a, x\rangle\}_{a\in I}
\]
is an isometric isomorphism from $H$ onto $\ell^2(I)$.

A natural question is whether different orthonormal basis of $H$ have the same cardinality, as we verified the separable case. The answer is positive.
\begin{theorem}\label{TheoEllTwoBasis}
Let $\{e_a:\ a\in I\}$ and $\{f_b:\ b\in J\}$ be o.n.b. of the same Hilbert space $H$. Then, $\sharp(I)=\sharp(J)$.
\end{theorem}
The cardinality in question is the {\it dimension} of the Hilbert space.
\begin{corollary}\label{CorEllTwoBasis}
All separable Hilbert spaces which are not finite dimensional have countable dimension, and are isomorphic to each other.
\end{corollary}
\begin{proof}
If both $I$ and $J$ are finite dimensional, the statement is a well know fact from linear algebra (go back to see its proof!).

If $I$ is finite dimensional then $J$ is finite dimensional as well. For, suppose $I$ is finite dimensional and let $S=\{x\in H;\ \|x\|=1\}$ be the unit sphere in $H$. Being $H$ a finite dimensional Euclidean space, and $S$ both closed and bounded, $S$ is compact. Suppose now $J$ is infinite and let $\{f_{b_n}:\ n\ge1\}$ an infinite countable subset of $\{f_b:\ b\in J\}$, which is contained in $S$ if a subsequence of it (which we might denote the same way) converges to some $x$ in $S$, by continuity of the inner product $\langle x, f_{b_n}\rangle\to \|x\|^2=1$ as $n\to\infty$. On the other hand, $x=\sum_{b\in J}\langle x, f_b\rangle f_b$, hence $\langle x, f_{b_n}\rangle\to 0$ as $n\to\infty$, because, for instance, $\sum_{n}|\langle x, f_{b_n}\rangle|^2\le \|x\|^2=1$.

The remaining case is that in which both $I$ and $J$ are infinite. To each $a\in I$ there corresponds an at most countable subset $J(a)$ of $J$ such that
\[
e_a=\sum_{b\in I(a)}\langle e_a, f_b\rangle f_b.
\]
The number of the involved $f_b$'s is
\[
\sharp\left(\cup_{a\in I}I(a)\right)\le\sharp(I).
\]
Suppose a basis element $f_{b_0}$ was not used. Then, $f_{b_0}$ is orthogonal to all $e_a$'s, hence to $H$: contradiction. 

We have then that $\sharp(J)\le\sharp(I)$, and the opposite inequality holds as well.
\end{proof}

%\begin{footnotesize}\begin{exercise}
%Let $H$ be a Hilbert space with infinite, countable dimension. Show that, as a vector space, it has infinite, uncountable dimension.
%\end{exercise}\end{footnotesize}
\end{footnotesize}

\section{The trigonometric system and Fourier series}\label{SSectTrigONB}\index{trigonometric system!completeness}
In item (IX) of section \ref{SSSectFourierseries} we saw that the trigonometric system $\left\{e_n\right\}_{n=-\infty}^\infty$, where
\begin{equation}\label{eqTrigSys}
    e_n(t)=\frac{e^{int}}{2\pi}
\end{equation}
is a complete orthonormal system (an orthonormal basis) for $L^2(\mathbb{T})$. 

It should be kept in mind that
a great amount of work goes on all over the world into finding, studying the special properties of, and applying, new orthornormal basis (ONB) for old and new Hilbert spaces of functions.
Broadly speaking, when the functions $e_n$ in an ONB share some qualitative feature, so do their linear combinations; hence, using that ONB to approximate a function $f$ means finding closer and closer approximations of $f$ which "look like" the functions in the basis. To wit, "if all you have in your toolbox is a hammer, then everything you look at seems a hammer". The trigonometric system, however, remains the most natural one, because of its interaction with the gropu structure of the torus,
\begin{equation}\label{eqExpHom}
e^{in(s+t)}=e^{ins}e^{int}.
\end{equation}
This is the reason why the trigonometric system {\it diagonalizes translations}. Let $\tau_af(t)=f(t-a)\ (\text{mod\ }2\pi)$: 
\[
\tau_ae_n=e^{-ina}e_n;
\]
hence the derivative operator, which is a linear combination of "infinitely close" translations,
\[
i^{-1}\frac{d\ }{dt}e_n(t)=\lim_{h\to0}\frac{\tau_{-h}e_n(t)-\tau_0e_n(t)}{ih}=\lim_{h\to0}\frac{e^{inh}-1}{ih}e_n(t)=ne_n(t).
\]
In fact, the trigonometric system diagonalizes all "reasonable" linear operators on $L^2(\mathbb{T})$ which commute with translations. The main example is that of the convolution operators $f\mapsto f\ast g$.

Next, we translate into the special context of the trigonometric system the general Hilbert spaces results we have seen.
\begin{itemize}
    \item[(I)] {\bf Fourier analysis.}\index{Fourier analysis!trigonometric system} For $f\in L^2(\mathbb{T})$ and $N\ge1$, let
    \[
    S_Nf(t)=\frac{1}{2\pi}\sum_{|n|\le N}\widehat{f}(n)e^{int},
    \]
    Then.
    \[
    \lim_{N\to\infty}\|S_Nf-f\|_{L^2(\mathbb{T})}=0.
    \]
    i.e.
    \begin{equation}\label{eqFourier}
    f(t)=S_\infty f(t):=\frac{1}{2\pi}\sum_{n=-\infty}^\infty\widehat{f}(n)e^{int},
    \end{equation}
    where the convergence is in $L^2(\mathbb{T})$.
    \item[(II)] {\bf Fourier synthesis.}\index{Fourier synthesis!trigonometric system} Let $\{c_n\}_{n=-\infty}^\infty\in\ell^2(\mathbb{Z})$,
    \[
    \|\{c_n\}_{n=-\infty}^\infty\|_{\ell^2}^2=\sum_{n\in\mathbb{Z}}|c_n|^2<\infty.
    \]
    Then, there exists $f\in L^2(\mathbb{T})$ such that
    \[
    \widehat{f}(n)=c_n.
    \]
    Explicitely, 
    \begin{equation}\label{eqFSynth}
    f(t)=\frac{1}{2\pi}\sum_{n=-\infty}^\infty c_ne^{int},
    \end{equation}
    where the series converges in the $L^2$ norm.
    \item[(III)] {\bf Plancherel formula.}\index{Plancherel formula!trigonometric system}
    \[
    \|f\|_{L^2(\mathbb{T})}^2=\frac{1}{2\pi}\sum_{n=-\infty}^\infty|c_n|^2.
    \]
\end{itemize}
The foundational problem of Fourier theory is understanding for which $f$, and in which sense, the series in \eqref{eqFourier} converges. For $f\in L^2(\mathbb{T})$, we have considered a quantitative answer, using the fact that $\{e_n\}_{n=-\infty}^{+\infty}$ is an orthonormal basis,
and the general theory of Hilbert spaces.

Fourier series provided, in the second half of the XIX century, the motivation for some of the most important innovations in mathematical analysis and in the foundations of mathematics. In S.M. Srivastava \href{https://www.ias.ac.in/article/fulltext/reso/019/11/0977-0999}{How did Cantor Discover Set Theory and Topology?} [Resonance, Vol. 19, p. 977-999, (2014)], for instance, you can find the history of how Cantor came to develop the theory of infinite sets to solve open problems concerning the convergence Fourier series.

It is noteworthy the role of Lebesgue integration. For the {\it Fourier analysis} of a function, the definition of integral provided by Riemann was considered wholly satisfactory, as it enabled to compute the Fourier coefficients of all functions mathmaticians of the age could think of (it was precisely to compute Fourier series that Riemann invented the definition). The problem is with {\it Fourier synthesis}, since it was not obvious (and in fact it is not true!) that the function $f$ designed by \eqref{eqFSynth} is Riemann integrable. The problem is that, as it came to be gradually realized, integrability in the sense of Riemann is not preserved under many limiting procedures, including those which are basilar in Fourier theory. Borel, then Lebesgue came, and this foundational problem was solved.

\chapter{Banach spaces}\label{chaptBanach}
\label{chap:Banach_spaces}

In this chapter we introduce some basic tools in Banach theory, {\it Hahn-Banach extension theorem} and {\it Baire's category theorem and it consequences}. The latter are, on the one hand, {\it (Banach-Steinhaus) uniform boundedness principle}, and on the other the sequence: {\it open mapping theorem}, {\it inverse mapping theorem}, {\it closed graph theorem}. We will also consider some consequences and applications of these cornerstones of functional analysis.

\section{Zorn's lemma and some of its consequences}\label{SectZorno}\index{Zorn lemma}

In the proof of Hahn-Banach Theorem we will use Zorn's lemma, which is a consequence of the axiom of choice, and it is in fact equivalent to it. In most courses Zorn's lemma is just stated as a principle, then used alongside the axiom of choice. One reason is that most instructors hope that the relationship between the two statements has been clarified in some other course. The other, probably, is that most instructors are research mathematicians with little taste for nonconstructive principles: they know them and use them, but feel that such principles are a sort of magic of last resort, and having a magic tool generating others, or having a supply of independent ones, does not really make to them much of a difference. A problem with highly nonconstructive existence theorems, in fact, is that we do not know much of the object we have proved the existence of. It is surprising, in view of this fact, that such theorems can be used to obtain more practical statements.

There is a price tag attached to it. Nonconstructive arguments allow us to prove the existence of counterexamples to very reasonable guesses. If we allow the logical basis of the nonconstructive arguments, we have to keep in mind that such monsters exist and have to be taken into account.

So much for the small talk. Before we move to Banach theory, in order to satisfy readers with a taste for logics we state Zorn's lemma and prove that it follows from the axiom of choice (but not the opposite implication).

\begin{theorem}[Zorn's lemma]\label{theoZorn}
    Let $(X,\le)$ be a partial order on a set $X$, and suppose that each totally ordered subset $C$ of $X$ has an upper bound: there is $u$ in $X$ such that $u\ge c$ for all $c\in C$. Then, $X$ has a maximal element $w$: $w\ge x$ for all $X$ in $X$.
\end{theorem}
\begin{proof} While reading the proof, it is very helpful drawing pictures. A partial order is a special instance of a directed graph, and a totally ordered subset of it might be thought of, in pictorial terms, as a subset of the real line (although it is clearly false that all totally ordered sets arise this way, for instance because there is no a priori restriction on their cardinality).

    Suppose by contradiction that $X$ has no maximal element. In particular, by the axiom of choice, for each totally ordered subset $C$, including $C=\emptyset$, we can choose an upper bound $u(C)$ in $X\setminus C$ (if the only upper bound, which exists by hypothesis, belonged to $C$, then it would be maximal for $X$). Now, to each subset $F$ of $X$ and $f\in F$ we associate its {\it tail set} $F_{<f}<=\{x\in F:x<f\}$.

    A subset $A$ is a $u$-set if
    \begin{itemize}
        \item[(i)] $A$ is totally ordered;
        \item[(ii)] $A$ contains no infinite descending sequence $a_1>a_2>\dots$ (e.g. $A$ can have the order type\footnote{By {\it order type} of a partially ordered set $(Z,<)$ we mean the equivalence class of all linearly ordered sets which are in an order-preserving bijection with $Z$.} of $\mathbb{N}$, but not that of $-\mathbb{N}$), hence, any subset $A'$ of $A$ has minumum;
        \item[(iii)] if $a\in A$, then $u(A_{<a})=a$.
    \end{itemize}
    The last property says that $u$-sets are linked in a very special way to the choice function $u$, and, as we will see below, they share many properties with the ordinals. 

    If $A$ is a $u$-set and $a\in A$, then $A_{<a}$ has properties (i-iii), hence, it is a $u$-set.
    
    If $A$ and $B$ are $u$-sets, by properties (i) and (ii) they have minimum, and $\{\min(A)\},\{\min(B)\}$ are $u$-sets. By property (iii),
    \[
    \min(A)=u(\{\min(A)\}_{<\min(A)})=u(\emptyset)=\dots=\min(B).
    \]
    The proof of Zorn's lemma mostly consists in a "transfinite iteration" of this argument, to show that all $u$-sets are tail sets of a largest $u$-set $E$, and by further applying $u$ to it we reach a 
    contradiction.
    
    {\bf Claim} Let $A\ne B$ be $u$-sets. Then, $A=B_{<b}$ for some $b\in B$, or $B=A_{<a}$ for some $a\in A$.

    Let
    \[
    C=\{c\in A\cap B:\ A_c=B_c\}.
    \]
    \noindent{\bf Sub-claim} $C=A$ or $C=A_{<a}$ for some $a\in A$.

    We show the sub-claim first. Suppose $C\ne A$ and let $a\in A\setminus C$ be minimal, so that $A_{<a}\subseteq C$. In the opposite direction, if $c\in C\setminus A_{<a}\subseteq A$, then $c>a$. Hence, 
    \[
    a\in A_{c}=B_{<c}\subseteq B,
    \]
    thus $a\in B$ and $A_{<a}=B_{<a}$, implying that $a\in C$, which is a contradiction. This implies that $C=A_{<}$.

    We now prove the claim. If $A\ne C\ne B$, by the sub-claim we have that $A_{<a}=C=B_{<b}$ for some $a\in A$, $b\in B$. By property (iii) of the $u$-sets, $a=u(A_{<a})=u(B_{<b})=b$, and we have that
    \[
    a=b\in A\cap B, \text{ and }A_{<a}=B_{<b}.
    \]
    By definition, then, $a\in C$, but this contradicts the fact that $A_{<a}=C$. The assumption $A\ne C\ne B$ can not hold, hence we have that, again by the sub-claim, either $A=C=B_{<b}$ or, the other way around, $B=C=A_{<a}$. (Hidden in the proof is that the cases $a=\min(A)$ or $b=\min(B)$, in which $C=\emptyset$, are covered because we defined $u(\emptyset)$).

    The claim allows us to glue together whatever $u$-sets we have into one. Let
    \[
    E=\cup_{A\text{ a $u$-set}}A.
    \]
    We first show that, if $A$ is a $u$-set and $a\in E$, then $A_{<a}=E_{<a}$. The direction $A_{<a}\subseteq E_{<a}$ follows from $A\subseteq E$. In the other direction, let $x\in E_{<a}$ and let $B\ni x$ be a $u$-set. If $B\subseteq A$, then $x\in A$. If not, by the claim $A=B_{<b}$ for some $b\in B$, in which case $x<a<b$, and so $x\in B_{<b}=A$. In both cases, $x\in A$, hence, $E_{<a}\subseteq A$, which clearly implies $E_{<a}\subseteq A_{<a}$.

It is now easy to verify that $E$ is a $u$-set.
\begin{itemize}
    \item[(i)] Let $a,b\in E$, and $A\ni a$, $B\ni b$ be $u$-sets. Then $A_{<a}=E_{<a}$ and $B_{<b}=E_{<b}$ are $u$-sets, hence, by the claim above, either they are equal, then $a=u(A_{<a})=u(B_{<b})=b$, or $A_{<a}=B_{<c}$ for some $c$ in $B_{<b}$, and $c=u(B_{<c})=u(A_{<a})=a$, or the other way around. In the second case, $a=c<b$, and in the remaining case the opposite inequality holds. Hence, $E$ is totally ordered. 
    \item[(ii)] Suppose $a_1>a>2>\dots$ is a descending sequence in $E$, and that $a_1\in A$, a $u$-set. For $i>1$, $a_i\in E_{a_1}=A_{a_1}$, hence the sequence, being contained in $A$, can not be infinite.
    \item[(iii)] If $a\in E$ and $A\ni a$ is a $u$-set, then $u(E_{<a})=u(A_{<a})=a$.
\end{itemize}
    Since $E$ is a $u$-set, $E\cup\{u(E)\}$ is still a $u$-set which is not contained in $E$, and we have reached a contradiction.
\end{proof}
In the spirit of the section, we use Zorn's lemma to prove a purely algebraic result.
\begin{theorem}\label{theoHamel}[Hamel's basis of a vector space]\index{Hamel basis}
Any vector space $X$ over a filed $\mathbb{F}$ has a basis. That is, there is a set $\{v_a\}_{a\in I}$ of linearly independent vectors of $X$ such that $X=\text{span}(\{v_a\}_{a\in I})$.
\end{theorem}
\begin{proof}
 Consider the set $\mathcal{V}$ of all families $V=\{v_b\}_{a\in J}$ of linearly independent vectors of $X$, partially ordered by inclusion. If $\{V_\alpha\}_{\alpha\in H}$ is a totally ordered subfamily in $\mathcal{V}$, then $\bigcup_{\alpha\in H}V_\alpha$ is an element of $\mathcal{V}$, and it is an upper bound for $\{V_\alpha\}_{\alpha\in H}$. By Zorn's lemma, there exists in $\mathcal{V}$ a maximal element $\{v_a\}_{a\in I}$, which is a basis for $X$, since a vector $u$ which is not in $\text{span}(\{v_a\}_{a\in I})$ could be added to the family, contradicting maximality.
\end{proof}
Below we use basic facts about cardinal numbers.
For properties of cardinal numbers, and of ordinal ones, you might look into \href{https://www2.math.uconn.edu/~solomon/math5026f18/OrdCard2.pdf}{Reed Solomon, Notes on ordinals and cardinals}.
\begin{theorem}\label{theoAlgebraicDimension}
    Let $X$ be a vector space over a field $\mathbb{F}$, and let $U=\{U_a\}_{a\in I}$ and $V=\{v_b\}_{a\in J}$ be two basis of $X$. Then, $\sharp(I)=\sharp(J)$, $I$ and $J$ have the same cardinality.
\end{theorem}
\begin{proof}
We know from linear algebra that $\sharp(I)$ is finite if and only if $\sharp(J)$ is, and they coincide. We can suppose then that both index sets are infinite.

    To each $a\in I$ associate the set $B(a):=\{b_1(a),\dots,b_{n(a)}(a)\}=\{b_1,\dots,b_n\}\subset J$, where
    \[
    U_a=\sum_{i=1}^n\lambda_i V_{b_i},
    \]
    with $0\ne\lambda_i\in\mathbb{F}$, is the unique expression of $U_a$ with respect to the basis $V$. The map $B$ is one to one from $I$ to the family of the finite subsets of $J$. Since the latter class has the same cardinality as $J$, we have $\sharp(I)\le\sharp(J)$, and the other inequality is proved in the same way.
\end{proof}
By theorem \ref{theoAlgebraicDimension}, we can attach to any vector space $X$ the cardinality of one, hence all, of its basis. Such number is called the {\it algebraic dimension}, or sometimes the {\it Hamel dimension}, of $X$. Let's denote it by $\text{dim}_{al}(X)$. The algebraic dimension is clearly a linear invariant: if $L:X\to Y$ is a linear bijection, then $\text{dim}_{al}(X)=\text{dim}_{al}(Y)$.

Infinite dimensional, but separable, Hilbert spaces have a countable orthonormal basis in the Hilbert sense, but their algebraic dimension is more than countable. This shows that Hilbert theory, and more general Banach theory, or even more generally topological vector spaces theory, is not reducible to linear algebra alone, even at the most basic level of dimension.
\begin{proposition}\label{propHilbHamel}
    As a vector space, any infinite dimensional Hilbert space $H$ does not have countable dimension.
\end{proposition}
\begin{proof}
    Let $\{v_n\}_{n=1}^\infty$ be any countable, linearly independent subset of $H$. Using Gram-Schmidt algorithm, find an orthonormal set $\{e_m\}_{m=1}^\infty$ such that $\text{span}(\{v_n\}_{n=1}^\infty)=\text{span}(\{e_m\}_{m=1}^\infty)$. Consider then
    \[
    x=\sum_{m=1}^\infty \frac{1}{m}e_m.
    \]
    Then, $x\in H\setminus \text{span}(\{e_m\}_{m=1}^\infty)=H\setminus \text{span}(\{v_n\}_{n=1}^\infty)$. This show that $\{v_n\}_{n=1}^\infty$ was not an algebraic basis for $H$. 
\end{proof}
The completeness of $H$ was used in that $x$ exists in $H$ since the series converges.

Using Hamel basis, we can construct such monstrous objects as everywhere defined, unbounded linear functionals on a Hilbert space.
\begin{proposition}\label{propUnboundedLinearEverywhere}
    Let $H$ be an infinite dimensional Hilbert space. Then, there exists a linear functional $l:H\to \mathbb{R}$ which is not bounded.
\end{proposition}
\begin{proof}
    Let $\{v_\alpha\}_{\alpha\in I}$ be a Hamel basis for $H$, let $\{v_{\alpha_n}\}_{n=1}^\infty$ a countable subset of it, which we might assume to be orthonormal after applying the Gram-Schmidt algorithm to it. We then define $l$ on the basis by
    \[
    \begin{cases}
        l(v_\alpha)=0\text{ if }\alpha\ne\alpha_n\text{ for all }n\ge1\crcr 
        l(v_{\alpha_n})=n.
    \end{cases}
    \]
    The operator extends to a linear operator $l$ on $H$ by linear algebra. On the other hand, $l$ is not bounded,
    \[
    \|l\|_{H^\ast}:=\sup_{\|x\|=1}|l(x)|\ge|l(v_{\alpha_n})|=n,
    \]
    hence, $\|l\|_{H^\ast}=\infty$.
\end{proof}
\section{The Hahn-Banach Theorem and some of its consequences}
\begin{footnotesize}\begin{exercise}
Let $H$ be a Hilbert space, $M\subset H$ a closed, linear space, and $S:M\rightarrow {\mathbb C}$ a (bounded) linear functional on $M$. Show that there exists an extension $T$ of $S$ to $H$ with $\|T\|=\|S\|$. Moreover, such an extension is unique.
\end{exercise}\end{footnotesize}
In general Banach spaces matters are more intricate.
\begin{theorem}[Real Hahn-Banach Theorem]\index{Hahn-Banach theorem!real}
Let $X$ be a real vector space, and $p:X\to{\mathbb R}$ a real-valued, convex functional:
\[
p(a x+(1-a) y)\le a p(x)+(1-a)p(y)
\]
if $x,y\in X$ and $a\in[0,1]$. Let $l:Y\to{\mathbb R}$ be a linear functional such that for $x\in Y$
\[
l(y)\le p(y).
\]
Then, there exists linear $L:X\to{\mathbb R}$ extending $l$ to $X$, and such that
\[
L(x)\le p(x)
\]
on $X$.
\end{theorem}
\begin{proof}
The main step consist in extending $l$ on $Y$ to $\lambda$ on ${\text span}(Y,z)$, where $z\notin Y$, which we do by a sort of separation of variables. Then we use Zorn's Lemma. We must have $\lambda(y+a z)=l(y)+a\lambda(z)$, and it suffices to determine $\lambda(z)$. We start with an inequality involving $z$ and $p$. Let $a,b>0$, $x,y\in Y$:
\begin{eqnarray*}
b l(x)+ a l(y)&=&(b+a) l\left(\frac{b}{b+a}x+\frac{a}{b+a} y\right)\crcr 
&\le& (b+a) p\left(\frac{b}{b+a}x+\frac{a}{b+a} y\right)\crcr
&=&(b+a) p\left(\frac{b}{b+a}(x-a z)+\frac{a}{b+a} (y+b z)\right)\crcr
&\le& b p(x-az)+a p(y+b z).
\end{eqnarray*}
That is,
\[
\frac{l(x)-p(x-a z)}{a}\le \frac{p(y+b z)-l(y)}{b}
\]
holds for all $a,b,x,y$ as above. There is then real $k$ such that $l(x)-p(x-a z)\le k a$ and $k b\le p(y+b z)-l(y)$, or
\[
l(x)-k a\le p(x-a z),\ l(y)+ k b\le p(y+b z).
\]
It follows that $\lambda(z)=k$ works.

Consider now the set $A$ having as element $(Z,\lambda)$ where $Z$ is a subspace of $X$ containing $Y$ and $\lambda$ is an extension of $l$ to $Z$ in such a way $\lambda(z)\le p(z)$ on $Z$. We partially order $A$ saying that $(Z,\lambda)\le(W,\mu)$ is $Z$ is a subspace of $W$ and $\mu$ extends $\lambda$.

Given a chain $\{(Z_c,\lambda_c):\ c\in C\}$ in $A$, we have that $Z=\cup_{c\in C}Z_c$ and $\lambda$ such that $\lambda|_{Z_c}=\lambda_c$ are a well defined element in $A$, and $(Z,\lambda)$ is maximal for $C$. By Zorn's Lemma, there is a maximal $(W,\nu)$ in $A$. Now, if $W\ne X$, we might apply the one-dimensional extension, and $(W,\nu)$ would not be maximal. Hence, $W=X$ and the theorem is proved.
\end{proof}
When proving Hahn-Banach's theorem for a separable Banach space $X$ with $p(x)=\|x\|$, as it is often the case, Zorn's lemma is not necessary. In fact, we can proceed in a way which is much similar to the Gram-Schmidt construction of a orthonormal basis for a Hilbert space.
\begin{footnotesize}\begin{exercise}\label{exeHBnoZorn}
    Let $(X,\|\cdot\|)$ be a separable Banach space and let $\{x_n\}_{n=1}^\infty$ be a dense set. Let $Z$ be a subspace of $X$ and let $l:X\to\mathbb{R}$ be a linear functional such that $|l(x)|\le K\|x\|$ for some constant $K\ge0$.
    \begin{enumerate}
        \item[(i)] Let $Z_0=\Bar{Z}$ be the closure of $Z$ in $X$. Show that $l$ has a unique extension $l_0$ to $Z_0$ with $|l_1(x)|\le K\|x\|$. 
        \item[(ii)] Let $x_{i_1}$ be the first $x_n$ in the sequence such that $x_n\notin Z_0$, and let $Z_1=\text{span}(Z_0,x_{i_1})$. Verify that the inductive step in the proof of the real Hahn-Banach theorem shows that $l_0$ extends to a linear functional $l_1$ on $Z_1$ such that $|l_1(x)|\le K\|x\|$.
        \item[(iii)] Write down how to iterate the procedure in (ii). What happens if it stops afetr finitely man steps?
        \item[(iv)] Suppose the procedure can be iterated indefinitely, exhausting $\{x_n\}$, and let $Z_\infty=\bigcup_{m=0}^\infty Z_m$, and $l_\infty:Z_\infty\to\mathbb{R}$ be such that $l_\infty|_{Z_m}=l_m$. Show that $Z_\infty$ is dense in $X$ and $|l_\infty(x)|\le K\|x\|$. Use then again (i) to show that $l_\infty$ has a unique extension $\lambda$ to $X$ which coincides with $l_\infty$ on $Z_\infty$ (hence, it coincides with $l$ on $Z$), and $|\lambda(x)|\le K\|x\|$.
    \end{enumerate}
\end{exercise}\end{footnotesize}

\begin{theorem}[Complex Hahn-Banach Theorem]\index{Hahn-Banach theorem!complex}
Let $X$ be a complex vector space and $p:X\to{\mathbb{R}}$ be such that
\[
p(a x+ b y)\le |a| p(x) + |b| p(y)
\]
whenever $a,b\in{\mathbb C}$, $|a|+|b|=1$, and $x,y\in X$. Let $l:Y\to{\mathbb C}$ be a linear functional defined on a subspace $Y$ on $X$, satisfying $|l(y)|\le p(y)$. Then, there exists a linear functional $L:X\to{\mathbb C}$ which extends $l$, and such that $|L(x)|\le p(x)$ on $X$.
\end{theorem}
\begin{proof}
We consider $X$ as a real space of "twice the dimension". Let $\lambda(x)={\text Re}(l(x))$, a real functional on $Y$, and observe that 
\[
\lambda(i x)={\text Re}(l(i x))={\text Re}(i l(x))=-{\text Im}(l(x)),
\]
so that $l(x)=\lambda(x)-i\lambda(i x)$ reconstructs $l$ from $\lambda$. By the real Hahn-Banach theorem, there exists a real functional $\Lambda:X\to{\mathbb R}$ which extends $\lambda$ to $X$, and such that $\Lambda(x)\le p(x)$. We define $L(x)=\Lambda(x)-i\Lambda(i x)$, which defines a complex linear functional on $X$:
\[
L(i x)=\Lambda(i x)-i\Lambda(-x)=i(\Lambda(x)-i \Lambda(i x))=i L(x).
\]
We also have, for some real $t$,
\[
|L(x)|=e^{i t}L(x)=L(e^{i t}x)=\Lambda(e^{i t}x)\le p(e^{i t}x)=p(x),
\]
where the third equality holds because $L(e^{i t}x)$ is real, and the last follows from the fact that $p(e^{i t}x)=p(x)$ for all $t$'s.
\end{proof}

\section{The dual of a Banach space} Let $X$ be a Banach space over ${\mathbb C}$ (virtually all we say holds for real Banach spaces). Its {\it dual space} $X^\ast$\index{Banach space!dual} is the space of the continuous linear functionals on $X$. It is a Banach space in itself. In fact more is true. Let $\mathcal{L}(X,Y)$  be the linear space of the bounded, linear operators $T:X\to Y$, normed with the operator norm,
\begin{footnotesize}\begin{exercise}
Show that the operator norm is in fact a norm on $\mathcal{L}(X,Y)$.\index{Banach space!of bounded operators}
\end{exercise}\end{footnotesize}
\begin{theorem}\label{TheoLinearMapsComplete}
Let $X,Y$ be Banach spaces.Then, $\mathcal{L}(X,Y)$ is complete.
\end{theorem}
\begin{proof}
Let $\{A_n:\ n\ge1\}$ be Cauchy in $\mathcal{L}(X,Y)$, so that $\{A_n x:\ n\ge1\}$ in Cauchy in $Y$ for each fixed $x$ in $X$, and $A x=\lim_{n\to\infty} A_n x$ is well defined, and linear. Also, $| \|A_m\|-\|A_n\| |\le\|A_m-A_n\|$, hence $\|A_n\|\to C$ has real limit as $n\to\infty$. We have:
\[
\|A x\|=\lim_{n\to\infty}\|A_n x\|\le \lim_{n\to\infty}\|A_n\| \|x\|\le C \|x\|,
\]
hence, $\|A\|\le C$.

We still have to show that $A_n\to A$ in $\mathcal{L}(X,Y)$ norm. But we have:
\[
\|(A-A_m) x\|= \lim_{n\to\infty}\|(A_n-A_m) x\|\le 
\lim_{n\to\infty}\|A_n-A_m\|\cdot\|x\|,
\]
hence,
\[
\|A-A_m\|\le \lim_{n\to\infty}\|A_n-A_m\|\le\epsilon
\]
if $m\ge m(\epsilon)$, so $\|A-A_m\|\to 0$ as $m\to\infty$.
\end{proof}

\begin{footnotesize}\begin{exercise}
Let $S\in \mathcal{L}(X,Y)$ and $T \in \mathcal{L}(Y,Z)$. Show that $T S:=T\circ S\in \mathcal{L}(X,Z)$ satisfies 
$\|T S\|\le\|T\| \|S\|$.
\end{exercise}\end{footnotesize}
The above inequality, when $X=Y=Z$, expresses the fact that $\mathcal{L}(X):=\mathcal{L}(X,X)$ is a {\it Banach algebra}\index{Banach algebra} with respect to the composition product. It is also a {\it unital} one, since the identity operator $I x=x$ satisfies $ I A=A I$ and $\|I\|=1$.

We will often use without mention the following fact.
\begin{footnotesize}
\begin{exercise}\label{exeUniqueExt}
    Let $f_0:X_0\to Y$ be a continuous map from a dense subspace $X_0$ of a metric space $X$ and a complete metric space $Y$. Then, $f_0$ can be uniquely extended to a continuous map $f:X\to Y$. If $X$ is a normed space, $Y$ is a Banach space, and $f_0$ is linear, then $f$ is linear and $\|f\|_{\mathcal{L}(X,Y)}=\|f_0\|_{\mathcal{L}(X_0,Y)}$.
\end{exercise}
\end{footnotesize}
On the space $\mathcal{L}(X,Y)$ we will make use of three different topologies, which naturally arise from the theory of operators and its applications.
\begin{enumerate}
    \item[($UOT$)] The {\it uniform operator topology (UOT)} is the one induced by the operator norm:\index{operator topology!uniform}
    \[
    T_n\xrightarrow[UOT]{}T\text{ if and only if }\|T_n-T\|\to0.
    \]
    \item[($SOT$)] The {\it strong operator topology (SOT)} is defined by:\index{operator topology!strong}
    \[
    T_n\xrightarrow[SOT]{}T\text{ if and only if }\|T_nx-Tx\|\to0
    \]
    for all $x$ in $X$.
    \item[($WOT$)] The {\it weak operator topology (WOT)} is defined by:\index{operator topology!weak}
    \[
    T_n\xrightarrow[WOT]{}T\text{ if and only if }l(T_nx-Tx)\to0
    \]
    for all $x$ in $X$ and all $l\in Y^\ast$.
\end{enumerate}

The Hahn-Banach Theorem has important consequences about $X^\ast$.
\begin{corollary}\label{CorHABOne}
Let $i:Y\to X$ be the inclusion in the Banach space $X$ of $Y$, a subspace of it. Then, the restriction map $l\mapsto l\circ i$ from $X^\ast$ to $Y^\ast$ is surjective: each $\lambda$ in $Y^\ast$ has an extension $l$ in $X^\ast$. Moreover, there exists one extension $l$ such that $\|l\|_{X^\ast}=\|\lambda\|_{Y^\ast}$.  
\end{corollary}
\begin{proof}
Apply the Hanh-Banach Theorem with $p(x)=\|\lambda\|_{Y^\ast}\cdot\|x\|$.
\end{proof}

\begin{corollary}\label{CorHABTwo}
Let $y\in X$, Banach. Then, there exists $l\in X^\ast$, $l\ne0$, such that $l(y)=\|l\|_{X^\ast}\|y\|$.
\end{corollary}
\begin{proof}
Set $\lambda(a y)=a\|y\|$ on ${\text span}(y)$, and extend it to $l:X\to{\mathbb C}$ by Hahn-Banach. Since $|\lambda(a y)|=\|a y\|$, we have $\|l\|=\|\lambda\|=1$, and $l(y)=\lambda(y)=\|y\|$.
\end{proof}

\begin{corollary}\label{CorHABThree}
Let $Z$ be a subspace of a normed linear space $X$, and let $d(y,Z):=\inf\{\|y-z\|:\ z\in Z\}$. Then, there is $L\in X^\ast$ such that $L(y)=d(y,Z)$, $\|L\|_{X^\ast}= 1$, and $L(z)=0$ for $z\in Z$.
\end{corollary}
\begin{proof}
We only have to define such $L$ on ${\text span}(Z,y)$, then use Hahn-Banach. We are forced to define $L(z+ a y)=a L(y) = a d(y,Z)$. This is a linear functional, and
\begin{eqnarray*}
\frac{|L(z+a y )|}{\|z+ a y\|}&=&|a|\frac{d(Z,y)}{\|z+ a y\|}\crcr
&=&|a|\frac{\inf\{\|w-y\|:\ w\in Z\}}{\|z+ a y\|}\crcr
&=&|a|\frac{\inf\{\|w-y\|:\ w\in Z\}}{|a|\cdot \|-\frac{z}{a}- y\|}\crcr
&\le&1,
\end{eqnarray*}
because $-z/a\in Z$. Hence, $\|L\|_{X^\ast}\le1$.

To show the opposite inequality, for $\epsilon>0$ find $z$ in $Z$ such that $\|y-z\|\le d(y,Z)+\epsilon$, so that
\[
\frac{|L(z-y)|}{\|z-y\|}\ge\frac{d(y,Z)}{d(y,Z)+\epsilon}.
\]
Since $\epsilon>0$ is arbitrary, $\|L\|_{X^\ast}\ge1$.
\end{proof}

\begin{footnotesize}\begin{exercise}
Give "Hilbertian" elementary proofs of the three corollaries above in the case when $X$ is a Hilbert space.
\end{exercise}\end{footnotesize}

The {\it bidual} of a Banach space $X$, $X^{\ast\ast}$, is the dual of $X^\ast$.
\begin{theorem}\label{TheoBidual}
Let $X$ be a Banach space and consider the map $i:X\rightarrow X^{\ast\ast}$ given by $[i(x)](\ell):=\ell(x)$ whenever $x\in X$ and $\ell\in X^\ast$. Then, $i$ is an isometric imbedding of $X$ into $X^{\ast\ast}$.
\end{theorem}
\begin{proof} We have:
\[
|[i(x)](\ell)|=|\ell(x)|\le \|\ell\|_{X^\ast}\|x\|_X,
\]
hence, $\|[i(x)]\|_{X^{\ast\ast}}\le \|x\|_X$.

In the other direction, by Corollary \ref{CorHABTwo} we can find $l$ in $X^\ast$ with $\|l\|_{X^\ast}=1$, and such that $l(x)=\|x\|$. i.e.,
\[
\|i(x)\|_{X^{\ast\ast}}\ge\frac{|[i(x)](l)|}{\|l\|_{X^\ast}}
=\frac{|l(x)|}{\|l\|_{X^\ast}}=\|x\|.
\]
\end{proof}

When the map $i:X\to X^{\ast\ast}$ is surjective (hence, a surjective isometry), we say that $X$ is {\it reflexive}.\index{Banach space!reflexive}
Not all Banach spaces are reflexive, as the following exercise shows.

\begin{footnotesize}\begin{exercise}
\begin{enumerate} Here $\ell^p=\ell^p({\mathbb N})$.
    \item[(i)] Show that for each $\varphi\in \ell^\infty$ the map $L_\varphi:\psi\mapsto \sum_{n=0}^\infty \varphi(n)\overline{\psi(n)}$, $L_\varphi:\ell^1\to{\mathbb C}$, is a continuous linear functional in $(\ell^1)^\ast$ and $\|L_\varphi\|_{(\ell^1)^\ast}=\|\varphi\|_\ell^\infty$. 
    \item[(ii)] Given $L\in(\ell^1)^\ast$, show that there exists $\varphi\in \ell^1$ s.t. $L=L_\varphi$. Moreover, $\|\varphi\|_{\ell^\infty}=\|L\|_{(\ell^1)^\ast}$. {\bf Hint.} Make use of $L(\delta_n)$ as "building blocks" to construct $\varphi$. This way, we have isometrically identified $(\ell^1)^\ast\equiv\ell^\infty$.
    \item[(iii)] Consider the subspace $C$ of $\ell^\infty$,
    \[
    C=\{\varphi\in\ell^\infty:\ \exists \lim_{n\to\infty}\varphi(n)\}.
    \]
    Show that is is a closed subspace $\ell^\infty$, and that $L:\varphi\mapsto \lim_{n\to\infty}\varphi(n)$ is a closed linear functional on $C$.
    \item[(iv)] Show that there is no $\psi\in \ell^1$ such that $C(\varphi)=\sum_{n=0}^\infty\psi(n)\varphi(n)$. Deduce from this that $\ell^1$ is not reflexive.
\end{enumerate}
\end{exercise}\end{footnotesize}

\begin{footnotesize}\begin{exercise} Let $\ell_c$ be the space of the sequences $h:{\mathbb N}\to{\mathbb C}$ which vanish outside a finite set, and $\ell$ be the space of all sequences $h:{\mathbb N}\to{\mathbb C}$.
\begin{enumerate}
\item[(i)]For $m:{\mathbb N}\to{\mathbb C}$, define the multiplication operator with symbol $m$, $T_m:\ell_c\to\ell$, by $T_m(h)(n)=m(n)h(m)$. Show that $\|T_m\|_{{\mathbb B}(\ell^2,\ell^2)}=\|m\|_{\ell^\infty}$.
    \item[(ii)] Find a similar statement with $L^2[0,1]$ instead of $\ell^2$, and prove it.
\end{enumerate}
\end{exercise}\end{footnotesize}

\begin{footnotesize}\begin{exercise}
Consider $\ell^1=\ell^1({\mathbb N})$, its $1$-dimensional, closed subspace ${\text span}(\delta_0)$, and $l(a \delta_0)=a$. Find all extensions $L$ of $l$ to a linear functional on $\ell^1$ satisfying $\|L\|_{(\ell^1)^\ast}=\|l\|=1$.

How is this different from the Hilbert case?
\end{exercise}\end{footnotesize}
Another reason why $\ell^1(\mathbb{N})$ is not reflexive is that, while $\ell^1(\mathbb{N})$ is separable, $\ell^\infty(\mathbb{N})$ is not.
\begin{theorem}\label{theoSepPreDual}[The predual of a separable space is separable]\index{theorem!predual of separable is separable}
    If $X$ is a Banach space and $X^\ast$ is separable, then $X$ is separable.
\end{theorem}
\begin{proof}
    We consider the case of a real Banach space, the complex case being identical.
    If $X^\ast$ is separable, it has a dense, countable subset $\{l_n\}_{n=1}^\infty$. For each $n$, let $x_n$ be such that $\|x_n\|=1$ and $l_n(x_n)\ge \|l_n\|_{X^\ast}/2$. Suppose $Y_0=\text{span}_{\mathbb{Q}}\{x_n\}_{n=1}^\infty$, the linear combinations of the $x_n$'s with rational coefficients, is not dense in $X$. So, $\Bar{Y_0}\ne X$, and we can find $z$ in $X\setminus Y$. By Corollary \ref{CorHABThree}, there exists $l\in X^\ast$ such that $l|_Y=0$, and $\|l\|_{X^\ast}=1$. We have that
    \[
    \|l-l_n\|_{X^\ast}\ge|(l-l_n)(x_n)|=|l_n(x_n)|\ge \|l_n\|_{X^\ast}/2.
    \]
    If $\{l_n\}$ were dense in $X^\ast$, we could find a subsequence $\{l_{n_k}\}$ such that
    \[
    0=\lim_{k\to\infty}\|l-l_{n_k}\|_{X^\ast}\ge\limsup_{k\to\infty}\|l_{n_k}\|_{X^\ast}/2,
    \]
    then,
    \[
    1=\|l\|_{X^\ast}=\lim_{k\to\infty}\|l_{n_k}\|_{X^\ast}=0,
    \]
    which is a contradiction.
\end{proof}
\begin{footnotesize}
\begin{exercise}\label{eqSepBanach}
\begin{enumerate}
    \item[(i)] Show that $\ell^\infty(\mathbb{N})$ is not separable, while $\ell^p(\mathbb{N})$ is separable for $1\le p<\infty$.
    \item[(ii)] Show that $L^\infty[0,1]$ is not separable, while $L^p[0,1]$ is separable for $1\le p<\infty$ (with respect to the Lebesgue measure).
    \item[(iii)] Show that $C[0,1]$ is separable (with respect to the uniform norm).
    \item[(iv)] Show that a Hilbert space $H$ is separable if and only if it has a countable basis.
    \item[(v)] Find Banach spaces $X,Y$ with $X\subset Y$, such that the imbedding map $x\mapsto x$ is bounded from $X$ to $Y$, $Y$ is separable, and $X$ is not separable.
\end{enumerate}
\end{exercise}
\end{footnotesize}
In the second part of the next exercise, we need a definition. If $A:X\to Y$ is a linear, bounded operator between Banach spaces, its {\it adjoint} $T':Y^\ast\to X^\ast$ is defined on a linear functional $l\in Y^\ast$ by $[T'l](x)=l(Tx)$.\index{operator!adjoint} 
\begin{footnotesize}
\begin{exercise}\label{exeMultOp}
   Let $1\le p\le\infty$ and let $g:\mathbb{R}\to\mathbb{C}$ be measurable. On $L^p(\mathbb{R})$, consider the multiplication operator $M_g:f\mapsto gf$.
   \begin{enumerate}
       \item[(i)] Show that
       \[
       \|M_g\|_{\mathcal{L}(L^p)}:=\sup_{0\ne f\in L^p[0,1]}\frac{\|M_gf\|_{L^p}}{\|f\|_{L^p}}=\|g\|_{L^\infty}.
       \]
       \item[(ii)] Let $1\le p<\infty$ and $g\in L^\infty(\mathbb{R})$. After identifying $[L^p(\mathbb{R})]^\ast=L^{p'}(\mathbb{R})$ ($1/p+1/p'=1$), what is the expression for $[M_g]'$, the adjoint of $M_g$? 
   \end{enumerate}
\end{exercise}
\end{footnotesize}

\section{Weak and weak$^\ast$ topologies, and the Banach-Alaoglu theorem}\label{SectBanachAlaoglu}
Our first experience is that we fix a topology on a set $Y$ (a notion of points "being close"), and this way we can tell which functions $f:Y\to\mathbb{R}$ are continuous (their values do not change "abruptly" from point to point). We might instead want to put in the forefront the functions (the "measurables"), rather than the points. We  fix a family $\mathcal{F}$ of functions $f:Y\to\mathbb{R}$, and require that those functions are continuous (that they do not change abruptly: that we consider two points to be close, that is, if they result to be close in all our observations). The {\it topology $\tau(\mathcal{F})$ generated by }$\mathcal{F}$ is the smallest one making all functions in $\mathcal{F}$ continuous.

It is easy to show that a basis of neighborhoods for $\mathcal{F}$ is given by finite intersections of basic sets having the form
\[
\mathcal{N}(f,a,r)=f^{-1}(a-r,a+r)=\{y\in Y:|f(x)-a|<r\}.
\]
One advantage of such {\it weak topologies} is that they are most economical in terms of open sets, hence they have the largest family of compact sets, which are good when convergence of (sub)sequences is concerned. 
We will consider below some important instances of this construction. 
\subsection{The weak and the weak$^\ast$ topologies}\label{SSectWeakTopology}
\subsubsection{The weak topology}\label{SSSectWeakTopology}\index{Banach space!weak topology}
Let $X$ be a Banach space, and let $X^\ast$ be its dual. The {\it weak topology} $w$ on $X$ is the coarsest (weaker) which makes all functionals $l\in X^\ast$ continuous. By general topology, a basis of neighborhoods for the weak topology is given by the family of the subsets $N(a;l_1,\dots,l_n;\epsilon)$ (with $a\in X$, $n\ge1$, $l_1,\dots,l_n\in X^\ast$, $\epsilon>0$), where:
\begin{eqnarray}\label{eqtWeakTopology}
    N(a;l_1,\dots,l_n;\epsilon)&=&\{x\in X:|l_1(x)-l_1(a)|<\epsilon,\dots,|l_n(x)-l_n(a)|<\epsilon\}\crcr 
    &=&N(a;l_1;\epsilon)\cap\dots\cap N(a;l_n;\epsilon)\crcr 
    &=&[N(0;l_1;\epsilon)+a]\cap\dots\cap[N(0;l_n;\epsilon)+a]\crcr 
    &=&[\epsilon N(0;l_1;1)+a]\cap\dots\cap[\epsilon N(0;l_n;1)+a]\crcr
    &=&\epsilon N(0;l_1,\dots,l_n;1)+a\crcr 
    &=&N(0;l_1,\dots,l_n;\epsilon)+a
\end{eqnarray}
are various ways to write and think of these basic neighborhoods. For instance, it is clear that the topology is invariant under translations and dilations, $x\mapsto \lambda x+a$ is a homeomorphism whenever $\lambda\ne 0$.

If $(X,w)$ is metrizable, and $\Omega$ is a metric space, then continuity of $f:X\to \Omega$ with respect to the weak topology is equivalent to continuity by sequences: $x_n\to x$ in $(X,w)$ implies $f(x_n)\to f(x)$ in $\Omega$. Unfortunately $(X.w)$ might not be metrizable, and we have to use {\it nets} instead. But some nice features remain.
\begin{proposition}\label{propWeakTopologyIsHausdorff}
    $(X,w)$ is a Hausdorff space.
\end{proposition}
\begin{proof}
By translation invariance, we just have to separate $a\ne0$ and $0$. By corollary \ref{CorHABTwo} to the Hahn-Banach theorem, there is $l\in X^\ast$ with $\|l\|_{X^\ast}=1$ and such that $l(a)=\|a\|$. If $N(0;l;\|a\|/2)\ni 0$ and $N(a;l;\|a\|/2)\ni a$ had a point $x$ in common, 
\[
\|a\|=|l(a)-l(0)|\le |l(a)-l(x)|+|l(x)-l(0)|< \|a\|,
\]
a contradiction.
\end{proof}
We say that the sequence $\{x_n\}$ in $X$ {\it converges weakly} to $a\in X$ if $l(x_n)\to l(a)$ for all $l$ in $X^\ast$. Although this notion is generally weaker than "weak net convergence", it is nonetheless useful in many applications, for instance in calculus of variations. We write $x_n\xrightarrow[\text{w}]{}a$, or $w-\lim_{n\to\infty}x_n=a$. By definition, weak convergence of $x_n$ to $a$ means that, for each $l$ in $X^\ast$ and each $\epsilon>0$, there exists $n(l,\epsilon)$ such that, if $n>n(l,\epsilon)$, then $x_n\in N(l;a;\epsilon)$.
\begin{proposition}\label{propWeakConvergence}
\begin{enumerate}
    \item[(1)] If $\|x_n-a\|\to0$, then $x_n\xrightarrow[w]{}a$.
    \item[(2)] If $x_n\xrightarrow[w]{}a$, then $\|a\|\le \liminf_{n\to\infty}\|x_n\|$.
\end{enumerate}
\end{proposition}
\begin{proof}
\noindent (1) For any $l$ in $X^\ast$, $|l(x_n)-l(a)|\le \|l\|_{X^\ast}\|x_n-a\|\to0$.

\noindent (2) Again by corollary \ref{CorHABTwo},  there is $l$ with $\|l\|_{X^\ast}=1$ and $\|a\|=l(a)$, hence:
\[
\|a\|=l(a)=\lim_{n\to\infty}|l(x_n)|\le \|l\|_{X^\ast}\liminf_{n\to\infty}\|x_n\|_X=\liminf_{n\to\infty}\|x_n\|_X.
\]
\end{proof}
We can not improve the statement in (2). Let $H$ be an infinite dimensional Hilbert space and let $\{e_n\}_{n=1}^\infty$ be a orthonormal system in it. Then, 
\begin{equation}\label{eqWnS}
e_n\xrightarrow[\text{w}]{}0,
\end{equation}
but $\|e_n\|=1$ for all $n$.
\subsubsection{The $\text{weak}^\ast$ topology}\label{SSSectWeakStarTopology}\index{Banach space!weak$^\ast$ topology}
Let again $X$ be a Banach space, and let $X^\ast$ be its dual. The {\it $\text{weak}^\ast$ topology} $w^\ast$ on $X^\ast$ is the coarsest which makes the functionals $l\mapsto l(x)$ continuous for all $x$ in $X$. Unless the natural identification $X\hookrightarrow X^{\ast\ast}$ is surjective, i.e. unless $X$ is {\it reflexive}, the  {\it $\text{weak}^\ast$ topology} on $X^\ast$ is weaker than the weak topology.

A basis of neighborhoods for the origin in $w^\ast$ is given by finite intersections of sets of the form
\begin{equation}\label{eqNeighborhoodsForWeakStar}
    N(0;x;\epsilon)=\{l\in X^\ast:|l(x)|<\epsilon\}.
\end{equation}
The topology is invariant under translations and a basis at any point is easy to write, as in the case of the weak topology. Its basic properties coincide with those of the weak topology.
\begin{proposition}\label{propWeakStarTopology}
    \begin{enumerate}
        \item[(1)]  $(X^\ast,w^\ast)$ is Hausdorff.
        \item[(2)]  If $\|x_n-a\|\to0$, then $x_n\xrightarrow[w^\ast]{}a$.
        \item[(3)] If $l_n\xrightarrow[w^\ast]{}l$, then $\|l\|_{X^\ast}\le \liminf_{n\to\infty}\|l_n\|_{X^\ast}$.
    \end{enumerate}
\end{proposition}
\begin{proof}
    Property (2) is weaker than the analogous statement for the weak topology. Property (1) is proved like the corresponding statement for the weak topology, using the fact that if $l\ne0$ is an element of $X^\ast$, there exists $x$ in $X$ with $\|x\|_X=1$ and $|l(x)|\ge \|l\|_{X^\ast}/2$. If $\lambda\in N(0;x;\|l\|_{X^\ast}/4)\cap N(l;x;\|l\|_{X^\ast}/4)$, then
    \[
    \|l\|_{X^\ast}/2\le|l(x)|\le|(l-\lambda)(x)|+|\lambda(x)|<\|l\|_{X^\ast}/4+\|l\|_{X^\ast}/4,
    \]
    a contradiction.
    
    The proof of (3) follows the same lines as that of the weak topology analog. For $\epsilon>0$ there is $x$ in $X$ with $\|x\|=1$ such that $|l(x)|\ge\|l\|_{X^\ast}-\epsilon$. Then,
    \[
    \|l\|_{X^\ast}-\epsilon\le |l(x)|=\lim_{n\to\infty}|l_n(x)|\le \liminf_{n\to\infty}\|l_n\|_{X^\ast}.
    \]
    Let then $\epsilon\to0$.
\end{proof}
\subsection{Two versions of the Banach-Alaoglu theorem}\label{SSectBanachAlaoglu}
The {\it raison d'être} of the $\text{weak}^\ast$ topology is related to compactness properties. To keep the exposition self.contained, below you find a proof of the Tychonoff theorem.
\subsubsection{Tychonoff's theorem}\label{SSectTychonoff}
If $\{X_i:i\in I\}$ is a family of sets, $X_I=\Pi_{i\in I}X_i$, their Cartesian product, is best interpreted as the set of the functions
\[
p:I\to\cup_{i\in I}X_i, \text{ with }p(i)\in X_i \text{ for each }i\in I.
\]
We call $I$ the {\it domain} of $p$.

If each $X_i$ is a topological space, the {\it product topology} on $\Pi_{i\in I}X_i$ is the weakest (coarsest, minimal) making all projections $\pi_i:X_I\to X_i$, $\pi_i(p)=p(i)$, continuous. A basis of open sets for it is provided by the sets $\pi_{i_1}^{-1}(U_{i_1})\cap\dots\cap \pi_{i_n}^{-1}(U_{i_n})$, where $n\ge1$, $i_1,\dots,i_n\in I$, and $U_{i_l}$ is open in $X_{i_l}$.
\begin{theorem}\label{theoTychonoff}\index{theorem!Tychonoff}
    If each $X_i$ is compact, $i\in I$, then $\Pi_{i\in I}X_i$ is compact.
\end{theorem}
We give the first of the three proofs surveyed in \href{https://hal.science/hal-03660150/document}{Three Proofs of Tychonoff's Theorem} by É. Matheron (2020), which the author labels the {\it Wisconsin proof}.

Some of the usual notions associated to functions carry over and are useful. If $J\subset H$, $\Pi_{H,J}p=p|_J$ is the restriction of $p\in X_H$ to $X_J$, $\pi_{H,J}:X_H\to X_J$, $\pi_{H,J}\circ\pi_{J,K}=\pi_{H,K}$. We set $\pi_J=\pi_{I,J}$. If $p\in X_H$ and $q\in X_J$ with $H$ and $J$ disjoint, then $p\vee q\in X_{H\cup J}$ is the function such that $(p\vee q)|_H= p$,  $(p\vee q)|_J= q$. The set $X_\emptyset$ reduces to the unique {\it empty function} from $\emptyset$ to itself, which we denote by $\emptyset$ (in agreement with the interpretation of functions as particular relations).

We fix a set $I$, and let $\mathbb{P}=\cup_{J\subseteq I}X_J$, which is partially ordered by the relation $p\le q$ if $q|_J=p$, where $J$ is the domain of $p$.
\begin{proof}
The proof is by contradiction. We suppose that there is a family $\mathcal{U}$ of open sets in $X_I$ such that no finite subfamily covers $X_I$, and we shall show that $\mathcal{U}$ itself does not cover $X_I$ by exhibiting an element $p\in X_I\setminus\cup_{U\mathcal{U}}U$.

We consider the set $B$ of the {\it bad} elements in $\mathbb{P}$: those $p\in X_J$ such that for all open $V\ni p$ open in $X_J$, $\pi_{J}^{-1}(V)$ can not be covered by finitely many sets in $\mathcal{U}$. The empty function lies in $B$, which is then nonempty. The only open set in $X_\emptyset$ containing $\emptyset$, in fact, is $X_\emptyset=\{\emptyset\}$ itself, and $\pi_\emptyset^{-1}(X_\emptyset)=X_I$, which can not be covered by a finite subfamily of $\mathcal{U}$ by assumption. We will prove three facts about $B$.

(i) {\it $B$ is downward close}: if $p\le q\in B$, then $p\in B$. In fact, if $H$ is the domain of $p$ and $J$ that of $q$, and $V\ni p$ is open in $X_H$, then $\pi_H^{-1}(V)=\pi_J^{-1}(\pi_{H,J}^{-1}(V))$ which can not be covered by a a finite subfamily of $\mathcal{U}$ because $\pi_{H,J}^{-1}(V)\ni q$ is open and $q$ is bad. We can assume $V_a=O_a\times W_a$

(ii) {\it If $p\in B$ has domain $J\subset I$ and $i_0\in I\setminus J$, then there is $a\in X_{i_0}$ such that $p\vee a\in B$.} Suppose by contradiction that for all $a$ in $X_{i_0}$ we have $p\vee a\notin B$; i.e. there exists $V_a\ni p\vee a$ open in $X_{J\cup\{i_0\}}$ such that $\pi_{J\cup\{i_0\}}^{-1}(V_a)$ can be covered by a finite subfamily of $\mathcal{U}$. We can assume that $V_a=O_a\times W_a$ with $O_a$ open in $X_J$ and $W_a\ni a$ open in $X_{i_0}$. Since $X_{i_0}$ is compact, it can be covered by finitely many $W_{a_1},\dots,W_{a_n}$. Let $O=O_{a_1}\cap\dots\cap O_{a_n}\ni p$, so that $O\times W_{a_j}\ni p\vee a_j$ and $\pi_{J\cup\{i_0\}}^{-1}(O\times W_{a_j})$ can be covered by a finite subfamily of $\mathcal{U}$. Now, 
\[
\pi_J^{-1}(O)=\bigcup_{j=1}^n\pi_{J\cup\{i_0\}}^{-1}\left(O\times W_{a_j}\right),
\]
and each set in the union can be covered by a finite subfamily of $\mathcal{U}$. Hence, $p\notin B$.

(iii) {\it $(B,\le)$ has a maximal element}. If any chain $C$ in $B$ has an upper bound in $B$,  by Zorn's lemma $B$ has a maximal element. Let $C$ be a chain in $B$, let $H$ be the union of the domains of all $q$'s in $C$, and define $p\in X_H$ to be such that $p|_J=q$ if $q\in C$ has domain $J$. The function $p$ is well defined because $C$ is totally ordered. Thus $p$ is an upper bound of $C$ in $\mathbb{P}$. It suffices to show that $p\in B$. Let $V\ni p$ be an open subset of $X_H$. By definition of the product topology, we can choose $V\supseteq\pi_{H,F}^{-1}(W)\ni p$ where $F\subset H$ is finite and $W\subset X_F$ is open with $W\ni p|_F$. Now, $C$ is a chain and each $i_l\in F$ belongs to some domain $J_l$ in the chain, so $F\subseteq J_1\cup\dots J_m=J_0\subseteq H$, where $J_0$ is a domain in the chain. Let $q_0\in X_{J_0}\cap C$. Since $q_0$ is bad, $p|_F=q_0|_F$ is bad as well by (i), thus $\pi_F^{-1}(W))$ can not be covered by a finite subfamily of $\mathcal{U}$. A fortiori, $\pi_H^{-1}(V)\supseteq \pi_F^{-1}(W)$ can not be covered by a finite subfamily of $\mathcal{U}$, hence $p$ is bad.

Summarizing, $B$ has a maximal element $p$ by (iii), having domain $I$ by (ii). In particular $X_I\ni p$ can not be covered by a finite subfamily of $\mathcal{U}$, but this is possible only if $p$ does not belong to any subset of $\mathcal{U}$, which henceforth does not cover the whole of $X_I$.
\end{proof}

\subsubsection{Banach-Alaoglu theorem: the topological form}\label{SSSectBAtop}\index{theorem!Banach-Alaoglu: topological form}
\begin{theorem}[Banach-Alaoglu]\label{theoBanachAlaoglu}
    Let $\overline{B}_1$ be the closed unit ball of $X^\ast$. Then, $\overline{B}_1$ is compact in the $\text{weak}^\ast$ topology.
\end{theorem}
\begin{proof}
    For each $x\in X$, consider $\overline{D(0,\|x\|)}$, the closed unit disc in the complex plane, and let
    \begin{equation}\label{eqBanachAlaoglu}
        K:=\Pi_{x\in X}\overline{D(0,\|x\|)},
    \end{equation}
    endowed with the product topology. Such is the topology having as basis at $\{f(x):x\in X\}$ with $|f(x)|\le\|x\|$, finite intersections of sets of the form
    \[
    M(f;a;\epsilon)=\{g:X\to\mathbb{C}\ \text{ s.t.}|g(x)|\le\|x\|\text{ and }|g(a)-f(a)|<\epsilon\}.
    \]
    Since each factor $\overline{D(0,\|x\|)}$ is compact, by Tychonoff theorem $K$ is compact, too.

    We imbed
    $\overline{B}_1$ into $K$, $\Phi:\overline{B}_1\to K$,
    \[
    \Phi(l)=\{\{l(x)\}_{x\in X}\}.
    \]
    Unraveling definitions, we have that
    \[
   \Phi(N(l;a;\epsilon))=M(l;a;\epsilon)\cap \Phi(\overline{B}_1),
    \]
    i.e. $\Phi$ is a homeomorphism of $\overline{B}_1$ onto its image into $K$. Since the latter is compact, if we show that $\Phi(\overline{B}_1)$ is closed in $K$, then it is compact, hence $\overline{B}_1$ is compact as well.

        Suppose $f$ lies in the closure of $\Phi(\overline{B}_1)$, and consider $x,y\in X$. For any $\epsilon>0$ there is $l\in \overline{B}_1$ such that $|l(x)-f(x)|<\epsilon/3$, $|l(y)-f(y)|<\epsilon/3$, and $|l(x+y)-f(x+y)|<\epsilon/3$,
    so that
    \[
    |f(x)+f(y)-f(x+y)|<\epsilon.
    \]
    Hence, $f(x)+f(y)=f(x+y)$. The same way, one shows that $f(\lambda x)=\lambda f(x)$ if $\lambda$ is a complex number. This shows that $f:X\to\mathbb{C}$ is a linear functional, which also lies in $\overline{B}_1$, since $|f(x)|\le\|x\|$. We have proved that $\Phi(\overline{B}_1)$ is closed, hence $K$ is compact.
\end{proof}
\subsubsection{Banach-Alaoglu theorem: the sequential form}\label{SSSectBAseq}\index{theorem!Banach-Alaoglu: sequential form}
A practical problem with this version of the Banach-Alaoglu theorem is that in general compactness is not in general equivalent to sequential compactness, which is very useful in many applications (measure theory and probability, claculus of variations, PDEs...). The version involving sequential compactness below could be deduced by the previous one, but we provide a direct proof which is similar to the original one for the Ascoli-Arzelà theorem, which is essentially constructive and does not require the Tychonoff theorem.
\begin{theorem}[Banach-Alaoglu, separable pre-dual]\label{theoBAseparable} Let $X$ be a separable Banach space, and let $\{l_n\}$ be a sequence in $\overline{B_1}$, the unit ball in $X^\ast$. then, there exists a subsequence $\{l_{n_k}\}$ and an element $l$ in $\overline{B_1}$ such that
    $l_{n_k}\xrightarrow[\text{w}^\ast]{}l$.
\end{theorem}
\begin{proof}
    Let $\{z_j\}_{j=1}^\infty$ be a dense sequence in $X$, let $X_0$ be the dense space it generates in $X$, and let $\{x_m\}_{m=1}^\infty$ be a maximal family of linearly independent vectors in $\{z_j\}_{j=1}^\infty$, so that $X_0=\text{span}(\{x_m\}_{m=1}^\infty)$. We will first proceed to define $l$ on $X_0$.
    
    Consider a subsequence $n^{(1)}=\{n^{(1)}_k\}_{k=1}^\infty$ such that
    \[
    \lim_{k\to\infty}l_{n^{(1)}_k}(x_1)=a_1
    \]
    exists in $\mathbb{C}$. Such sequence exists because $|l_n(x_1)|\le\|l_n\|_{X^\ast}\|x_1\|_X\le\|x_1\|_X $ is bounded in $\mathbb{C}$. Suppose that subsequences $\{n\}_{n=1}^\infty\supset\{n^{(1)}_k\}_{k=1}^\infty\supset\dots\supset\{n^{(m)}_k\}_{k=1}^\infty$ have been chosen in such a way (i) $\lim_{k\to\infty}l_{n^{(j)}_k}(x_j)=a_j\in\mathbb{C}$ when $1\le j\le m$, and (ii) $n^{(1)}_1<n^{(2)}_1<\dots<n^{(m)}_1$. As above, we can choose a subsequence $\{n^{(m+1)}_k\}_{k=1}^\infty$ of $\{n^{(m)}_k\}_{k=1}^\infty$ such that $\lim_{k\to\infty}l_{n^{(m+1)}_k}(x_{m+1})=a_{m+1}\in\mathbb{C}$, and $n^{(m+1)}_1>n^{(m)}_1$. 
    The sequence $\{x^{(m)}_1\}_{m=1}^\infty$ is a subsequence of all subsequences $\{n^{(m)}_k\}_{k=1}^\infty$, hence,
    \[
    \lim_{m\to\infty}l_{n^{(m)}_1}(x_j)=a_j
    \]
    for all $x_j$'s. Define $l(x_j)=a_j$. 
    
    The linear extension of $l$ to $X_0$ is forced by its values on the basis elements,
    \[
    l\left(\sum_{j=1}^nc_jx_j\right):=\sum_{j=1}^nc_jl(x_j)=\sum_{j=1}^nc_j\lim_{k\to\infty}l_{n_k}(x_j)=\lim_{k\to\infty}l_{n_k}\left(\sum_{j=1}^nc_jx_j\right).
    \]
    As a consequence,
    \begin{eqnarray}\label{eqBnddOnXo}
        \left|l\left(\sum_{j=1}^nc_jx_j\right)\right|&=&\lim_{k\to\infty}\left|l_{n_k}\left(\sum_{j=1}^nc_jx_j\right)\right|\\
        &\le&\liminf_{k\to\infty}\|l_{n_k}\|_{X^\ast}\left\|\sum_{j=1}^nc_jx_j\right\|_X\\ 
        &\le& \left\|\sum_{j=1}^nc_jx_j\right\|_X.
    \end{eqnarray}
    The unique continuous extension of $l$ to $X=\overline{X_0}$ (see Exercise \ref{exeUniqueExt}) satisfies $\|l\|_{X^\ast}\le1$.

    We have to verify that $\lim_{k\to\infty}l_{n_k}(x)=l(x)$ for all $x$ in $X$. This is a simple $3\epsilon$ argument. For $x\in X$ and $x_0\in X_0$,
    \begin{eqnarray*}
        |l(x)-l_{n_k}(x)|&\le&|l(x)-l(x_0)|+|l(x_0)-l_{n_k}(x_0)|+|l_{n_k}(x)-l(x)|\\ 
        &\le&2\|x-x_0\|_X+|l(x_0)-l_{n_k}(x_0)|.
    \end{eqnarray*}
    For given $\epsilon>0$, choose $x_0$ with $\|x-x_0\|_X\le\epsilon$, then $k(\epsilon)>0$ such that for $k>k(\epsilon)$ one has $|l(x_0)-l_{n_k}(x_0)|<\epsilon$.   
\end{proof}
\section{Baire's Theorem and the uniform boundedness principle}\label{SectBaire}\index{theorem!Baire}
Let $(X,d)$ be a metric space. A subset $A$ of $X$ is {\it nowhere dense} in $X$ if $\overline{A}$ has empty interior (it does not contain nonempty open subsets).
\begin{theorem}\label{TheoBaire}[Baire's Theorem]
If $(X,d)$ is a complete metric space, and ${A_n:\ n\ge1}$ is a countable union of nowhere dense sets, then $\cup_{n=1}^\infty A_n\ne X$.
\end{theorem}
\begin{proof}
Since $\overline{A_1}$ does not contain a nonempty open set, we can find $\overline{B(x_1,r_1)}\subset X\setminus \overline{A_1}$ with $r_1<1/2$. 
Since $\overline{A_2}$ does not contain a nonempty open set, $B(x_1,r_1)\setminus\overline{A_2}$ is a nonempty open set, hence, it contains $\overline{B(x_2,r_2)}$ with $r_2<1/2^2$.
By iteration, we find $\overline{B(x_{n},r_{n})}\supseteq
\overline{B(x_{n-1},r_{n-1})}$ with $r_n<1/2^n$, and $\overline{B(x_{n},r_{n})}\cap\left(\overline{A_1}\cup\dots\overline{A_n}\right)=\emptyset$.

Since $X$ is complete, the intersection of the balls $\overline{B(x_{n},r_{n})}$ contains (a unique) point $z\in X\setminus\left(\cup_{n=1}^\infty A_n\right)$.
\end{proof}

The next, important theorem is due to Banach and Steinhaus.\index{theorem!Banach-Steinhaus, or uniform boundedness principle}
\begin{theorem}\label{TheoUBP} [Uniform Boundedness Principle] Let $\mathcal{F}=\{T\}$ be a family of bounded, linear operators $T:X\to Y$ from a Banach space $X$ to a normed, linear space $Y$. Suppose that, for each $x$ in $X$, $\|T x\|\le C(x)$ independent of $T\in \mathcal{F}$, although possibly dependent on $x$. Then, $\mathcal{F}$ is bounded,
\[
\|T\|\le C
\]
for some $C>0$ independent of $T\in\mathcal{F}$.
\end{theorem}
\begin{proof}
For $m\ge1$, let $A_m=\{x\in X:\ \sup\{\|T x\|:\ T\in\mathcal{F}\}\le m\}$. 
By hypothesis, $\cup_{m=1}^\infty A_m=X$, hence, some $A_m$ contains a closed open ball $\overline{B(z,r)}$. Since $z\in A_m$,
\[
\overline{B(0,r)}\subseteq A_m-z\subseteq A_m+A_m\subseteq A_{2m}.
\]
This is what we need, since, for $\|x\|\le1$ and $T\in\mathcal{F}$,
\begin{eqnarray*}
\|T x\|&=&\frac{1}{r}\|T rx\|\le\frac{2m}{r}.
\end{eqnarray*}
\end{proof}

Here is an application of Banach-Steinhaus.
\begin{corollary}\label{CorBSLimits}
Let $T_n\in\mathcal{L}(X,Y)$ be a sequence of bounded, linear maps from $X$ to $Y$, Banach spaces. If for all $x$ in $X$ there exists
\[
T x:=\lim_{n\to\infty}T_n x,
\]
then $T\in\mathcal{L}(X,Y)$ is bounded.
\end{corollary}
\begin{proof}
By hypothesis, $\{\|T_n x\|,\ n\ge1\}$ is bounded for each $x$, hence $\{\|T_n\|,\ n\ge1\}$ is bounded by some finite $C>0$. We have, then, if $\|x\|\le1$,
\[
\|T x\|=\lim_{n\to\infty}\|T_n x\|\le \limsup_{n\to\infty}\|T_n\|\le C.
\]
Thus, $\|T\|\le C$.
\end{proof}
Banach-Steinhaus allows us to transform (without control of the constants) weak {\it quantitative} information, into strong quantitative information. The following exercise provides an example.
\begin{exercise}\label{exeBSlip}
    Let $f:[a,b]\to X$ a function with values in a Banach space $X$, and suppose $f$ is {\bf weakly Lipschitz},
    \[
    |l(f(t+h))-l(f(t))|\le C(l)|h|
    \]
    for all $t, t+h\in [a,b]$ and $l\in X^\ast$. Show that $f$ is Lipschitz: there is a constant $C>0$ such that
    \[
    \|f(t+h)-f(t)\|\le C|h|,
    \]
\end{exercise}
We will see another important example in the next subsection, where we treat Banach space valued holomorphic functions.

{\it Qualitative} information, however, behaves differently. For instance, it is possible to exhibit functions $f:[0,1]\to\ell^2$ which are {\it weakly continuous} (for all $h$ in $\ell^2$, $x\mapsto\langle h,f(x)\rangle_{\ell^2}$ is continuous), but which are not continuous. The following example is modeled on the weak, but not strongly convergent sequence in \eqref{eqWnS}. Let
\begin{equation}\label{eqWCbutnotSC}
    [f(x)](n)=\begin{cases}
        \frac{1}{1+|n-1/x|^2}\text{ if }0<x\le 1,\\ 
        0\text{ if }x=0.
    \end{cases}
\end{equation}
Since for $x\ne0$ we have $\|f(x)\|_{\ell^2}^2\ge\sum_{n=0}^\infty\frac{1}{(1+n^2)^2}$, $f$ is not continuous at $x=0$. On the other hand,
\begin{equation}\label{eqWCholds}
    \lim_{x\to0}\langle h,f(x)\rangle_{\ell^2}=0.
\end{equation}
After changing variables to $y=1/x\to\infty$, and considering WLOG $h\ge0$, we can estimate (using an elementary estimate and Cauchy-Schwarz)
\begin{eqnarray*}
    \langle h,f(1/y)\rangle_{\ell^2}&=&\sum_{0\le n\le y/2}\frac{h(n)}{1+(n-y)^2}+\sum_{n> y/2}\frac{h(n)}{1+(n-y)^2}\\ 
    &\le&\frac{C}{y}\sum_{0\le n\le y/2}\frac{h(n)}{(1+(n-y)^2)^{1/2}}\\ 
    &+&2\left(\sum_{n>y/2}h(n)^2\right)^{1/2}\left(\sum_{n\ge0}\frac{1}{(1+n^2)^2}\right)^{1/2}\\ 
    &\le&\frac{C}{y}\|h\|_{\ell^2}\left(\sum_{n\ge0}\frac{1}{(1+n^2)^2}\right)^{1/2}\\ 
    &+&2\left(\sum_{n>y/2}h(n)^2\right)^{1/2}\left(\sum_{n\ge0}\frac{1}{(1+n^2)^2}\right)^{1/2}.
\end{eqnarray*}
The last expression tends to $0$ as $y\to\infty$ by monotone convergence. It is easy to see that $f$ is continuous on $(0,1]$ (for instance, by dominated convergence).
\subsection{Banach space-valued holomorphic functions}\label{SSectHolBanach}
Complex valued power series of a complex variable define holomorphic functions, which are at the heart of an elegant and powerful theory. At the root of it, is the fact that holomorphic functions can be defined in several, equivalent ways: through power series, through complex integrals, via the Cauchy-Riemann equations, and others. These different viewpoints can be translated to the world of functions with values in Banach spaces $X$; or even better in Banach spaces of the form ${\mathcal{L}}(X)$, where a product is part of the structure; or in Banach algebras, which is the most general Banach structure with products. This way, a number of results, tools, and techniques from holomorphic function theory become available to Functional Analysis, a fact of the greatest importance.

A function $f:{\mathbb C}\supseteq\Omega\to X$ defined from a region of the complex plane with values in a Banach space is {\it holomorphic} if
\[
f'(z):=\lim_{h\to0\text{ in }{\mathbb C}}\frac{f(z+h)-f(z)}{h}\in X
\]
exists for all $z$ in $\Omega$. The function $f$ is {\it weakly holomorphic} if for all $l$ in $X^\ast$ the function
\[
z\mapsto l(f(z))\in{\mathbb C}
\]
is holomorphic in the usual sense.
\begin{footnotesize}\begin{exercise}\label{exeVeryEasy}
Show that a holomorphic function $f:\Omega\to X$ is weakly holomorphic.\index{holomorphic!Banach space valued function}
\end{exercise}\end{footnotesize}
\begin{theorem}\label{TheoBanachSeries}
Let $f:\Omega\to X$ be a map defined from an open subset of ${\mathbb C}$ with values in a Banach space $X$. The following are equivalent.
\begin{enumerate}
    \item[(i)] The function $f$ is holomorphic in $\Omega$. 
    \item[(ii)] The function $f$ is weakly holomorphic in $\Omega$.
    \item[(iii)] For each $z_0$ in $\Omega$, there is $r>0$ such that, for $|z-z_0|<r$,
\[
f(z)=\sum_{n=0}a_n (z-z_0)^n,
\]
where $a_n\in X$ and the series converges absolutely uniformly for $|z-z_0|\le \rho< r$.
\end{enumerate}
The value of $r$ can be taken to be that of the larger disc centered at $z_0$ and contained in $\Omega$.

In fact, the series converges to an $X$-valued holomorphic function in $B(z_0,R)$ where
\begin{equation}\label{eqRoCbanach}
R=\frac{1}{\limsup_{n\to\infty}\|a_n\|^{1/n}}
\end{equation}
is the radius of the largest disc centered at $z_0$ and contained in $\Omega$.
\end{theorem}
The computational part of the proof is contained in the following.
\begin{lemma}\label{lemmaBanachSeries}
    Consider the power series 
    \begin{equation}\label{eqHolStarStarBis}
g(z)=\sum_{n=0}^\infty a_n (z-z_0)^n
\end{equation}
with coefficients $a_n\in X$, a Banach space, and with radius of convergence as in \eqref{eqRoCbanach}.

Then, $g:B(0,R)\to X$ is strongly holomorphic, and 
\begin{equation}\label{eqGderivative}
    g'(z)=\sum_{n=1}^\infty na_n (z-z_0)^{n-1},
\end{equation}
which has the same radius of convergence as $g$.

In particular, $g$ is infinitely differentiable, and
\begin{equation}\label{eqDerSeries}
    g^{(n)}(z_0)=n!a_n.
\end{equation}
\end{lemma}
\begin{proof}[Proof of the lemma]
    The usual proof from holomorphic theory makes use of tools we do not have, and we do not want to develop. We provide instead an XVIII century style proof which does not require them. We can suppose $z_0=0$. We start with the estimate
\begin{equation}\label{eqEstimateDerivativePower}
    \left|\frac{(z+h)^n-z^n}{h}-nz^{n-1}\right|\le|h|\frac{n(n-1)}{2}(|z|+|h|)^{n-2},
\end{equation}
which holds for $n\ge1$ and $z,h\in\mathbb{C}$.
The proof is just a calculation,
\begin{eqnarray*}
    \left|\frac{(z+h)^n-z^n}{h}-nz^{n-1}\right|&=&\left|\sum_{j=2}^{n}\binom{n}{j}z^{n-j}h^{j-1}\right|\crcr 
    &\le&|h|\sum_{l=0}^{n-2}\binom{n}{l+2}|z|^{n-2-l}|h|^{l}\crcr 
    &=&|h|\sum_{l=0}^{n-2}\frac{n(n-1)}{(l+2)(l+1)}\binom{n-2}{l}|z|^{n-2-l}|h|^{l}\crcr
    &\le&|h|\frac{n(n-1)}{2}\sum_{l=0}^{n-2}\binom{n-2}{l}|z|^{n-2-l}|h|^{l}\crcr
    &=&|h|\frac{n(n-1)}{2}(|z|+|h|)^{n-2}.
\end{eqnarray*}
 Let now
\[
    \psi(z)=\sum_{n=1}^\infty n a_n z^{n-1},
\]
which has the same radius of convergence.
For $|z|<R$ and $|h|<R-|z|$, we have:
\begin{eqnarray*}
    \left\|\frac{g(z+h)-g(z)}{h}-\psi(z)\right\|&=&\left\|\sum_{n=1}^\infty a_n\left(\frac{(z+h)^n-z^n}{h}-nz^{n-1}\right)\right\|\crcr 
    &\le&|h|\sum_{n=2}^\infty\frac{n(n-1)}{2}\|a_n\|(|z|+|h|)^{n-2}\crcr 
    &\to&0
\end{eqnarray*}
as $h\to0$ in $\mathbb{C}$, since the series in the line before the last converges. Hence, $\psi=g'$. Since $g'$ has the same radius of convergence as $g$, we can iterate the calculation,
\begin{equation}\label{eqPowSder}
    g^{(m)}(z)=\sum_{n=m}^\infty \frac{n!}{(n-m)!}a_n(z-z_0)^{n-m},\ g^{(m)}(z_0)=m!a_m.
\end{equation}
\end{proof}
\begin{proof}[Proof of the theorem] (i) implies (ii) by exercise \ref{exeVeryEasy}, and (iii) implies (i) is the lemma above. 

We show that (ii) implies (iii). We fix some notation. If $x\in X$, $\widehat{x}\in X^{\ast\ast}$ is the functional on $X^\ast$ for which $\widehat{x}(l)=l(x)$. We set $\widehat{X}\subseteq X^{\ast\ast}$ the set of such functionals.

We start by showing that if $f:\Omega\to X$ a weakly holomorphic function, and $\gamma$ is a closed curve in $\Omega$, then $z\mapsto f(z)$ is bounded on $\gamma$. In fact, and for $z$ on $\gamma$:
\begin{eqnarray*}
    |\widehat{f(z)}(l)|&=&|l(f(z))|\crcr 
    &\le&\sup_{z\in\gamma}|l(f(z))|\crcr 
    &=&C(l),
\end{eqnarray*}
which is finite because $z\mapsto l(f(z))$ is continuous on $\gamma$. By Banach-Steinhaus theorem,
\[
\sup_{z\in\gamma}\|f(z)\|_X=\sup_{z\in\gamma}\|\widehat{f(z)}\|_{X^{\ast\ast}}<\infty.
\]
Let $B(z_0,r)$ be a disc contained in $\Omega$, and let $\gamma$ be a circle of radius $\rho<r$ centered at $z_0$ and contained in $\Omega$. Let again $l$ be in $X^\ast$.  Then, $z\mapsto l(f(z))$ can be expanded as a power series,
\[
l(f(z))=\sum_{n=0}^\infty a_n(l) (z-z_0)^n,
\]
where
\[
a_n(l)=\frac{1}{2\pi i}\int_\gamma \frac{l(f(z))}{(z-z_0)^{n+1}}d z.
\]
The functional $l\mapsto a_n(l)$ is linear and
\[
|a_n(l)|\le \sup_{z\in\gamma}\|f(z)\|_X\frac{\|l\|_{X^\ast}}{\rho^{n}}.
\]
That is, $a_n\in X^{\ast\ast}$ with 
\begin{equation}\label{eqHolRadiusConv}
\|a_n\|_{X^{\ast\ast}}\le \frac{\sup_{z\in\gamma}\|f(z)\|_X}{\rho^n}.
\end{equation}
Thus,
\begin{equation}\label{eqHolStarStar}
g(z)=\sum_{n=0}^\infty a_n (z-z_0)^n
\end{equation}
converges for $|z-z_0|<\rho$ to a function $g$ with values in $X^{\ast\ast}$. You just have to check that the usual proof in holomorphic function theory works when you have power series with coefficients in a Banach space (geometric series only are involved). The Hadamard formula for the radius of convergence of the series is proved like in holomorphic theory. By the lemma, $g$ is an $X^{\ast\ast}$ valued holomorphic function. We next prove that $g(z)=\widehat{f(z)}$ is the image in $X^{\ast\ast}$ of $\widehat{f(z)}$, i.e. that $[g(z)](l)=l(f(z))$ for all $l\in X^\ast$.

In fact, all $l$ in $X^\ast$ we have
\[
[g(z)](l)=\sum_{n=0}^\infty a_n(l) (z-z_0)^n=l(f(z)).
\]
By the lemma, $g$ is infinitely differentiable and
\[
a_n=\frac{g^{(n)}(z_0)}{n!}=\frac{\widehat{f^{(n)}(z_0)}}{n!}\in X,
\]
because $X$ is closed in $X^{\ast\ast}$, and the derivatives, which are limits of $X$ valued functions, belong to $X$. Thus,
$a_n=\widehat{\alpha_n}$, with $\alpha_n\in X$, and 
\[
f(z)=\sum_{n=0}^\infty \alpha_n (z-z_0)^n,
\]
as wished.
\end{proof}
\begin{footnotesize}\begin{exercise}\label{exeHolBanachProduct}
    Let $\Omega\subseteq\mathbb{C}$ be open, and let $\lambda:\Omega\to\mathbb{C}$ be holomorphic, and $f:\Omega\to X$ be a Banach space-valued holomorphic function. Show that their product $\lambda\cdot f:\Omega\to X$ is holomorphic.
\end{exercise}\end{footnotesize}

\section{The Open Mapping Theorem and the Closed Graph Theorem}\label{SectOMandCG}
A map $F:M\to N$ between metric spaces is {\it open} if the image of an open set in $M$ is open in $N$. 
\begin{footnotesize}\begin{exercise}\label{exeOpenEquiv}
If $T:X\to Y$ is a linear map between normed, linear spaces, then the following are equivalent:
\begin{enumerate}
    \item[(i)] $T$ is open;
    \item[(ii)] there is a ball $B_X(0,r)$ in $X$ such that $T(B_X(0,r))$ contains a ball in $Y$;
    \item[(iii)] there is a ball $B_X(0,r)$ in $X$ such that $T(B_X(0,r))$ contains a ball centered at $0$ in $Y$.
\end{enumerate}
Moreover, if $T$ is open, then it is onto.

\noindent{\bf Hint.} $(i)\implies(ii)$ is clear. For $(ii)\implies(iii)$, you can show that if $T(B(0,r))\supset B(y,R)$, then $T(B(0,2r))\supset B(0,R)$. The proof that $(iii)\implies(i)$ is easy. Also, (iii) implies that $T$ is onto by the homogeneity of $T$.
\end{exercise}\end{footnotesize}

\begin{theorem}\label{TheoOpenMapping}[Open Mapping Theorem]\index{theorem!open mapping}
Let $X,Y$ be Banach spaces, and let $T:X\to Y$ be a linear, bounded map from $X$ onto $Y$. Then, $T$ is open.
\end{theorem}
\begin{proof}
Let $B_n:=B(0,n)\subset X$. Since $\cup_n T(B_n)=Y$, there is $n$ such that $\overline{T(B_n)}\supset B(y_0,\epsilon)$. If $y\in Y$, then $y=(y+y_0)- y_0\in \overline{T(B_n)}-\overline{T(B_n)} \subset \overline{T(B_{2n})}$ provided $\|y\|\epsilon$. i.e. $B(0,\epsilon/(2n))=B(0,\eta)\subseteq \overline{T(B_{1})}$.

It suffices then to show that $\overline{T(B_1)}\subset T(B_2)$.

Let $y\in \overline{T(B_1)}$, and pick $x_1\in B_1$ such that $\|y-T x_1\|<\eta/2$, so that $T x_1-y\in B(0,\eta/2)\subseteq \overline{T(B_{1/2})}$. We can then find
pick $x_2\in B_{1/2}$ such that $\|y-T x_1-T x_2\|<\eta/2^2$.

Iterating, we have $x_n\in B_{1/2^{n-1}}$ such that $\|y-T(x_1+\dots+x_n)\|<\eta/2^n$. We have that $x_1+1dots+x_n\to x\in B_2$ as $n\to\infty$, and $T(x)=\lim_{n\to\infty} T(x_1+\dots+x_n)=y$, as wished.
\end{proof}
\begin{footnotesize}\begin{exercise} Let $X={\mathbb R}^2$ with the Euclidean metric. Show that for all $\epsilon>0$ there is a linear bijection $T:X\to X$ such that $\|T\|=1$, yet $T(B(0,1))$ does not contain $B(0,\epsilon)$. 
That is, the Open Mapping Theorem is not quantitative. 
\end{exercise}\end{footnotesize}

\begin{theorem}\label{TheoInVt}[Inverse Mapping Theorem]\index{theorem!inverse mapping}
Let $T:X\to Y$ be a bounded bijection of Banach spaces. Then, $T^{-1}:Y\to X$ is bounded.
\end{theorem}
\begin{proof}
The fact that $T$ is open means that for some $\epsilon>0$, if $\|y\|<\epsilon$ in $Y$, then there is $x\in X$ with $\|x\|<1$ and $T x=y$, i.e. $x=T^{-1}y$. Said it differently, $T^{-1}$ maps $B(o,\epsilon)$ into $B(0,1)$, so that $\|T^{-1}\|\le 1/\epsilon$.
\end{proof}

\begin{footnotesize}\begin{exercise}
Show that if $\|\cdot\|_1$ and $\|\cdot\|_2$ are two norms on $X$, if $\|x\|_2\le\|x\|_1$ on $X$ and $X$ is Banach with respect to $\|\cdot\|_1$, then there is $C>0$ such that $\|x\|_1\le C\|x\|_2$ on $X$. 
\end{exercise}\end{footnotesize}
\begin{footnotesize}\begin{exercise}
Show that (even!) in ${\mathbb R}^2$ one can find norms $\|\cdot\|_1$ and $\|\cdot\|_2$ with $\|x\|_2\le\|x\|_1$, yet 
$\|x\|_1\le C\|x\|_2$ only holds if $C$ is (arbitrarily) large. That is, the result in the previous exercise is not quantitative.
\end{exercise}\end{footnotesize}

The {\it graph} of a function $f:M\to N$ is the set $\Gamma(f)=\{(x,y\in M\times N:\ y=f(x))\}$. If $T:X\to Y$ is a linear operator between linear spaces, then $\Gamma(T)$ is a linear subspace of $X\times Y$. If $X$ and $Y$ are normed, then $\|(x,y)\|:=\|x\|+\|y\|$ defines a norm on $X\times Y$, hence on $\Gamma(T)$.

\begin{footnotesize}\begin{exercise}
If $X$ and $Y$ are normed linear space and $T\in\mathcal{L}(X,Y)$, then $\Gamma(T)$ is closed.
\end{exercise}\end{footnotesize}

\begin{theorem}\label{TheoCGT}[Closed Graph Theorem] \index{theorem!closed graph}
Let $T:X\to Y$ be a linear operator defined from a Banach space $X$ to a Banach space $Y$. If $\Gamma(T)$ is closed, then $T$ is bounded.
\end{theorem}
\begin{proof}
By assumption, $\Gamma(T)$ is closed in a Banach space, hence Banach itself: if $(x_n,T x_n)\to(x,y)$ is Cauchy, then $y=T x$, so $(x,y)\in\Gamma(T)$.

Consider the projections $\pi_X:(x,T x)\mapsto x$ and $\pi_Y:(x,T x)\mapsto T x$. Both projections are bounded, $\pi_X$ is invertible (hence, its inverse is continuous), and $T=\pi_Y\circ \pi_X^{-1}$. Hence, $T$ is bounded.
\end{proof}

\section{Integrals of continuous, Banach space valued functions}\label{SectIntValX}\index{Banach space!valued functions: integral}
The integral of functions with values on a Banach space can be defined in several ways, choosing which depends on the application we have in mind\footnote{There are different definitions for the integral of a Banach space valued function. A reasonably general, Lebesgue style one, is provided by {\it Bochner integrals}. See e.g. \href{http://home.ustc.edu.cn/~anprin/Bochnerintegral.pdf}{The Bochner integral} by Wenjing Wu. The weak version of the Bochner integral is the {\it Pettis integral}. Here, however, we are integrating continuous functions, and more elementary definitions of integral can be used.}. Think of Banach space valued holomorphic functions: extending to them the notion of Cauchy integrals requires integrating vector valued functions. Having in mind holomorphic theory, we sketch here a construction of the integral for a continuous function $f:[a,b]\to X$, where $X$ is Banach. The definition of integral we discuss here does not require any foundational result concerning Banach spaces. We only make use of linearity and completeness.

Let $f:[a,b]\to X$ be a continuous function with values in the Banach space $X$. Since $[a,b]$ is compact, $f$ is uniformly continuous. For fixed $n\ge1$, let
\begin{equation}\label{eqIntValXn}
    S_n(f)=\sum_{j=1}^{2^n}f\left(\frac{j}{2^n}\right)\frac{b-a}{2^n}\in X.
\end{equation}
The sequence $\{S_n(f)\}_{n=1}^\infty$ is Cauchy in $X$. The calculation is similar to the one we met when defining the Lebesgue measure. For any fixed $\epsilon>0$, there exists $\delta>0$ such that $\|f(s)-f(t)\|\le\epsilon$ if $|s-t|\le\delta$. If $1/2^n\le\delta$, then
\begin{eqnarray*}
    \|S_{n+m}(f)-S_n(f)\|&=&\left\|\sum_{j=1}^{2^n}\frac{b-a}{2^n}\frac{1}{2^m}\sum_{l=1}^{2^m}\left[f\left(\frac{j}{2^n}\right)-f\left(\frac{j-1}{2^n}+\frac{l}{2^{n+m}}\right)\right]\right\|\crcr 
    &\le&\sum_{j=1}^{2^n}\frac{b-a}{2^n}\frac{1}{2^m}\sum_{l=1}^{2^m}\left\|f\left(\frac{j}{2^n}\right)-f\left(\frac{j-1}{2^n}+\frac{l}{2^{n+m}}\right)\right\|\crcr 
    &\le&(b-a)\epsilon.
\end{eqnarray*}
By definition,
\begin{equation}\label{eqIntValXdef}
    \int_a^bf(t)dt=\lim_{n\to\infty}S_n(f).
\end{equation}
You might complain that this definition depends on the choice of very special partitions of $[a,b]$ (this makes additivity of the integral problematic, for instance). Fortunately, the dependence is apparent. Given a partition $a=t_0<t_1<\dots<t_{m-1}<t_m=b$ of $[a,b]$ and sampling points $t_{j-1}\le t_j^\ast\le t_j$, let
\[
S(\{t_j\}_{j=0}^m,\{t_j^\ast\}_{j=1}^m)=\sum_{j=1}^mf(t_j^\ast)(t_j-t_{j-1}).
\]
The {\it resolution} of the partition is $\min_{j=1,\dots,m}(t_j-t_{j-1})=\delta\left(\{t_j\}_{j=0}^m\right)$.
\begin{proposition}\label{propIntValX}
Let $X$ be Banach and $f:[a,b]\to X$ continuous. Then, for all $\epsilon>0$ there is $\delta>0$ such that if $\delta\left(\{t_j\}_{j=0}^m\right)\le\delta$, then $\left\|\int_a^bf(t)dt-S(\{t_j\}_{j=0}^m,\{t_j^\ast\}_{j=1}^m)\right\|\le\epsilon$.
\end{proposition}
\begin{proof}
    Let $n\ge1$ be such that $\left\|\int_a^bf(t)dt-S_n(f)\right\|\le\epsilon$ and that $\|f(s)-f(t)\|\le \epsilon$ if $|s-t|\le 1/2^{n}$. Consider a partition $\{t_j\}_{j=0}^m$ with resolution less than $\delta=1/2^n$, and a sampling set $\{t_j^\ast\}_{j=1}^m$ as above. 
    Below, we denote by $|I|$ the Lebesgue measure of an interval $I$.
    \begin{eqnarray*}
       &\ & \left\|S_n(f)-S(\{t_j\}_{j=0}^m,\{t_j^\ast\}_{j=1}^m)\right\|\crcr 
       &=&\left\|\sum_{l=1}^{2^n}\sum_{j:[t_{j-1},t_j]\cap[(l-1)/2^n,l/2^n]\ne\emptyset}(f(l/2^n)-f(t_j^\ast))|[t_{j-1},t_j]\cap[(l-1)/2^n,l/2^n]|\right\|\crcr 
        &\le&\sum_{l=1}^{2^n}\sum_{j:[t_{j-1},t_j]\cap[(l-1)/2^n,l/2^n]\ne\emptyset}\left\|f(l/2^n)-f(t_j^\ast)\right\|\cdot|[t_{j-1},t_j]\cap[(l-1)/2^n,l/2^n]|\crcr 
        &\le&(b-a)\epsilon.
    \end{eqnarray*}
\end{proof}
\begin{exercise}\label{exeIntLSvalX}
    Let $\alpha:[a,b]\to\mathbb{R}$ be increasing (or, more generally, let $\alpha:[a,b]\to\mathbb{C}$ be of bounded variation). For $f:[a,b]\to X$ continuous (where $X$ is Banach) provide a definition of the vector valued Stjelties integral
    \[
    \int_a^bf(t)d\alpha(t),
    \]
    and show that it is well defined.
\end{exercise}

\begin{corollary}\label{corIntValX}
    Let $f:[a,b]\to X$ be a continuous function with values in a Banach space $X$, and let $T:X\to Y$ be a bounded operator between $X$ and another Banach space $Y$. Then,
    \[
    T\left(\int_a^bf(t)dt\right)=\int_a^b T(f(t))dt.
    \]
\end{corollary}
\begin{proof}
    It suffices to pass in the limit for $n\to\infty$ the two sides of the equality
    \[
    T\left(\sum_{j=1}^mf(t_j^\ast)(t_j-t_{j-1})\right)=\sum_{j=1}^mT\left[f(t_j^\ast)\right](t_j-t_{j-1}).
    \]
\end{proof}
We give a sample application to Banach space valued holomorphic functions. The reader might find it interesting and rewarding to translate in the Banach valued world the chapter of a book of complex analysis concerning complex integrals.
\begin{corollary}\label{corMoreraBanach}[Morera's theorem for Banach space valued functions] Let $f:\Omega\to X$ be a continuous, Banach space valued function. Then $f$ is holomorphic if and only if $\int_{\partial T}f(z)dz=0$ for all triangles $T$ contained in $\Omega$.
\end{corollary}
\begin{proof}
    The function $f$ is holomorphic if and only if it is weakly holomorphic, and by Morera's theorem this holds if and only if
    \[
    0=\int_{\partial T}l(f(z))dz=l\left(\int_{\partial T}f(z)dz\right)
    \]
    for all triangles and $l\in X^\ast$, the second equality following from corollary \ref{corIntValX}. Thus, $\int_{\partial T}f(z)dz=0$.
\end{proof}
\section{Sesquilinear forms on Hilbert spaces}\label{SectHilbApp}\index{sesquilinear form}
In this section we consider a few applications of the Banach theory considered above to Hilbert space theory.

Let $H$ be a Hilbert space. A map $[\cdot,\cdot]:H\times H\to \mathbb{C}$, $[x,y]=\langle x,y\rangle$, is a {\it sesquilinear form} if it is conjugate linear in the first component and linear in the second.
\[
[x,\lambda y+\mu z]=\lambda[ x,y]+\mu[x,z],\ [\lambda y+\mu z,x]=\overline{\lambda}[ y,x]+\overline{\mu}[z,x].
\]
For instance, the inner product $\langle\cdot,\cdot\rangle$ is sesquilinear.

Sesquilinear forms are determined by their action on the diagonal.
\begin{proposition}\label{propPolarization}[Polarization identity\index{polarization identity!sesquilinear forms}] Let $H$ be a complex Hilbert space. If $[\cdot,\cdot]$ is a sesquilinear form, then the expression $[x,y]$ can be expressed as a linear combination of expressions of the form $[z,z]$:
    \begin{equation}\label{eqSesqui}
       [x,y]=\frac{1}{4}([x+y,x+y]-[x-y,x-y]-i[x+iy,x+iy]+i[x-iy,x-iy]). 
    \end{equation}
\end{proposition}
The proof consists in working out the right hand side of \eqref{eqSesqui}. The polarization identity holds in Hilbert spaces over $\mathbb{R}$ under the extra hypothesis that $[x,y]=[y,x]$ is commutative. The form
\[
\left[\begin{pmatrix}x\\ y
\end{pmatrix},\begin{pmatrix}u\\ v
\end{pmatrix}\right]=(x,y)\begin{pmatrix}1&1\\ 0&1
\end{pmatrix}\begin{pmatrix}u\\ v
\end{pmatrix}=xu+xv+yv
\]
provides an example of a bi-linear form on $\mathbb{R}^2$ for which the polarization identity does not hold.

A sesquilinear form $H\times H\xrightarrow{[\cdot,\cdot]}\mathbb{C}$ is {\it bounded} if there is a positive constant $C$ such that
\[
|[x,y]|\le C\|x\|\cdot\|y\|
\]
for all $x,y$ in $H$. Using homogeneity in each component, this is equivalent to requiring that there is a (possibly different) constant $C>0$ such that
\[
\sup_{\|x\|,\|y\|\le1}|[x,y]|\le C.
\]
All bounded sesquilinear forms are obtained by "sandwiching" a bounded linear operator between two vectors. 
\begin{proposition}\label{propSesquiBLO}
    \begin{itemize}
        \item[(i)] If $A\in\mathcal{L}(H)$, then
\begin{equation}
    [x,y]:=\langle x,Ay\rangle
\end{equation}
is a bounded, sesquilinear form: $|[x,y]|\le\|A\|\|x\|\|y\|$.
        \item[(ii)] Viceversa, for any bounded sesquilinear form $[\cdot,\cdot]$ there exists $A\in\mathcal{L}(H)$.
    \end{itemize}
\end{proposition}
\begin{proof}
    (i) is clear. (ii) The functional $y\mapsto[x,y]$ is bounded, then, by Riesz theorem for Hilbert spaces, there is $Ax\in H$ such that $[x,y]=\langle Ax,y\rangle$. It is easy to show that $x\mapsto Ax$ is a linear map from $H$ to $H$. Clearly $A$ is bounded,
    \[
    \|A\|=\sup_{\|y\|\le1}\sup_{\|x\|\le1}|[x,y]|.
    \]
\end{proof}
Boundedness can be checked separately for each variable.
\begin{proposition}\label{propSesquiOperator} Let $[\cdot,\cdot]$ be a sesquilinear form on $H$.
 If $x\mapsto [x,y]$ is a bounded conjugate-linear functional for each $y$, and $y\mapsto [x,y]$ is a bounded linear functional for each $x$, then $[\cdot,\cdot]$ is bounded.
\end{proposition}
\begin{proof}
    Since $x\mapsto\overline{[x,y]}=\lambda_y(x)$ is a bounded linear functional, there is $Ay\in H$ such that $[x,y]=\langle x,Ay\rangle$. The map $A$ is clearly linear, $A:H\to H$, and $\|\lambda_y\|_{H^\ast}=\|A y\|$. For each $x$ we have that
    \[
    \sup_{\|y\|\le1}|[x,y]|\le C(x),
    \]
    a constant depending on $x$. By Banach-Steinhaus,
    \[
    \infty>\sup_{\|y\|\le1}\|\lambda_y\|_{H^\ast}=\sup_{\|y\|\le1}\|A y\|=\|A\|,
    \]
    as wished.
\end{proof}

\vskip0.5cm

The sesquilinear form $[\cdot,\cdot]$ is {\it symmetric} if $[y,x]=\overline{[x,y]}$. If there is a linear operator $A:H\to H$ such that $[x,y]=[x,Ay]$, symmetry translates into:
\[
\langle x,Ay\rangle=[x,y]=\overline{[y,x]}=\overline{\langle y,Ax\rangle}=\langle Ax,y\rangle. 
\]
That is, $A=A^\ast$ is self-adjoint.

It is remarkable that all such symmetric sesquilinear forms on a Hilbert space are bounded.\index{theorem!Hellinger-Toeplitz}
\begin{theorem}\label{TheoHellToep}[Hellinger-Toeplitz] Let $A:H\to H$ be an everywhere defined self-adjoint operator on a Hilbert space $H$:
\[
\langle A x , y\rangle=\langle x , A y\rangle\ {\text if }\ x,y\in H.
\]
Then, $A$ is bounded.
\end{theorem}
\begin{proof}
Let $(x_n,A x_n)$ be a sequence in $\Gamma(A)$ with $x_n\to x$ and $A x_n\to y$. For all $z\in H$:
\begin{eqnarray*}
\langle z,y \rangle&=&\lim_{n\to\infty}\langle z , A x_n \rangle=\lim_{n\to\infty}\langle A z , x_n \rangle\crcr
&=&\langle A z , x \rangle=\lim_{n\to\infty}\langle z , A x \rangle,
\end{eqnarray*}
hence, $A x=y$, so $\Gamma(A)$ is closed, and we can use the closed graph theorem.
\end{proof}

A sesquilinear $[\cdot,\cdot]$ form is {\it coercive} if
\begin{equation}\label{eqCoerciveDefinition}
    [x,x]\ge c\|x\|^2
\end{equation}
for some $c>0$.

\index{theorem!Lax-Milgram}
\begin{theorem}[Lax-Milgram]\label{theoLaxMilgram}
    Let $[\cdot,\cdot]$ be a bounded sesquilinear form with bound $C$, which is coercive with constant $c>0$. Then, for all $w$ in $H$ there is a unique $x$ such that:
    \begin{equation}\label{eqLaxMilgram}
        [x,y]=\langle w,y\rangle. 
    \end{equation}
\end{theorem}
\begin{proof}
    We saw earlier that $[x,y]=\langle Ax,y\rangle$, with $A\in\mathcal{L}(H)$. We have to prove that $A$ is a bijection.
    \begin{enumerate}
        \item[(i)] {\bf $A$ is injective}. In fact, 
        \[
        c\|x\|^2\le [x,x]=\langle Ax,x\rangle\le \|x\|\cdot\|Ax\|,
        \]
        which imples that $x=0$ when $Ax=0$.
        \item[(ii)] {\bf $\text{Ran}(A)$ is closed in $H$}. Let $Ax_n\to w$ in $H$. In particular, $\{Ax_n\}$ is Cauchy and, using coercivity as in (i),
        \[
        \|x_n-x_{n+j}\|\le 1/c\|Ax_n-Ax_{n+j}\|,
        \]
        so that $\{x_n\}$ is Cauchy, too, hence $x_n\to x$. Since $A$ is bounded, $w=\lim_n Ax_n=Ax$.
        \item[(iii)] {\bf $\text{Ran}(A)=H$}. By (ii) and the general theory, we can write $H=\text{Ran}(A)\oplus(H\ominus\text{Ran}(A))$. Let $k\in H\ominus\text{Ran}(A)$:
        \[
        c\|w\|^2\le [w,w]=\langle Aw|w\rangle=0, 
        \]
        thus $w=0$. Hence, $\text{Ran}(A)=H$.
    \end{enumerate}
\end{proof}

\section{The Hilbert space adjoint of an operator}\label{SectAdjoint}
In this section we define the Hilbert space adjoint of a bounded operator, and prove some of its properties. This material is not prerequisite for the chapter on spectral theory.
\subsection{Basic properties of the adjoint}\label{SSectBasicAdjointProp}
Let $T\in\mathcal{L}(H)$. Its {\it adjoint} $T^\ast\in\mathcal{L}(H)$ is defined on $y\in H$ by:
\begin{equation}\label{eqAdjointOperator}
    \langle x,T^\ast y\rangle=\langle Tx,y\rangle.
\end{equation}
Equation \eqref{eqAdjointOperator} defines a linear operator. For $x,y,z\in H$ and $a,b\in\mathbb{C}$,
\begin{eqnarray*}
    \langle x,T^\ast (ay+bz)\rangle&=&\langle Tx, ay+bz\rangle\crcr 
    &=&a\langle Tx, y\rangle+b\langle Tx,z\rangle\crcr 
    &=&a\langle x, T^\ast y\rangle+b\langle x,T^\ast z\rangle\crcr
    &=&\langle x,aT^\ast y+bT^\ast z\rangle.
\end{eqnarray*}
The map $T\mapsto T^\ast$ is conjugate linear,
\begin{eqnarray*}
    \langle x,(aS+bT)^\ast y\rangle&=&\langle (aS+bT)x, y\rangle\crcr 
    &=&\Bar{a}\langle Sx, y\rangle+\Bar{b}\langle Tx, y\rangle\crcr
    &=&\langle x,(\Bar{a}S^\ast+\Bar{b}T^\ast) y\rangle.
\end{eqnarray*}
The operator $T$ is {\it self-adjoint} if $T^\ast=T$, and, more generally, it is {\it normal} if $TT^\ast=T^\ast T$.
\begin{theorem}\label{theoAdjointBasics} Let $S$, $T$ be bounded operators on $H$.
    \begin{enumerate}
        \item[(i)] $(T^\ast)^\ast=T$.
        \item[(ii)] $\|T^\ast\|=\|T\|$.
        \item[(iii)] $(TS)^\ast=S^\ast T^\ast$.
        \item[(iv)] $T$ is invertible if and only if $T^\ast$ is invertible, and $(T^{-1})^\ast=(T^\ast)^{-1}$.
        \item[(v)] $T\mapsto T^\ast$ is continuous with respect to the uniform and the weak operator topology (but not with respect to the strong operator topology). 
        \item[(vi)] $\|T^\ast T\|=\|T\|^2$.
    \end{enumerate}
\end{theorem}
\begin{proof}
    \begin{enumerate}
        \item[(i)]We have 
        \[\langle x, (T^\ast)^\ast y\rangle=\langle T^\ast x, y\rangle=\overline{\langle y,T^\ast x\rangle}=\overline{\langle Ty, x\rangle}=\langle x, T y\rangle. \]
        \item[(ii)] We have
        \begin{eqnarray*}
            \|T^\ast\|&=&\sup_{\|y\|=1}\|Ty\|=\sup_{\|x\|=1}\sup_{\|y\|=1}|\langle x,T^\ast y\rangle|\crcr 
            &=&\sup_{\|x\|=1}\sup_{\|y\|=1}|\langle Tx, y\rangle|\crcr 
            &\le&\|T\|,
        \end{eqnarray*}
        showing that $\|T^\ast\|\le\|T\|$, and so, by (i), $\|T\|=\|(T^\ast)^\ast\|\le \|T^\ast\|$.
        \item[(iii)] Simple algebra.
        \item[(iv)] By definition,
        \[
            \langle x, (T^{-1})^\ast T^\ast y\rangle=\langle T T^{-1}x,  y\rangle =\langle x,  y\rangle=\langle T^{-1}Tx,  y\rangle=\langle x,  T^\ast (T^{-1})^\ast y\rangle.
        \]
        \item[(v)] For the uniform topology, see (ii). For the weak operator topology,
        \[
            \langle x,(T_n^\ast-T^\ast) y\rangle=\langle (T_n-T)x, y\rangle\to0
        \]
        if $T_n\xrightarrow[WOT]{} T$.
        \item[(vi)] The inequality $\|T^\ast T\|\le\|T\|^2$ follows from (ii). In the other direction,
        \[
        \|Tx\|^2=\langle Tx,Tx\rangle=\langle x,T^\ast Tx\rangle\le \|T^\ast T\|\cdot\|x\|^2.
        \]
    \end{enumerate}
\end{proof}
About the negative statement on SOT in (iii), see the examples after the definition of the operator topologies.
\begin{corollary}\label{corSelfAdjointSpectralRadius}
\begin{enumerate}
    \item[(i)] If $S$ and $T$ are self-adjoint and commute, $ST=TS$, then $ST$ is self-adjoint.
    \item[(ii)] If $T$ is self-adjoint, then $\|T^2\|=\|T\|^2$, and $r(T)$, the spectral radius of $T$, is
\begin{equation}\label{eqSelfAdjointSpectralRadius}
    r(T)=\|T\|.
\end{equation}
\end{enumerate}
\end{corollary}
\begin{proof}
(i) is a simple algebraic calculation,
\begin{eqnarray*}
   \langle x,TSy\rangle&=& \langle x,STy\rangle=\langle T^\ast S^\ast x,y\rangle\crcr 
   &=&\langle TS x,y\rangle,
\end{eqnarray*}
shows that $(TS)^\ast=TS$.

    The first assertion in (ii) follows from (vi) in the previous theorem, so that
    \[
    r(T)=\lim_{n\to\infty}\|T^n\|^{1/n}=\lim_{m\to\infty}\|T^{2^m}\|^{1/2^m}=\|T\|,
    \]
    since $\|T^{2^m}\|=\|T^{2^{m-1}}\|^2$.
\end{proof}
\subsection{Spectrum and adjoint}\label{SSectSAspectra}
\begin{lemma}\label{lemmaAdjointSpectrum}
    The spectrum of $T^\ast$ is $\sigma(T^\ast)=\overline{\sigma(T)}$.
\end{lemma}
\begin{proof}
    For complex $\lambda$, $\lambda I-T$ is invertible if and only if $(\lambda I-T)^\ast=\overline{\lambda}I-T^\ast$ is invertible.
\end{proof}
When we look more closely inside the spectrum, the picture is more complex. For instance, the shift $\tau_1$ on $\ell^2(\mathbb{N})$ does not have eigenvalues, while any $\lambda\in\mathbb{D}$ is an eigenvalue for the the back-shift $\tau_1^\ast$. We split the spectrum of $T$ into three components.
\begin{enumerate}
    \item[(i)] the {\it point spectrum} $\sigma_p(T)$ is the set of the {\it eigenvalues} $\lambda$ of $T$.
    \begin{equation}\label{eqEigenvalue}
        Tv=\lambda v
    \end{equation}
    for some $v\ne0$ in $H$. Any such $v$ is called an {\it eigenvector} relative to $\lambda$. The closed, linear subspace $V_\lambda$ of the eigenvectors relative to $\lambda$ is the {\it eigenspace} relative to $\lambda$. In particular, if $\ker(T)\ne0$, $0$ is an eigenvalue with eigenspace $\ker(T)$. 
    \item[(ii)] The {\it continuous spectrum} $\sigma_c(T)$ is the set of those $\lambda\in\sigma(T)$ such that $\lambda I-T$ has dense range in $H$.
    \item[(iii)] The {\it residual spectrum} $\sigma_{res}(T)$ of $T$ is the set of those $\lambda$ which are not eigenvalues, and such that $\lambda I-T$ does not have dense range in $H$.
\end{enumerate}
Items (ii) and (iii) only happen in infinite dimension. When $H=\mathbb{C}^N$, the spectrum is the set of the eigenvalues, and if $\lambda$ is an eigenvalue of $A$, then $\dim(\text{Ran}(\lambda I-A))=N-\dim(\ker(\lambda I-A))<N$, hence the range of $\lambda I-A$ is not dense.

\begin{proposition}\label{propSpectraDCR}
    \begin{enumerate}
        \item[(i)] If $\lambda\in\sigma_{res}(T)$, then $\overline{\lambda}\in\sigma_d(T)$.
        \item[(ii)] If $\lambda\in\sigma_{p}(T)$, then $\overline{\lambda}\in\sigma_d(T)\cup\sigma_{res}(T)$.
    \end{enumerate}
\end{proposition}
\begin{proof}
    \noindent (i) Suppose $\lambda\in\sigma_{res}(T)$ and consider $0\ne k\in H\ominus\text{Cl}(\text{Ran}(\lambda I-T))$. For all $x$ in $H$,
    \[
    \langle x,(\overline{\lambda} I-T^\ast)k\rangle= \langle (\lambda I-T)x,k\rangle=0,
    \]
    hence, $T^\ast k=\overline{\lambda}k$.

    \noindent (ii) Let $k\ne 0$ be such that $Tk=\lambda k$. Then, for all $x$ in $H$.
    \[
    \langle (\overline{\lambda} I-T^\ast)x,k\rangle=\langle x,(\lambda I-T)k\rangle=0,
    \]
    showing that $k$ is orthogonal to $\text{Ran}(\overline{\lambda} I-T^\ast)$, which is then not dense in $H$.
\end{proof}
Consider for instance the spectra of the shift $M_z$ on $H^2(\mathbb{D})$, and of its adjoint $M_z^\ast$,
\[
M_z^\ast f(z)=\frac{f(z)-f(0)}{z}.
\]
If $\lambda\in\mathbb{D}$, then it is an eigenvalue of $M_z^\ast$. Since $\sigma_p(M_z)=\emptyset$, by (ii) $\overline{\lambda}\in\sigma_{res}(M_z)$. 

If $|\lambda|=1$ and $\lambda\in \sigma_{res}(M_z)$, by (i) we would have $\overline{\lambda}\in \sigma_p(M_z^\ast)$, which is not true because the non-vanishing solutions of
\[
\frac{f(z)-f(0)}{z}=\overline{\lambda}f(z)
\]
are multiples of $f(z)=\frac{1}{1-\overline{\lambda}z}=\sum_{n=0}^\infty\overline{\lambda}^nz^n$, which does not belong to $H^2$. Hence, $\lambda\in \sigma_c(M_z)$.

Again by (i), if $|\lambda|=1$, $\lambda\notin \sigma_{res}(M_z^\ast)$, because $\overline{\lambda}$ is not an eigenvalue of $M_z$. Summarizing:
\begin{enumerate}
    \item[(i)] $\sigma_p(M_z^\ast)=\mathbb{D}$, $\sigma_c(M_z^\ast)=\mathbb{T}$, $\sigma_{res}(M_z^\ast)=\emptyset$.
    \item[(ii)] $\sigma_p(M_z)=\emptyset$, $\sigma_c(M_z)=\mathbb{T}$, $\sigma_{res}(M_z)=\mathbb{D}$.
\end{enumerate}
\subsection{The spectrum of a bounded, self-adjoint operator}\label{SSectSpectrumSA}
We have next a foundational result on self-adjoint operators.
\begin{theorem}\label{theoSAspectrumIsReal}
    If $T=T^\ast$, then $\sigma(T)\subset\mathbb{R}$.
\end{theorem}
\begin{proof}
    We show that if $a,b\in\mathbb{R}$ and $b\ne0$, then $a+ib\in\rho(T)$. Consider the sesquilinear form
    \[
    B(x,y)=\langle(T-(a+ib))x,(T-(a+ib))y\rangle=\langle(T-(a-ib))(T-(a+ib))x,y\rangle, 
    \]
    which is coercive since, using again $T=T^\ast$,
    \begin{eqnarray*}
        B(x,x)&=&\langle(T-(a+ib))x,(T-(a+ib))x\rangle\crcr 
        &=&\|(T-a)x\|^2+b^2\|x\|^2+ib(-\langle Tx,x \rangle+ \langle x,Tx \rangle)\crcr 
        &=&\|(T-a)x\|^2+b^2\|x\|^2\crcr 
        &\ge&b^2\|x\|^2.
    \end{eqnarray*}
    By Lax-Milgram lemma, $(T-(a-ib))(T-(a+ib))$ is invertible, hence $T-(a+ib)$ is $1-1$ and $T-(a-ib)$ is onto. Applying the same reasoning $T-(a-ib)$ instead of $T-(a+ib)$, we have that $T-(a+ib)$ is $1-1$ and $T-(a-ib)$ is onto. Hence, $T-(a+ib)$ is invertible (incidentally, as well as $T-(a+-ib)$).
\end{proof}
    We record here a useful corollary of Lax-Milgram lemma.
        \begin{corollary}\label{corOfLaxMilgram}
        Let $\mu\in\mathbb{C}$. Then, $\mu\in\sigma(A)$ if and only if there exist vectors $\{h_j\}_{j=1}^\infty$ in $H$ with $\|h_j\|=1$, such that
\begin{equation}\label{eqOfLaxMilgram}
\|(A-\mu I)h_j\|\to0.
\end{equation}
    \end{corollary}
    \begin{proof} Suppose such vectors do not exist.
        Then, there is $\epsilon>0$ such that $\|(A-\mu I)h\|\ge \epsilon\|h\|$, hence the sesquilinear form
        \[
        B(h,k)=\langle (A-\mu I)h,(A-\mu I)k\rangle= \langle (A-\mu I)(A-\mu I)h,k\rangle
        \]
        is coercive, and by Lax-Milgram $(A-\mu I)^2$ is invertible, hence $A-\mu I$ is invertible, $\mu\in\rho(A)$.

        In the other direction, suppose $A-\mu I$ is invertible. For $\|h\|=1$,
        \[
        \infty>\|(A-\mu I)^{-1}\|\ge\frac{\|(A-\mu I)^{-1}(A-\mu I)h\|}{\|(A-\mu I)h\|}=\frac{1}{\|(A-\mu I)h\|},
        \]
        showing that \eqref{eqOfLaxMilgram} does not hold.
    \end{proof}
    The spectrum of a self-adjoint operator has no residual part.
\begin{corollary}\label{corSAhasNoResidual}
    If $T=T^\ast$, then $\sigma_{res}(T)=\emptyset$.
\end{corollary}
\begin{proof}
    If $\lambda\in\sigma_{res}(T)$, then $\lambda=\overline{\lambda}\in\sigma_p(T^\ast)=\sigma_p(T)$, which is a contradiction.
\end{proof}
In infinite dimensions, it can happen that $A:H\to H$ has an eigenvalue $\lambda$, yet $\text{Ran}(\lambda I-A)$ is the whole of $H$: think of the back-shift $\tau_1^\ast$, for which $\ker(\tau_1^\ast)\ne0$ ($\lambda=0$ is an eigenvalue), but $\text{Ran}(\tau_1^\ast)=\ell^2(\mathbb{N})$. If $A=A^\ast$ is self-adjoint, however, things are similar to the finite dimensional case.
\begin{corollary}\label{corSAdiscSpectrumAsUsual}
If $T=T^\ast$ and  $\lambda\in\sigma_{p}(T)$, then $\text{Ran}(\lambda I-T)$ is not dense in $H$.
\end{corollary}
\begin{proof}
    Since $\lambda$ is real,  $\lambda I-T$ is self-adjoint, and it suffices to show the statement for $\lambda=0$. We show, more precisely, that if $0\ne x\in\ker(T)$, then $x\perp\text{Ran}(T)$, so that the latter is not dense in $H$. In fact,
    \[
    \langle Ty,x\rangle=\langle y ,Tx\rangle=0.
    \]
\end{proof}
Actually, we can be more precise in locating the spectrum of $T=T^\ast$.
\begin{theorem}\label{theoSpectrumSaLocation}
If $T$ is selfadjoint, and
\begin{equation}
   a(T):=\inf_{\|x\|=1}\langle x,Tx\rangle \text{ and } b(T):=\sup_{\|x\|=1}\langle x,Tx\rangle,
\end{equation}
then $a(T),b(T)\in\sigma(T)$ and $[a(T),b(T)]\supseteq\sigma(T)$.
\end{theorem}
\begin{proof}
    Observe first that $a(T+cI)=a(T)+c$, $b(T+cI)=b(T)+c$, and
    \[
    \sigma(T+cI)=\{\lambda:\ \lambda I-(T+cI)\text{ is not invertible}\}=\sigma(T)+c,
    \]
    and we can suppose $a(T)=0$. If $\lambda>b(T)$, then
    \[
    \langle x,(\lambda I-T)x\rangle \ge(\lambda-b(T))\|x\|^2
    \]
    shows that $B(x,y):=\langle x,(\lambda I-T)y\rangle$ is bounded and coercive, and we have seen before that this implies that $\lambda I-T$ is invertible, hence that $\lambda\in \rho(T)$. The same reasoning applies, changing sign, with $\lambda<a(T)$. This shows $\sigma(T)\subseteq[a(T),b(T)]$.
    
On the other hand, $b(T)\le\|T\|=r(T)$, the spectral radius of $T$. Hence, $b(T)\in \sigma(T)$. The same reasoning, but choosing the normalization $b(T)=0$, shows that $a(T)\in\sigma(T)$.   
\end{proof}
\subsection{The Hausdorff distance between spectra}\label{SSectHausDistSpectra}
Given two closed sets $A,B$ in $\mathbb{R}$ (or in any metric space), their {\it Hausdorff distance} $d_H(A,B)$ is defined as
\begin{equation}\label{eqHausdorffDistance}
    d_H(A,B)=\max\left(\sup_{x\in A}\inf_{y\in B}|x-y|,\sup_{y\in B}\inf_{x\in A}|x-y|\right).
\end{equation}
To have a better picture of this quantity, let $d(x,B)=\inf_{y\in B}|x-y|$. Asking $\sup_{s\in A}\inf_{y\in B}|x-y|\le\epsilon$ means that all points of $A$ are within distance $\epsilon$ from $B$, and asking $d_H(A,B)\le\epsilon$ means that this requirement also holds with the roles of $A$ and $B$ interchanged. When $\epsilon=0$, we are requiring that $A\subseteq B$ and $B\subseteq A$, i.e. $A=B$.
\begin{exercise}\label{exeHausdorffDist}
    Let $A,B$ be subsets, not necessarily closed, of a metric space. Show that $d_H(A,B)=0$ if and only if $\text{Cl}(A)=\text{Cl}(B)$.

    Also, show that $d_H$ defines a distance, when restricted to nonempty closed sets.
\end{exercise}
\begin{theorem}\label{theoHausdDistSpectra}
    If $S,T$ are self-adjoint, then 
    \begin{equation}\label{eqHausdDistSpectra}
        d_H(\sigma(S),\sigma(T))\le\|S-T\|.
    \end{equation}
\end{theorem}
\begin{proof} The case $S=T$ is trivial. By compactness of spectra, we can replace $\sup/\inf$ by $\max/\min$.
    Suppose by contradiction that $\max_{s\in \sigma(S)}\min_{t\in \sigma(T)}|s-t|>\|S-T\|$, i.e. that for some $s\in\sigma(S)$,
    \[
    \min_{t\in \sigma(T)}|s-t|>\|S-T\|.
    \]
    Thus $s\in\rho(T)$. i.e. $T-sI$ is invertible, and by the spectral mapping theorem for Banach algebras,
    \[
    \sigma((T-sI)^{-1})=(\sigma(T)-s)^{-1}.
    \]
    The spectral radius $r((\sigma(T)-sI)^{-1})$ satisfies, since $(\sigma(T)-sI)^{-1}$ is self-adjoint:
    \begin{eqnarray}\label{eqHausdDistSpectraIntermediate}
        \|(T-sI)^{-1}\|&=&r((\sigma(T)-sI)^{-1})\crcr 
        &=&\max_{t\in\sigma(T)}|t-s|^{-1}<\|S-T\|^{-1}.
    \end{eqnarray}
    We next decompose $S-sI$ in a way which allows us to use the information we have gathered so far:
    \[
    S-sI=T-sI+S-T=(T-sI)\left(I+(T-sI)^{-1}(S-T)\right).
    \]
    The first factor is invertible, and the second is, too, because by \eqref{eqHausdDistSpectraIntermediate}
    \[
    \|(T-sI)^{-1}(S-T)\|\le\|(T-sI)^{-1}\|\cdot\|S-T\|<1.
    \]
    This shows that $s\notin\sigma(S)$, a contradiction.
\end{proof}
We finish this section with some observations on the adjoint.
\begin{exercise}\label{exeAdjProperties}
Prove the following properties.
    \begin{enumerate}
    \item[(i)] Let $T\in\mathcal{L}(H)$. Then, $T_R=\frac{T+T^\ast}{2}$ and $T_I=\frac{T-T^\ast}{2i}$ are self-adjoint, and
    \[
    T=T_R+iT_I.
    \]
    This decomposition is somehow analogous to the algebraic expression of a complex number.
    \item[(ii)] The decomposition in (i) is unique. If $T=A+iB$ with $A,B$ self-adjoint, then $A=T_R$ and $B=T_I$.
    \item[(iii)] If $A$ is self-adjoint, then $\langle x,Ax\rangle$ is real,
    \[
    \langle x,Ax\rangle=\overline{\langle Ax,x\rangle}=\overline{\langle x,Ax\rangle}.
    \]
    \item[(iv)] $A$ is selfadjoint if and only if $B=iA$ is {\bf skew-adjoint}, $B^\ast=-B$. In this case,  $\langle x,Bx\rangle$ is purely imaginary.
    \item[(v)] $T_IT_R=T_RT_I$ if and only if $T^\star T=TT^\ast$, i.e. if $T$ is {\bf normal}.
    \item[(vi)] If $T$ is normal, then $\|T\|=\|T^\ast T\|^{1/2}=\sqrt{\|T_I^2+T_R^2\|}$.
\end{enumerate}
\end{exercise}

\chapter{Compact operators and their spectra}\label{chapcCmpctSpctr}

If you are just interested in compact, self-adjoint operators on a Hilbert space, you find a different narrative of the same theory in chapter \ref{chaptUSAandUO}.

\section{Preliminaries on compactness}
\subsection{Compact subsets of Banach spaces}\label{SSectCompactSetsInBanach}
The following lemma shows that closed subspaces of Banach spaces exhibit some orthogonality properties.
\begin{lemma}[Riesz lemma in Banach spaces]\label{lemmaRieszLemmaBanach}\index{lemma!of Riesz in Banach spaces}
Let $X$ be a Banach space and $Y\subset X$ a proper, closed subspace. Let $0<\alpha<1$. Then, there exists $z$ in $X$ with $\|z\|=1$ such that $\alpha\le d(z,Y)\le 1$.     
\end{lemma}
\begin{proof}
    Let $v\in X\setminus Y$, $a:=d(v,Y)>0$, and choose $y_0\in Y$ such that $a\le\|v-y_0\|\le\frac{a}{\alpha}$ (draw a picture to have an intuition of what's going on). Set $z=\frac{v-y_0}{\|v-y_0\|}$, so that $\|z\|=1$. For all $y$ in $Y$ we have
    \begin{eqnarray*}
        \|z-y\|&=&\frac{\|v-y_0-\|v-y_0\|y\|}{\|v-y_0\|}\crcr 
        &=&\frac{\|v-y_1\|}{\|v-y_0\|}\text{ with }y_1\in Y\crcr 
        &\ge&\frac{a}{\|v-y_0\|}\ge\alpha
    \end{eqnarray*}
    by our choice of $y_0$.
\end{proof}
\begin{corollary}\label{corUnitBallIsNotCompact}
    Let $X$ be a Banach space. The unit ball in $X$ is compact in the norm topology if and only if $X$ is finite dimensional.
\end{corollary}
\begin{proof}
    The {\it if} direction easily follows from Weierstrass theorem on extremals of continuous functions.

    In the other direction, suppose $X$ is infinite dimensional, and construct a sequence of vectors and subspaces as follows.
    \begin{enumerate}
        \item[(i)] Pick any vector $x_1$ with $\|x_1\|=1$, and let $Y_1=\text{span}\{x_1\}\subset X$.
        \item[(ii)] By Riesz lemma, you can pick a second vector $x_2$ with $\|x_2\|=1$ and $d(x_2,Y_1)\ge1-1/2$, and let  $Y_2=\text{span}\{x_1,x_2\}\subset X$.
        \item[(iii)] Iterate, so that $d(x_n,Y_{n-1})\ge 1-1/n$, and let $Y_n=\text{span}\{x_1,\dots,x_n\}\subset X$.
    \end{enumerate}
    This way, $\|x_m-x_n\|>1/3$ for all $m\ne n$. Cover the unit ball by open balls of radius $1/6$, so that no one can contain more than one point from the sequence $\{x_n\}$; hence, the unit ball is not compact.
\end{proof}

\subsection{Compact operators: generalities}\label{SSectCompactonBanach}
Linear functionals on a Banach space $X$ are maps $l:X\to{\mathbb C}$, with values in the $1$-dimensional Banach space ${\mathbb C}$. Next level of complexity are bounded, linear maps $T:X\to Y$, where $Y$ is a finite dimensional Banach space. If $\{l_n\}_{n=1}^N$ is a basis for $Y^\ast$, then $C:y\mapsto (l_1(y),\dots,l_N(y))=\sum_{n=1}^N l_j(y)e_j$ is an isomorphism between $Y$ and ${\mathbb C}^N$ ($\{e_j\}$ is here the standard basis of ${\mathbb C}^N$), and $C\circ T(x)=\sum_{n=1}^N l_j(T x)e_j$ identifies $T$ with a map from $X$ to ${\mathbb C}^N$. All this makes sense with little changes if $T(X)$ is finite dimensional in $Y$, a Banach space. In this case we say that $T$ is {\it finite rank}.

Now the closed unit ball $\overline{B_1}$ in $X$ is not compact if $X$ is infinite dimensional and Banach; but if $T$ is finite rank, $T(\overline{B_1})$ is both bounded and closed in $Y$, hence it is compact.

We say that an operator $T:X\to Y$ between Banach spaces is {\it compact} if $\overline{T(B_1)}$ is compact in $Y$.\index{operator!compact}

\begin{footnotesize}
\begin{exercise}
Let $m:{\mathbb N}\to{\mathbb C}$ be bounded, and consider the bounded operator $M_b\varphi=m \varphi$, $M_m:\ell^2\to\ell^2$. Show that $M_m$ is finite rank if and only if $m(j)=0$ for all, but a finite set of $j$'s in ${\mathbb N}$; and that is is compact if and only if $\lim_{j\to\infty}m(j)=0$.
\end{exercise}
\end{footnotesize}

Recall that a sequence $\{x_n\}$ in a Banach space $X$ {\it converges to $x$ weakly}, $x_n\xrightarrow[w]{}x$, if
\[
\lim_{n\to\infty}l(x_n)=l(x)
\]
for all $l\in X^\ast$.

The theorem below is crucial in many applications, e.g. to Calculus of Variations.\index{weak topology}
\begin{theorem}\label{TheoWeakCompactOp}\index{weak topology!and compact operators}
If $T$ is compact, it maps weakly convergent sequences to norm convergent sequences. 
\end{theorem}
\begin{proof}
By Banach-Steinhaus, $x_n\xrightarrow[w]{}x$ implies that $\|x_n\|\le C$ is norm bounded, hence
\[
l(T x_n)-l(T x)=(T' l)(x_n-x)\to 0
\]
as $n\to\infty$ for all $l\in Y^\ast$. This means that $T(x_n)\xrightarrow[w]{} T(x)$ weakly in $Y$, hence that $\|T x_n\|\le K$ is bounded, too.

On the other hand, if by contradiction $\|T x_n-T x\|\nrightarrow 0$, then $\|T x_{n_j}-T x\|\ge\epsilon$ for a subsequence, and by compactness of $T$ there exists a sub-subsequence $\{T x_{n_{j_k}}\}$ strongly convergent to some $y\ne T x$ in $Y$. But then $T x_{n_{j_k}}\xrightarrow[w]{}y\ne T x$, contradicting the assumption.
\end{proof}
The oppposite implication holds, for instance, if $X$ is a reflexive space.
\begin{theorem}\label{theoWeakCompactOpInverse}
    Let $X\equiv X^{\ast\ast}$ be a reflexive Banach space, and $T:X\to Y$ be an operator with values in a Banach space $Y$. If $T$ maps weakly convergent sequences to norm convergent sequences, then $T$ is compact.
\end{theorem}
\begin{proof}
    Let $\Bar{B}$ be the unit ball in $X$ and suppose that 
\end{proof}
The spectrum of a compact operator on a infinite dimensional Banach space always contains $0$.
\begin{proposition}\label{propSpectrumOfCompactHasZero}
    Let $T:X\to X$ be a compact operator on an infinite dimensional Banach space $X$. Then, $T$ is not surjective. In particular, $0\in\sigma(T)$. 
\end{proposition}
\begin{proof}
    Let $B$ be the unit ball in $X$. Suppose $T$ is surjective. By the open mapping theorem, $T(\Bar{B})$ contains a closed ball $\epsilon \Bar{B}$, which is compact, since $T(\Bar{B})$ is compact. By corollary \ref{corUnitBallIsNotCompact}, $X$ is finite dimensional.
\end{proof}

\section{Compact operators on Hilbert spaces}\label{SectCompactOnHilbert}
\subsection{Compactness vs. finite rank}\label{SSectCmptVsFiniteRank}
The picture in the Hilbert space case is much cleaner.
\begin{theorem}\label{TheoCompactHilbert}\index{operator!compact vs. finite rank}
Let $T\in{\mathcal L}(H)$ be a compact operator on a separable Hilbert space. Then, there are finite rank operators $T_n$ on $H$ such that $\|T_n-T\|\to0$.
\end{theorem}
\begin{proof}
Let $\{e_n\}$ be a o.n.b. for $H$ and set
\[
\lambda_n:=\sup_{\|x\|=1,x\in \text{span}\{e_1,\dots,e_n\}^\perp} \|T x\|\ge\lambda_{n+1}.
\]
Let $\lambda=\lim_n\lambda_n$ We start by proving that $\lambda=0$.
Pick first $x_n\in \text{span}\{e_1,\dots,e_n\}^\perp$ such that $\|T x_n\|\ge\lambda_n/2\ge\lambda/2$.
Observe that $x_n\xrightarrow[w]{}0$:
\begin{eqnarray*}
\left|\langle x_n,y \rangle\right|&=&\left|\langle x_n,\sum_{j=1}^\infty\langle y,e_j \rangle e_j\rangle\right|\crcr 
&=&\left|\sum_{j=n+1}^\infty\langle x_n, e_j\rangle \overline{\langle y,e_j \rangle}\right|\crcr
&\le&\left(\sum_{j=n+1}^\infty|\langle x_n, e_j\rangle |^2\right)^{1/2}\left(\sum_{j=n+1}^\infty|\langle y,e_j \rangle|^2\right)^{1/2}\crcr
&\le&\|x_n\| \left(\sum_{j=n+1}^\infty|\langle y,e_j \rangle|^2\right)^{1/2}\crcr 
&\to&0
\end{eqnarray*}
as $n\to\infty$. By the previous theorem, $\|T x_n\|\to0$ because $T$ is compact, hence, $\lambda=0$.

Let now $\pi_n$ be the projection onto $\text{span}(e_1,\dots,e_n)$, and
\begin{eqnarray*}
T_n h&=&(T \pi_n)h\crcr 
&=&T\left(\sum_{j=1}^n\langle e_j,h\rangle e_j\right)\crcr 
&=&\sum_{j=1}^n \langle e_j,h\rangle T e_j,
\end{eqnarray*}
which has finite rank, and satisfies, with $\pi_n^\perp$ being the projection onto $\text{span}(e_1,\dots,e_n)^\perp$:
\begin{eqnarray*}
\|T-T_n\|&=&\|T(I-\pi_n)\|\crcr 
&=&\|T \pi_n^\perp\|\crcr 
&=&\lambda_n\to0
\end{eqnarray*}
as $n\to\infty$.
\end{proof}
\subsection{Fredholm Theorem and its consequences}\label{SSectFredholmAndSoOn}\index{theorem!Fredholm}
\begin{theorem}\label{TheoFredholm}
Let $T$ be a compact operator on a separable Hilbert space $H$. Then, $(I-z t)^{-1}$ exists for $z\in{\mathbb C}\setminus D$, where $D$ is a discrete set in ${\mathbb C}$. 

Moreover, $z\in D$ if and only if the equation $x=z T x$ has a solution $x\ne0$.
\end{theorem}
The idea of the proof consists in first localizing the problem, then reducing it, in substance, to a similar problem concerning finite rank operators. The latter can be rephrased in terms of zeros of a holomorphic determinant, and this is where the discrete set makes its appearence.

\begin{proof}
We start with a formal calculation also involving a finite rank operator $F$,
\begin{eqnarray}
I-z T&=&(I-z T+F) -F\crcr 
&=&(I-F(I-z T+F)^{-1})(I-z T+ F),
\end{eqnarray}
provided $(I-z T+F)^-1$ exists. The advantage of the first factor in the last expression over $I -z T$ is that the former has the form $I-[\textbf{finite rank operator}]$, on which we can use tools from ordinary linear algebra.

We start by proving the theorem on discs.
Let $z_0$ a complex number, and consider $z$'s such that
\[
\|z T-z_0 T\|=|z-z_0|\cdot\|T\|<1/2,
\]
and a finite rank $F$ such that
\[
\|z_0 T-F\|<1/2,
\]
so that $(I-z T+F)^{-1}=\sum_{n=0}^\infty(z T-F)^n$ exists by the convergence of the Neumann series.\index{Neumann series} Define then $g(z)=F(I-z T+F)^{-1}$, a finite rank operator which holomorphically depends on $z$, and note that above we have computed:
\[
I-z T=(I-g(z)) (I-z T+ F),
\]
the second factor on the right being invertible. Then,
\begin{itemize}
    \item $x=g(z)x$ has a solution $x\ne0$ if and only if $y= zT y$ has a solution $y\ne 0$;
    \item $I-g(z)$ is invertible if and only $I-z T$ is invertible.
\end{itemize}
Let $e_1,\dots,e_n$ be a bases for ${\text Ran}(F)$, and observe that a solution to $x=g(z)x$ must have the form $x=\sum_{j=0}^n a_j e_j$, since it belongs to $\text{Ran}(g(z))\subseteq \text{Ran}(F)$. On the other hand, the operator $F$ must have the form:
\[
F x=\sum_{j=1}^n\langle f_j,x \rangle e_j
\]
for some $f_1,\dots,f_n\in H$.
Plug this into the equation,
\begin{eqnarray*}
\sum_{j=0}^n a_j e_j&=&x=g(z)x\crcr
&=&\sum_{l=0}^n a_l g(z)e_l\crcr 
&=& F((I-z T+F)^{-1}e_l)\crcr 
&=&\sum_{l=0}^n a_l\sum_{j=1}^n\langle f_j, (I-z T+F)^{-1}e_l  \rangle e_j\crcr 
&=&\sum_{j=1}^n\left(
\sum_{l=0}^n a_l\langle [(I-z T+F)^{-1}]^\ast f_j,e_l \rangle
\right)e_j.
\end{eqnarray*}
i.e. we look for a nonzero solution to
\begin{equation}\label{eqDeterminant}
a_j=\sum_{l=0}^n a_l\langle [(I-z T+F)^{-1}]^\ast f_j,e_l \rangle
=\sum_{l=0}^n a_l\langle  f_j(z),e_l \rangle,\ j=1,\dots,n.
\end{equation}
Such solution exists if and only if $d(z)\ne0$, where
\[
d(z)=\det\left[\delta_{j,l}-\langle  f_j(z),e_l \rangle\right]_{j,l=1,\dots,n}.
\]
The matrix has holomorphic entries, hence, one of the two can happen: $d(z)\equiv0$ in $B(z_0,\frac{1}{2\|T\|})$, or, there is a discrete subset $D_{z_0}$ of $B(z_0,\frac{1}{2\|T\|})$ such that $d(z)=0$ for $z\in D_{z_0}$.

We claim that if $d(z)\ne0$, then $I-g(z)$ is invertible, and so is $I-z T$. Fix $u\in H$ and look for solution $x=u+\sum_{j=0}^n a_j e_j$ to the equation
\[
x-g(z)x=u,
\]
i.e.
\begin{eqnarray*}
g(z)u&=&\sum_{j=0}^n a_j e_j-g(z)\left(\sum_{l=0}^n a_l e_l\right)\crcr 
&=&\sum_{j=0}^n a_j e_j
-\sum_{j=1}^n\left(
\sum_{l=0}^n a_l\langle [(I-z T+F)^{-1}]^\ast f_j,e_l \rangle
\right)e_j,
\end{eqnarray*}
which is a nonhomogeneous version of  (\ref{eqDeterminant}), having solution if $d(z)\ne0$.

At his point, we have also shown that $(I-z T)^{-1}$ does not exist if and only if $x=z T x$ has a nonzero solution.

We have now to pass from local to global. The subset $D'$ of the limit (accumulation) points of $D$ is closed by general topology, and it is open by what we have just seen, hence it is either empty or the whole complex plane. Since $\{1/z:z\in D\}\subseteq\sigma(T)$ is nonempty and compact, $D'=\emptyset$. Hence, $D$ is discrete.
\end{proof}

\begin{corollary}\label{CorFredAlt}[Fredholm alternative]\index{theorem!Fredholm alternative}
Let $T$ be a compact operator on a separable Hilbert space. If $I-T$ is not invertible, then $x=T x$ has a nonzero solution.
\end{corollary}
\begin{proof}
Apply the theorem to $z=1$.
\end{proof}

\index{theorem!Riesz-Schauder}
\begin{theorem}\label{TheoRiszSchauder}[Riesz-Schauder Theorem] Let $T$ be a compact operator on a separable Hilbert space $H$. Then, $\sigma(T)$ can only have $0$ as accumulation point (equivalently: $\sigma(T)\setminus\{0\}$ is discrete in ${\mathbb C}\setminus\{0\}$).

Moreover, the eigenspace $E_\lambda$ relative to a nonzero $\lambda\in\sigma(T)$ is finite dimensional. 
\end{theorem}
\begin{proof}
For $\lambda\ne0$, we have $\lambda I-T=\lambda\left(I-\frac{1}{\lambda}T\right)$, which is invertible for $1/\lambda$ in a discrete set in ${\mathbb C}$, i.e. for $\lambda$ in a bounded, discrete set in ${\mathbb C}\setminus\{0\}$, having $0$ as accumulation point if the set of the $1/\lambda$'s was infinite.

Suppose that $E_\lambda$ is infinite dimensional for an eigenvalue $\lambda\ne0$, observe that $E_\lambda$ is closed, and let ${\mathcal E}=\{e_n\}_{n=1}^\infty$ be an o.n.b. for it. Then, ${\mathcal E}$ is bounded, but $T{\mathcal E}=\lambda{\mathcal E}$ is not precompact (as we have seen in an exercise).
\end{proof}
\subsection{Self-adjoint compact operators}\label{SSectCompactSA}
If we put together what we have seen so far on the topic, we have the fundamental:\index{theorem!Hilbert-Schmidt}
\begin{theorem}\label{TheoHilbertSchmidt} [Hilbert-Schmidt Theorem] Let $A$ be a compact self-adjoint operator\index{operator!compact, self-adjoint} on a (infinite dimensional) separable Hilbert space $H$. Then  $\sigma(A)\subset{\mathbb R}$. Moreover, there is a o.n.b. $\{e_n\}_{n=1}^\infty$ and a real sequence $\{\lambda_n\}_{n=1}^\infty$ of $H$ (possibly with repetitions) such that
\[
A e_n=\lambda_n e_n
\]
and
\begin{equation}\label{eqEigenvaluesGoToZero}
\lim_{n\to\infty}\lambda_n=0.
\end{equation}
If $\lambda_n\ne0$, the eigenspaces $E_{\lambda_n}=\{x\in H:Ax=\lambda_n x\}$ are finite dimensional.
\end{theorem}
\begin{proof}
Since $A=A^\ast$, $\sigma(A)\subset\mathbb{R}$.
For an eigenvalue $\mu\ne 0$ of $A$, the image of the closed unit ball $B$ in $H$ contains $|\mu|$ times the unit ball in $E_{\mu}$, thus the latter has to be compact, which only is if $\dim(E_{\mu})<\infty$.

 By Fredholm theorem, $\sigma(A)\setminus\{0\}=\{\mu_m\}$, where $1/\mu_m$ has no point of accumulation. If the number of such $\mu_m$'s is finite, then $A$ is finite rank. If it is infinite, then $\mu_m\to0$.

Arrange the elements in decreasing order $|\mu_1|\ge|\mu_2|\ge\dots>0$, where it might happen that $\mu_i=-\mu_{i-1}>0$ for some $i$'s.  Let $M_i=\sharp(E_{\mu_i})$. Consider then a orthonormal basis $e^i_1,\dots,e^i_{M_i}$ for each $E_{\mu_i}$, and a possibly countable orthonormal basis $f_1,f_2,\dots$ of $\ker(A)$. For both the cases where the number of the nonzero eigenvalues if finite, or infinite, we can arrange all these vectors in a sequence $\{e_n\}$ with $Ae_n=\lambda_n$ (and $\lambda_n$ is one of the $\mu_i$'s, or $\lambda_n=0$), and $\lambda_n\to0$ as $n\to\infty$.

Finally, if $\lambda_m\ne\lambda_n$, and $Ax=\lambda_m x$, $Ay=\lambda_ny$, then
\[
\lambda_m \langle x,y\rangle= \langle Ax,y\rangle=\langle x,Ay\rangle=\lambda_n\langle x,y\rangle,
\]
hence, $\langle x,y\rangle=0$. In particular, $\{e_n\}$ is a orthonormal system. If it were not complete, then
\[
K:=\overline{\text{span}\{e_n\}}\subset H, \text{ with }A(K)= K.
\]
We have then $A:H\ominus K\to H\ominus K$. In fact, $\langle e_n,x\rangle=0$ implies that $\langle e_n,Ax\rangle=\langle Ae_n,x\rangle=\lambda_n\langle e_n,x\rangle=0$. Clearly $A|_{H\ominus K}:H\ominus K\to H\ominus K$ is self-adjoint and compact. Since $\ker(A)\subseteq K$. $\ker(A|_{H\ominus K})=0$. Then, $A_{H\setminus K}$ has an eigenvalue  $\mu\in\mathbb{R}\setminus\{0\}$. The eigenvalue $\mu$ is not even one of the previous nonzero eigenvalues, otherwise its eigenspace in $H\ominus K$ would lie in $E_{\mu_i}\subseteq K$ for some $i$. It is an eigenvalue for $A$, hence $\{\mu_i\}\cup\{0\}$ did not exhaust $\sigma(A)$, and we have reached a contradiction.  
\end{proof}
\begin{corollary}\label{corOfHS}
    Let $A$ be a compact self-adjoint operator on a (infinite dimensional) separable Hilbert space $H$. Then,
    \begin{equation}\label{eqCorOfHS}
        H=\overline{\text{Ran}(A)}\oplus\ker(A).
    \end{equation}
\end{corollary}
\begin{proof} Let $\{\lambda_n\}$ be the sequence of the nonzero eigenvalues of $A$.
    By Hilbert-Schmidt theorem, each $x$ in $H$ can be written as
    \[
    x=\sum_n\Pi_{E_{\lambda_n}}x+\Pi_{\ker(A)}x, 
    \]
    and (i) $Ax=\sum_n\lambda_n\Pi_{E_{\lambda_n}}x$; (ii) $\sum_n\Pi_{E_{\lambda_n}}x\perp \Pi_{\ker(A)}x$. Hence, 
    \[
    \text{Ran}(A)=\overline{\bigoplus_nE_{\lambda_n}}\perp\ker(A).
    \]
\end{proof}
\section{The matrix representation of bounded operators}\label{SectOPeratirsAsMatrices}
Any separable, infinite dimensional Hilbert space is unitarily isomorphic to $\ell^2=\ell^2(\mathbb{N}_\ast)$ (where $\mathbb{N}_\ast=\{1,2,\dots\}$ might be replaced by any infinite, countable set). To any linear operator $A:H\to A$ we can then associate the infinite matrix $[A_{m,n}]_{m,n=1}^\infty$, where $A_{m,n}=\langle e_m,A e_n\rangle=\overline{A^\ast_{n,m}}$. We can perform on infinite matrices representing bounded operators many of the calculations we are used to do in the finite dimensional case, which is often useful. In order to make calculations into proofs, it has to be clear in which sense the involved infinite sums converge, and which manipulations are licit. In this section we will see some results, examples, and counterexamples.

This material is not needed in the sequel, but it can be useful in working with concrete operators.
\subsection{Bounded operators on separable Hilbert spaces as infinite matrices}\label{SSectOPeratirsAsMatrices}
If $A\in\mathcal{L}(H)$ is bounded, and $x=\sum x_ne_n$, then, by the continuity of $A$ and of the inner product, 
\begin{eqnarray}\label{eqMatOpOne}
    \langle e_m,Ax\rangle&=&\left\langle e_m,A \sum_n x_ne_n\right\rangle=\left\langle e_m,\sum_n x_n Ae_n\right\rangle\\ 
    &=&\sum_n\left\langle e_m,A  e_n\right\rangle x_n=\sum_nA_{m,n}x_n.
\end{eqnarray}
Since $x\mapsto \langle e_m,Ax\rangle$ is a bounded linear functional on $H$, the row sequence $\{A_{m,n}\}_{n=1}^\infty$ belongs to $\ell^2$. This also implies that the numerical series above converge absolutely.
We will write $A$ in terms of its rows and columns,
\[
A=\begin{pmatrix}
  A_1\crcr 
  A_2\crcr 
  \dots
\end{pmatrix}
=(A^1|A^2|\dots).
\]
Basically, we have proved the following.
\begin{proposition}\label{propMatOpOne}
The infinite matrix $A=[A_{m,n}]$ represents a bounded operator on $\ell^2$ if and only if:
\begin{enumerate}
    \item[(i)] each row of $A$ belongs to $\ell^2$, and
    \item[(ii)] for all $x\in\ell^2$, $(A_{m,n}x_n)_{m=1}^\infty=(\langle A_m^t,x\rangle)_{m=1}^\infty$ (where $^t$ denotes the transpose) belongs to $\ell^2$, and the operator $(x_n)_{n=1}^\infty\mapsto (A_{m,n}x_n)_{m=1}^\infty$.
\end{enumerate}
Moreover, (i) and (ii) might be equivalently required to hold for the columns instead of the rows.
\end{proposition}
An immediate consequence of proposition \ref{propMatOpOne} is the following version of Fubini theorem for infinite sums, and the usual row-times-product rule for the product of matrices.
\begin{corollary}[Change of order of summation in infinite sums]\label{corMatOpOne}
    Let $A$ be a bounded operator on $\ell^2$ and $x,y\in\ell^2$. Then,
    \begin{equation}\label{eqChangeOfOrderInInfiniteSums}
        \sum_n\left(\sum_m \overline{y_m} A_{m,n}\right)x_n=\sum_m \overline{y_m}\left(\sum_n A_{m,n}x_n\right).
    \end{equation}
    If $A$ and $B$ are bounded operators on $\ell^2$, then
    \begin{equation}\label{eqChangeOfOrderInInfiniteSumsTwo}
        (BA)_{m,n}=\sum_j B_{m,j}A_{j,n}.
    \end{equation}
\end{corollary}
\begin{proof}
    The identity \eqref{eqChangeOfOrderInInfiniteSums} follows from
    \begin{eqnarray*}
        \sum_m \overline{y_m}\left(\sum_n A_{m,n}x_n\right)&=&\langle y,Ax\rangle\crcr
        &=&\langle A^\ast y,x\rangle=\sum_n\left(\sum_m\overline{A^\ast_{n,m}y_m}\right)x_n\crcr
        &=&\sum_n\left(\sum_m \overline{y_m} A_{m,n}\right)x_n.
    \end{eqnarray*}
    The row-times-product rule for the product of operators follows from Plancherel formula (used in the equality from first to second line):
    \begin{eqnarray*}
        (BA)_{m,n}&=&\langle e_m,BAe_n\rangle=\langle B^\ast e_m,Ae_n\rangle\crcr 
        &=&\sum_j\overline{(B^\ast e_m)_j}(Ae_n)_j=\sum_j\overline{\langle e_j,B^\ast e_m\rangle}\langle e_j,Ae_n\rangle\crcr
        &=&\sum_j\langle e_m,B e_j\rangle \langle e_j,Ae_n\rangle\crcr 
        &=&\sum_j B_{m,j}A_{j,n}.
    \end{eqnarray*}
\end{proof}
The relation \eqref{eqChangeOfOrderInInfiniteSums}  can be used to change order of summation when several operators are composed. For instance, if $A$ and $B$ are bounded and $x\in \ell^2$, then
\begin{equation}\label{eqChangeOfOrderInInfiniteSumsThree}
    \sum_n\left(\sum_j B_{m,j}A_{j,n}\right)x_n=\sum_j B_{m,j}\left(\sum_n A_{j,n} x_n\right),
\end{equation}
since $(B_{m,j})_{j=1}^\infty$ lies in $\ell^2$.

We might prove \eqref{eqChangeOfOrderInInfiniteSumsThree} by a different argument:
\begin{eqnarray*}
    \sum_n\left(\sum_j B_{m,j}A_{j,n}\right)x_n&=&\sum_n(BA)_{m,n}x_n=(BAx)_m\crcr 
    &=&[B(Ax)]_m=\sum_j B_{m,j}(Ax)_j\crcr 
    &=&\sum_j B_{m,j}\left(\sum_n A_{j,n} x_n\right).
\end{eqnarray*}
\begin{exercise}\label{exeUnitaryMatrix}
    Let $U$ be a bounded operator on $\ell^2$, $U=(U^1|U^2|\dots)=\begin{pmatrix}
       U_1\crcr 
       U_2\crcr 
       \dots 
    \end{pmatrix}$. Show that $U$ is unitary if and only if $U_1^t,U_2^t,\dots$ are an orthonormal basis for $\ell^2$, if and only if $U^1,U^2,\dots$ are an orthonormal basis for $\ell^2$.
\end{exercise}

One might think that \eqref{eqChangeOfOrderInInfiniteSums} follows from Fubini theorem for series. The following exercise shows that this is not the case. In fact, Fubini theorem requires that
\begin{equation}\label{eqGattoInArrivo}
\sum_{m,n}|y_m||A_{m,n}||x_n|<\infty,
\end{equation}
which is a condition involving the absolute values of the matrix entries. The requirement \eqref{eqGattoInArrivo} is that $\Tilde{A}=[|A_{m,n}|]_{m,n}$ defines a bounded operator on $\ell^2$.
\begin{proposition}\label{propAbsValInfiniteMatrix}
\begin{enumerate}
    \item[(i)] Let $\{x_n\}_{n=1}^\infty$ be a sequence of positive numbers, $x(n)\ge0$. Then, $x\in\ell^2$ if and only if
     \begin{equation}\label{eqGattoInArrivoTer}
\sum_{n}x(n)y(n)<\infty
\end{equation}
for all positive $y\in\ell^2$.   
    \item[(ii)] Let $A=[A(m,n)]_{m,n=1}^\infty$ be an infinite matrix with positive entries, $A(m,n)\ge0$. Then, $A$ defines a bounded operator on $\ell^2$ if and only if 
    \begin{equation}\label{eqGattoInArrivoter}
\sum_{m,n}x(m)A(m,n)y(n)<\infty
\end{equation}
for all positive $x,y\in\ell^2$.
\end{enumerate}
\end{proposition}
\begin{proof}
\noindent (i) The only if direction follows from Cauchy-Schwarz. Suppose \eqref{eqGattoInArrivoTer} holds. Then, it is easy to show that $x(n)\to0$ as $n\to\infty$. Suppose by contradiction that $\sum_nx(n)^2=\infty$: we will find a sequence $\epsilon(n)\searrow 0$ such that (a) $\sum_nx(n)^2\epsilon(n)=\infty$, but (b) $\sum_n x(n)^2\epsilon(n)^2<\infty$. If we do this, the sequence $y(n)=x(n)\epsilon(n)$ belongs to $\ell^2$, but fails to satisfy \eqref{eqGattoInArrivoTer}. Let's do the construction, assuming with no loss of generality that $x(n)\le1$ for all $n$. Find an increasing sequence $n_j$ ($n_1=0$) such that
\[
1\le \sum_{n=n_{j-1}+1}^{n_j}x(n)^2\le 2,
\]
and let $\epsilon(n)=\frac{1}{j}$ for $n_{j-1}+1\le n\le n_j$. Then, both (a) and (b) hold.

(ii) Again, the only if direction is elementary. Suppose now \eqref{eqGattoInArrivoter} holds. By (i), $m\mapsto \sum_nA(m,n)y(n)$ and $n\mapsto \sum_mA(m,n)x(m)$  are $\ell^2$ sequences,
\[
\sum_m\left(\sum_nA(m,n)y(n)\right)^2<\infty,\ \sum_n\left(\sum_mA(m,n)x(m)\right)^2<\infty.
\]
Thus, both $A$ and $A^\ast$ are everywhere defined, and so is $A^\ast A$, which is self-adjoint. By the Hellinger-Toeplitz theorem, $A^\ast A$ is bounded, hence
\[
\|Ax\|_{\ell^2}^2=\langle x,A^\ast Ax\rangle\le\|A^\ast A\|_{\mathcal{L}(\ell^2)}\|x\|_{\ell^2}^2,
\]
which shows that $A$ is bounded.
\end{proof}
\begin{exercise}\label{exeAbsoValDoesNotGetIn} Construct a bounded operator on $\ell^2(\mathbb{N})$ such that \eqref{eqGattoInArrivo} does not hold. Here is possible pathway.
    \begin{enumerate}
        \item[(i)] Consider the {\bf discrete Fourier transform} $\mathcal{F}_N:\ell^2(\{0,1,\dots,N-1\})$,
        \begin{equation}\label{eqDiscreteFTdomenica}
            (\mathcal{F}_Nx)_m:=\frac{1}{\sqrt{N}}\sum_{n=0}^{N-1}e^{2\pi i mn/N}x_n.
        \end{equation}
    \end{enumerate}
That is, the matrix elements of $\mathcal{F}_N$ are $(\mathcal{F}_N)_{m,n}=\frac{e^{2\pi i mn/N}}{\sqrt{N}}$.
    
    Show that $\mathcal{F}_N$ is a unitary operator.
    \item[(ii)] Consider the operator $A_N$ having as entries the moduli of the entries of $\mathcal{F}_N$, $(A_n)_{m,n}=\frac{1}{\sqrt{N}}$. Show that $\|A_n\|=\sqrt{N}$.
    \item[(iii)] Construct a unitary operator $A$ on $\ell^2(\mathbb{N})$ in block diagonal form, where the blocks are $\mathcal{F}_N$ with increasing $N$. Show that $A$ does not satisfy \eqref{eqGattoInArrivo}.
    \item[(iv)] Modify the construction to construct an operator $A$ whjich is bounded and self-adjoint, yet \eqref{eqGattoInArrivo} fails.
\end{exercise}
\subsection{Hilbert-Schmidt operators on Hilbert spaces}\label{SSectHSop}
Let $A$ be a bounded operator on a separable Hilbert space $H$, $\{e_n\}$ be a orthonormal basis for $H$, and $\{A_{m,n}\}_{m,n}$ be the matrix coefficients of $A$ with respect to the basis. The {\it Hilbert-Schmidt norm} of $A$ is\index{operator!Hilbert-Schmidt}
\begin{equation}\label{eqHSnorm}
    \|A\|_{HS}^2:=\sum_{m,n}|A_{m,n}|^2.
\end{equation}
\begin{proposition}\label{propHSnormIndepOfBasis}
    The Hilbert-Schmidt norm of $A$ is independent of the basis.
\end{proposition}
\begin{proof}
    The statement is equivalent to show that, if $U:\ell^2\to\ell^2$ is unitary, then
    \[
    \|U^\ast A U\|_{HS}=\|A\|_{HS},
    \]
    which is what we do below. The matrix coefficients of $U^\ast A U$ are
    \[
    (U^\ast A U)_{m,n}=\sum_{i,j}U^\ast_{m,i}A_{i,j}U_{j,n}=\sum_{i,j}\overline{U_{i,m}}A_{i,j}U_{j,n}.
    \]
    We can freely change the order of summation by corollary \ref{corMatOpOne}:
    \begin{eqnarray*}
      \|U^\ast A U\|_{HS}^2&=&\sum_{m,n,i,j,k,l}U_{i,m}\overline{A_{i,j}}\overline{U_{j,n}}\overline{U_{l,m}}A_{l,k}U_{k,n}\crcr 
      &=&\sum_{i,j,k,l}\left(\sum_m U_{i,m}\overline{U_{l,m}}\right)\overline{A_{i,j}}\left(\sum_n \overline{U_{j,n}}U_{k,n}\right)A_{l,k}\crcr 
      &=&\sum_{i,j,k,l}\delta_{i,l}\overline{A_{i,j}}\delta_{j,k}A_{l,k}\crcr 
      &=&\sum_{i,k}\overline{A_{i,k}}A_{i,k}\crcr 
      &=&\|A\|_{HS}^2.
    \end{eqnarray*}
\end{proof}
A bounded operator $A$ is in the {\it Hilbert-Schmidt class}, $A\in HS(H)$, if $\|A\|_{HS}<\infty$. The map $A\mapsto[A_{m,n}]$ is unitary identification of $HS$ with $\ell^2(\mathbb{N}_\ast\times \mathbb{N}_\ast)$. In particular, $HS(H)$ has a Hilbert space structure. We also have
\begin{equation}\label{eqHSgreaterNorm}
    \|A\|\le\|A\|_{HS}.
\end{equation}
Also observe that membership in $HS(H)$ only depends on the moduli $|A_{m,n}|$ of the matrix coefficients (compare with exercise \ref{exeAbsoValDoesNotGetIn}).
We have the following important result.
\begin{theorem}\label{theoHSisCompact}
    Hilbert-Schmidt operators are compact.
\end{theorem}
\begin{proof}
    For $M\ge1$ fixed, let $A^{(M)}$ be the infinite matrix with
    \[
    A^{(M)}_{m,n}=\begin{cases}
        A_{m,n}\text{ if }1\le m\le M\crcr 
        0 \text{ if }m>M.
    \end{cases}
    \]
    Then $A^{(M)}$ represents a finite rank operator, and
    \[
    \|A-A^{(M)}\|^2\le \|A-A^{(M)}\|_{HS}^2=\sum_{m>M}\sum_{n\ge1}|A_{m,n}|^2\to0
    \]
    as $M\to\infty$ by dominated convergence.
\end{proof}
\begin{exercise}\label{exeNotAllCompactAreHS}
    \begin{enumerate}
    \item[(i)] Characterize self-adjoint Hilbert-Schmidt operators in terms of their spectrum.
    \item[(ii)] Find a compact operator which is not Hilbert-Schmidt.         
    \end{enumerate}
\end{exercise}
The Hilbert-Schmidt class can be introduced in a coordinate-free way, but we will not do it here.
\subsection{The trace of a compact operator}\label{SSectOperatorTrace}\index{operator!trace}
Let $A\ge0$ be a positive, compact operator on a separable Hilbert space $H$, and let $\{e_n\}$ be a orthonormal basis of $H$. The {\it trace} of $A$ is
\begin{equation}\label{eqOPeratorTrace}
    \text{Trace}(A)=\sum_{n=1}^\infty \langle e_n,A e_n\rangle\ge0.
\end{equation}
It is readily verifies that the trace is a unitary invariant. We use matrix notation $\langle e_n,A e_n\rangle=A_{nn}$. If $U=[U_{ij}]$ is unitary, and $B=U^\ast AU$, then
\begin{eqnarray*}
\sum_nB_{nn}&=&\sum_{njk} U^\ast_{nj}A_{jk}U_{kn}=\sum_{jk}A_{jk}\sum_n U^\ast_{nj}U_{kn}\crcr 
&=&\sum_{jk}A_{jk}\delta_{jk}\crcr 
&=&\sum_nA_{nn}.
\end{eqnarray*}
\section{Compact integral operators}\label{sectCompactIO}
The study of compact operators had its inception in the study of integral operators, especially those arising from the analysis of differential equations. In this section, we consider some important examples.
We are interested in operators having the form
\begin{equation}\label{eqIntegralOperator}
    T_kf(x)=\int_Y k(x,y)f(y)d\nu(y),
\end{equation}
where the input is $f:Y\to\mathbb{C}$ is defined on a measure space $(Y,\nu)$; $k:X\times Y\to\mathbb{C}$ has some regularity properties; and the output is $T_kf:X\to\mathbb{C}$. The function $k$ is the {\it kernel} of the integral operator $K$.\index{operator!integral kernel} 
\subsection{Integral operators with continuous kernel}\label{ssectIOcontinuousK}
We consider here the case where the kernel $k$ is continuous. Statement and proof are on closed intervals, and the extension to a more general framework is left as an exercise in critical analysis of proofs.
\begin{theorem}\label{theoCompactIntOpBasic}
    Let $k:[a,b]\times[c,d]\to\mathbb{C}$ be continuous, and define
    \[
    T_kf(x)=\int_c^dk(x,y)f(y)dy.
    \]
    Then,
    \begin{enumerate}
        \item[(i)] $T_k:C[a,b]\to C[a,b]$ is compact.
        \item[(ii)] For $1\le p\le\infty$, $T_k:L^p[a,b]\to C[a,b]$ is compact.
        \item[(iii)] For $1\le p<\infty$, $T_k:L^p[a,b]\to L^p[a,b]$ is compact.
    \end{enumerate}
\end{theorem}
\begin{proof}
The main tool here is Ascoli-Arzelà theorem\index{theorem!Ascoli-Arzelà} and we use that kernel $k$ is uniformly continuous on $[a,b]\times[c,d]$ by compactness. We start with (i), proving at the same time that $K$ maps $C[a,b]$ into itself, and that it is a compact operator.  For real $\delta$ and $f$ continuous with $\|f\|_{L^\infty}\le1$, in fact,
\begin{eqnarray*}
   |T_kf(x+\delta)-T_kf(x)|&\le&\int_c^d|k(x+\delta,y)-k(x,y)|\cdot|f(y)|dy\crcr 
   &\le&\epsilon(d-c)\|f\|_{L^\infty}\le \epsilon[d-c],
\end{eqnarray*}
for any $\epsilon>0$, provided $|\delta|\le\eta(\epsilon)$ is such that $|k(x+\delta,y)-k(x,y)|\le\epsilon$ for all $(x,y)\in[a,b]\times[c,d]$. Hence, the family $\{T_kf:\|f\|_{L^\infty}\le1\}$ is \index{equicontinuous}equicontinuous, and also uniformly bounded because
\[
|T_kf(x)|\le (d-c)\max_{(x,y)\in[a,b]\times[c,d]}|k(x,y)|\cdot \|f\|_{L^\infty}.
\]
By the Ascoli-Arzelà theorem, the family has compact closure in $C[a,b]$ (with respect to the norm topology).

The proof of (ii) is similar. Using H\"older inequality with $1/p+1/p'=1$, if $\|f\|_{L^p}\le1$,
\begin{eqnarray*}
    |T_kf(x+\delta)-T_kf(x)|&\le&\int_c^d|k(x+\delta,y)-k(x,y)|\cdot|f(y)|dy\crcr
    &\le&\left(\int_c^d|k(x+\delta,y)-k(x,y)|^{p'}dy\right)^{1/p'}\|f\|_{L^p}\crcr 
    &\le&(d-c)^{1/p'}\epsilon \|f\|_{L^p},
\end{eqnarray*}
by choosing $\delta$ as in (i). We have equicontinuity, and uniform boundedness is proved similarly.

(iii) follows from (ii) and the fact that the identity $I:C[a,b]\to L^[a,b]$ is bounded, and the composition of a compact map and a bounded map is compact.
\end{proof}
\subsection{Hilbert-Schmidt integral operators}\label{ssectIOhilbertSchmidt}\index{operator!Hilbert-Schmidt}
Let $k:[c,d]\times[a,b]\to\mathbb{C}$, $k\in L^2([c,d]\times[a,b])$. To each function $f\in L^2[a,b]$ we associate, as before,
\begin{equation}\label{eqHilbertSchmidt}
    T_kf(x)=\int_a^bk(x,y)f(y)dy,\ T_kf:[c,d]\to\mathbb{C}.
\end{equation}
\begin{theorem}\label{theoHSintegralOp}
    The operator $T_k:L^2[a,b]\to L^2[c,d]$ is a Hilbert-Schmidt operator.
\end{theorem}
\begin{proof}
    Boundedness follows, as earlier, from Cauchy-Schwarz and Fubini:
    \begin{eqnarray*}
        \|T_kf\|_{L^2}^2&=&\int_c^d\left|\int_a^b k(x,y)f(y)dy\right|^2dx\crcr 
        &\le&\int_c^d\left(\int_a^b |k(x,y)|^2dy\cdot \int_a^b |f(z)|^2dz\right)dx\crcr
        &=&\iint_{[c,d]\times[a,b]}|k(x,y)|^2dxdy\cdot \int_a^b |f(z)|^2dz.
    \end{eqnarray*}
    Let now $\{e_n\}$ be a orthonormal basis for $L^2[a,b]$ and $\{\epsilon_m\}$ be one for $L^2[c,d]$, so that $\{e_n(y)\epsilon_m(x)\}_{m,n}$ is a orthonormal basis for $L^2([a,b]\times[c,d])$, and we can expand
    \begin{align*}
    k(x,y)&=\sum_{m,n}k_{n,m}e_n(y)\epsilon_m(x),\crcr 
    \|k\|_{L^2([a,b]\times[c,d])}^2&=\sum_{m,n}|k_{n,m}|^2<\infty.
    \end{align*}
    Then, if $f(y)=\sum_n f_n\cdot e_n(y)$,
    \[
        Kf(x)=\sum_{m} \epsilon_m(x)\sum_n k_{n,m}f_n.
    \]
    With respect to the basis $\{e_n\}$ and $\{\epsilon_m\}$, the operator $T_k$ is represented by the Hilbert-Schmidt matrix $[k_{n,m}]_{m,n=1}^\infty$, hence it is a Hilbert-Schmidt operator.
\end{proof}
Statement and proof hold for more general measure spaces.
\subsection{The Volterra operator}\label{ssectVolterraOp}
The {\it Volterra operator} $V$ is defined as
\begin{equation}\label{eqVolterra}
    Vf(x)=\int_0^xf(t)dt,
\end{equation}
where $f\in L^1[0,1]$. It has a holomorphic avatar, with $f:\mathbb{D}\to\mathbb{C}$, which we will consider below. We will see some basic properties of $V$ as operator on $L^2[0,1]$.
\begin{theorem}\label{theoVolterraReal}
    \begin{enumerate}
        \item[(i)] The operator $V$ is compact on $L^2[0,1]$ (in fact, it a Hilbert-Schmidt integral operator). Its range is
        \[
        \{F\in AC[0,1]:F(0)=0\text{ and }F'\in L^2[0,1]\},
        \]
        where $AC[0,1]$ is the class of the absolutely continuous functions.
        \item[(ii)] $V$ does not have eigenvalues, and $\sigma(V)=\{0\}$.
        \item[(iii)] The adjoint $V^\ast$ of $V$ is
        \begin{equation}\label{eqVolterraAst}
            V^\ast g(x)=\int_x^1g(y)dy,
        \end{equation}
        and
        \begin{equation}\label{eqVolterraAstVolterra}
            V^\ast Vf(x)=\int_0^1(1-\max(x,y))f(y)dy.
        \end{equation}
    \end{enumerate}
\end{theorem}
\begin{proof} We show (i), (iii) first.
    Since we can write $Vf(x)=\int_0^1k(x,y)f(y)dy$ with $k(x,y)=\chi(0\le y\le x\le 1)$, and $k\in L^2[0,1]^2$, $V$ is a Hilbert-Schmidt operator. Equations \eqref{eqVolterraAst} and \eqref{eqVolterraAstVolterra} follow by, respectively, observing that $V^\ast$ has kernel $k^\ast(x,y)=k(y,x)=\chi(0\le x\le y\le1)$, and computing the kernel
    \[
    h(x,y)=\int_0^1k^\ast(x,u)k(u,y)du.
    \]
    About (ii), $\lambda$ is an eigenvalues of $V$ if and only if the equation
    \[
    F=\lambda F'
    \]
    has a nonzero solution $F$ with $F(0)=0$. You might worry about the fact that a priori derivatives exist just almost everywhere; but $Vf=\lambda f$ implies that $f$ is continuous, hence it is $C^1$. For $\lambda=0$, $F=0$, hence $ker(F)=\{0\}$. For the other $\lambda$'s, $F(x)=C e^{x/\lambda}$, which satisfies $F(0)=0$ only if $F=0$.
\end{proof}
The {\it holomorphic Volterra operator} is defined on holomorphic functions $f:\mathbb{D}\to\mathbb{C}$,  
\begin{equation}\label{eqVolterraHol}
    Vf(z)=\int_0^zf(w)dw.
\end{equation}
We use the same symbol $V$, although the operator is not the one we had considered earlier.
We might consider it as an operator defined on the Taylor coefficients of $f$,
\[
f(z)=\sum_{n=0}^\infty a_nz^n\implies Vf(n)=\sum_{n=0}\frac{a_n}{n+1}z^{n+1},
\]
which we might write as
\begin{equation}\label{eqVolterraHolSeq}
    \widehat{V}(\{a_n\}_{n=0}^\infty)=\left\{\frac{a_{n-1}}{n}\right\}_{n=0}^\infty,
\end{equation}
or $\widehat{V}=\tau_1\circ S$, where $S(\{a_n\})=\{a_n/(n+1)\}$ is compact.
\begin{theorem}\label{theoVolterraHol}
    \begin{enumerate}
        \item[(i)] The operator $V$ is a Hilbert-Schmidt (hence, compact) on $H^2(\mathbb{D})$.
        \item[(ii)] $V$ does not have eigenvalues, and $\sigma(V)=\{0\}$.
        \item[(iii)] The adjoint $V^\ast$ of $V$ is
        \begin{equation}\label{eqVolterraAstO}
            \widehat{V}^\ast(\{b_n\}_{n=0}^\infty) =\left\{\frac{b_{n+1}}{n+1}\right\}_{n=0}^\infty.
            \end{equation}
        \end{enumerate}
\end{theorem}
\begin{exercise}\label{exeVolterraHol}
    Prove theorem \ref{theoVolterraHol}. Write explicit expressions for $\widehat{V}^\ast\widehat{V}$, and for $\widehat{V}\widehat{V}^\ast$.
\end{exercise}
\section{Intermezzo: Sturm-Liouville theory}\label{SectSturmLiouville}\footnote{For this section I basically follow the exposition in \href{https://www.uio.no/studier/emner/matnat/math/MAT4400/v21/sturml.pdf}{AN INTRODUCTION TO STURM-LIOUVILLE THEORY
} by ERIK BÉDOS}.\index{Sturm-Liouville theory}
In this section we see an application of the spectral theory for self-adjoint compact operators to ODE's. The differential operator appearing in the ODE turns out to be the inverse of a (compact) integral operator. The eigenvalues of the latter will have a role in the analysis of the former.

\subsection{The Sturm-Liouville problem, spaces and operators}\label{SSectSturmLiouSpOp}
In its general form, Sturm-Liouville theory deals with the differential equations
\begin{equation}\label{eqStLiou}
-(p y')'+qy=\lambda\rho y,    
\end{equation}
with $y\in C^2[a,b]$; $\lambda\in\mathbb{C}$; $p\in C^1[a,b]$ real valued with $p(x)\ne0$ in $[a,b]$; $q,\rho\in C[a,b]$ real valued with $\rho(x)\ne0$ in $[a,b]$; and non-degenerate boundary conditions
\begin{equation}\label{eqStLiouBC}
\alpha_0 y(a)+\alpha_1 y'(a)=0=\beta_0 y(b)+\beta_1 y'(b),
\end{equation}
where $(\alpha_0,\alpha_1)\ne(0,0)\ne(\beta_0,\beta_1)$ have real entries. Such conditions might be read as orthogonality relations in $\mathbb{R}^2\subset\mathbb{C}^2$. The problem is finding information about {\it eigenvalues}: the $\lambda$'s such that \eqref{eqStLiou}, \eqref{eqStLiouBC} has a nontrivial solution $y$.

We consider here the case $p\equiv1\equiv\rho$, which contains the main features of the theory:
\begin{align}\label{eqStLiouSimpler}
    \lambda y&=-y''+qy,\crcr 
    0&=\alpha_0 y(a)+\alpha_1 y'(a)=\beta_0 y(b)+\beta_1 y'(b).
\end{align}
Let $F\subset C^2[a,b]$ be the space of the functions satisfying the boundary conditions \eqref{eqStLiouBC}, and $D:F\to C[a,b]$ be the operator in \eqref{eqStLiouSimpler},
\begin{equation}\label{eqOperatorD}
Dy=-y''+qy.
\end{equation}
The Sturm-Liouville problem, then, consists in finding for which $\lambda\in\mathbb{C}$ the equation
\begin{equation}\label{eqPasseggiataNotturna}
    Dy=\lambda y,\text{ where }y\in F
\end{equation}
has a solution $y\ne 0$.

From elementary results on linear ODEs we have the following.
\begin{lemma}\label{lemmaSLodeFact}
    For $\lambda\in\mathbb{C}$, let $S_\lambda=\{y\in C^2[a,b]:D(y)=\lambda y\}$. Then, $\dim_\mathbb{C}(S_\lambda)=2$.

    In fact, for $c\in[a,b]$ and $\lambda\in\mathbb{C}$ fixed, the map $L_{\lambda,c}:S_\lambda\to\mathbb{C}^2$, $L_{\lambda,c}y=(y(c),y'(c))$, is an isomorphism.
\end{lemma}
Given $y_1,y_2\in C^2[a,b]$, the {\it Wronsky determinant} $W(y_1,y_2)\in C^1[a,b]$ is
\[
W(y_1,y_2)=
\begin{vmatrix}
   y_1&y_2\crcr 
   y_1'&y_2'
\end{vmatrix}.
\]
The following {\it Lagrange identity} is elementary.
\begin{lemma}\label{lemmaLagrangeSvWronsky}
    For $f,g\in C^2[a,b]$,
    \begin{equation}\label{eqLagrangeSvWronsky}
        W(f,g)'=D(f)g-fD(g).
    \end{equation}
    If $y_1,y_2$ are solution of $D(y)=0$, then, $W(y_1,y_2)=W$ is constant.
\end{lemma}
Having two linearly independent solutions of $D(y)=0$, which is \eqref{eqStLiouSimpler} with $\lambda=0$, using variation of parameters it is easy to write the general solution of
\begin{equation}\label{eqStLiouSimplerf}
    -y''+qy=f,
\end{equation}
where $f\in C[a,b]$. The solution can be written in a way which is symmetric with respect to the interval's endpoints.
\begin{lemma}\label{lemmaSLsolutionSymmetric}
    Let $u,v$ be two linearly independent solutions of $-y''+qy=0$, and let $f$ be continuous. Let $W$ be the Wronski determinant of $u,v$. Then, a solutions of \eqref{eqStLiouSimplerf} can be written in the form
    \begin{equation}\label{eqGreenFunctionForSL}
        y(x)=\int_a^b G(x,t)dt,\text{ where }G(x,t)=\frac{1}{W}\begin{cases}
            u(x)v(t)\text{ if }a\le t\le x\le b,\crcr 
            u(t)v(x)\text{ if }a\le x\le t\le b
        \end{cases}
    \end{equation}
    is symmetric, $G(x,t)=G(t,x)$, and continuous on $[a,b]\times[a,b]$.
\end{lemma}
\begin{proof}
  We look for a solution of \eqref{eqStLiouSimplerf} of the form $y=cu+dv$ for some functions $c,d$ to be chosen. We have:
 \begin{align}\label{eqVariazioniCostanti}
     y&=cu+dv,\text{ so that}\crcr 
     y'&=cu'+dv'\text{ if }c'u+d'v=0,\text{ and}\crcr 
     y''&=cu''+dv''-f\text{ if }u'c'+v'd'=-f.
 \end{align}
 It is immediately verified that, with these requirements, $y$ in fact solves \eqref{eqStLiouSimplerf}.
What we need is
\[
\begin{cases}
    c'u+d'v&=0\crcr 
    u'c'+v'd'&=-f
\end{cases}, \text{ i.e. }\begin{pmatrix}
    u&v\crcr 
    u'&v'
\end{pmatrix}\begin{pmatrix}
    c'\crcr 
    d'
\end{pmatrix}=\begin{pmatrix}
    0\crcr 
    -f
\end{pmatrix}.
\]
Solving, we obtain
\[
\begin{cases}
    c'&=\frac{fv}{W}\crcr 
    d'&=-\frac{fu}{W}
\end{cases}
\]
Usually, the next step would be writing the general solution in the form
\[
\Tilde{y}=\frac{1}{W}\int_{x_0}^x(u(x)v(t)-v(x)u(t))f(t)dt+\rho u(x)+\sigma v(x), 
\]
with real $\rho,\sigma$, but this expression is unnatural for our boundary conditions, since we have just one endpoint $x_0$ to play with. We consider instead the two "halves" of $\Tilde{y}$ separately, and we see what can be done with them.
\begin{align}\label{eqAandBseparately}
    A(x)&=\frac{u(x)}{W}\int_{a}^xv(t)f(t)dt\text{ and }B(x)=\frac{v(x)}{W}\int_{b}^xu(t)f(t)dt,\crcr 
    A'(x)&=\frac{u'(x)}{W}\int_{a}^xv(t)f(t)dt+\frac{(uvf)(x)}{W}\text{ and }\crcr 
    B'(x)&=\frac{v'(x)}{W}\int_{b}^xu(t)f(t)dt+\frac{(uvf)(x)}{W}\crcr 
    A''(x)&=\frac{u''(x)}{W}\int_{a}^xv(t)f(t)dt+\frac{[u'vf+(uvf)'](x)}{W}\text{ and }\crcr 
    B''(x)&=\frac{v''(x)}{W}\int_{b}^xu(t)f(t)dt+\frac{[v'uf+(uvf)'](x)}{W}.
\end{align}
The function $y$ in \eqref{eqGreenFunctionForSL} is $y(x)=A(x)-B(x)$, and it satisfies
\begin{eqnarray*}
    -y''+qy&=&\frac{(-u''+qu)(x)}{W}\int_{a}^xv(t)f(t)dt-\frac{(-v''+qv)(x)}{W}\int_{b}^xu(t)f(t)dt\crcr 
    &\ &+\frac{(v'u-uv')(x)f(x)}{W}\crcr 
    &=&f,
\end{eqnarray*}
since $v'u-uv'=W$.
\end{proof}
\subsection{The solution to the Sturm-Liouville problem}\label{SSectSTsolved}
\subsubsection{The solution of the Sturm-Liouville problem when $D$ is $1-1$}\label{SSSectSLifDisOneOne}
Here we suppose that $D:F\to C[a,b]$ is $1-1$. In general this is not the case, as the example $D(y)=-y''+y$ with conditions $y(0)=y(\pi)=0$ clearly illustrates: $D(\sin x)=0$ and $\sin$ verifies the boundary conditions. Observe that the same operator is $1-1$ on $F$ if the conditions are $y(0)=y(\pi/2)=0$. Injectivity of $D$ on $F$ means that the boundary problem 
\[
D(y)=0,\ \alpha_0y(a)+\alpha_1y'(a)=0=\beta_0y(b)+\beta_1y'(b)
\]
only has the trivial solution $y=0$.

We introduce spaces where either of the boundary conditions hold for the equation $D(y)=\lambda y$ (but here we are interested in $\lambda=0$),
\begin{align*}
L_\lambda&=\{y\in C^2[a,b]: D(y)=\lambda y\text{ and } \alpha_0y(a)+\alpha_1y'(a)=0\},\crcr  
R_\lambda&=\{y\in C^2[a,b]: D(y)=\lambda y\text{ and } \beta_0y(b)+\beta_1y'(b)=0\},
\end{align*}
so that $E_\lambda=L_\lambda\cap R_\lambda$ is the space of the solutions of \eqref{eqStLiouSimpler}. By lemma \ref{lemmaSLodeFact}, $\dim(L_\lambda)=\dim(R_\lambda)=1$, hence $\dim(E_\lambda)\le1$.  We are here assuming that $\dim(E_0)=0$, hence, there exist $v\in L_0$ and $u\in R_0$ which are linearly independent.
\begin{proposition}\label{propSLwhenDisOneOne}
    Suppose $D$ is $1-1$ and let $v\in L_0$, $u\in R_0$. Then, the solution of $D(y)=f$ in lemma \ref{lemmaSLsolutionSymmetric} is the only one satisfying the boundary conditions.
\end{proposition}
\begin{proof}
   Uniqueness is clear, since $D$ is $1-1$ on $F$. We have to verify that the boundary conditions are satisfied.
   \begin{align*}
    A(a)&=0\text{ and }B(a)=\frac{v(a)}{W}\int_{b}^au(t)f(t)dt,\crcr 
    A(b)&=\frac{u(b)}{W}\int_{a}^bv(t)f(t)dt\text{ and }B(b)=0,\crcr 
    A'(a)&=\frac{(uvf)(a)}{W}\text{ and }B'(a)=\frac{v'(a)}{W}\int_{b}^au(t)f(t)dt+\frac{(uvf)(a)}{W}\crcr 
    A'(b)&=\frac{u'(b)}{W}\int_{a}^bv(t)f(t)dt+\frac{(uvf)(b)}{W}\text{ and }B'(b)=\frac{(uvf)(b)}{W},
   \end{align*}
   hence
   \begin{eqnarray*}
      \alpha_0y(a)+\alpha_1y'(a)&=&-\frac{\alpha_0v(a)}{W}\int_{b}^au(t)f(t)dt-\frac{\alpha_1 v'(a)}{W}\int_{b}^au(t)f(t)dt\crcr 
      &=&0
   \end{eqnarray*}
   because $v\in L_0$, and similarly
   \begin{eqnarray*}
       \beta_0y(b)+\beta_1y'(b)&=&\frac{\beta_0 u(b)}{W}\int_{a}^bv(t)f(t)dt+\frac{\beta_1 u'(b)}{W}\int_{a}^bv(t)f(t)dt\crcr 
       &=&0
   \end{eqnarray*}
   because $v\in R_0$.
\end{proof}
The result can be rephrased in a more functional form.
\begin{corollary}\label{corSLwhenDisOneOne}
    If $D:F\to C[a,b]$ is $1-1$, then it is a bijection, and its inverse $T_G:C[a,b]\to F$ has the expression
    \begin{equation}\label{eqSLwhenDisOneOne}
        T_G(f)=\int_a^bG(x,t)f(t)dt,
    \end{equation}
    where $G:[a,b]\times[a,b]\to\mathbb{R}$ is that given in \eqref{eqGreenFunctionForSL}, and $u\in R_0$, $v\in L_0$. 
\end{corollary}
\begin{proof}
   We have verified above that $D\circ T_G(f)=f$. Since $D$ is injective, $T_G$ is surjective, and we also have $T_G\circ D(y)=y$ for $y\in F$. 
\end{proof}
Since the operator $T_G$ extends to a compact, self-adjoint operator on $L^2[a,b]$, we are now in a position where the Hilbert-Schmidt theorem can be applied.
\begin{theorem}\label{theoSLwhenOneOne}
    Suppose the operator $D$ is $1-1$ on $F$, and consider the Sturm-Liouville problem in \eqref{eqPasseggiataNotturna},
    \begin{equation}\label{eqDimenticata}
    D(y)=\lambda y,\ y\in F.
    \end{equation}
    Then, the following hold.
    \begin{enumerate}
        \item[(i)] The set $\{\lambda_n\}$ of the $\lambda$'s for which there is a solution $y\ne0$ is countable, each eigenvalue is real and nonzero, and $|\lambda_n|\to\infty$.
        \item[(ii)] For each $n$, the eigenspace $E_{\lambda_n}=\{y:\eqref{eqDimenticata}\text{ holds}\}$ is $1$-dimensional.
        \item[(iii)] Let $e_n$ be a unit vector in $E_{\lambda_n}$. Then, $\{e_n\}$ is a orthonormal basis for $L^2[a,b]$.
    \end{enumerate}
\end{theorem}
\begin{proof}
    The space $F$ is dense in $L^2[a,b]$ (see exercise below) and $F=T_G(C[a,b])$, hence,
    \[
    L^2[a,b]=\overline{F}\subseteq\overline{T:G(L^2[a,b])}\subseteq L^2[a,b]. 
    \]
    Thus, $\overline{T_:G(L^2[a,b])}$, and corollary \ref{corOfHS} implies that $\ker(T_G)=L^2[a,b]\ominus \overline{T_G(L^2[a,b])}=0$.

    We apply the spectral theorem to $T_G$, and find eigenvalues $\{\mu_n\}$ with $\mu_n\ne0$, and $\mu_n\to0$ as $n\to\infty$. Let $f_n\in L^2[a,b]$ be an eigenfunction for $\mu_n$,
    \[
    f_n=1/\mu_n T_G(f_n),\text{ hence, }f_n\in C[a,b], \text{ and }f_n=1/\lambda_n^2T_G(T_G f_n)\in F.
    \]
    Applying $D$ to both sides of the first equality,
    \[
    D(f_n)=1/\mu_n f_n.
    \]
    Set $\lambda=\lambda_n=1/\mu_n$. Then \eqref{eqDimenticata} holds for $y=f_n$. The eigenspace $E_{\lambda_n}$ is one dimensional, as observed earlier.

    We have to make sure that $\{\lambda_n\}$ exhausts the eigenvalues of $D$. If there were another $\lambda\ne0$ for which $D(y)=\lambda y$ has a nonzero solution in $F$, then $y=\lambda T_G(y)$, and $1/\lambda$ would be an other eigenvalue for $T_G$, which is not possible.

    Finally, by normalizing the $f_n$'s, $e_n=f_n/\|f_n\|_{L^2}$, we find a orthonormal basis of (real valued) eigenfunctions of $D$.
\end{proof}
\subsubsection{The solution of the Sturm-Liouville problem in general}\label{SSSectSLinGeneral}
Suppose $D(y)=-y''+qy$ does not define a $1-1$ operator on $F$. The results of the previous section, though, still hols if there is $\mu\in\mathbb{R}$ such that $D_\mu(y)=-y''+(q-\mu)y$ is injective (we just have to shift the old $\lambda$'s by $\mu$). Hence, we are done if we prove the following.
\begin{lemma}\label{lemmaShiftByMu}
    There exists $\mu\in\mathbb{R}$ which is not an eigenvalue of $D$ on $F$.
\end{lemma}
\begin{proof}
    If there are uncountably many eigenvalues, then there exists an uncountable orthonormal system in $L^2[a,b]$, which is separable. Contradiction.
\end{proof}
We have then the complete theorem.
\begin{theorem}\label{theoSLabOnlyNaevoVindicatus}
    Consider the Sturm-Liouville problem in \eqref{eqPasseggiataNotturna},
    \begin{equation}\label{eqDimenticataO}
    D(y)=\lambda y,\ y\in F.
    \end{equation}
    Then, the following hold.
    \begin{enumerate}
        \item[(i)] The set $\{\lambda_n\}$ of the $\lambda$'s for which there is a solution $y\ne0$ is countable, each eigenvalue is real and nonzero, and $|\lambda_n|\to\infty$.
        \item[(ii)] For each $n$, the eigenspace $E_{\lambda_n}=\{y:\eqref{eqDimenticataO}\text{ holds}\}$ is $1$-dimensional.
        \item[(iii)] Let $e_n$ be a unit vector in $E_{\lambda_n}$. Then, $\{e_n\}$ is a orthonormal basis for $L^2[a,b]$.
    \end{enumerate}
\end{theorem}

\chapter{Intermezzo: reproducing kernel Hilbert spaces}\label{chaptRKHS}

The {\it structure} of a complex Hilbert space $(H,\langle\cdot,\cdot\rangle)$ includes a vector space $H$ and an inner product $\langle\cdot,\cdot\rangle$ on $H$, plus the requirement that $H$ is complete with respect to the norm associated to the inner product (this is not especially demanding, since, if $H$ is not complete, its completion is a Hilbert space). All orthonormal basis of $H$ have  the same cardinality, the {\it (Hilbert) dimension} of $H$; and if $H_1,H_2$ are Hilbert, there exists a bounded, linear bijection $T:H_1\to H_2$ if and only if $H_1$ and $H_2$ have the same dimension. Moreover, $T$ can be taken to be {\it unitary},
\begin{equation}\label{eqUnitaryHilbert}
   \langle Tx,Ty\rangle_{H_2}=\langle x,y\rangle_{H_2}. 
\end{equation}
The dimension, that is, is a {\it unitary invariant} of the Hilbert space, and it completely classifies the elements of the Hilbert class modulo unitary equivalence.

In mathematics and its applications, however, the elements of the Hilbert spaces we work with are typically {\it functions} $f$, which are defined on some space $X$, taking values at its {\it points}. The expression {\it function Hilbert space} is often used. The points of $X$ (and, more generally, subsets of $X$) provide extra structure. The structure, that is, has the form $(H,\langle\cdot,\cdot\rangle,X)$.

The points of $X$ can enter the picture in many different ways. Consider the (unitarily equivalent) Hilbert spaces $\ell^2(\mathbb{N})$ and $L^2[0,1]$. The elements of the first space are functions $n\mapsto f(n)$ defined on $\mathbb{N}$, while the elements of the second, we learn in Lebesgue theory, are {\it equivalence classes} of functions $x\mapsto f(x)$. In the second case, the difficulty is that $\|f\|_{L^2}=0$ does not imply that $f=0$, but just that $f(x)=0\ a.e.x$.

This fact is not just an artifact depending on our construction of $L^2[0,1]$: it is a phenomenon which is intrinsic to $L^2[0,1]$ as a function Hilbert space: the value of a function $f\in L^2[0,1]$ at a point $x\in[0,1]$ can not be detected using the Hilbert structure. Suppose $f:[0,1]\to\mathbb{C}$ is continuous, so that $f(1)$ can be defined unambiguously. In fact, $\eta_1:f\mapsto f(1)$ is a linear functional defined on $C[0,1]$. Consider now the functions $f_n(x)=x^n$. We have that:
\[
\eta_1(f_n)=1, \text{ yet }\lim_{n\to\infty}\|f_n-0\|_{L^2}=0, \text{ and } \eta_1(0)=0.
\]
The {\it evaluation functional }$\eta_1$, that is, is not bounded in the $L^2$-norm on $C[0,1]$, and in no way can be extended to a bounded, linear functional on $L^2[0,1]$, no matter how $L^2[0,1]$ is thought of. (We might use Zorn's lemma to extend $\eta_1$ to an {\it unbounded} linear functional on $L^2[0,1]$, with no reasonable relation to evaluation at $x=1$ of functions which are not continuous).
\begin{exercise}\label{exeLtwoIsNotRKHS}
    Let $E\subseteq[0,1]$ be measurable with positive measure, and define $\eta_E:L^2[0,1]\to\mathbb{C}$ be the mean value of $f$ on $E$,
    \[
    \eta_E(f)=\frac{\int_E fdm}{m(E)},
    \]
    which might be thought of as an "average evaluation" on the set $E$. Show that $\|\eta_E\|_{H^\ast}=m(E)^{-1/2}\to\infty$ as $m(E)\to0$.
\end{exercise}

In $\ell^2(\mathbb{N})$ we have a completely different story, since, for any positive integer $m$, $|\eta_m(f)|=|f(m)|\le\|f\|_{\ell^2}$: evaluating $f$ at the point $m$ is a functional having finite energy, so to speak.

These examples are rather trivial, and do not provide sufficient motivation to develop a theory. In the 1930' Hilbert spaces where still a novelty, they provided a unified approach to a variety of problems, and unified in a simple way a great deal of apparently distant mathematical objects, phenomena, and techniques. Stefan Bergman, while developing ways to efficiently approximate conformal mappings, introduced a Hilbert space of holomorphic functions, which nowadays goes under his name: the {\it Bergman space}. Soon after, Nachman Aronszajn framed Bergman's ideas in a general theory, presented in \href{https://www.ams.org/journals/tran/1950-068-03/S0002-9947-1950-0051437-7/S0002-9947-1950-0051437-7.pdf}{Theory of reproducing kernels} (1950), which is still the one of best sources to start studying the subject.

From its inception, RKHS theory followed several paths: holomorphic function theory in one, then several variables; and statistics. At the end of the past century there was an explosion of new interest in both directions: in statistics in view of applications to machine learning, and in mathematical analysis after several breakthroughs concerning specific problems, and the realization that underlying there were structural properties of general families of RHHS's (Pick spaces, De Branges spaces...).

The renewed interest, and the speed with which results are nowadays communicated, opened new areas, perspectives, and cross-contamination between different fields.

\section{Equivalent definitions of RKHS}\label{SectDefRKHS}
In this subsection we see different, equivalent definitions of a reproducing kernel Hilbert space.

Let $X$ be a set and $H$ a Hilbert space of functions $f:X\to{\mathbb C}$. We say that $H$ has {\it bounded point evaluation} if for each $x$ in $X$ there the evaluation functional $\eta_x:f\mapsto f(x)$ is bounded on $H$.

By Riesz Representation Theorem, this holds if and only if for each $x$ there is a (unique) vector $k_x$ in $H$ such that
\[
f(x)=\langle  k_x, f\rangle.
\]
Moreover, choosing $f=k_x$ we see that $\|\eta_x\|_{H^\ast}=\|k_x\|_H$.

\vskip0.5cm

\noindent{\bf An example where point evaluation is not bounded} We use Zorn's lemma. Consider in $\ell^2(\mathbb{N})$ the orthonormal basis $\{\delta_n:n\in\mathbb{N}\}$ given by the Dirac's deltas, and use Zorn's lemma to extend it to a Hamel (algebraic) basis $\{v_\alpha\}_{\alpha\in I}$. Define then the non-bounded, linear functional $\lambda:\ell^2\to\mathbb{C}$ by its action on the Hamel basis,
\[
\lambda(v_\alpha)=\begin{cases}
    n\text{ if }v_\alpha=e_n\text{ for some }n,\crcr 
    0\text{ otherwise}.
\end{cases}
\]
Let then $X=\mathbb{N}\cup\{\lambda\}$ and to each $h\in\ell^2(\mathbb{N})$ associate $f_h:X\to\mathbb{C}$,
\[
f_h(x)=\begin{cases}
    \langle \delta_n,h\rangle_{\ell^2}\text{ if }x=n,\crcr 
    \lambda(h)\text{ if }x=\lambda.
\end{cases}
\]
The space $H$ of such functions is unitarily equivalent to $\ell^2(\mathbb{N})$ with respect to the inner product $\langle f_h,f_k\rangle_H:=\langle h,k\rangle_{\ell^2(\mathbb{N})}$, and evaluation at $\lambda$, $f_h(\lambda)=\lambda(h)$ is not bounded.

\vskip0.5cm

The functions $\{k_x\}_{x\in X}$ are the {\it kernel functions} for $H$, which is then called a {\it Reproducing kernel Hilbert space}, with kernel $k:X\times X\to{\mathbb C}$ defined by
\[
k(x,y):=\langle k_x, k_y\rangle=k_y(x).
\]
\begin{proposition}\label{PropRKHS}
The following properties hold.
\begin{enumerate}
    \item[(i)] $k(x,x)=\|k_x\|^2\ge0$, with equality if and only if $f(x)=0$ for all $f$ in $H$ (and such points can be removed).
    \item[(ii)] $k_x(y)=k(y,x)=\overline{k(x,y)}=\overline{k_y(x)}$.
    \item[(iii)] If $x_1,\dots,x_n\in X$ and $c_1,\dots,c_n\in{\mathbb C}$, then
    \[
    \sum_{i,j=1}^n c_i \overline{c_j}k(x_j,k_i)=\left\|\sum_{i=1}^n c_i k_{x_i}\right\|^2\ge0.
    \]
\end{enumerate}
Property (iii), by definition, say that $k:X\times X\to\mathbb{C}$ is a {\it positive semi-definite function}.
\end{proposition}
\begin{exercise}
Prove the proposition.
\end{exercise}
Property (iii) obviously implies that $k(x,x)\ge0$, hence that $\overline{k(x_1,x_2)}=k(x_2,x_1)$. In fact, setting $c_1=1$ and $c_2=i$:
\[
0\le k(x_1,x_1)+ik(x_1,x_2)-ik(x_2,x_1)+k(x_2,x_2)
\]
which implies that $\text{Re}(k(x_1,x_2))=\text{Re}(k(x_2,x_1))$, and setting $c_1=c_2=1$,
\[
0\le k(x_1,x_1)+k(x_1,x_2)+k(x_2,x_1)+k(x_2,x_2)
\]
which implies that  $\text{Im}(k(x_1,x_2))+\text{Im}(k(x_2,x_1))=0$. If in the definition we had just considered real scalars $c_1,\dots,c_n$, as it is sometimes done in the real theory, dealing with real kernels $k$, then the condition $k(x,y)=k(y,x)$ does not follow from the real version of (iii), as illustrated by the matrix
\[
K=\begin{pmatrix}
    1&1/2\crcr 
    0&1
\end{pmatrix},
\]
which can be seen as a non-symmetric kernel on a space $X$ with two points, having an associated quadratic form which is positive definite.

It is noteworthy that to each positive semi-definite function $k:X\times X\to\mathbb{C}$ we can associate a RKHS of which $k$ is the reproducing kernel.
\begin{theorem}[Moore–Aronszajn theorem]\label{theoRKHSfromPD}
    Let $k:X\times X\to\mathbb{C}$ be a positive semi-definite function on some set $X$, and define on $\text{span}(k_x:\ x\in X)$ by
    \begin{equation}\label{eqRKHSfromPD}
        \left\langle\sum_{i=1}^m a_ik_{x_i}\left\vert\sum_{j=1}^n b_jk_{y_j}\right.\right\rangle:=\sum_{i=1}^m\sum_{j=1}^n\Bar{a_i}b_j k(x_i,y_j).
    \end{equation}
    Then, the completion $H$ of $\text{span}(k_x:\ x\in X)$ with respect to the inner product is a reproducing kernel Hilbert space having kernel $k$.
\end{theorem}
\begin{proof}

\noindent{\bf Step I.}
Since $k$ is positive semi-definite, $\langle\cdot,\cdot\rangle$ defines a inner product on the function space $\text{span}(k_x:\ x\in X)$ provided we show that, for $f\in \text{span}(k_x:\ x\in X)$, $\langle f,f\rangle=0$ implies that $f(x)=0$ for all $x$ in $X$. Let $f=\sum_{i=1}^n a_ik_{x_i}$ and suppose that $\langle f,f\rangle=0$. If $x=x_{n+1}$ is any other point in $X$, we let $a_{n+1}=0$.

The matrix $K=[k(x_i,x_j)]_{i,j=0}^{n+1}$ is Hermitian and positive, then it can be written as $K=U^\ast \Delta U$, where $\Delta$ is the diagonal matrix having diagonal entries 
\[
\lambda_1\ge\lambda_2\ge\dots\ge\lambda_m>0=\lambda_{m+1}=\dots=\lambda_{n+1}.
\]
With $a\in\mathbb{C}^{n+1}$ and $a^\ast=\Bar{a}^t$, We are assuming that
\[
0=a^\ast Ka=(Ua)^\ast\Delta(Ua)=b^\ast \Delta b=\sum_{l=1}^m\lambda_l|b_l|^2,
\]
which implies that $b_1=\dots=b_m=0$, hence $\Delta b=0$, so
\[
0=U^\ast\Delta b=Ka.
\]
We have, then, for $j=1,\dots,x_{n+1}$:
\begin{eqnarray*}
    f(x_j)&=&\sum_{i=1}^n a_ik_{x_i}(x_j)=\sum_{i=1}^n a_ik(x_j,x_i)=\sum_{i=1}^n a_iK(j,i)\crcr 
    &=&(Ka)_j=0.
\end{eqnarray*}
Hence, $f=0$ as a function.

\noindent{\bf Step II.} The completion $H$ of $\text{span}(k_x:\ x\in X)$ under the norm induced by the inner product is a Hilbert space, having as generators the functions $k_x$. Each element $f$ of it can be interpreted as a function on $X$,
\[
f(x):=\langle k_x,f\rangle.
\]
The function $k_x\in H$ is then reproducing for the value of $f$ at $x$, and, as we saw above, the corresponding reproducing kernel is $\langle k_x,k_y\rangle=k(x,y)$, as stated.
\end{proof}
The theorem above is somehow deceptive. Its proof, in fact, gives little hint on how to find a "natural" expression for the inner product associated with a positive semi-definite function. We will see below that $k(z,w)=\frac{1}{1-\Bar{w}z}$ defines a positive semi-definite function, hence a kernel, on the unit disc $\mathbb{D}$ of the complex plane, but at this moment is not clear how the inner product really looks like. We have, however, a characterization of reproducing Hilbert spaces based on orthonormal basis, which is sometimes useful in this respect.
\begin{theorem}\label{theoRKHSfromONB}
    Let $H$ be a reproducing kernel Hilbert space on some set $X$, and let $\{e_\alpha\}_{\alpha\in I}$ be a orthonormal basis of it. Then, for $y\in X$, the series
    \begin{equation}\label{eqRKHSfromONB}
        k_y:=\sum_{\alpha\in I}\overline{e_\alpha(y)}e_\alpha 
    \end{equation}
    converges in $H$, and it converges to the reproducing function at $y$. 
    
    Viceversa, suppose $H$ is a Hilbert function space on $X$ and that the series \eqref{eqRKHSfromONB} converges in $H$. Then $k(x,y)=k_y(x)$ reproduces the values of functions in $\text{span}(e_\beta:\beta\in I)$: 
    \begin{equation}\label{eqRKfromONB}
        k(x,y)=\sum_{\alpha\in I}\overline{e_\alpha(y)}e_\alpha(x),\ \langle k_x,e_\alpha\rangle=e_\alpha(x) \text{ for }\alpha\in I. 
    \end{equation}    
    Moreover, the series in \eqref{eqRKfromONB} converges pointwise for all $x,y\in X$. After associating to each function $f$ in $H$ its modification $\Tilde{f}(x):=\langle k_x,f\rangle$, the space $\Tilde{H}$ of such $\Tilde{f}$'s is a reproducing Hilbert space with kernel $k=k(x,y)$.
\end{theorem}
The reproducing kernel is uniquely determined by the structure of RKHS, hence it does not depend on the orthonormal basis under consideration. Also, we have further evidence that $L^2[0,1]$ is not a RKHS. The series
\[
y\mapsto\sum_{n=-\infty}^{+\infty}e^{-2\pi i ny}e^{2\pi i n x}
\]
is not, in fact, convergent in $L^2[0,1]$.
\begin{proof}
Observe that, by basic properties of the orthonormal basis, the series on the right of \eqref{eqRKfromONB} converges if and only if
\[
\sum_{\alpha\in I},e_\alpha(y),^2 \text{ converges in }\mathbb{R}_+,
\]
or, equivalently, if for all $x,y\in X$
\[
\sum_{x,y}\overline{e_\alpha(x)}e_\alpha(y) \text{ converges absolutely in }\mathbb{C}.
\]
Suppose $H$ is a RKHS with kernel $k$. Then,
\[
k_y=\sum_{\alpha\in I}\langle e_\alpha,k_y\rangle e_\alpha=\sum_{\alpha\in I}\overline{\langle k_y e_\alpha\rangle} e_\alpha
=\sum_{\alpha\in I}\overline{e_\alpha(y)}e_\alpha
\]
converges in $H$.

Suppose, viceversa, that the series in \eqref{eqRKHSfromONB} converges in $H$. Then,
    \begin{equation}\label{eqTildeNoTilde}
    \langle k_y,e_\beta\rangle=\left\langle\sum_{\alpha\in I}\overline{e_\alpha(y)}e_\alpha \left|e_\beta\right.\right\rangle
    =\sum_{\alpha\in I}{e_\alpha(y)}\left\langle e_\alpha \left|e_\beta\right.\right\rangle=e_\beta(y).
    \end{equation}
    Thus, $k_y$ is the reproducing function for the basis elements, hence for all elements in $\text{span}(e_\beta:\beta\in I)$. Also, the function $k(x,y)=\langle k_x,k_y\rangle$ is positive semi-definite:
    \begin{eqnarray*}
        \left\langle\sum_{i=1}^nc_i k_{x_i}\left|\sum_{i=j}^nc_i k_{x_j}\right.\right\rangle&=&\sum_{i,j=1}^n\sum_{\alpha\in I}\overline{c_i}c_j e_\alpha(x_i)\overline{e_\alpha(x_j)}\crcr 
        &=&\sum_{\alpha\in I}\left|\sum_{i=1}^nc_i \overline{e_\alpha(x_i)}\right|^2\crcr 
        &\ge0&.
    \end{eqnarray*}
    By Moore-Aronszajn theorem, $k$ is the reproducing kernel of $\Tilde{H}$ and $\Tilde{e_\alpha}=e_\alpha$ by \eqref{eqTildeNoTilde}.
\end{proof}
Finally, we can see reproducing kernel Hilbert spaces "from the inside" of a given, possibly abstract Hilbert space. This viewpoint is especially useful in applications to machine learning. Observe that to a RKHS $H$ on $X$ with reproducing functions $\{k_x:x\in X\}$ we can associate the map $x\mapsto k_x$ (the {\it feature map} of machine learning) which associates to each point in $X$ a vector in $H$. Below, we see that the viceversa holds as well.
\begin{theorem}\label{theoRKHSfromVectorSets}
    Let $H$ be a Hilbert space, and let $\{v_x:x\in X\}$ be a family of vectors in $H$, indexed by some set $X$, such that $\text{span}(v_x:x\in X)$ is dense in $H$. To each $h$ in $H$ associate $f_h:x\mapsto \langle v_x,h\rangle_H$, let $H_X$ be the vector space of such functions, and set $\langle f_h,f_k\rangle_{H_X}:=\langle h,k\rangle_H$. Then, $H_X$ is a RKHS having $\{f_{v_x}:x\in X\}$ as kernel functions.
\end{theorem}
\begin{proof} By definition, $f_{v_x}(y)=\langle v_y,v_x\rangle_H$, and $(x,y)\mapsto k(x,y):=\langle v_x,v_y\rangle_H$ is positive definite on $X$, as is easily verified. Still by definition we have that
\[
f_h(x)=\langle v_x,h\rangle_H=\langle f_{v_x},f_h\rangle_{H_X},
\]
hence, $f_{v_x}$ is the reproducing function for $x$ in $H_X$, as wished.
\end{proof}
At the end of the day, we have several, equivalent viewpoints on RKHS's.
\begin{enumerate}
    \item[(i)] Hilbert function spaces with bounded point evaluation.
    \item[(ii)]  Hilbert function spaces with a reproducing kernel.
    \item[(iii)] Positive semi-definite functions on sets.
    \item[(iv)] Hilbert function spaces where the values at each point of an orthonormal basis' elements are square integrable. 
    \item[(v)] Families of vectors in a Hilbert space.
\end{enumerate}
\section{Some RKHS's}\label{SectRKHSexamples}
\subsection{The Hardy space $H^2$ and its real counterpart}\label{SSectHardyRKHS} 
We denote by $\mathbb{D}=\{z\in\mathbb{C}:\ |z|<1\}$ the unit disc in the complex plane. Let $H^2=H^2(\mathbb{D})$ be the space of the power series $f(z)=\sum_{n=0}^\infty a_n z^n$, $z\in{\mathbb{C}}$, for which the norm
\[
\left\|\sum_{n=0}^\infty a_n z^n\right\|_{H^2}^2=
\left\|\{a_n\}_{n=0}^\infty\right\|_{\ell^2}^2=\sum_{n=0}^\infty|a_n|^2
\]
is finite. In analogy with Fourier theory, in fact consistently with it, we write $a_n=\widehat{f}(n)$.
\begin{lemma}\label{LemmaHTwoConverge}
If $f\in H^2$, then the series defining $f$ converges on ${\mathbb D}$.
\end{lemma}
\begin{proof}
If $m,j\ge1$,
\begin{eqnarray*}
\left|\sum_{n=m+1}^{m+j}a_n z^n\right|^2&\le&\sum_{n=m+1}^{m+j}\left|a_n\right|^2 \cdot
\sum_{n=m+1}^{m+j}|z|^{2n}\crcr
&\le&\sum_{n=m+1}^{m+j}\left|a_n\right|^2 \cdot \frac{|z|^{2(m+1)}}{1-|z|^2}
\end{eqnarray*}
if $|z|<1$, hence Cauchy criterion is verified.
\end{proof}
By the theory of series in a complex variable, $f(z)$ defines a holomorphic function in ${\mathbb D}$.
\begin{proposition}\label{PropHardyKernel}
The reproducing kernel of the Hardy space is
\[
k_z(w)=\frac{1}{1-\overline{z}w}.
\]
In particular, the best estimate for the value of a Hardy function at $z\in{\mathbb D}$ is:
\[
|f(z)|\le\frac{\|f\|_{H^2}}{(1-|z|^2)^{1/2}}.
\]
\end{proposition}
\begin{proof}
\begin{eqnarray*}
\sum_{n=0}^\infty \widehat{f}(n) z^n&=&f(z)\crcr
&=&\left\langle f|k_z \right\rangle_{H^2}= 
\left\langle \sum_{n=0}^\infty \widehat{f}(n) w^n | 
\sum_{n=0}^\infty \widehat{k_z}(n) w^n\right\rangle_{H^2}\crcr
&=&\sum_{n=0}^\infty \widehat{f}(n)\overline{\widehat{k_z}(n)},
\end{eqnarray*}
which is equivalent to $\overline{\widehat{k_z}(n)}=z^n$,
hence,
\[
k_z(w)=\sum_{n=0}^\infty \widehat{k_z}(n)w^n
=\sum_{n=0}^\infty \overline{z}^n w^n=\frac{1}{1-\overline{z}w}.
\]
\end{proof}
An interesting variation on $H^2$ is  $h^2$, the {\it real (or harmonic) Hardy space}: the space populated by the series $u(r e^{it})=\sum_{n=-\infty}^{+\infty}a_n r^{|n|}e^{i n t}$, $0\le r<1$, such that $\|u\|_{h^2}^2=\sum_{n=-\infty}^{+\infty}|a_n|^2<\infty$. Another way to write a function in $h^2$ is
\begin{equation}\label{eqRealHardyFunction}
    u(z)=\sum_{n=0}^\infty a_n z^n+\sum_{n=1}^\infty a_{-n}\Bar{z}^n.
\end{equation}
We argue as above in order to find the reproducing kernel:
\[
u(w)=\sum_{n=0}^\infty a_n \overline{\Bar{w}^n}+\sum_{n=1}^\infty a_{-n}\overline{w^n}=\langle k_w,u\rangle_{h^2},
\]
where
\begin{eqnarray*}
    k_w(z)&=&\sum_{n=0}^\infty z^n\Bar{w}^n+\sum_{n=1}^\infty a_{-n}\Bar{z}^n w^n\crcr 
    &=&\frac{1}{1-\Bar{w}z}+\frac{1}{1-w\Bar{z}}-1\crcr 
    &=&\frac{1-|z|^2|w|^2}{|1-\Bar{w}z|^2}.
\end{eqnarray*}
The reason why $h^2$ is also called {\it harmonic} Hardy space is that its elements are functions which are harmonic on $\mathbb{D}$.
Recall that a function $u:\Omega\to\mathbb{R}$ defined on an open set $\Omega\subseteq\mathbb{C}$ is {\it harmonic} it is $C^2(\Omega)$ (it has continuous, partial second order derivatives), and it satisfies the {\it Laplace equation},
\begin{equation}\label{eqLaplaceDimTwo}
0=\Delta u:=\partial_{xx}u+\partial_{yy}u.
\end{equation}
A complex valued function is harmonic if its real and imaginary parts are. By Cauchy-Riemann equations, for instance, a holomorphic function is harmonic. In particular, $z\mapsto z^n=u+iv$ is harmonic on $\mathbb{C}$, and for the same reason $z\mapsto\Bar{z}^n=u-iv$ is harmonic, too. Any function $u$ of the form \eqref{eqRealHardyFunction}, if the series converges in $\mathbb{D}$ (hence, by the theory of power series, it converges absolutely and uniformly on any compact in $\mathbb{D}$), is harmonic, because we can indefinitely differentiate the series under the sum.

Let's verify this for the derivative w.r.t. $x$. We have that $\partial_x(z^n)=\partial_x(x+iy)^n=nz^{n-1}$, and $\partial_x(\Bar{z}^n)=n\Bar{z}^{n-1}$. The power series of the derivatives,  
\[
w(z):=\sum_{n=1}^\infty na_n z^{n-1}+\sum_{n=1}^\infty na_{-n}\Bar{z}^{n-1},
\]
has the same radius of convergence as that representing $u$, hence, by a well know results of advanced calculus, $w(z)=\partial_xu(z)$. The reasoning is the same with $y$ instead of $x$, and
\[
\partial_yu(z)=i\sum_{n=1}^\infty na_n z^{n-1}-i\sum_{n=1}^\infty na_{-n}\Bar{z}^{n-1}.
\]
Now, $\partial_xu$ and $\partial_yu$ satisfy the same hypothesis of $u$, and we can take further partial derivatives.

With this at hands,
\[
\Delta u(z)=\sum_{n=0}^\infty a_n \Delta(z^n)+\sum_{n=1}^\infty a_{-n}\Delta(\Bar{z}^n)=0.
\]

\vskip0.5cm

We introduced the Hardy space from the viewpoint of the "Fourier side" (which engineers consider the {\it time domain}: $\mathbb{N}$ for the holomorphic, and $\mathbb{Z}$ for the harmonic case), imposing a size condition on the coefficients of the power series. The {\it raison d'être} of representing the sequence $\{a_n\}$ by means of the holomorphic (or harmonic) function $f(z)$ of which it provides the Taylor coefficients at the origin (its {\it generating function}, as it is named in combinatorics, or its {\it $z$-transform}, as it is called in engineering, where the set of the admissible $z$'s is called the {\it frequency domain}), is that many algebraic and, more important, quantitative properties of the sequence, and of the operators acting on it, emerge in the complex domain, as we will see later on. To start, we give a characterization of the $H^2$ norm in "frequency space".

We denote by $\mathbb{T}:\partial\mathbb{D}=\{e^{it}:t\in[0,2\pi)\}$ the {\it torus}, i,e, the unit circle in the complex plane. Sometimes we identify $t\leftrightarrow e^{it}$, and $[0,2\pi)\leftrightarrow\mathbb{T}$.
\begin{theorem}\label{theoHtwoAsInBooks}
    Let $f:\mathbb{D}\to\mathbb{C}$ be a holomorphic function. Then,
    \begin{equation}\label{eqHtwoAsInBooks}
        \|f\|_{H^2}^2=\sup_{0\le r<1}\frac{1}{2\pi}\int_{0}^{2\pi}|f(re^{it})|^2dt=\lim_{r\to1}\frac{1}{2\pi}\int_{0}^{2\pi}|f(re^{it})|^2dt.
    \end{equation}
    Moreover, if $f_\ast(e^{it})=\sum_{n=0}^\infty a_ne^{int}$, where $a_n$ is the $n^{th}$ coefficients in the power expansion of $f$ at the origin, and $f_r(e^{it})=f(re^{it})$, then
    \begin{equation}\label{eqHtwoCauchyCondition}
        \lim_{r\to1}\|f_r-f_\ast\|_{L^2(0,2\pi)}=0.
    \end{equation}
    In particular,
    \begin{equation}\label{eqHtwoFromBoundary}
\|f\|_{H^2}^2=\frac{1}{2\pi}\int_{0}^{2\pi}|f_\ast(e^{it})|^2dt.
    \end{equation}
\end{theorem}
In \eqref{eqHtwoCauchyCondition} we have that the boundary values $f_\ast$ of $f(r\cdot)$ as $r\to1$ exist in the $L^2$ sense. A cornerstone of $H^2$-theory is that $\lim_{r\to1}f(r^{it})$ exists for $a.e.\ t\in[0,2\pi)$. We are not concerned, here, with these "hard analysis" developments of the theory, which are, however, important also for the more functional analytic side.
\begin{proof}
   For each $0\le r<1$, $f_r$ is continuous on $\mathbb{T}$, and by Plancherel isometry we have that
   \[
   \frac{1}{2\pi}\int_{0}^{2\pi}|f(re^{it})|^2dt=\sum_{n=0}^\infty |a_n|^2r^{2n},
   \]
   so that \eqref{eqHtwoAsInBooks} follows by dominated convergence, and $f\in H^2(\mathbb{D})$ if and only if $f_\ast\in L^2(\mathbb{T})$. Still by Plancherel,
   \[
   \|f_r-f\|_{L^2(0,2\pi)}^2=\sum_{n=0}^\infty|a_n|^2(1-r^n)^2\to0,\text{ as }r\to1
   \]
by dominated convergence, which shows \eqref{eqHtwoCauchyCondition}.
\end{proof}
For the harmonic Hardy space we have a similar statement, which immediately follows from the holomorphic one provided you know the following.
\begin{theorem}\label{theoHolHarm}
    Let $\Omega$ be a simply connected set in $\mathbb{C}$, and $u:\Omega\to\mathbb{R}$ be harmonic. Then, there exists $v$ harmonic in $\Omega$ such that $f=u+iv$ is holomorphic in $\Omega$. Moreover, there is a unique such $v$ with $v(z_0)=0$, where $z_0$ is any distinguished point in $\Omega$. 

    In particular, $u=\frac{f+\Bar{f}}{2}$ can be written as the sum of a holomorphic and an anti-holomorphic function. 
\end{theorem}
Applying the theorem to real and negative parts of $u$, and using known properties of holomorphic functions, we have the following.
\begin{corollary}\label{CorHolHarm}
    Let $\Omega$ be a simply connected set in $\mathbb{C}$. The functions $u:\Omega\to\mathbb{C}$ which are harmonic in $\Omega$ are exactly those having the form $u(z)=f(z)+\Bar{g(z)}$, with $f,g$ holomorphic in $\Omega$. Moreover, $f$ and $g$ are uniquely determined once we require $g(z_0)=0$ at a distinguished point $z_0$ in $\Omega$.

    If $\Omega=\mathbb{D}$, then
    \begin{equation}\label{eqHarmIsSeries}
        u(re^{it})=\sum_{n=-\infty}^{+\infty}a_n r^{|n|}e^{int},
    \end{equation}
    where both the series with negative $n$'s and positive $n$'s converge in $\mathbb{D}$. The coefficients themselves are given by
    \[
    a_n=\frac{r^{-|n|}}{2\pi}\int_0^{2\pi}u(re^{it})e^{-int}dt,
    \]
    and are independent of $0\le r<1$.
\end{corollary}
\begin{proof}[Proof of the theorem]
   We want $u+iv$ to satisfy Cauchy-Riemann equations, i.e. $\partial_xv=-\partial_yu$ and $\partial_yv=\partial_xu$, or
   \begin{equation}\label{eqCRforPoincare}
       dv=-\partial_yudx+\partial_xudy.
   \end{equation}
   Now, the $1$-form on the right hand side of \eqref{eqCRforPoincare} is closed because $\partial_y(-\partial_yu)=\partial_x(\partial_xu)$ by harmonicity of $u$. Poincaré theorem ensures that the form is exact, hence
   \[
   v(z):=\int_{\gamma(z_0,z)}(-\partial_yudx+\partial_xudy)
   \]
   is well defined, independently of the chosen integration path $\gamma(z_0,z)$, as long as it starts in $z_0$ and ends in $z$.

Moreover, $u+iv$ satisfies the Cauchy-Riemann equation, since $\nabla v=(-\partial_yu,\partial_xu)$, hence $u+iv$ is holomorphic.

Uniqueness holds because, in a connected domain, a function is determined by its gradient, modulo an additive constant.
\end{proof}
We have then the counterpart of theorem \ref{theoHtwoAsInBooks}.
\begin{theorem}\label{theohtwoAsInBooks}
    Let $u:\mathbb{D}\to\mathbb{C}$ be a harmonic function, with series expansion like in \eqref{eqHarmIsSeries}. Then,
    \begin{equation}\label{eqhhtwoAsInBooks}
        \|u\|_{H^2}^2=\sup_{0\le r<1}\frac{1}{2\pi}\int_{0}^{2\pi}|u(re^{it})|^2dt=\lim_{r\to1}\frac{1}{2\pi}\int_{0}^{2\pi}|u(re^{it})|^2dt.
    \end{equation}
    Moreover, if
    \[
    u_\ast(e^{it})=\sum_{n=-\infty}^\infty a_ne^{int}
    \]
     and $u_r(e^{it})=u(re^{it})$, then
    \begin{equation}\label{eqhhtwoCauchyCondition}
        \lim_{r\to1}\|u_r-u_\ast\|_{L^2(0,2\pi)}=0.
    \end{equation}
    In particular,
    \begin{equation}\label{eqhhtwoFromBoundary}
\|h\|_{H^2}^2=\frac{1}{2\pi}\int_{0}^{2\pi}|u_\ast(e^{it})|^2dt.
    \end{equation}
\end{theorem}
The proof of the theorem easily follows from theorem \ref{theoHolHarm}, corollary \ref{CorHolHarm}, and theorem \ref{eqHtwoAsInBooks}.

\subsection{The Dirichlet space $\mathcal{D}(\mathbb{D})$}\label{SSectDirichletSpace}
Consider the subspace $\mathcal{D}=\mathcal{D}(\mathbb{D})$ of the Hardy space $H^2(\mathbb{D})$ having as elements functions $f(z)=\sum_{n=0}^\infty$ such that
\begin{equation}\label{eqDirichletNorm}
    \|f\|_{\mathcal{D}}^2=\sum_{n=0}^\infty(n+1)|a_n|^2=\|f\|_{H^2}^2+[f]<\infty.
\end{equation}
The expression $[f]^{1/2}$, which is the "largest part" of the norm, is just a seminorm, since it vanishes on constant functions. We will see below that is is a {\it conformal invariant}.
\begin{exercise}\label{exeDirichletRK}
    Show that $\mathcal{D}$ has reproducing kernel $k(z,w)=\frac{1}{\Bar{w}z}\log\frac{1}{1-\Bar{w}z}$, where $\log$ is the determination of the complex logarithm which takes real values on the positive axis. 
\end{exercise}
The study of the Dirichlet space started in 1940 with Arne Beurling's\footnote{Arne Beurling ( 1905-1986) was one of the great analysts of the middle part of the XX century, and he gave fundamental contributions to a number of areas. This \href{https://people.kth.se/~haakanh/publications/hed-crypto-beurling.pdf}{article} summarizes his role in the Swedish cryptanalysis program during WWII.} \href{https://projecteuclid.org/journals/acta-mathematica/volume-72/issue-none/Ensembles-exceptionnels/10.1007/BF02546325.full}{Ensembles exceptionnels}.
\begin{lemma}\label{lemmaDirichletNorm}
    Let $f$ be a function which is holomorphic in the unit disc. Then, 
    \begin{equation}\label{eqDirichletNormReal}
        [f]=\frac{1}{\pi}\int_{\mathbb{D}}|f'(z)|^2dxdy,\ z=x+iy.
    \end{equation}
\end{lemma}
\begin{proof}
    For $f(z)=\sum_{n=0}^\infty a_nz^n=\sum_{n=0}^\infty a_n r^ne^{int}$, using the orthogonality of the trigonometric system,
    \begin{eqnarray*}
        \int_{\mathbb{D}}|f'(z)|^2dxdy&=&\int_{0}^{2\pi}dt\int_0^1\left|\sum_{n=1}^\infty na_n r^{n-1}e^{i(n-1)t}\right|^2rdr\crcr
        &=&\int_0^1\left(\int_{0}^{2\pi}\sum_{m=1}^\infty\sum_{n=1}^\infty mn a_m\overline{a_n} r^{m-n} e^{i(m+n-2)t}dt\right)rdr\crcr 
        &=&2\pi\sum_{n=1}^\infty n^2|a_n|^2 \int_0^1r^{2n-2}rdr\crcr 
        &=&\pi \sum_{n=0}^\infty n|a_n|^2\crcr 
        &=&\pi[f].
    \end{eqnarray*}    
\end{proof}
The first geometric interpretation of $[f]$ is in terms of areas.
\begin{proposition}\label{profDirichletArea}
    We have that
    \begin{equation}\label{eqDirichletArea}
    [f]=\frac{1}{\pi}\text{Area}(f(\mathbb{D}))
    \end{equation}
    is $1/\pi$ times the area of $f(\mathbb{D})$, computed with multiplicities (that is, if $f$ is a $n-\text{to}-1$ map from open $A\subset\mathbb{D}$ to $f(A)$, then the area of $f(A)$ appears multiplied times $n$ in the right hand side of \eqref{eqDirichletArea}).
\end{proposition}
\begin{proof}
   The function $f=u+iv$ is holomorphic, hence $f'=\partial_xf=\partial_xu+i\partial_xv$ and, when $f$ is thought of as a map from $\mathbb{R}^2$ to itself, its Jacobian is, by Cauchy-Riemann equations,
   \[
   Jf=\begin{pmatrix}
       \partial_xu&\partial_yu\crcr 
       \partial_xv&\partial_yv
   \end{pmatrix}
   =\begin{pmatrix}
       \partial_xu&-\partial_xv\crcr 
       \partial_xv&\partial_xu
   \end{pmatrix},
   \]
   hence,
   \[
   \det(Jf)=(\partial_xu)^2+(\partial_xv)^2=|f'|^2.
   \]
   By the change of variables formula for double integrals,
   \[
   \int_{\mathbb{D}}|f'(z)|^2dxdy=\int_{\mathbb{D}}\det(Jf)(z)dxdy=\int_{f(\mathbb{D})}dudv=\text{Area}(f(\mathbb{D})),
   \]
   as wished.
\end{proof}
In order to draw some further geometric information on the seminorm $[f]^{1/2}$, we review some facts concerning the unit disc. Any biholomorphic map $\varphi:\mathbb{D}\to\mathbb{D}$ (in literature such maps are often called {\it automorphisms} of $\mathbb{D}$, or {\it Moebius maps}) can be uniquely expressed in the form
\begin{equation}\label{eqMoebiusMap}
    \varphi(z)=v\frac{a-z}{1-\Bar{a}z},
\end{equation}
where $v,a$ are complex numbers, $|a|<1$ and $|v|=1$. Since $f\circ\varphi$ is a map having the same image as $f$ (counting multiplicities), we have that
\begin{equation}\label{eqDirichletConformalInvariant}
    [f\circ\varphi]=[f],
\end{equation}
i.e. the quantity $[f]$ is a {\it conformal invariant}.
\begin{exercise}\label{exeDirichletConformalInvariant}
    Let $p_n(z)=z^n$, and $f\in\mathcal{D}$. Show that $[f\circ p_n]=n[f]$.
\end{exercise}
Another interpretation of the Dirichlet seminorm comes from hyperbolic geometry, and it is based of a couple of "magic relations" concerning automorphisms of the disc. If $\varphi$ is the map in \eqref{eqMoebiusMap}
\begin{align}\label{eqMagicHyperbolic}
    1-|\varphi(z)|^2&=\frac{(1-|a|^2)(1-|z|^2)}{|1-\Bar{a}z|^2}\crcr 
    |\varphi'(z)|^2&=\frac{1-|a|^2}{|1-\Bar{a}z|^2}=\frac{1-|\varphi(z)|^2}{1-|z|^2}.
\end{align}
\begin{exercise}\label{exeMagicHyperbolic}
    Verify the equations \eqref{eqMagicHyperbolic}.
\end{exercise}
Consider on $\mathbb{D}$ the Riemannian metric
\[
ds^2=\frac{4|dz|^2}{(1-|z|^2)^2},
\]
which is called the {\it hyperbolic metric}, or {\it Poincaré metric} on $\mathbb{D}$. The area form associated to $ds^2$ is
\[
dA_{hyp}(z)=\frac{4dxdy}{(1-|z|^2)^2}.
\]
Let $d$ be the Riemannian distance associated to the metric,
\[
d(z,w)=\log\frac{1+\left|\frac{z-w}{1-\Bar{w}z}\right|}{1-\left|\frac{z-w}{1-\Bar{w}z}\right|}.
\]
(Look for a proof on the web, or produce your own).

The relations \eqref{eqMagicHyperbolic} imply that the automorphisms of $\mathbb{D}$ are isometries of $ds^2$:
\[
\frac{4|d\varphi(z)|^2}{(1-|\varphi(z)|^2)^2}=\frac{4|dz|^2}{(1-|z|^2)^2}.
\]
Now, for a holomorphic map $f:\mathbb{D}\to\mathbb{C}$, the {\it hyperbolic distorsion} is (written in a somehow old fashioned way) the ratio between the infinitesimal increment of $f$ in the Euclidean plane and that of the independent variable in the disc endowed with the hyperbolic distance:
\[
\delta_{hyp}(f):=\frac{|f(z+dz)-f(z)|}{d(z+dz,z)}=\frac{|f'(z)|\cdot|dz|}{\frac{2|dz|}{1-|z|^2}}=(1-|z|^2)|f'(z)|/2.
\]
Since they express quantities of hyperbolic geometry, the hyperbolic area element and the hyperbolic distorsion of a holomorphic function $f:\mathbb{D}\to\mathbb{C}$ are invariant under any automorphism $\varphi$, 
\[
dA_{hyp}(\varphi(z))=dA_{hyp}(z), \text{ and }\delta_{hyp}^2(f\circ\varphi)(z)=\delta_{hyp}^2(f)(z).
\]
Hence, the same holds for
\[
\pi[f]=\int_{\mathbb{D}}|f'(z)|^2dxdy=\int_{\mathbb{D}}\delta_{hyp}^2(f)dA_{hyp}.
\]
That is,
\[
\pi[f\circ\varphi]=\pi[f].
\]
This new, and different, argument proving the invariance of the Dirichlet seminorm, has the advantage that it allows us to write down a number of different, conformally invariant quantities. For instance,
\begin{equation}\label{eqBlochNorm}
    \|f\|_{\mathcal{B},\ast}:=\sup_{z\in \mathbb{D}}(1-|z|^2)|f'(z)|=2\sup_{z\in \mathbb{D}}\delta_{hyp}(f)(z),
\end{equation}
the {\it Bloch seminorm} of the holomorphic function $f$, is conformally invariant,
\begin{equation}\label{eqBlochInvariance}
\|f\circ\varphi\|_{\mathcal{B},\ast}=\|f\|_{\mathcal{B},\ast}.
\end{equation}
The {\it Bloch space} $\mathcal{B}=\mathcal{B}(\mathcal{D})$ is the space of functions for which $\|f\|_{\mathcal{B},\ast}<\infty$, and it is a Banach space under the norm $\|f\|_{\mathcal{B}}=\|f\|_{\mathcal{B},\ast}+|f(0)|$ (which is {\it not} conformally invariant), where the term $|f(0)|$ was added to avoid that the norm of a constant vanishes.

\section{Multipliers of a RKHS}\label{SectMultipliersRKHS}
The translation of convolution on the "time side" into products on the "frequency side" suggests that multiplication operators should play a prominent role in reproducing Hilbert space theory. This happens, in fact, even outside Fourier theory proper. In this section we see the definition and the first properties of such multiplication operators in general, then we characterize the multiplicators of the (holomorphic) Hardy space.

\subsection{Multipliers in RKHS}\label{SSectMultRKHS}
Let $H$ be a reproducing kernel Hilbert space on a set $X$, with reproducing kernel $k(x,y)=k_y(x)$. A function $b:X\mapsto\mathbb{C}$ is a {\it multiplier } of $H$ if the linear operator $M_b:f\mapsto bf$ maps $H$ into $H$.
\begin{lemma}\label{lemmaMultRHHSbndd}
    If $b$ is a multiplier, then $M_b$ is bounded on $H$.
\end{lemma}
\begin{proof}
    We show that the graph of $M_b$ is closed in $H\times H$.
    Suppose $f_n\to f$ and $bf_n\to g$ in $H$. By Cauchy-Schwarz and the reproducing property,
    \[
    |f_n(x)-f(x)|=|\langle k_x,f_n\rangle-\langle k_x,f\rangle|\le \|k_x\|\cdot \|f_n-f\|\to0\text{ as }n\to\infty.
    \]
    Then,
    \begin{eqnarray*}
        g(x)&=&M_b^\ast k_x=\lim_{n\to\infty}\langle k_x,bf_n\rangle\crcr 
        &=&\lim_{n\to\infty}b(x)f_n(x)\crcr 
        &=&b(x)f(x).
    \end{eqnarray*}
    Hence, $bf\in H$, proving the closedness of $M_b$.
\end{proof}
We denote by $\mathcal{M}(H)$ the {\it multiplier space} of $H$, and by $\left\|b\right\|_{\mathcal{M}(H)}:=\|M_b\|_{\mathcal{B}(H)}$ the {\it multiplier norm }of $b$.
\begin{theorem}\label{theoMultiplierDiagonalized}
    Let $b\in \mathcal{M}(h)$ and let $M_b^\ast$ be the Hilbert space adjoint of the corresponding multiplication operator. Then, for all $x$ in $X$,
    \begin{equation}\label{eqMultiplierDiagonalized}
        M_b^\ast k_x=\overline{b(x)}k_x.
    \end{equation}
    Viceversa, if $T:H\to H$ is a bounded operator having each kernel function $k_x$ as eigenfunctions, with eigenvalue $\overline{b(x)}$, then $T^\ast=M_b$ is a multiplication operator.
\end{theorem}
\begin{proof}
    About the direct statement,
    \begin{eqnarray*}
        M_b^\ast k_x(y)&=&\langle k_y,M_b^\ast k_x\rangle\crcr 
        &=&\langle M_b^\ast k_y, k_x\rangle\crcr 
        &=&\overline{\langle k_x, M_b k_y\rangle}\crcr 
        &=&\overline{\langle k_x, bk_y\rangle}=\overline{b(x)k_y(x)}\crcr 
        &=&\overline{b(x)}k_x(y).
    \end{eqnarray*}
    For the converse one, the operator $T$ satisfies
    \begin{eqnarray*}
        (T^\ast f)(x)&=&\langle k_x|T^\ast f\rangle\crcr 
        &=&\langle T k_x| f\rangle\crcr
        &=&\langle \overline{b(x)} k_x| f\rangle\crcr
        &=&b(x)\langle k_x, f\rangle\crcr 
        &=&b(x)f(x),
    \end{eqnarray*}
    it is then multiplication times $b$.
\end{proof}
\begin{corollary}\label{corMultEstimateGeneric}
    We have the estimate $\|b\|_{\mathcal{M}(H)}\ge\sup_{x\in X}|b(x)|$.
\end{corollary}
\begin{proof} In fact,
\[
\|M_b\|_{\mathcal{B}(H)}=\|M_b^\ast\|_{\mathcal{B}(H)}\ge \sup_{x\in X}\frac{\|M_b^\ast k_x\|}{\|k_x\|}=\sup_{x\in X}\left(|\overline{b(x)}|\frac{\|k_x\|}{\|k_x\|}\right)=\sup_{x\in X}|b(x)|.
\]  
\end{proof}
The fact that the reproducing functions are eigenvalues of $M_b^\ast$, and that $\langle k_x|k_y\rangle=k(x,y)$ does not vanish identically, but for trivial examples, says that multiplication operators in RKHS are very far from being self-adjoint. 
\begin{exercise}\label{exeMultComplete}
    The multiplier space $\mathcal{M}(H)$ of a RKHS $H$ is complete with respect to the multiplier norm.
\end{exercise}

It is reasonable asking if, in general, $\|b\|_{\mathcal{M}(H)}=\sup_{x\in X}|b(x)|$, or at least if we can provide some estimate $\|b\|_{\mathcal{M}(H)}\approx\sup_{x\in X}|b(x)|$. The answer is negative, and much research goes into finding sharp estimates of the multiplier norm for specific RKHS's in terms of concrete "size" information concerning the multiplier. In the Hardy space $H^2(\mathbb{D})$, however multiplier norm and sup norm coincide, as we show in the next subsection.

\subsection{The multiplier space of $H^2(\mathbb{D})$}\label{SSectHtwoHinfty}
We denote by $H^\infty(\mathbb{D})$ the space of the bounded, holomorphic functions on the unit disc. We endow it with the sup norm,
\[
\|b\|_{H^\infty(\mathbb{D})}=\sup_{|z|<1}|b(z)|.
\]
\begin{theorem}\label{theoH2Hinfty}
A function $b:\mathbb{D}\to\mathbb{C}$ is a multiplier of $H^2(\mathbb{D})$ if and only if it is holomorphic and bounded. Moreover,
\begin{equation}\label{eqH2Hinfty}
    \|b\|_{\mathcal{M}(H^2(\mathbb{D}))}=\|b\|_{H^\infty(\mathbb{D})}.
\end{equation}
\end{theorem}
\begin{proof}
Since $1\in H^2(\mathbb{D})$, if $b$ is a multiplier, then $b=b\cdot 1\in H^2$, hence it is holomorphic. We already know from the general theory that $\|b\|_{\mathcal{M}(H^2(\mathbb{D}))}\ge\|b\|_{H^\infty(\mathbb{D})}$. In the other direction,
\begin{eqnarray*}
    \|bf\|_{H^2}^2&=&\sup_{r\to1}\frac{1}{2\pi}\int_0^{2\pi}|b(re^{it})f(re^{it})|^2dt\crcr 
    &\le&\|b\|_{H^\infty(\mathbb{D})}^2\sup_{r\to1}\frac{1}{2\pi}\int_0^{2\pi}|f(re^{it})|^2dt\crcr 
    &=&\|b\|_{H^\infty(\mathbb{D})}^2\|f\|_{H^2}^2.
\end{eqnarray*}
Passing to sup over $f$'s with $\|f\|_{H^2}$, we see that $\|b\|_{\mathcal{M}(H^2(\mathbb{D}))}\le\|b\|_{H^\infty(\mathbb{D})}$.
\end{proof}

\chapter{Banach algebras}\label{chaptBaAl}\index{Banach algebra}

\section{Banach algebras}\label{SectBanachALgebras}
\noindent {\bf Notation} In this chapter we denote by $Cl(E)$ the closure of a subset $E$ in $\mathbb{C}$. We write $\overline{E}=\{\overline{z}:z\in\mathbb{C}\}$ to denote the set of the complex conjugates of elements in $E$.
\subsection{Definition and basic properties}\label{SSectDefBAandBasica}
An {\it associative algebra} is a vector space $\mathcal{A}$ endowed with a product $\cdot:\mathcal{A}\times\mathcal{A}\to\mathcal{A}$ such that, for $A,B,C\in\mathcal{A}$ and $\lambda\in\mathbb{C}$,
\begin{enumerate}
    \item[(a)] $(AB)C=A(BC)$,
    \item[(b)] $\lambda(AB)=(\lambda A)B=A(\lambda B)$,
    \item[(c)] $(A+B)C=AB+BC$ and $A(B+C)=AB+AC$.
\end{enumerate}
An element $1\in \mathcal{A}$ is a {\it unit} if $1A=A1=A$ for $A\in\mathcal{A}$.

The algebra $\mathcal{A}$ is {\it normed} when endowed with a norm $A\mapsto\|A\|$ such that
\begin{equation}\label{eqNormAlgebra}
    |AB|\le |A|\cdot|B|.
\end{equation}
It is a {\it Banach algebra} if it is complete with respect to the given norm. If the Banach algebra has a unit $I$, we require that
\[
|I|=1.
\]
As we have seen before, examples of Banach algebras are given by $\mathcal{L}(X)$, where $X$ is a Banach space. Even the finite dimensional Banach algebra of the $2\times 2$ complex matrices is instructive in guessing which phenomena might occur.

If a Banach algebra $\mathcal{A}$ does not have a unit, we can add it in a canonical way. Let $\Tilde{A}=\mathcal{A}\times\mathbb{C}$ containing elements $(A,z)\equiv A+1z$ with the product
\[
(A+Iz)(B+1w):=AB+zB+wA+Izw.
\]
Then, $\Tilde{\mathcal{A}}$ is an associative algebra with unit $I$, $\mathcal{A}$ has codimension one in $\Tilde{\mathcal{A}}$, and it becomes a normed algebra under the norm
\[
\|A+Iz\|:=|A|+|z|.
\]
\begin{exercise}
Show that $(\Tilde{\mathcal{A}},\|\cdot\|)$ is a Banach algebra if $(\mathcal{A},|\cdot|)$ is.
\end{exercise}
From now on, all our algebras will have a unit, unless otherwise stated.

Here are some notable examples of Banach algebras.
\begin{enumerate}
    \item[(1)] Let $X$ be a Banach space. Then, $\mathcal{L}(X)$, the algebra of the bounded operators on $X$, is a Banach algebra. 
    \item[(2)] In particular, $\mathcal{L}(H)$ is a Banach algebra when $H$ is a Hilbert space.
    \item[(3)] By Young's inequality, $L^1(\mathbb{R})$ with the convolution product is a (commutative) Banach algebra. The same holds for $\ell^1(\mathbb{Z})$ and, more generally, for $L^1(G)$ (with respect to Haar measure), when $G$ is a locally compact group (but commutativity ceases to hold in general, unless $G$ is a commutative group).
    \item[(4)] The $n\times n$ matrices with complex entries form a Banach algebra with respect to the operator norm.
    \item[(5)]  Let $\Omega$ be a locally compact topological space. Then, $C_0(\Omega)$ (continuous functions vanishing at infinity) and $C_b(\Omega)$ (bounded, continuous functions) are (commutative) Banach algebras. $C_0(\Omega)$ has unit if and only if $\Omega$ is compact.
    \item[(6)] The multiplier space $\mathcal{M}(H)$ of a reproducing kernel Hilbert space $H$ is a (commutative) Banach algebra.
    \item[(7)] In particular, $H^\infty(\mathbb{D})$ is a (commutative) Banach algebra.
\end{enumerate}

An element $A$ in $\mathcal{A}$ is {\it invertible} if there is $E$ in $\mathcal{A}$ such that $AE=EA=I$. We denote it by $A^{-1}$.
It might happen that $A$ has just a {\it left inverse} $E$, $EA=I$, or a {\it right inverse} $F$, $AF=I$. If they both exists, $E=F$ is an inverse for $A$, since
\[
E=EI=E(AF)=(EA)F=IF=F.
\] 
The shift $\tau_1$ on $\ell^{2}(\mathbb{N})$, $\tau_1f(n)=\begin{cases}0\text{ if }n=0\crcr f(n-1)\text{ if }n\ge1\end{cases}$ has $\tau_1^\ast f(n)=f(n-1)$, the {\it back-shift}, as left inverse, but it does not have a right inverse:
\[
\tau_1^\ast\tau_1=I, \text{ but }\tau_1\tau_1^\ast=I-\delta_0.
\]
On the algebraic side, we have some useful properties.
\begin{theorem}\label{theoBanachAlgebraic}
    \begin{enumerate}
        \item[(1)] If $A,B\in\mathcal{A}$ are invertible, then $AB$ is invertible. If $AB=BA$, then $B^{-1}A^{-1}=A^{-1}B^{-1}$.
        \item[(2)] If $AB=BA$ and $AB$ is invertible, then $A$ and $B$ are invertible.
    \end{enumerate}
\end{theorem}
The example above shows that $AB=BA$ is necessary in (2), since $\tau_1^\ast\tau_1=I$, but neither factor is invertible.
\begin{proof}
\noindent(1) We have
\[
(AB)B^{-1}A^{-1}=I=B^{-1}A^{-1}(AB).
\]
The second assertion is an easy exercise.
\noindent(2) The elements $(AB)^{-1}A$ and $A(AB)^{-1}$ are left and right inverse of $B$,
\[
[(AB)^{-1}A]B=(AB)^{-1}(AB)=I=(AB)(AB)^{-1}=(BA)(AB)^{-1}=B[A(AB)^{-1}],
\]
hence, they coincide and are a bilateral inverse for $B$.
\end{proof}
\subsection{Invertibility, resolvent, and spectrum}\label{SSectBAspectrum}
The next result has both an algebraic and a topological nature.
\begin{theorem}\label{theoInvertibleIsOpen}
The set of the invertible elements in $\mathcal{A}$ is open. More precisely, if $A_0$ is invertible and $|A-A_0|<|A_0^{-1}|^{-1}$, then $A$ is invertible.
\end{theorem}
\begin{proof} We start by elements which are close to the identity, $A=I-H$ with $|H|<1$.
    The idea is mimicking the usual geometric series argument in the complex plane,
    \[
    \frac{1}{1-h}=\sum_{n=0}^\infty h^n.
    \]
    The series $\sum_{n=0}^\infty H^n$ converges for $|H|<1$, and using the telescopic property of the sum,
    \begin{equation}\label{eqNeumannSeries}
    (I-H)\sum_{n=0}^\infty H^n=\sum_{n=0}^\infty (H^n-H^{n+1})=I-H^{m+1}+\sum_{n=m+1}^\infty (H^n-H^{n+1})=I,
    \end{equation}
and similarly $\sum_{n=0}^\infty H^n(I-H)=I$.

    Suppose $A_0$ is invertible. Then, $A=A_0-H=A_0(I-A_0^{-1}H)$ is invertible if $|A-A_0|=|H|<|A_0^{-1}|^{-1}$, with inverse
    \[
    A^{-1}=(I-A_0^{-1}H)^{-1}A_0^{-1}=\sum_{n=0}^\infty \left(A_0^{-1}H\right)^{n}A_0^{-1},
    \]
    by what we have seen above.
\end{proof}
The series in \eqref{eqNeumannSeries} is called the {\it Neumann series}\index{Neumann series} of $H$.

The {\it resolvent } of an element $A$ of $\mathcal{A}$ is\index{resolvent!set}
\begin{equation}\label{eqResolventBA}
    \rho(A):=\{\lambda\in\mathbb{C}:\ \lambda I-A \text{ is invertible}\}\subset\mathbb{C},
\end{equation}
and the {\it spectrum } of $A$ is 
\begin{equation}\label{eqSpectrumBA}
    \sigma(A):=\{\lambda\in\mathbb{C}:\ \lambda I-A \text{ is not invertible}\}=\mathbb{C}\setminus \rho(A).
\end{equation}
\begin{theorem}\label{theoSpecOpen}
    \begin{enumerate}
        \item[(1)] For any $A$, $\rho(A)$ is open and $\sigma(A)$ is compact.
        \item[(2)] The {\bf resolvent function} of $A$, $R_\lambda(A):\rho(A)\to\mathbb{C}$,\index{resolvent!function}
        \begin{equation}\label{eqResolventFunction}
            R_\lambda(A):=(\lambda I-A)^{-1}
        \end{equation}
        is holomorphic on $\rho(A)$.
    \end{enumerate}
\end{theorem}
\begin{proof}
    \noindent(1) $\lambda\in\rho(A)$ if and only if $\lambda I-A$ is invertible. If $h\in\mathbb{C}$,
    \[
    (\lambda-h) I-A=(\lambda I-A)-hI,
    \]
    which is invertible by theorem \ref{theoInvertibleIsOpen}, provided $|h|<|(\lambda I-A)^{-1}|^{-1}$. We so have that $\rho(A)$ is open, hence that $\sigma(A)$ is closed. On the other hand, $\lambda I-A$ is invertible if $|\lambda|>|A|$, hence $\rho(A)\supseteq \mathbb{C}\setminus \overline{B(0,|A|)}$, and $\sigma(A)\subseteq\overline{B(0,|A|)}$ is compact.

    \noindent(2) With $A$ and $\lambda\in \rho(A)$ as in (1), and complex $h$,
    \begin{eqnarray*}
        R_{\lambda-h}(A)&=&[(\lambda-h) I-A]^{-1}\crcr 
        &=&[(\lambda I-A)-hI]^{-1}=(\lambda I-A)^{-1}[I-h(\lambda I-A)^{-1}]^{-1}\crcr 
        &=&(\lambda I-A)^{-1}\sum_{n=0}^\infty (\lambda I-A)^{-n}h^n,
    \end{eqnarray*}
    which defines a holomorphic function of the variable $h$ for $|h|<|(\lambda I-A)^{-1}|^{-1}$.
\end{proof}
\begin{corollary}\label{corSpecOpen}
    For the distance of $\lambda\in\rho(A)$ to $\sigma(A)$ we have the estimate
    \[
    \min\{|\lambda-z|:z\in\sigma(A)\}\ge |(\lambda I-A)^{-1}|^{-1}.
    \]
\end{corollary}

\begin{theorem}\label{theoSpctrumNonempty}
    For any $A$ in $\mathcal{A}$, $\sigma(A)\ne\emptyset$.
\end{theorem}
\begin{proof}
    We have
    \begin{equation}\label{eqLaurentResolvent}
        R_\lambda(A)=(\lambda I-A)^{-1}=\lambda^{-1}(I-A/\lambda)^{-1}=\sum_{n=0}^\infty A^n\lambda^{-n-1},
    \end{equation}
    which converges for $|\lambda|>|A|$ and it is a Laurent series. Its coefficients can be computed via complex integrals on countours $\gamma$ in $|\lambda|>|A|$ and, in particular,
    \begin{equation}\label{eqLaurentCoefficientOne}
        I=\frac{1}{2\pi i}\int_\gamma(\lambda I-A)^{-1}d\lambda.
    \end{equation}
    If the spectrum were empty, $\lambda\mapsto R_\lambda(A)$ would be holomorphic on the whole plane and bounded by \eqref{eqLaurentResolvent}, hence constant, and again by \eqref{eqLaurentResolvent} we have $\lim_{\lambda\to\infty}(\lambda I-A)^{-1}=0$,
    \[
    R_\lambda(A)=\frac{\sum_{n=0}^\infty A^n\lambda^{-n}}{\lambda}\to0 \text{ as }\lambda\to\infty.
    \]
    Hence $(\lambda I-A)^{-1}=0$ for all $\lambda$, which contradicts \eqref{eqLaurentCoefficientOne}.
\end{proof}
\subsection{Example: the spectrum of a multiplication operator on $L^2(\nu)$}\label{SSectSpectrumOfMult} Many spectral theorems for operators $L$ defined on a Hilbert space have the form: {\it $L$ is unitarily equivalent to a multiplication operator on some $L^2$ space}. Here we consider the bounded operators having this form.

Let $(X,\nu)$ be a measure space in which for each measurable $E$, 
\[
\nu(E)=\sup_{F\subseteq E,\nu(F)<\infty}\nu(F).
\]
Let $\varphi:X\to\mathbb{C}$ be a measurable function. Denote by $M_\varphi:f\mapsto \varphi f$ the corresponding multiplication operator. The {\it (natural) domain} of $M_\varphi$ is $D(M_\varphi)=\{g\in L^2(\nu):M_\varphi g\in L^2(\nu)\}$.
\begin{lemma}\label{lemmaMultBndd}
    \begin{enumerate}
        \item[(i)] $D(M_\varphi)=\{g:\int_X(1+|\varphi|^2)|g|^2d\nu\}<\infty$. Endowed with the norm
        \[
        \|g\|_{D(M_\varphi)}^2=\int_X(1+|\varphi|^2)|g|^2d\nu,
        \]
        $D(M_\varphi)$ is a Hilbert space, which is dense in $L^2(\nu)$ with respect to the $L^2(\nu)$-norm.
        \item[(ii)] $M_\varphi\in\mathcal{L}(L^2(\nu))$ if and only if $\varphi\in L^\infty(\nu)$. In fact,
        \[
        \|M_\varphi\|_{\mathcal{L}(L^2(\nu))}=\|\varphi\|_{L^\infty(\nu)}.
        \]
        \item[(iii)] $\varphi\mapsto M_\varphi$ is an isometric algebra homomorphism from $L^\infty(\nu)$ into $\mathcal{L}(L^2(\nu))$; in particular, $M_\varphi M_\psi=M_{\varphi\psi}$. Moreover, it is a $\ast$-homomorphism: $M_\varphi^\ast=M_{\overline{\varphi}}$. Also, the operator $M_\varphi$ is {\bf normal}, $M_\varphi M_\varphi^\ast =M_\varphi^\ast M_\varphi$.
    \end{enumerate}
\end{lemma}
\begin{proof}
(i) $\int_X(1+|\varphi|^2)|g|^2d\mu<\infty$ if and only if $\int_X|g|^2d\mu<\infty$ and $\int_X\varphi g|^2d\mu<\infty$, if and only if $g\in L^2(\nu)$ and $M_\varphi g\in L^2(\nu)$. Hence, the identification of $D(M_\varphi)$. The norm $\|g\|_{D(M_\varphi)}$ is the $L^2$ norm with respect to the measure  $(1+|\varphi|^2)d\nu$, which defines a Hilbert space.

(ii) Since $\|\varphi g\|_{L^2(\nu)}\le\|\varphi\|_{L^\infty(\nu)}\|g\|_{L^2(\nu)}$, we have $\|M_\varphi\|_{\mathcal{L}(L^2(\nu))}\le\|\varphi\|_{L^\infty(\nu)}$. In the other direction, we can suppose $\|\varphi\|_{L^\infty(\nu)}>0$. Let $0<\lambda<\|\varphi\|_{L^\infty(\nu)}$, and let $E\subseteq\{x\in X:|\varphi(x)|>\lambda\}$ such that $0<\nu(E)<\infty$. Then $\nu(E)>0$ and
\[
\frac{\|\varphi\chi_E\|_{L^2(\nu)}^2}{\|\chi_E\|_{L^2(\nu)}^2}\ge \lambda^2,
\]
hence, $\|M_\varphi\|_{\mathcal{L}(L^2(\nu))}\ge\|\varphi\|_{L^\infty(\nu)}$.

(iii) The proof is an easy exercise in bringing abstract definitions down to earthly objects.
\end{proof}
Let $\varphi:X\to\mathbb{C}$ be measurable. Its {\it essential range}  is 
\[
\text{essRan}(\varphi)=\{w\in\mathbb{C}:\nu(\varphi^{-1}(D(w,\epsilon)))>0\text{ for all }\epsilon>0\}.
\]
An arbitrarily small perturbation of a point in the essential range of $\varphi$, that is, belongs to the range of $\varphi$ with "positive probability".
\begin{theorem}\label{theoSpecMultLtwo}
    If $\varphi\in L^\infty(\nu)$, then
    \begin{equation}\label{eqSpecMultLtwo}
        \sigma(M_\varphi)=\text{essRan}(\varphi).
    \end{equation}
\end{theorem}
\begin{proof}
    If $w\notin \text{essRan}(\varphi)$, then there is $\epsilon>0$ such that $|\varphi(x)-w|\ge\epsilon$ for $a.e.$ $x$. Then, $(\varphi-w)^{-1}\in L^\infty(\nu)$, $\|(\varphi-w)^{-1}\|_{L^\infty(\nu)}\le 1/\epsilon$, hence, $M_{(\varphi-w)^{-1}}=(M_\varphi-wI)^{-1}$ is the inverse of $M_\varphi-wI$. This shows that $w\in\rho(M_\varphi)$.

    In the opposite direction, let $w\in \text{essRan}(\varphi)$ and fix any $\epsilon>0$. Then, there $E_\epsilon\subseteq X$ such that $0<\nu(E_\epsilon)<\infty$ and $|\varphi(x)-w|\le\epsilon$ for all $x\in E_\epsilon$. Then, $\|(\varphi-w)^{-1}\|_{L^\infty(\nu)}\ge1/\epsilon$. This shows that $M_\varphi-wI$ does not have a bounded inverse, hence that $w\in\sigma(M_\varphi)$.
\end{proof}
Some purely operator theoretic properties of multiplication operators are easily spotted and proved.

A complex number $\lambda$ is an {\it eigenvalue } of a bounded operator $L:H\to H$ ($H$ Hilbert) is there exists $0\ne h\in H$ such that
\[
Lh=\lambda h.
\]
We say that $h$ is an {\it eigenvector} relative to $\lambda$, and the linear space $E_\lambda$ of such eigenvectors is the {\it eigenspace} relative to the eigenvector $\lambda$. Clearly, $\lambda\in\sigma(L)$.
\begin{exercise}\label{exeMultOpSpecProp}
Let $\varphi\in L^\infty(\nu)$. Then, $\lambda$ is an eigenvector for $\lambda$ if and only if there is a measurable set $A$ with $0<\nu(A)$ such that $\varphi(x)=\lambda$ for all $x\in A$. 
\end{exercise}
An operator $U:H\to H$ is {\it unitary} if $UU^\ast=U^\ast U=I$. A bounded operator $A:H\to H$ is {\it self-adjoint} if and only if $A=A^\ast$.
\begin{exercise}\label{exeMultOpSpecPropBis}
    Let $\varphi\in L^\infty(\nu)$. Show that $M_\varphi$ is unitary if and only if $|\varphi(x)|=1$ $a.e.$ in $X$. Show that $M_\varphi$ is self-adjoint if and only if $\varphi(x)\in\mathbb{R}$ $a.e.$ in $X$. Deduce that if $M_\varphi$ is self-adjoint, then $M_{e^{i\varphi}}$ is self-adjoint.
\end{exercise}
\subsection{Example: the spectrum of a multiplier on $H^2(\mathbb{D})$}\label{SSectBSspectrum} If instead of $L^2(\nu)$ we consider one of its subspaces, the spectrum of a multiplication operator might be very different from what we have seen in \S  \ref{SSectSpectrumOfMult}.

Let's compute spectrum and resolvent of the back-shift operator $\tau_1^\ast$ on $\ell^2(\mathbb{N})$. We start from the most obvious enemies of invertibility, which are {\it eigenvalues}, i.e. $\lambda$'s for which $\tau_1^\ast \varphi=\lambda\varphi$ has a nonvanishing solution $\varphi\in \ell^2$. Such $\varphi$ are the solution to
\[
\varphi(n+1)=\tau_1^\ast\varphi(n)=\lambda\varphi(n),
\]
that is $\varphi(n)=c\lambda^n$ for some $c\in\mathbb{C}$, which lies in $\ell^2$ if and only if $|\lambda|<1$. The same argument shows that $\lambda$ is not an eigenvalue if $|\lambda|\ge1$. The {\it point spectrum} $\sigma_p(\tau_1^\ast)$ of the back-shift, the set of its eigenvalues, is the open disc $\mathbb{D}$. But the spectrum is closed, hence $\sigma(\tau_1^\ast)$ contains the closed unit disc $\overline{\mathbb{D}}$.

Consider now $|\lambda|>1$. Since $\|\tau_1^\ast\|\le1<|\lambda|$, the inverse of $I-\frac{\tau_1^\ast}{\lambda}$ exists, hence $\lambda\in\rho(\tau_1^\ast)$. We have proved:
\begin{proposition}\label{propSpecBackShift}
    The spectrum of the back-shift is $\sigma(\tau_1^\ast)=\overline{\mathbb{D}}$. Its point spectrum is $\sigma_p(\tau_1^\ast)=\mathbb{D}$.
\end{proposition}
We will say more on the role of $\mathbb{T}$ in the spectrum when we do the spectral theory of operators on a Hilbert space.

In frequency space, the eigenvectors corresponding to $\lambda$ are constant multiples of
\[
f_\lambda(z)=\sum_{n=0}^\infty \lambda^n z^n=\frac{1}{1-\lambda z}=k_{\overline{\lambda}}(z),
\]
is the reproducing kernel at $\overline{\lambda}$! In fact, we showed that $M_z^\ast k_{\overline{\lambda}}=\lambda k_{\overline{\lambda}}$. Our result was for general multipliers, and we obtain a general result as a consequence of it.
\begin{theorem}\label{theoSpectMultipliers}
    \begin{enumerate}
        \item[(1)] Let $H$ be a RKHS on a set $X$, and let $b$ be a multiplier. Then, 
        \[
        \sigma(M_b^\ast)\supseteq \text{Cl}(\overline{b(X)}).
        \]
        \item[(2)] If $b\in H^\infty(\mathbb{D})$, then
        \[
        \sigma(M_b^\ast)=\text{Cl}(\overline{b(\mathbb{D})}).
        \]
    \end{enumerate}
\end{theorem}
\begin{proof}
(1) holds because each $\overline{b(x)}$ is an eigenvalue of $M_b^\ast$, corresponding to the eigenvector $k_x$. (2) In $H^2(\mathbb{D})$ it is easy to see that there are no other elements in the spectrum. We use the fact that $\sigma(M_b^\ast)=\overline{\sigma(M_b)}$. Suppose $\lambda\notin b(\mathbb{D})$. Then, $\frac{1}{\lambda-b}\in H^\infty(\mathbb{D})$. In fact we can say more:
\[
0<c\le \left|\frac{1}{\lambda-b(z)}\right|\le C<\infty.
\]
This implies that $M_{(\lambda-b)^{-1}}=M_{b-\lambda}^{-1}$ is an inverse for $M_{\lambda-\lambda}$, hence, $\lambda\in\rho(M_b)$.
\end{proof}

\footnote{We can do the same calculations in frequency domain. The shift operator corresponds to $M_z$, multiplication times the $z$ variable on $H^2(\mathbb{D})$, and the back-shift to
\begin{eqnarray*}
   \langle g|M_z^\ast f\rangle_{H^2}&=&\langle M_z g| f\rangle_{H^2}=\lim_{r\to1}\frac{1}{2\pi}\int_{|z|=r}\overline{z}\overline{g(z)}f(z)\frac{dz}{iz}\crcr 
   &=&\lim_{r\to1}\frac{1}{2\pi}\int_{|z|=1}\overline{z}\overline{g(z)}(f(z)-f(0))\frac{dz}{iz}\crcr 
   &=&\lim_{r\to1}\frac{r^2}{2\pi}\int_{|z|=1}\overline{g(z)}\frac{f(z)-f(0)}{z}\frac{dz}{iz}\crcr 
   &=&\left\langle g| \frac{f-f(0)}{z}\right\rangle_{H^2},
\end{eqnarray*}
then
\begin{equation}\label{eqBackshiftZtransf}
    M_z^\ast f(z)=\frac{f(z)-f(0)}{z}.
\end{equation}
\begin{exercise}\label{exeBackshiftZtransf}
    Deduce \eqref{eqBackshiftZtransf} from the $\ell^2(\mathbb{N})$ expression of $\tau_1^\ast$.
\end{exercise}
The eigenvalues $\lambda$ of $M_z^\ast$ are those for which
\begin{equation}\label{eqBackShiftEigenZ}
    \frac{f(z)-f(0)}{z}=\lambda f(z)
\end{equation}
has a solution $f\in H^2$. The equation has solution $f(z)=\frac{f(0)}{1-\lambda z}$, which is holomorphic in $\mathbb{D}$ if and only $|\lambda|\le1$, and it belongs to $H^2(\mathbb{D})$ if and only if $|\lambda|<1$.}

\subsection{The spectral radius}\label{SSectSpectralRadius}
The {\it spectral radius} $r(T)$ of $T$ is
\begin{equation}\label{eqSpectralRadius}
    r(T)=\max\{|\lambda|:\lambda\in\sigma(T)\}\le|T|.
\end{equation}
To compute it, we have a beautiful formula.
\begin{theorem}[Gelfand]\label{theoSpectralRadius}
    We have
    \begin{equation}\label{eqSpectralRadiusCompute}
        r(T)=\lim_{n\to\infty}|T^n|^{1/n}.
    \end{equation}
\end{theorem}
\begin{proof}
The Laurent expansion of $R_\lambda(T)$ at infinity is, as we have seen before,
\begin{equation}\label{eqPerhapsALreadyThere}
(\lambda I-T)^{-1}=\sum_{n=0}^\infty \frac{T^n}{\lambda^{n+1}},
\end{equation}
As in the case of holomorphic functions of one variable, one shows that the radius of convergence of the series with respect to $1/\lambda$ is
\[
\frac{1}{\limsup_{n\to\infty}|T^n|^{1/n}}.
\]
Hence, $R_\lambda(T)$ is holomorphic for $|\lambda|>\limsup_{n\to\infty}|T^n|^{1/n}$. If it were holomorphic for $|\lambda|>Q$ with $Q<\limsup_{n\to\infty}|T^n|^{1/n}$, arguing as is in holomorphic function theory we would deduce (using Cauchy formula) that the series in \eqref{eqPerhapsALreadyThere} would converge for $|\lambda|>Q$, contradicting formula for the radius of convergence. Hence,
\[
r(T)=\limsup_{n\to\infty}|T^n|^{1/n}.
\]
    We then show that $\lim_{n\to\infty}|T^n|^{1/n}$ exists. Observe first that $|T^{m+n}|\le|T^m| |T^n|$, hence that $a_n:=\log|T^n|$ satisfies $a_{m+n}\le a_m+a_n$.
For $n=mq+r$, $r=0,\dots,m-1$, we have
\[
\frac{a_n}{n}\le \frac{q a_m}{m q+r}+\frac{a_r}{m q+r}=\frac{ a_m}{m +r/r}+\frac{a_r}{m q+r},
\]
thus,
\[
\limsup_{n\to\infty}\frac{a_n}{n}\le\frac{a_m}{m},
\]
for each fixed $m$, hence,
\[
\limsup_{n\to\infty}\frac{a_n}{n}\le\liminf_{m\to\infty}\frac{a_m}{m}.
\]
Then, $lim_{n\to\infty}\frac{a_n}{n}=\lim_{n\to\infty}\log(|T^n|^{1/n})$ exists.
\end{proof}
\section{Holomorphic calculus}\label{SectHolCalc}\index{holomorphic!calculus in a Banach algebra}
If $A$ is an element of  an algebra $\mathcal{A}$ on  $\mathbb{C}$, $q\in\mathbb{C}[z]$ is a polynomial with complex coefficients, the definition of $q(A)\in\mathcal{A}$ is given in the obvious way, 
\begin{equation}\label{eqCalculusWithPolynomials}
q(z)=\sum_{j=0}^n\alpha_j z^j\mapsto \Phi_A(p):=q(A)=\sum_{j=0}^n\alpha_j A^j.
\end{equation}
Interesting facts happen when $A\ne0$, yet $p(A)=0$ for some $p\ne0$, as we see below.
In this section, we extend the definition from polynomials to holomorphic functions.  We first review some basic algebra in order to have a basic theory to make comparisons with.
\subsection{The soon-to-be-lost paradise of polynomial calculus}\label{SSectBasicAlgebraOfAlgebras}
Given an element $A$ in an algebra $\mathcal{A}$ and a polynomial $p(z)$ with complex coefficients, we can compute $p(A)$ in the obvious way. This operation gives a algebra homomorphism $\Phi_A:p\mapsto p(A)$ from $\mathbb{C}[z]$ to $\mathcal{A}$. Actually, we can be more precise. Let $\mathbb{C}[A]$ be the smallest sub-algebra of $\mathcal{A}$ which contains $A$ and $I$ (which is a commutative algebra). Then, $\Phi_A$ is a surjective homomorphism from $\mathbb{C}[z]$ to $\mathbb{C}[A]$.

In general, $\Phi_A$ is far from being injective. Let $\mathcal{A}$ be the algebra of the $n\times n$ matrices with complex coefficients. Since $\mathcal{A}$ has dimension $n^2$, while $\mathbb{C}[z]$ has infinite dimension, $\Phi_A$ has an infinite dimensional kernel. By contrast, let $\mathcal{A}=H^\infty(\mathbb{D})$ and let $A(z)=z$, the shift. Then, $p(A)(z)=p(z)$, which we might think of as multiplication times $p(z)$, hence $\Phi_A=Id$ is the identity. Similarly, if $\mathcal{A}=C[0,1]$ is the algebra of the continuous functions on $[0,1]$, and $A(x)=x$, $\Phi_A(p)(x)=p(x)$.

 Let $\mathcal{I}=\ker\Phi_A$, which is an {\it ideal} in $\mathbb{C}[z]$. Recall that all ideals in the algebra of polynomials are {\it principal ideals}, i.e. there is a {\it generating polynomial} $p$, uniquely determined but for a multiplicative constant, such that $\mathcal{I}=\langle p\rangle:=\{p\cdot q:q\in\mathbb{C}[z]\}$. To see this, it suffices to show that for $a,b\in\mathcal{I}$, $GCD(a,b)\in\mathcal{I}$, but this follows from the Euclidean algorithm to compute the greatest common divisor of two polynomials. If $\ker\Phi_A=\langle p\rangle$. then $p(A)=0$ and $q(A)\ne0$ for all polynomials which are not multiple of $p$ (in particular, for all proper factors of $p$).

We have then a complete classification of the algebras $\mathbb{C}[A]$, since they are all isomorphic to $\mathbb{C}[z]/\langle p\rangle$ for a uniquely determined, monic $p$, or (when $\ker\Phi_A=0$) $\mathbb{C}[A]$ is isomorphic to $\mathbb{C}[z]$. The isomorphism is canonical: $q(A)\mapsto [q]_{\mod p}$.

Summarizing, we have the first two items of the following statement. 
\begin{theorem}\label{theoParadiseLostOne}
    \begin{enumerate}
        \item[(i)] Let $A\in\mathcal{A}$. Then, there is a unique minimal monic polynomial $p(z)$ such that $p(A)=0$. Let $[q]_{\mod p}=\langle p\rangle +q$ denote the equivalence class of $q$ in $\mathbb{C}[z]/\langle p\rangle$. Then, the map
        \begin{equation}\label{eqParadiseLostOne}
            [q]\mapsto q(A)
        \end{equation}
        is an isomorphism of $\mathbb{C}[z]/\langle p\rangle$ onto $\mathbb{C}[A]$.
        \item[(ii)] In (i), $\dim(\mathbb{C}[A])=\deg(p)$, if $\deg(p)\ne 0$, and $\dim(\mathbb{C}[A])=\infty$ if $\deg(p)=0$.
        \item[(iii)] Viceversa, let $p$ be a monic polynomial with $\deg(p)\ge1$, let $\mathcal{A}=\mathbb{C}[z]/\langle p\rangle$, and set $A=[z]_{\mod p}$. Then, $q(A)=[q]_{\mod p}$.
    \end{enumerate}
\end{theorem}
\begin{proof}
    The third one follows from the fact that $\pi_p:q\mapsto [q]_{\mod p}$ is a homomorphism of algebras, $\pi_p:\mathbb{C}[z]\to \mathbb{C}[z]/\langle p\rangle$.
\end{proof}
The nonconstant polynomial $p(z)$ can be written as
\begin{equation}\label{eqPolDec}
    p(z)=(z-\lambda_1)^{m_1}\dots(z-\lambda_l)^{m_l},
\end{equation}
where $\lambda_1,\dots,\lambda_l$ are the distinct roots of $p$, and $m_j\ge1$ is the multeplicity of $\lambda_j$.

Before we proceed, let's consider two instructive examples. If all multiplicities are one, then we can consider the diagonal $l\times l$ matrix
\[
\text{Diag}(\lambda_1,\dots,\lambda_l)=
\begin{pmatrix}
\lambda_1&0&\dots\crcr 
0&\lambda_2&\dots\crcr 
0&0&\dots\crcr 
\dots 
\end{pmatrix}
\]
for which
\[
p(\text{Diag}(\lambda_1,\dots,\lambda_l))=
\begin{pmatrix}
p(\lambda_1)&0&\dots\crcr 
0&p(\lambda_2)&\dots\crcr 
0&0&\dots\crcr 
\dots 
\end{pmatrix}=0,
\]
and $p$ is clearly the minimal polynomial with this property.

Consider $p(z)=(z-1)^2$. A "concrete" model of an algebra $\mathcal{A}$ and of an element $A$ in it, such that $\mathcal{A}$ is isomorphic to $\mathbb{C}[z]/\langle p\rangle$ is the algebra of the $2\times 2$ matrices generated by the identity and
\[
A=\begin{pmatrix}
    1&1\crcr 
    0&1
\end{pmatrix},
\]
which satisfies $(A-I)^2=0$, but $I\ne A$. In general, a representation of this kind can be obtained in terms of the Jordan decomposition of a matrix.

An element in $\mathbb{C}[z]/\langle p\rangle$ can be thought of as a polynomial with $\deg(r)<\deg(p)$. For any polynomial $a\in \mathbb{C}[z]$, $a=qp+r$, and what $a$ and $r$ share (a set of invariants for $[a]_{\mod p}$) are the $m_1+\dots+m_l$ complex numbers:
\begin{eqnarray}\label{eqSpectrumIsNotEnough}
    &\ &a(\lambda_1)=r(\lambda_1),a'(\lambda_1)=r'(\lambda_1),\dots,a^{(m_1-1)}(\lambda_1)=r^{(m_1-1)}(\lambda_1),\crcr 
    &\ &\dots\crcr 
    &\ &a(\lambda_l)=r(\lambda_l),a'(\lambda_l)=r'(\lambda_l),\dots,a^{(m_l-1)}(\lambda_l)=r^{(m_l-1)}(\lambda_l).
\end{eqnarray}
This numbers, in turn, are sufficient to determine $r$. There are interpolating formulas dating back to Lagrange, reconstructing $r$ by its values on the zeros of $p$, if all zeros have multiplicity one, or by  values of $r$ and its derivatives in the case of higher multiplicity (see \href{https://www.fields.utoronto.ca/talk-media/1/43/01/slides.pdf}{this} lecture).

The set $\{\lambda_1,\dots,\lambda_l\}$ is the {\it (algebraic) spectrum} of the algebra $\mathbb{C}[z]/\langle p\rangle=\mathbb{C}[A]$. We have just seen that the spectrum alone is not sufficient to determine the algebra $\mathbb{C}[A]$, but that we also need the multiplicities attached to each of the points in it. But for the constant $p$, that is, the algebras $\mathbb{C}[A]$ are parametrized by sets of the form
\[
\{(\lambda_1,m_1),\dots,(\lambda_l,m_l)\},
\]
with the $\lambda_j$'s complex, and the $m_j\ge1$'s positive integers.

\vskip0.5cm 

An interesting special case is that of the {\it Hermitian matrices}, those matrices $A$ such that $\langle z,Aw\rangle_{\mathbb{C}^N}=\langle Az,w\rangle_{\mathbb{C}^N}$. The infinite dimensional version of them are the self-adjoint operators, to be discussed in the chapter on spectral theory. A basic fact fact of linear algebra is that, after a suitable unitary change of coordinates (i.e. passing from $A$ to $U^\ast AU$, where $UU^\ast=U^\ast U=I$), all such matrices can be written as diagonal matrices with real entries, the elements $\mu_1\le\mu_2\le\dots\le\mu_N$ on the diagonal being the eigenvalues. If $\lambda_1<\dots<\lambda_L$ are the {\it distinct} eigenvalues, then 
\begin{equation}\label{eqPVMfinite}
    A=\Delta((\lambda_1,m_1),\dots,(\lambda_L,m_L)),
\end{equation}
where $m_i$ is the multiplicity of $\lambda_i$. Let $E(\lambda_i)$ be the eigenspace corresponding to the eigenvalue $\lambda_i$, and $\Pi_{E(\lambda_i)}$ be the (orthogonal) projection onto it. Then,\eqref{eqPVMfinite} can be written as:
\begin{equation}\label{eqPVMspectralResolution}
    A=\sum_{i=1}^L\lambda_i\Pi_{E(\lambda_i)},
\end{equation}
which is the {\it spectral resolution} of $A$. The {\it projection valued measure} (p.v.m.) corresponding to $A$ is the family $\{\Pi_{\lambda_i}:1\le i\le L\}$, which might be written as a singular measure  on $\mathbb{R}$, having values in the space of projection,
\begin{equation}\label{eqPVMspectralResolutionTwo}
    \Pi=\sum_{i=1}^L\Pi_{E(\lambda_i)}\delta_{\lambda_i}, \text{ i.e. } \Pi(B)=\sum_{\lambda_i\in B}\Pi_{E(\lambda_i)},
\end{equation}
where $\delta_{\lambda_i}$ is the Dirac mass at $\lambda_i$, and $B$ is a Borel set in $\mathbb{R}$.

\subsection{Holomorphic calculus}\label{SSectHolomorphicCalculus}
For a power series like
\begin{equation}\label{eqSeriesWithNoAlgebra}
    f(z)=\sum_{n=0}^\infty \alpha_n z^n,
\end{equation}
the natural definition of $f(A)$ is:
\begin{equation}\label{eqSeriesWithAlgebras}
    f(A)=\sum_{n=0}^\infty \alpha_n A^n.
\end{equation}
Let $R$ be the radius of convergence of $f(z)$.  By comparison with power series, \eqref{eqSeriesWithAlgebras} converges if $R>|A|$, and, in fact, for $R>r(A)$, the spectral radius of $A$. Power series are not the only way to represent holomorphic functions. For instance, we have Cauchy formula, which is more flexible since it deals with regions which are more general than discs.

Let $A$ be an element of a Banach algebra $\mathcal{A}$, $\Omega\supset\sigma(A)$ an open set in the complex plane, and $f:\Omega\to\mathbb{C}$ a function which is holomorphic in $\Omega$. Let $\gamma$ be a loop (or a formal sum of loops) in $\Omega$ such that $\text{Ind}(\gamma,\zeta)=1$ for all $\zeta$ in $\sigma(A)$, i.e. such that $\gamma$ circles around each $\zeta$ in the spectrum of $A$ exactly once, counterclockwise. We define, then
\begin{equation}\label{eqHolomorphicCalculusViaCauchyIntegral}
    f(A):=\frac{1}{2\pi i}\int_\gamma (zI-A)^{-1}f(z)dz.
\end{equation}
We have already seen that the expression $z\mapsto (zI-A)^{-1}$ defines an $\mathcal{A}$-valued holomorphic function in $\mathbb{C}\setminus\sigma(A)$, and the integral in \eqref{eqHolomorphicCalculusViaCauchyIntegral} does not depend on the curve $\gamma$ by Cauchy integral theorem. The element $f(A)$, that is, is defined if $f$ is holomorphic on (any open set containing) $\sigma(A)$.

If we replace $A$ by $w$, the formula gives back $f(w)$, and for the moment this is justification enough for the definition.

We need a relation which has an independent interest.
\begin{lemma}\label{lemmaResolventIdentity}[Resolvent identity]\index{lemma!resolvent identity}
    If $z,w\in\rho(A)$, then
    \begin{equation}\label{eqResolventIdentity}
        (zI-A)^{-1}(wI-A)^{-1}=(w-z)^{-1}[(zI-A)^{-1}-(wI-A)^{-1}].
    \end{equation}
\end{lemma}
\begin{proof}
    Multiply both sides of \eqref{eqResolventIdentity} by all inverses in sight:
    \[
    (w-z)I=(wI-A)-(zI-A),
    \]
    which in fact holds.
\end{proof}
\begin{theorem}\label{theoHolomorphicCalculusViaCauchyIntegral}
\begin{enumerate}
    \item[(i)] If $q$ is a polyomial, or more generally a function  of the form \eqref{eqSeriesWithNoAlgebra}, the definition of $q(A)$ coincides with that given by \eqref{eqCalculusWithPolynomials}, or \eqref{eqSeriesWithAlgebras} in the general case. 
    \item[(ii)]  The map $f\mapsto f(A)$ is an algebra homomorphism.
    \item[(iii)]  {\bf (Spectral mapping theorem)} $\sigma(f(A))=f(\sigma(A))$.
    \item[(iv)]  $(g\circ f)(A)=g(f(A))$, where $g:E\to\mathbb{C}$ is holomorphic on an open set $E\supseteq f(\Omega)\supset f(\sigma(A))$.
\end{enumerate}
\end{theorem}
\begin{proof}
\noindent (i) Let $L_f(A)$ be the algebra element in \eqref{eqHolomorphicCalculusViaCauchyIntegral}: we have to show that $L_q(A)=q(A)$.
We can assume $\gamma$ to be a circle containing $\sigma(q)$ in its interior. By linearity, it suffices to show it for powers $q(z)=z^m$. By developing a Neuman series, which converges of $\gamma$, we have
\begin{eqnarray*}
    l_q(A)&=&\frac{1}{2\pi i}\int_\gamma (zI-A)^{-1}z^mdz\crcr 
    &=&\sum_{n=0}^\infty \frac{1}{2\pi i}\int_\gamma z^{m-1-n}dz A^n\crcr 
    &=&A^m,
\end{eqnarray*}
as wished.

\noindent (ii) We have to show that $L_{fg}=L_fL_g$. Cosider a curve $\gamma$ as in the definition, and another curve $\delta$ with the same properties, but outside $\gamma$, so that $\text{Ind}(\delta,z)=1$ for all $z$ on $\gamma$, and $\text{Ind}(\gamma,w)=0$ for all $w$ on $\delta$. Using the resolvent identity from first to second line,
\begin{eqnarray*}
    L_f(A)L_g(A)&=&\frac{1}{(2\pi i)^2}\int_\gamma\int_\delta (zI-A)^{-1}(wI-A)^{-1}f(z)g(w)dwdz\crcr 
    &=&\frac{1}{(2\pi i)^2}\int_\gamma\int_\delta (w-z)^{-1}[(zI-A)^{-1}-(wI-A)^{-1}] f(z)g(w)dwdz\crcr 
    &=&\frac{1}{2\pi i}\int_\gamma(zI-A)^{-1}f(z)\left(\frac{1}{2\pi i}\int_\delta (w-z)^{-1}g(w)dw\right)dz\crcr 
    &\ &-\frac{1}{2\pi i}\int_\delta\left(\frac{1}{2\pi i}\int_\gamma (w-z)^{-1}f(z)dz\right)(wI-A)^{-1}g(w)dw.
\end{eqnarray*}
The inner integral in the second summand vanishes, while the inner integral in the first is $g(z)$, hence,
\begin{eqnarray*}
    L_f(A)L_g(A)&=&\frac{1}{2\pi i}\int_\gamma (zI-A)^{-1}f(z)g(z)dz\crcr 
    &=&L_{fg}(A).
\end{eqnarray*}
\noindent (iii) If $\mu\notin f(\sigma(A))$, then $f-\mu$ does not vanish on $\sigma(A)$, hence $g=(f-\mu)^{-1}$ defines a function which is holomorphic in a neighborhood containing $\sigma(A)$. We can then define $g(A)$, and by (ii) 
\[
g(A)(f(A)-\mu I)=I=(f(A)-\mu I)g(A),
\]
i.e. $\mu\notin\sigma(f(A))$. This shows $\sigma(f(A))\subseteq f(\sigma(A))$.

In the other direction, suppose that $\lambda\in\sigma(A)$ and observe that
\[
h(z)=\frac{f(z)-f(\lambda)}{z-\lambda}
\]
defines a function which is holomorphic on an open set containing $\lambda$, with $h(\lambda)=f(\lambda)$. Again by (ii) we can define $h(A)$ and
\[
h(A)(A-\lambda I)=f(A)-f(\lambda)I.
\]
The first term is not invertible because $\lambda\in\sigma(A)$, thus the second term is not invertible either, showing that $f(\lambda)\in\sigma(f(A))$, hence $f(\sigma(A))\subseteq\sigma(f(A))$.

\noindent (iv) Our hypothesis guarantee that $g$ is holomorphic on an open set $\Omega$ containing $f(\sigma(A))=\sigma(f(A))$. Let $\gamma$ be a simple curve in $\Omega$ such that $\text{Ind}(\gamma,\xi)=1$ for all $\xi$ in $f(\sigma(f(A)))$, and let $E=\text{Interior}(\gamma)\supset\sigma(f(A))$ be the open region inside $\gamma$. Then, $f^{-1}(E)$ is an open set containing $A$ in its interior. Let $\delta$ be a curve in $f^{-1}(B)$ such that $\text{Ind}(\delta,\zeta)=1$ for all $\zeta$ in $A$, and observe that $f(z)\ne w$ for $z$ in $f^{-1}(E)$ and $w$ on $\delta$, hance, $z\mapsto (w-f(z))^{-1}$ is holomorphic on $f^{-1}(E)$.

Using the definition of the holomorphic function of an operator in the first two and in the last equality:
\begin{eqnarray*}
    g(f(A))&:=&\frac{1}{2\pi i}\int_\gamma (wI-f(A))^{-1}g(w)dw\crcr 
    &=&\frac{1}{2\pi i}\int_\gamma\left(\frac{1}{2\pi i}\int_\delta (w-f(z))^{-1} (Iz-A)^{-1}dz\right)g(w)dw\crcr 
    &=&\frac{1}{2\pi i}\int_\delta\left(\frac{1}{2\pi i}\int_\gamma (w-f(z))^{-1} g(w)dw\right)(Iz-A)^{-1}dz\crcr 
    &=&\frac{1}{2\pi i}\int_\delta g(f(z))(Iz-A)^{-1}dz\crcr 
    &\ &\text{by Cauchy formula}\crcr 
    &=&(g\circ f)(A).
\end{eqnarray*}
\end{proof}
\section{Intermezzo: holomorphic calculus and Laurent series}\label{SectHolCalcLauSer}
Let $A$ be an element of a Banach algebra 
$\mathcal{A}$ having spectrum 
$\sigma(A)\subseteq \text{Cl}B(0,r)=\{z\in\mathbb{C}:|z|\le R\}$, 
let $\Omega=\Omega(r,R)=\{z\in\mathbb{C}:0\le r<|z|\le R\}$, and let $f$ be holomorphic in $\Omega$. If $\gamma$ is a loop in $\Omega$ such that $\text{Ind}(\gamma,\zeta)=1$ for all $\zeta\in \text{Cl}B(0,r)$, then the expression
\begin{equation}\label{eqMerCalc}
    L_A(f,\gamma)=\frac{1}{2\pi i}\int_\gamma (zI-A)^{-1}f(z)dz
\end{equation}
certainly makes sense and it defines an element $L_A(f,\gamma)\in\mathcal{A}$. Can we identify in some way $L_A(f,\gamma)$ with $f(A)$? We provide here an example of a function $f\ne0$ such that $L_A(f,\gamma)=0$.

Let $\mathcal{A}=\mathbb{C}^{2\times 2}$ be the space of the $2\times2$ matrices with complex entries, let 
$A=\begin{pmatrix}
    1&0\\ 
    0&0
\end{pmatrix}$ be the projection onto the first coordinate, so that $\sigma(A)=\{0,1\}$, and let $f(z)=\frac{1}{z}$, which is holomorphic in $\Omega(1,\infty)$. If $\gamma(t)=p\cdot e^{it}$ ($t\in[0,2\pi]$, $p>1$ fixed) is as above, then $L_A(f,\gamma)$ is well defined, although it can not be $A^{-1}$, since $A$ is not invertible.
We claim that $L_A(f,\gamma)=0$.

The function $g$,
\[
g(w)=\frac{1}{2\pi i}\int_\gamma \frac{f(z)}{z-w}dz=\frac{1}{2\pi i}\int_\gamma \frac{dz}{(z-w)z}
\]
vanishes in the interior of $\gamma$. In fact, by the residue theorem,
\[
g(w)=\text{Res}\left(\frac{1}{(z-w)z},z=0\right)+\text{Res}\left(\frac{1}{(z-w)z},z=w\right)=\frac{1}{-w}+\frac{1}{w}=0.
\]
We can now compute
\begin{eqnarray*}
    \frac{1}{2\pi i}\int_\gamma (zI-A)^{-1}f(z)dz&=&\frac{1}{2\pi i}\int_\gamma \left[z\begin{pmatrix}
    1&0\\ 
    0&1
\end{pmatrix}-\begin{pmatrix}
    1&0\\ 
    0&0
\end{pmatrix}\right]^{-1}\frac{dz}{z}\\ 
&=&\frac{1}{2\pi i}\int_\gamma\begin{pmatrix}
    z-1&0\\ 
    0&z
\end{pmatrix}^{-1}\frac{dz}{z}\\ 
&=&\frac{1}{2\pi i}\int_\gamma\begin{pmatrix}
    \frac{1}{(z-1)z}&0\\ 
    0&\frac{1}{z^2}
\end{pmatrix}dz\\ 
&=&\begin{pmatrix}
    g(1)&0\\ 
    0&g(0)
\end{pmatrix}\\ 
&=&\begin{pmatrix}
    0&0\\ 
    0&0
\end{pmatrix}.
\end{eqnarray*}
All this seems rather {\it ad hoc}, but it can be rephrased in a more conceptual way. Recall that any function $f$ which is holomorphic in the annulus $\Omega(r,R)$ can be expanded in a Laurent series there,
\begin{equation}\label{eqLaurent}
    f(z)=\sum_{n=-\infty}^{+\infty}a_nz^n,
\end{equation}
where $f_+(z)=\sum_{n=0}^{+\infty}a_nz^n$ converges in $B(0,R)$ and $f_-(z)=\sum_{n=-\infty}^{-1}a_nz^n$ converges in $\mathbb{C}\setminus[\text{Cl}(B(0,r))]$. Then, with the hypothesis on the objects contained in \eqref{eqMerCalc},
\begin{equation}\label{eqMerCalcTwo}
    L_A(f,\gamma)=f_+(A)
\end{equation}
is the element provided by the "projection" of $f$ onto the functions holomorphic in $B(0,R)$. Here is the proof (I leave it to you to verify all needed convergences).
\begin{eqnarray*}
    L_A(f,\gamma)&=&\frac{1}{2\pi i}\int_\gamma (zI-A)^{-1}f(z)dz\\
    &=&\sum_{n=-\infty}^{+\infty}a_n\frac{1}{2\pi i}\int_\gamma (zI-A)^{-1}z^ndz\\
    &=&\sum_{n=-\infty}^{-1}a_n\frac{1}{2\pi i}\int_\gamma (zI-A)^{-1}z^ndz+\sum_{n=0}^{+\infty}a_n A^N\\
    &=&\sum_{n=-\infty}^{-1}a_n\frac{1}{2\pi i}\int_\gamma (zI-A)^{-1}z^ndz+f_+(A).
\end{eqnarray*}
In the third equality we have used the holomorphic calculus.
To finish, observe that for $n<0$
\begin{eqnarray*}
      \frac{1}{2\pi i}\int_\gamma (zI-A)^{-1}z^ndz&=&\frac{1}{2\pi i}\int_\gamma (I-z^{-1}A)^{-1}z^{n-1}dz\\ 
      &=&\sum_{k=0}^\infty\frac{1}{2\pi i}\int_\gamma A^k z^{-k}z^{n-1}dz\\ 
      &=&\sum_{k=0}^\infty A^k\frac{1}{2\pi i}\int_\gamma z^{n-1-k}dz\\
      &=&0,
\end{eqnarray*}
because $n-1-k<0$. The second equality seems to require $|z|>|A|$, but if you look more closely, you realize that $|z|>r(A)$ suffices ($r(A)$ is the spectral radius of $A$), and since $r(A)\le r<p$ by hypothesis, we are done.

\section{Gelfand theory of commutative Banach algebras}\label{SectGelfand}
There are many interesting examples of commutative Banach algebras with unit. Here are some examples, not always independent of each other.
\begin{enumerate}
    \item[(i)] $C_b(X)$, the algebra of the bounded and continuous functions on a topological space $X$ with respect to the uniform norm.
    \item[(ii)] $L^\infty(\mu)$ on a measure space $(X,\mu)$ with respect to the $L^\infty$ norm.
    \item[(iii)] $C^1[0,1]$, the space of the continuously differentiable functions on $[0,1]$, with respect to the norm $\|f\|_{C^1}=\|f\|_C+\|f'\|_C$.
    \item[(iv)] $\mathcal{M}(H)$, the multiplier space of a reproducing kernel Hilbert space $H$, is a commutative Banach algebra.
    \item[(v)] In particular, $H^\infty(\mathbb{D})$, the multiplier space of $H^2(\mathbb{D})$, is a commutative Banach algebra.
\end{enumerate}
In example (iv) you might worry about completeness. If you use both completeness of $\mathcal{L}(H)$ and reproducing functions, worries evaporate.
\begin{exercise}\label{exeCompleteMultiplierSPace}
    $\mathcal{M}(H)$ is complete.
\end{exercise}
Gelfand theory identifies a commutative Banach algebra $\mathcal{A}$ with a subalgebra of $C(X)$, the algebra of the continuous fnctions on a compact set $X$. The points of such set are the maximal ideals of $\mathcal{A}$. This section is devoted to precisely stating, and proving, this fact.

\subsection{Multiplicative functionals and ideals}\label{SSSectMFandIdeals}
We start with some definitions. Let $\mathcal{A}$ be an algebra. a {\it multiplicative functional }$p$ on $\mathcal{A}$ is a homomorphism of algebras $p:\mathcal{A}\to\mathbb{C}$. We include in the definition that $p\ne0$ is not the trivial homomorphism. The set of the nonzero algebra homomorphisms $\Phi:\mathcal{A}\to\mathcal{L}$ is denoted by $\text{Hom}(\mathcal{A},\mathcal{L})$, so that $\text{Hom}(\mathcal{A},\mathbb{C})$ is the set of the multiplicative functionals.

A {\it left ideal} $\mathcal{I}$ of an algebra $\mathcal{A}$ is a vector subspace of $\mathcal{A}$ such that $A\mathcal{I}\subseteq\mathcal{I}$ for all $A$ in the algebra; it is a {\it right ideal} if $\mathcal{I}A\subseteq\mathcal{I}$ instead; it is a {\it bilateral ideal } if it is both right and left. Since we deal with commutative Banach algebras, we simply call $\mathcal{I}$ and {\it ideal}. A bilateral ideal $\mathcal{I}$ is {\it maximal} if it is not contained in any other proper, maximal ideal.

An element $A$ in a algebra $\mathcal{A}$ with unit $I$ is {\it invertible} if $A^{-1}$ exists in $\mathcal{A}$: $AA^{-1}=A^{-1}A=I$. 

We list in an exercise some basic, purely algebraic facts about algebras.
\begin{exercise}\label{exeAlgebraProperties} Let $\mathcal{A}$ be an algebra. Prove the following. 
\begin{enumerate}
\item[(i)] Show that a left ideal that contains $I$, or an invertible element, is $\mathcal{A}$.
\item[(ii)] Let $\Phi:\mathcal{A}\to\mathcal{L}$ a homomorphism of algebras. Show that $\ker\Phi$ is a bilateral ideal. 
\item[(iii)] If $p$ is a multiplicative functional and $\mathcal{A}$ has unity, then $p(I)=1$.
\item[(iv)] If $A$ is invertible in the algebra with unity $\mathcal{A}$ and $p\in{\text{Hom}(\mathcal{A},\mathbb{C})}$, then $p(A)\ne0$.
\item[(v)] Conversely to (ii), if $\mathcal{I}$ is a bilateral ideal of $\mathcal{A}$, then $\mathcal{I}$ is the kernel of the algebra homomorphism
\[
\pi:\mathcal{A}\to\mathcal{A}/\mathcal{I}:=\{\pi(A)=A+\mathcal{I}:\ A\in\mathcal{A}\}.
\]
On the quotient we define the product $(A+\mathcal{I})(B+\mathcal{I})=AB+\mathcal{I}$, which is well defined (prove it).
\item[(vi)] If $\Phi\in\text{Hom}(\mathcal{A},\mathcal{L})$ and $\mathcal{I}$ is a bilateral ideal in $\mathcal{L}$, then $\Phi^{-1}(\mathcal{I})$ is a bilateral ideal in $\mathcal{A}$.
\item[(vii)] By Zorn's lemma, in a commutative algebra with identity a proper ideal $\mathcal{I}$ is always contained in a maximal ideal.
\end{enumerate}
\end{exercise}
At the end of the section we will see that (iv) has an inverse in commutative Banach algebras: if $p(A)\ne0$ for all $p\in\text{Hom}(\mathcal{A},\mathbb{C})$, then $A$ is invertible.

We close this algebraic introduction with a lemma.
\begin{lemma}\label{lemmaDivisionAlgebra}
    Let $\mathcal{M}$ be a maximal, bilateral ideal in $\mathcal{A}$, an algebra with unity. Then, $\mathcal{L}=\mathcal{A}/\mathcal{M}$ is a {\bf division algebra}: all nonzero elements in it have an inverse.
\end{lemma}
\begin{proof} Let $\pi:\mathcal{A}\to\mathcal{L}$ be the projection. By contradition, 
    let $0\ne \pi(B)$, with $\pi(B)$ noninvertible in $\mathcal{L}$, and let $\mathcal{L}\supset\mathcal{I}$ be the ideal generated by $B$, which is proper since $I+\mathcal{M}\notin\mathcal{I}$. Now, $\pi^{-1}(\mathcal{I})$ is a proper ideal of $\mathcal{A}$, $I\notin \pi^{-1}(\mathcal{I})$, and it properly contains $\mathcal{M}$: $B\notin \mathcal{M}$, but $B\in \pi^{-1}(\mathcal{I})$.
    Hence, $\mathcal{M}$ was not maximal.
\end{proof}

We now come to analytic-algebraic results.
\begin{lemma}\label{lemmaInvInBanachAlg}
    Let $\mathcal{A}$ be a commutative Banach algebra, and $p\in\text{Hom}(\mathcal{A},\mathbb{C})$. Then, $|p(A)|\le|A|$.
\end{lemma}
\begin{proof}
    Suppose that $|p(A)|>1$ for some $A$, and let $B=\frac{A}{p(A)}$, so that $|B|<1$. Then, $I-B$ is invertible (by a Neumann series), hence,
    \[
    0\ne p(I-B)=1-\frac{p(A)}{p(A)}=0,
    \]
    contradiction.
\end{proof}
The lemma self-improves to the following inequality.
\begin{corollary}\label{corpspectral}
Let $\mathcal{A}$ be a commutative Banach algebra, and $p\in\text{Hom}(\mathcal{A},\mathbb{C})$. Then, $|p(A)|\le r(A)$.
\end{corollary}
\begin{proof}
    We have
    \[
    |p(A)|=|p(a^n)|^{1/n}\le |A^n|^{1/n}\to r(A)
    \]
    as $n\to\infty$.
\end{proof}
\begin{lemma}\label{lemmaClosedIdeals}
Let $\mathcal{I}$ be a proper ideal in a commutative Banach algebra $\mathcal{A}$ with unity. Then, its closure $\overline{\mathcal{I}}$ is a proper ideal.
\end{lemma}
\begin{proof}
    The continuity of the algebra operations ensure that $\overline{\mathcal{I}}$ is an ideal, hence it suffices to show that $I\notin \overline{\mathcal{I}}$. Since elements $A$ with $|I-A|<1$ are invertible, $\{A:|I-A|<\}\cap \mathcal{I}=\emptyset$, and this shows that $I\notin \overline{\mathcal{I}}$.
\end{proof}
\begin{corollary}\label{corClosedIdeals}
    A maximal ideal in $\mathcal{A}$ is closed.
\end{corollary}

\subsection{The algebra of continuous functions}\label{SSSectALgebraCofX}
We consider here the special case where $\mathcal{A}=C(X)$ is the algebra of the continuous functions on a Hausdorff, compact space $X$, endowed with the uniform norm. In this case, the evaluations at point,
\begin{equation}\label{eqPointEvaluationOnCofX}
    \eta_x:f\mapsto f(x)
\end{equation}
is a multiplicative functional, and $\ker(\eta_x)=\mathcal{M}_x=\{f\in C(X):f(x)=0\}$ is a maximal ideal in $C(X)$. Maximality, which we obtained by abstract reasoning, can be rephrased as follows. If $a\in X$ and $h\in C(X)$ is such that $h(a)\ne0$, then there are $\lambda\in\mathbb{C}$ and $f\in C(X)$ with $f(a)=0$ such that $1=\lambda h+f$. In fact, we obtain this by letting $\lambda=1/h(a)$ and $f=1-h/h(a)$.
\begin{theorem}\label{theoPreSToneWeierstrass}
All maximal ideals in $C(X)$ have the form $\mathcal{M}_a$ for some $a\in X$.
\end{theorem}
\begin{proof}
    By contradiction, let $\mathcal{M}$ be a maximal ideal which is is not of the form $\mathcal{M}_a$ for any $a$. Then, for each $x$ in $X$, $\mathcal{M}$ contains a function $f_x$ with $f_x(x)\ne0$ (otherwise $\mathcal{M}_x\subseteq\mathcal{M}$), and we can assume that $f_x(x)>0$. After multiplying $f_x$ times a cut-off function $h_x$ in $C(X)$, $h_xf_x\in\mathcal{M}$, we can also assume  that $h_xf_x\ge0$ on $X$ and $h_x(x)f_x(x)>0$. Let $U_x=\{y\in X:h_xf_x(y)>0\}$. By compactness, $X=\cup_{i=1}^nU_{x_i}$, and so, by Weierstrass theorem, 
    \[
    c\le \sum_{i=1}^n h_{x_i}f_{x_i}\le C.
    \]
    That is, $(\sum_{i=1}^n h_{x_i}f_{x_i})^{-1}\in C(X)$, hence,
    \[
    \mathcal{M}\ni\sum_{i=1}^n \frac{h_{x_i}}{\sum_{j=1}^n h_{x_j}f_{x_j}}f_{x_i}=1. 
    \]
    This shows that $\mathcal{M}=C(X)$.
\end{proof}
\begin{exercise}\label{exeContinuousALgebra}
    Show that the maximal ideal space of the algebra $C^1[0,1]$ is formed by the elements $\mathcal{M}_x\cap C^1[0,1]$, as $x$ ranges in $[0,1]$.
\end{exercise}
%The basic result concerning algebras of continuous functions is provided by the Stone-Weierstrass theorem in its various avatars.
%\begin{theorem}[Stone-Weierstrass (real case)]\label{theoStoneWeierstrassReal}
%    Let $X$ be a Hausdorff, compact space, and let $\mathcal{A}$ be an algebra of functions in $C(X,\mathbb{R})$. Then, $\mathcal{A}$ is dense in $C(X)$ if and only if
%    \begin{enumerate}
%        \item[(i)] $\mathcal{A}$ {\bf is nowhere vanishing}: for each $x\in X$ there is $f\in\mathcal{A}$ s.t. $f(x)\ne0$;
%        \item[(ii)] $\mathcal{A}$ {\bf separates points}: for all $x\ne y$ there is $f\in\mathcal{A}$ s.t. $f(x)\ne f(y)$.
%    \end{enumerate}
%\end{theorem}
%\begin{proof}
%    
%\end{proof}
\subsection{Intermezzo: the algebra of multipliers on a reproducing kernel Hilbert space}\label{SSSectALgebraMultRKHS}
Other notable examples of ideals in algebras are found in multiplier algebras $\mathcal{M}(H)$ of RKHSs $H$ on a set $X$:
\begin{equation}
    \mathcal{M}_x=\{f\in\mathcal{M}(H):f(x)=0\},
\end{equation}
where $x\in X$. Observe that $\mathcal{M}_x$ is proper since $1\ne \mathcal{M}_x$.

More generally, for $x_1,\dots,x_n\in X$, we can define 
\[
\mathcal{M}_{\{x_1,\dots,x_n\}}=\{f\in\mathcal{M}:f(x_1)=\dots=f(x_n)=0\}.
\]
If $E\subset F\subset X$, then $\mathcal{M}_E\supseteq \mathcal{M}_E$.

We recall a result from basic linear algebra.
\begin{lemma}\label{lemmaKernelLinearFunctional}
Let $V$ be a vector space over $\mathbb{C}$ and $l:V\to\mathbb{C}$ be linear, $l\ne0$. Then, $\ker(l)$ is a maximal proper subspace of $V$. 
\end{lemma}
\begin{proof}
    Let $u,v\in V\setminus \ker(l)$. Then, $l(u)v-l(v)u=w\in\ker(l)$, i.e. $v=\frac{l(v)}{l(u)}u+\frac{z}{l(u)}\in\text{span}(\ker(l),u)$, showing that $V=\text{span}(\ker(l),u)$.
\end{proof}
\begin{proposition}\label{propMaximIdealMultiplierRJHS}
The evaluation map $\eta_x:f\mapsto f(x)$ is a multiplicative functional, hence it is bounded.
The ideal $\mathcal{M}_x$ is maximal, hence closed.
\end{proposition}
\begin{proof} The functional $f\mapsto f(x)$ is well defined by definition of multiplier ($f$ is a function) and it is multiplicative. Hence, it is a contraction by lemma \ref{lemmaInvInBanachAlg}.
Its kernel, $\mathcal{M}_x$, is maximal by lemma \ref{lemmaKernelLinearFunctional}.
\end{proof}
In the chapter on reproducing kernel Hilbert spaces we had a different proof that $\eta_x$ is a contraction on $\mathcal{M}(H)$, which followed from the fact that the conjugates of values of $f$, if $f$ is a multiplier, are eigenvalues of the multiplier operator's adjoint.

\subsection{Gelfand's representation}\label{SSSectGelfandRepresentation}

We return to the general theory. 
\begin{theorem}\label{theoQuotientOfBanachALgebra}
    Let $\mathcal{I}$ be a proper, closed ideal in a commutative Banach algebra $\mathcal{A}$. Then, $\mathcal{A}/\mathcal{I}$ is a Banach algebra with respect to its quotient norm as Banach space,
    \[
    |A+\mathcal{I}|:=\inf\{|A+N|:N\in\mathcal{I}\}.
    \]
\end{theorem}
\begin{proof}
    \begin{eqnarray*}
        |AB+\mathcal{I}|&=&\inf\{|AB+N|:N\in\mathcal{I}\}\crcr 
        &\le&\inf\{|AB+AN+BM+MN|:M,N\in\mathcal{I}\}\crcr
        &\le&\inf\{|A+M|\cdot |B+N|:M,N\in\mathcal{I}\}\crcr 
        &=&\inf\{|A+M|:M\in\mathcal{I}\}\inf\{|B+N|:N\in\mathcal{I}\}\crcr 
        &=&|A+\mathcal{I}|\cdot |B+\mathcal{I}|.
    \end{eqnarray*}
\end{proof}
    Next, we have a useful characterization of the complex field.
    \begin{theorem}[Mazur]\label{theoMazur}
        Let $\mathcal{A}$ be a Banach algebra which unit and which is a division algebra. Then, $\mathcal{A}$ is isometrically isomorphic to $\mathbb{C}$.
    \end{theorem}
    \begin{proof}
        Let $A\ne0$ be an element in $\mathcal{A}$. The spectrum $\sigma(A)$ is nonempty, hence there is $\alpha$ such that $\alpha I-A$ is not invertible. Since $\mathcal{A}$ is a division algebra, $A=\alpha I$, and $\sigma(A)=\{A\}$.  We then have a map $\Phi:\mathcal{A}\to\mathbb{C}$, $\Phi(A)=\Phi(\alpha I)=\alpha$, which is an isometric algebra isomorphism.
    \end{proof}

    Finally can state and prove the inverse of property (iv) in exercise \ref{exeAlgebraProperties}.
\begin{theorem}\label{theoMaximalIdealsVsMultiplicativeFunctionals} Let $\mathcal{A}$ be a commutative Banach algebra with unit.
    \begin{enumerate}
        \item[(i)] If $A$ is noninvertible, then there exists $p\in\text{Hom}(\mathcal{A},\mathbb{C})$ such that $p(A)=0$.
        \item[(ii)] To each maximal ideal $\mathcal{M}$ in $\mathcal{A}$ there corresponds a unique $p\in\text{Hom}(\mathcal{A},\mathbb{C})$ such that $\ker(p)=\mathcal{M}$.
    \end{enumerate}
\end{theorem}
\begin{proof}
    \noindent (i) By (vii) in exercise \ref{exeAlgebraProperties}, $A$ is contained in a maximal ideal $\mathcal{M}$. Then,
    \[
    p_\mathcal{M}:\mathcal{A}\to\mathcal{A}/\mathcal{M}\equiv\mathbb{C}
    \]
    provides an element in $\text{Hom}(\mathcal{A},\mathbb{C})$ with $p_\mathcal{M}(A)=0$.

    \noindent (ii) To each maximal $\mathcal{M}$ we associate $p_\mathcal{M}$, and the functional $p$ is uniquely determined since $\mathcal{A}=\text{span}(\mathcal{M},I)$ and $p(I)=1$. 
\end{proof}
Let $\Omega=\Omega(\mathcal{A})$ be the set of the maximal ideal is $\mathcal{A}$, a commutative Banach algebra with unity. By theorem \ref{theoMaximalIdealsVsMultiplicativeFunctionals} we have a well defined map
\begin{equation}\label{eqGelfandRepresentation}
    p:\Omega\times\mathcal{A}\to\mathbb{C},\ p(\mathcal{M},A)=p_\mathcal{M}(A).
\end{equation}
We might look at it as a map $\mathcal{G}:\mathcal{A}\mapsto p(\cdot,A)=p^A$,
\begin{eqnarray}\label{eqGelfandRepresentationName}
    \mathcal{G}:A\to \Omega^{\mathbb{C}},\ \mathcal{G}(A):\Omega\to\mathbb{C},\ [\mathcal{G}(A)](\mathcal{M})=p_\mathcal{M}(A).
\end{eqnarray}
The map $\mathcal{G}$, representing elements of the Banach algebra $\mathcal{A}$ as functions on its maximal ideal space $\Omega$, is the {\it Gelfand representation} of $\mathcal{A}$. We summarize the properties we have so far seen concerning it.
\begin{theorem}\label{theoGelfandALgebraic}
    \begin{enumerate}
        \item[(i)] The Gelfand representation $\mathcal{G}$ is a homomorphism of algebras, $\mathcal{G}\in\text{Hom}(\mathcal{A},\Omega^\mathbb{C})$.
        \item[(ii)] The representation is contractive in the sense that $\sup_{\mathcal{M}\in \Omega}|[\mathcal{G}(A)](\mathcal{M})|\le|A|$.
        \item[(iii)] $\sigma(A)=[\mathcal{G}(A)](\Omega)=\{p(\mathcal{M},A):\mathcal{M}\in \Omega\}$.
        \item[(iv)] $\mathcal{G}(I)=1$ is the constant function. 
        \item[(v)] The functions $\mathcal{G}(A)$ separate the points of $\Omega$: if $\mathcal{M}\ne\mathcal{N}$ are maximal ideals, then $[\mathcal{G}(A)](\mathcal{N})=p(\mathcal{M},A)\ne p(\mathcal{N},A)=[\mathcal{G}(A)](\mathcal{N})$.
    \end{enumerate}
\end{theorem}
\begin{proof}
    
\end{proof}
The heart of Gelfand theory consists in making the representation sharp. In order to do this, we have to make $\Omega$ into a topological space.

Consider on $X$ the {\it Gelfand topology} $\tau$, the coarsest making all maps $p_A:\Omega\to\mathbb{C}$ continuous. A basis of neighborhoods for $\mathcal{M}_0\in \Omega$ in $\tau$ are finite intersections of the sets
\[
    N(A;U)=\{\mathcal{M}\in \Omega:\ p(\mathcal{M},A)\in U\},
\]
where $A\in\mathcal{A}$ and $U\subseteq\overline{D(0,|A|)}$ is open in the Euclidean topology of $\overline{D(0,|A|)}=\{z\in\mathbb{C}:|z|\le|A|\}$.
\begin{theorem}\label{theoGelfandCompact}
    The maximal ideal space $\Omega$ is compact with respect to the Gelfand topology.
\end{theorem}
The proof is similar to the one of the Banach-Alaoglu theorem.
\begin{proof}
    Consider the space    
    \[
    Z=\Pi_{A\in\mathcal{A}}\overline{D(0,|A|)}
    \]
    endowed with the product topology, which is the coarsest making the projections $\pi_{A_o}:Z\to\overline{D(0,|A_0|)}$, 
    $\pi_{A_0}(\{z_A\}_{A\in\mathcal{A}})=z_{A_0}$. A basis of neighborhoods for the product topology are finite intersections of sets of the form
    \[
    \mathcal{N}(A;U)=\{\{z_A\}_{A\in\mathcal{A}}:z_{A_0}\in U\},
    \]
    with $A_0$ and $U$ as above. By Tikhonov theorem, $Z$ is compact.

    Consider then the map $\Phi:\Omega\to Z$
    \[
    \Phi(\mathcal{M})=\{p(\mathcal{M},A):A\in\mathcal{A}\}.
    \]
    By the way the topologies on $\Omega$ and $Z$ are defined, $\Phi$ is a homomorphism onto its image. To show that $\Omega$ is compact, then, it suffices to show that $\Phi(\Omega)$ is closed in $Z$.

    Suppose $\{z_A\}_{A\in\mathcal{A}}$ lies in the closure of $\Phi(X)$, and consider $A,B\in\mathcal{A}$. For any $\epsilon>0$ there is $\mathcal{M}\in X$ such that $|p(\mathcal{M},A)-z_A|<\epsilon/3$, $|p(\mathcal{M},B)-z_B|<\epsilon/3$, and $|p(\mathcal{M},A+B)-z_{A+B}|<\epsilon/3$,
    so that
    \[
    |z_A+z_B-z_{A+B}|<\epsilon.
    \]
    Hence, $z_A+z_B=z_{A+B}$. The same way, one shows that $z_{cA}=cz_A$ if $c$ is a complex number, and $z_{AB}=z_Az_B$. This shows that $p:A\mapsto z_A$ is a multiplicative functional from $\mathcal{A}$ to $\mathbb{C}$. We know that $p=p(\mathcal{M},\cdot)$ for some $\mathcal{M}$ in $\Omega$, hence that $\{z_A\}_{A\in\mathcal{A}}=\Phi(\mathcal{M})$.
\end{proof}
\subsection{The Gelfand representation of $C(X)$ and other examples}\label{SSSectGelfandCofX} When $\mathcal{A}=C(X)$, the space of the continuous functions on a compact, Hausdorff space $X$ with topology $\sigma$, the Gelfand representation is just the identity.

We saw that there a bijection between $X$ and the maximal ideal space, $a\mapsto \mathcal{M}_a$. The multiplicative functional associated to $\mathcal{M}_a$ is evaluation at $a$, $\eta_a(f)=f(a)$. the Gelfand transform, then 
acts as
\[
[\mathcal{G}(f)](\mathcal{M}_a)=\eta_a(f)=f(a)
\]
is nothing but $f$ itself. The topology $\tau$ on $X$ interpreted as a maximal ideal space is the coarsest making each $f\in C(X)$ continuous, hence $\tau\subseteq\sigma$ and $Id:(X,\sigma)\to(X,\tau)$ is continuous. But both $\sigma$ and $\tau$ make $X$ compact, hence sets which are closed for $\sigma$ are also compact for $\sigma$, then they are compact, hence closed for $\tau$. This shows that $\sigma\subseteq\tau$, hence that $\sigma=\tau$.

\begin{exercise}\label{exeGelfandCone}
    Show that the Gelfand representation applied to $C^1[0,1]$ just identifies $C^1[0,1]$ with a particular subalgebra of $C[0,1]$.
\end{exercise}

The space $H^\infty[0,1]$ is a Banach algebra, hence it has a Gelfand representation. Its maximal ideal space is huge, but there seems to be no concrete, known example of maximal ideal.

\section{Intermezzo: corona theorems}\label{SectCoronaTheoremsInGeneral}

\subsection{The corona problem in reproducing kernel Hilbert spaces}\label{SSectCoronaRKHS}

Let $H$ be a reproducing Hilbert space of functions on a set $X$, and let $\mathcal{M}(H)$ be its multiplier algebra. In many important cases, the multiplier algebra is too small even to separate points on $X$. In fact, it might very well be that $\mathcal{M}(H)$ reduces to the constant functions. In many cases of interest, however, including the Hardy space $H=H^2(\mathbb{D})$, the multiplier algebra is large. For instance, $f(z)=z$ already separates points in the disc.

Suppose $\mathcal{M}(H)$ separates points of $X$, and consider its maximal ideal space $\Omega=\Omega(\mathcal{M}(H))$. Since point evaluations are multiplicative functionals, we have a map $i:X\to\Omega$ associating to $x$ the maximal ideal $\mathcal{M}_x=\{f\in \mathcal{M}(H):f(x)=0\}$. The map is injective because $\mathcal{M}(H)$ separates points, and we can thus view $X$ a subset of $\Omega$, which is compact with respect to the Gelfand topology.

Let $\overline{X}$ be the closure of $X$ in $\Omega$. The {\it corona} of $\Omega(\mathcal{M}(H))$ is $\Omega(\mathcal{M}(H))\setminus \overline{X}$. We say that $\Omega(\mathcal{M}(H))$ {\it has a corona} if $\Omega(\mathcal{M}(H))\setminus \overline{X}$ is nonempty.

In 1941 Gelfand published his article on the representation of commutative Banach algebras, and in the same year Kakutani asked if $H^\infty(\mathbb{D})$, the multiplier space of $H^2(\mathbb{D})$, has a corona. The question remained open for twenty years, during which several reformulations of it where found. In 1962, Carleson showed that, as conjectured by Kakutani, $H^\infty(\mathbb{D})$ has no corona.

Carleson proof was hard, and it involved a number of other foundational results. In 1979, Tom Wolff came with a much simplified proof, which he never published, but which has become the standard one in monographs\footnote{See e.g. \href{https://www.cambridge.org/core/books/introduction-to-hp-spaces/13A053E7F427D0F31E94B3DC26DA520A}{Paul Koosis, Introduction to $H^p$ Spaces}, $2^{nd}$ edition (1999) Cambridge University Press.}. This is not the end of the story, because the maximal ideal space of $H^\infty(\mathbb{D})$ {\it is not} homeomorphic to the closed unit disc, as one might naively guess by comparison with the algebra of the continuous functions on a compact set. A concrete understanding of it came later.

After Carleson's result, the corona theorem was extended to many other multiplier spaces, and examples of multiplier spaces where it does not hold where found. The area is an active one, and new results and open problems abound. In this section, we limit ourselves to show that the corona theorem is equivalent to an analytic statement, which is the the starting point of most of the "corona proofs" in literature. We do it for $H^\infty$, but the strategy has general value.

\subsection{An equivalent, analytic formulation of the corona problem}\label{SSectCoronaAndAnalysis}
Let $\mathcal{G}(H^\infty(\mathcal{D}))$ be the maximal ideal space of $H^\infty(\mathbb{D})$. Each $\mathcal{M}$ in it is the kernel of a unique multiplicative functional $p_\mathcal{M}:H^\infty(\mathbb{D})\to\mathbb{C}$. For $f$ in $H^\infty(\mathbb{D})$, we write
\[
f(\mathcal{M}):=p_\mathcal{M}(f),
\]
thus thinking of $f$ as a function on $\mathcal{G}(H^\infty(\mathcal{D}))$.
\begin{lemma}[An equivalent statement of the corona theorem]\label{lemmaCoronaAnalytic}
    The following are equivalent.
    \begin{enumerate}
        \item[(i)] The corona theorem holds in $H^\infty(\mathbb{D})$. 
        \item[(ii)] Suppose that for all $n\ge1$, if $f_1,\dots,f_n\in H^\infty(\mathbb{D})$ there is $j$ such that $|f_j(z)|\ge\delta>0$ for some $\delta>0$ independent of $z\in\mathbb{D}$. Then there exist $g_1,\dots,g_n\in H^\infty(\mathbb{D})$ such that 
        \begin{equation}\label{eqCoronaAnalytic}
            f_1g_1+\dots+f_ng_n=1,
        \end{equation}
        i.e. the ideal generated by $f_1,\dots,f_n$ is the whole of $H^\infty(\mathbb{D})$, $\mathcal{M}(f_1,\dots,f_n)=H^\infty(\mathbb{D})$.
    \end{enumerate}
\end{lemma}
The implication we need for proving the corona theorem is (ii)$\implies$(i). It is useful, however, knowing that the two properties are equivalent. Property (ii) is obvious if we drop the requirement that $g_1,\dots,g_n$ be holomorphic:
\begin{eqnarray*}
    1&=&\frac{f_1\overline{f_1}+\dots+f_n\overline{f_n} }{f_1\overline{f_1}+\dots+f_n\overline{f_n}}\crcr 
    &=&f_1 \frac{\overline{f_1}}{\sum_{j=1}|f_j|^2}+\dots+f_n \frac{\overline{f_n}}{\sum_{j=1}|f_j|^2},
\end{eqnarray*}
and $g_i=\frac{\overline{f_i}}{\sum_{j=1}|f_j|^2}$ is $C^\infty$ and bounded on $\mathbb{D}$ by our assumptions on the $f_j$'s.
\begin{proof}
    \noindent (ii)$\implies$(i). Suppose $H^\infty(\mathbb{D})$ is not dense in $\mathcal{G}(H^\infty(\mathbb{D}))$ with respect to the Gelfand topology. Then, there are a maximal ideal $\mathcal{M}$, $\delta>0$, and $h_1,\dots,h_n\in H^\infty(\mathbb{D})$ such that, for all $z\in\mathbb{D}$, $\mathcal{M}_z\notin N(\mathcal{M};h_1,\dots,h_n;\delta)$, i.e.
    \[
    |h_1(\mathcal{M})-h_1(z)|\ge\delta, \text{ or }|h_2(\mathcal{M})-h_2(z)|\ge\delta, \text{ or }\dots, \text{ or }|h_n(\mathcal{M})-h_n(z)|\ge\delta.
    \]
    So, $f_j:=h_j-h_j(\mathcal{M})\in H^\infty(\mathbb{D})$, $j=1,\dots,n$, satisfy 
    \[f_j(\mathcal{M})=0,\]
but for each $z$ in the disc there is $f_j$ such that $|f_j(z)|\ge\delta$. By (ii), we find $g_1,\dots,g_n\in H^\infty(\mathbb{D})$ such that $f_1g_1+\dots+f_ng_n=1$, hence
    \[
    1=p_\mathcal{M}(1)=\sum_{j=1}^n p_\mathcal{M}(f_j)p_\mathcal{M}(g_j)=0
    \]
    and we have found a contradiction.

\noindent (i)$\implies$(ii). Suppose the corona theorem holds and let $f_1,\dots,f_n$, and $\delta>0$, be as in (ii). Suppose by contradiction that for no $g_1,\dots,g_n\in H^\infty(\mathbb{D})$ we have that $f_1g_1+\dots-f_ng_n=1$, i.e. that 
\[
\mathcal{M}(f_1,\dots,f_n)=\{f_1g_1+\dots+f_ng_n:g_1,\dots,g_n\in H^\infty(\mathbb{D})\}
\]
is a proper ideal in $H^\infty(\mathbb{D})$. Let $\mathcal{M}$ be a maximal ideal in $H^\infty(\mathcal{D})$ which contains $\mathcal{M}(f_1,\dots,f_n)$, and let $p_\mathcal{M}$ be the corresponding multiplicative functional, so that $p_\mathcal{M}(f_j)=0$ for $j=1,\dots,n$. Since the corona theorem holds, there exists $z\in\mathbb{D}$ such that $\mathcal{M}_z\in N(\mathcal{M};f_1,\dots,f_n;\delta)$, i.e. such that for $j=1,\dots,n$,
\[
\delta>\left|p_\mathcal{M}(f_j)-p_{\mathcal{M}_z}(f_j)\right|=|f_j(z)|,
\]
which contradicts our assumption on $f_1,\dots,f_n$.
\end{proof}
A critical analysis of the proof shows that we used very little of the structure of $H^\infty(\mathcal{D})$. 
\begin{theorem}\label{theoCoronaAnalytic}
    Let $\mathcal{A}$ be a Banach algebra of functions $f:X\to\mathbb{C}$, defined on some set $X$, and suppose that the {\bf evaluation functionals} $\eta_x:\mathcal{A}\to \mathbb{C}$,
    \[
    \eta_x(f)=f(x)
    \]
    are bounded on $\mathcal{A}$ for all $x$ in $X$. Then, the following are equivalent.
        \begin{enumerate}
        \item[(i)] The set $\{\ker(\eta_x):x\in X\}$ is dense in $\mathcal{G}(\mathcal{A})$. 
        \item[(ii)] Suppose that for all $n\ge1$, if $f_1,\dots,f_n\in \mathcal{A}$ there is $j$ such that $|f_j(x)|\ge\delta>0$ for some $\delta>0$ independent of $x\in X$. Then there exist $g_1,\dots,g_n\in \mathcal{A}$ such that 
        \begin{equation}\label{eqCoronaAnalyticGeneral}
            f_1g_1+\dots+f_ng_n=1,
        \end{equation}
        i.e. the ideal generated by $f_1,\dots,f_n$ is the whole of $\mathcal{A}$, $\mathcal{M}(f_1,\dots,f_n)=\mathcal{A}$.
    \end{enumerate}
\end{theorem}
The theorem can be in particular be applied to the multiplier space of a reproducing kernel Hilbert space. In fact, most proofs of corona theorems try to show (ii), which is an analytic statement.
\begin{exercise}\label{exeCoronaAnalytic}
    Prove theorem \ref{theoCoronaAnalytic}.
\end{exercise}
The case $n=1$ in \eqref{eqCoronaAnalyticGeneral} is already an interesting one. It says that if $f$ lies in the multiplier space, and $\mathcal{M}(H)$ and $|f|$ is bounded from below, then $1/f\in\mathcal{M}(H)$.

\chapter{Spectral theory of unitary, and of self-adjoint operators}\label{chaptUSAandUO}
Most of this chapter is based on the expository  \href{https://terrytao.wordpress.com/2011/12/20/the-spectral-theorem-and-its-converses-for-unbounded-symmetric-operators/}{The spectral theorem and its converses for unbounded symmetric operators (2011)} by Terence Tao. The advantage of this approach is that we do not need any previous version of the spectral theorem: not that for bounded self-adjoint operators, or more generally that for the normal operators, and not even that for the finite dimensional case. In fact, the self-adjoint bounded case, and the self-adjoint compact case, can be recovered as special cases, the peculiarities of which can be easily spotted and proved using rather basic measure theory. The pedagogical risk, on the other hand, is that the unbounded case presents many subtle points which have to be untangled before reaching the heart of the matter.

A good source for a more comprehensive treatment of spectral theory using Herglotz representation is \href{https://www.mat.univie.ac.at/~gerald/ftp/book-schroe/index.html}{Mathematical Methods in Quantum Mechanics
With Applications to Schrödinger Operators (2014)} by Gerald Teschl (Graduate Studies in Mathematics, Volume 157, Amer. Math. Soc.). A new approach to the whole matter is in \href{https://arxiv.org/abs/2003.06130}{The Functional Calculus Approach to the Spectral Theorem (2020)} by Markus Haase, arXiv:2003.06130. Another recommended reading is Michael Taylor's \href{https://mtaylor.web.unc.edu/wp-content/uploads/sites/16915/2018/04/chap8.pdf}{Spectral Theory}.

A spectral theorem with a simpler proof is the one for unitary operators. I have included it because the result is important by itself, and it can serve as a starter, its proof having several ideas in common with the one for self-adjoint operators. The variant I present here is from Michael Taylor's lecture notes \href{https://mtaylor.web.unc.edu/wp-content/uploads/sites/16915/2018/04/specthm.pdf}{The Spectral Theorem for Self-Adjoint and Unitary Operators}. 

\vskip0.5cm

The spectral theorem which is considered here tells the life and works of {\it one} self-adjoint operator $L$, and of a commutative Banach algebra it generates (in an appropriate way). In fact, the generated algebra is a {\it von Neumann algebra}. Here is a short list of topics we do not cover.
\begin{enumerate}
    \item If $L,M$ are {\it commuting} bounded, self-adjoint operators, $LM=ML$, can they be simultaneously represented in the same model (can they be simultaneously diagonalized)? The answer is yes, but you have to find it elsewhere.
    \item What about normal operators? See Michael Reed, Barry Simon,  \href{https://books.google.it/books/about/Functional_Analysis.html?id=16jFzQEACAAJ&redir_esc=y}{Methods of Modern Mathematical Physics: Functional analysis}, Academic Press, 1972, chapter VII, exercise 5 p. 246.
    \item On the way towards our spectral theorem we will make a number of arbitrary choices. Is there a {\it canonical way} to represent a self-adjoint operator? Yes, up to a point, but you have to learn it somewhere else (e.g. Reed and Simon's Functional Analysis).
    \item An $N\times N$ Hermitian matrix is not in general unitarily determined by its eigenvalues: we also need to know their multiplicity. Is there a multiplicity theory even in the more general Hilbert universe? Again, the answer is positive, but for the details you have to find other sources. e.g. \href{https://books.google.it/books/about/Introduction_to_Hilbert_Space_and_the_Th.html?id=FXBADwAAQBAJ&redir_esc=y}{Introduction to Hilbert Space and the Theory of Spectral Multiplicity},  New York, Chelsea Pub. Co. 1957, by Paul R. Halmos.
\end{enumerate}

\section{Intermezzo: spectral theory and quantum mechanics}\label{SectSTandQM}
In this "motivational" section I try to sketch the reason why the notions of spectral theory are foundational in  quantum mechanics. I do everything in a finite dimensional Hilbert space, so what is needed to follow is just some notion of linear algebra, and you do not need to know the theory which lies ahead. You also might read the section after you have read the rest of the chapter. Or you might skip altogether, if you so wish.

Spectral theory had its origins in the study of integral equations, with foundational contributions of Fredholm, then Hilbert, and it first appeared in its infinite dimensional form. In the course of the 1920's, physicists found that it provided the mathematical framework for quantum mechanics, which was being developed in those years, and both finite and infinite dimensional theories play there a role. It is a really amazing coincidence that the word {\it spectrum} introduced by Hilbert in his work turned out to represent, in the light of the quantum mechanical interpretation, the spectrum of electromagnetic radiation physicists were accustomed to. It was von Neumann who wrote the foundations of quantum mechanics in the language of abstract Hilbert spaces, which has since become standard.

I know very little physics, and the interested student is invited to compare what is to follow with some authoritative source. I suggest the video lectures \href{https://theoreticalminimum.com/courses/quantum-mechanics/2012/winter}{The theoretical minimum: quantum mechanics} by Leonard Susskind, or, even better, the book with the same title by Susskind and Art Friedman (2014). I like the book because:  (i) the mathematics involved are kept to a minimum, so the math student can quickly browse through the preliminaries of each chapter and speedily get to (ii) the conceptual exposition of the physics, which are presented in rigorous terms, (iii) in a way which is not encumbered by a long discussion of experimental results and history of concepts. Susskind's text, that is, is a fast vehicle for those who want to proceed from operator theory to its motivations and applications in  theoretical physics. 

Of course, at some point experiments and history of concepts have to be considered, as always in natural sciences. For instance, the \href{https://arxiv.org/abs/2208.10151#}{Lecture notes on operator algebras and their application in physics} (2022) by Alexander Kegeles provide a more advanced, but nice and readable, introduction to the quantum mechanical interpretation of operator theory, and to operator theory itself, for Master students in mathematics.

I thank Giacomo De Palma for correcting a conceptual mistake in the definition of "complete system of commuting observables".

\subsection{States and observables}\label{SSectStatesQM}
A {\it physical system} is described by a (complex!) Hilbert space $H$, whose vectors represent the {\it states} of the system (or, better: states are vectors of unitary norm). Here we consider $H$ to be finite dimensional. Any self-adjoint operator $A$ represents an {\it observable} of such system, and its eigenvalues $\lambda_i$'s are the possible {\it measures} of the observable: $\lambda_i$ is the value of the observation represented by $A$ when the system is in a state $e\in E(\lambda_i)$, where $E(\lambda_i)$ is, recall, the eigenspace relative to $\lambda_i$.

Consider two states $x,y$. They are {\it distinguishable} if there is an observable $A$ which returns (with certainty: see below) different values for the two states. i.e. there is a self-adjoint operator $A$ such that $Ax=\lambda x$, $Ay=\mu y$, $\lambda\ne\mu$. This happens if and only if $x$ and $y$ are orthogonal: $\langle x|y\rangle=0$. 

More generally, if the system lies in the state $h\in H$, $\|h\|=1$, then the {\it probability} that the measure $X$ (a random variable) of the observable $A$ has the value $\lambda_i$ is
\begin{equation}\label{eqProbabilityQM}
\mathbb{P}_h(X=\lambda_i)=\|\Pi_{E(\lambda_i)}h\|^2.
\end{equation}
The states, that is, encode probabilities of observables. If $h=\alpha k$ with $\alpha\in\mathbb{C}\setminus\{0\}$, $h$ and $k$ return the same probabilities, and they can not be distinguished. The state space really is the projectivization $(H\setminus\{0\})/\mathbb{C}$ of $H$. Many properties of the system, however, depend on the linear structure of $H$. Some examples of this will be presented below. 

With this probabilistic interpretation, the {\it expected value} (or {\it mean}) of the measure in the state $h$ is
\begin{eqnarray}\label{eqExpectationQM}
\mathbb{E}_h(X)&=&\sum_i\lambda_i\mathbb{P}_h(X=\lambda_i)=\sum_i\lambda_i\|\Pi_{E(\lambda_i)}h\|^2\crcr
&=&\left\langle h, \sum_{i=1}^L\lambda_i\Pi_{E(\lambda_i)}h\right\rangle   \crcr 
&=&\langle h,Ah\rangle. 
\end{eqnarray}
The variance of $X$ is called the {\it uncertainty} of $A$, and it can be expressed in a number of ways,
\begin{eqnarray}\label{eqExpectationQMtwo}
    \text{Var}_h(X)&=&\mathbb{E}_h[X-\mathbb{E}_hX]^2\crcr 
    &=&\sum_i[\lambda_i-\sum_j\lambda_i\mathbb{P}_h(X=\lambda_i)]^2\mathbb{P}_h(X=\lambda_i)\crcr
    &=&\mathbb{E}_hX^2-(\mathbb{E}_hX)^2\crcr
    &=&\langle h,A^2h-\langle h,Ah\rangle^2\crcr 
    &=&\langle h,[A-\langle h,Ah\rangle I]^2h\rangle\crcr 
    &=&\langle h,A^2h\rangle,
\end{eqnarray}
where the last expression holds if we put the mean of $A$ to zero.

For vectors $h$ in $E(\lambda_i)$ there is no uncertainty: $\mathbb{P}_h(X=\lambda_i)=1$, $\mathbb{P}_h(X=\lambda_j)=0$ for $j\ne i$. Such vectors are called {\it pure states} for the observable $A$. A consequence of the spectral theorem for self-adjoint matrices is that there is a orthonormal bases of pure states; any state is {\it super-position} (linear combination) of pure states.

\subsection{Commuting observables}\label{SSectCommutingObs}
After the measure of $A$ in the state $h$ is performed, and the eigenvalue $\lambda_i$ is its outcome, the system {\it collapses} to the state $\Pi_{E(\lambda_i)}h$ (which has to be properly normalized by considering $\Pi_{E(\lambda_i)}h/\|\Pi_{E(\lambda_i)}h\|$)\footnote{You might be worried about the denominator: where does $h$ collapse when $\Pi_{E(\lambda_i)}h=0$? That case has null probability, $\mathbb{P}_h(X=\lambda_i)=\|\Pi_{E(\lambda_i)}h\|^2=0$, and we can simply ignore it.}. Suppose we have a second observable $B$, with associated random variable $Y$, eigenvalues $\mu_k$, and eigenspaces $F(\mu_k)$. If the initial state is $h$, and the measurement $A$ is performed with outcome $\lambda_i$, then the observation of $B$ is performed in the state $\Pi_{E(\lambda_j)}h$. By \eqref{eqProbabilityQM}, this means that:
\begin{eqnarray}\label{eqSeconObservation}
    \mathbb{P}_h(Y=\mu_k|X=\lambda_i)&=&\mathbb{P}_{\Pi_{E(\lambda_i)h}}(Y=\mu_k)\crcr 
    &=&\frac{\|\Pi_{F(\mu_k)}\Pi_{E(\lambda_i)}h\|^2}{\|\Pi_{E(\lambda_i)}h\|^2}, \text{ and }\crcr 
    \mathbb{P}_h(Y=\mu_k \text{ after } X=\lambda_i)&=&\|\Pi_{F(\mu_k)}\Pi_{E(\lambda_i)}h\|^2.
\end{eqnarray}
Here we can see a peculiar feature of quantum theory. If the measurement of $A$ did not affect the system in those parts measured by $B$, we would have equality from first to second line in:
\begin{eqnarray}\label{eqCompatibleQMone}
    \|\Pi_{F(\mu_k)}h\|^2&=&\mathbb{P}_h(Y=\mu_k)\crcr 
    &=&\sum_i \mathbb{P}_h(Y=\mu_k \text{ after } X=\lambda_i)\crcr 
    &=&\sum_i \|\Pi_{F(\mu_k)}\Pi_{E(\lambda_i)}h\|^2,
\end{eqnarray}
simply because the events $\{X=\lambda_i\}$ exhaust all possible outcomes of the measure of $A$.
This is certainly the case if $[A,B]:=AB-BA=0$, since in that case 
\begin{equation}\label{eqCompatibleQMthree}
\Pi_{F(\mu_k)}\Pi_{E(\lambda_i)}=\Pi_{E(\lambda_i)}\Pi_{F(\mu_k)},
\end{equation}
and the equality between the first and the last terms in \eqref{eqCompatibleQMone} would follow from the fact that $\{\Pi_{E(\lambda_i)}\}_i$ is pvm, so that
\begin{equation}\label{eqCompatibleQMtwo}
    \mathbb{P}_h(Y=\mu_k)=\sum_i \mathbb{P}_h(Y=\mu_k \text{ after } X=\lambda_i).
\end{equation}
We say that the the observables represented by $A$ and $B$ are {\it compatible} if \eqref{eqCompatibleQMtwo} holds. We wrote "after", not "and", because the system, "after" measuring $X$, is in a different state. The event $\{Y=\mu_k,\text{ in the state }h\}$, that is, is not a priori (and typically is not, as we shall see below), the union $\cup_i\{Y=\mu_k \text{ after }X=\lambda_i,\text{ in the state }h\}$.

In the finite dimensional case, compatibility is easily characterized in terms of spectral theory. The extension of this and other notion to infinite dimensions requires much more technicalities, although one can guess what the statements are.
\begin{theorem}[Compatibility theorem]\label{theoCompatibilityQM}
    The following are equivalent.
    \begin{enumerate}
        \item[(i)] The observables represented by $A$ and $B$ are compatible.
        \item[(ii)] $A$ and $B$ commute.
        \item[(iii)]  There is a common system of eigenvectors for $A$ and $B$ (i.e. $A$ and $B$ can be simultaneously diagonalized).
    \end{enumerate}
\end{theorem}
\begin{proof}
    \noindent (ii)$\implies$(i). By the spectral theorem, $[A,B]=0$ if and only if $[f(A),B]=0$ for all $f:\sigma(A)\to\mathbb{C}$, if and only if $[f(A),g(B)]=0$ for all $f:\sigma(A)\to\mathbb{C}$ and $g:\sigma(B)\to\mathbb{C}$, which implies that $[\Pi_{F(\mu_k)},\Pi_{E(\lambda_i)}]=0$ for all $i,j$ (and by inserting this in the spectral resolutions of $A$ and $B$, we see that in fact the last condition is equivalent to (ii)). We have observed above that (i) follows.

    \noindent (iii)$\implies$(ii) is clear: with the respect to the common eigenbasis, both $A$ and $B$ are diagonal matrices. 

    \noindent (i)$\implies$(iii). Let $x$ be a vector in $F(\mu_k)$. By \eqref{eqCompatibleQMone},
    \begin{eqnarray*}
       \sum_i \|\Pi_{F(\mu_k)}\Pi_{E(\lambda_i)}x\|^2&=&\|\Pi_{F(\mu_k)}x\|^2\crcr 
       &=&\|x\|^2=\sum_i \|\Pi_{E(\lambda_i)}\Pi_{F(\mu_k)}x\|^2\crcr 
       &=&\sum_i \|\Pi_{E(\lambda_i)}x\|^2.
    \end{eqnarray*}
    Each summand of the first term is less or equal to each sum of the last term, then 
    \[
    \|\Pi_{F(\mu_k)}\Pi_{E(\lambda_i)}h\|^2=\|\Pi_{E(\lambda_i)}x\|^2.
    \]
    This can only happen if $\Pi_{E(\lambda_i)}x\in F(\mu_k)$. i.e. $\Pi_{E(\lambda_i)}:F(\mu_k)\to F(\mu_k)$ is the projection onto $F(\mu_k)\cap E(\lambda_i)$. Since the latter spaces form an orthogonal decomposition of $F(\mu_k)$, there is an orthonormal basis $\{e_j^k\}_{j=1}^{\dim(F(\mu_k))}$ of $F(\mu_k)$, with respect to which each $\Pi_{E(\lambda_i)}|_{F(\mu_k)}$ is diagonal. Hence, with respect to the basis $\bigcup_k\{e_j^k\}_{j=1}^{\dim(F(\mu_k))}$ both $A$ and $B$ are diagonal.
\end{proof}
\begin{corollary}\label{corCompatibilityQM}
    Two observables $A$ and $B$ are compatible if and only if for all states $h$:
    \[
    \mathbb{P}_h(Y=\mu_k \text{ after } X=\lambda_i)=\mathbb{P}_h( X=\lambda_i\text{ after } Y=\mu_k).
    \]
\end{corollary}
\begin{exercise}\label{exeCompatibilityQM}
    Prove corollary \ref{corCompatibilityQM}.
\end{exercise}
\subsection{Complete systems of commuting observables}\label{SSectCompleteCommuting}
The space $H$ expresses all possible states of the system, and two states are physically distinct if it is possible to devise a measurement which assumes different values (or, at least, different expectations of such values) if the system is in one or the other state. In modelling physical systems, it is important that there are no redundancies: the state space $H$ must be large enough to accommodate all physical states of the system, and no more.

In the finite dimensional case, a single observable $A$ does the job if and only if its eigenvalues are simple, $m_1=\dots=m_L=1$. If $m_i\ge 2$, in fact, there are distinct pure states $h,k\in E(\lambda_i)$ which return the same values for the observable $A$. 
The problem arises, then, of understanding when a set $\{A_1,\dots,A_M\}$ of commuting observables, posssibly having eigenvalues with higher multiplicities (that is, some eigenvalues are{\it degenerate}) has the property that, for each eigenvalue $\lambda$, the intersection $E_1(\lambda)\cap\dots\cap E_M(\lambda)$ of the eigenspaces $E_j(\lambda)$'s corresponding to the $A_j$'s is at most one dimensional.  Such set $\{A_1,\dots,A_M\}$ is called a {\it complete system of commuting operators}. In the finite dimensional theory these problems are elementary. In the infinite dimensional case they are less elementary and they are the subject of the {\it multiplicity theory}, which is best seen having at hands the multiplication form of the spectral theorem, which we will see later on.

We see here how the multiplication form of the spectral theorem looks in the finite dimensional case. We saw in \eqref{eqPVMspectralResolution} and \eqref{eqPVMspectralResolutionTwo} the pvm associate to a self-adjoint operator on a finite dimensional $H$. If $\dim(E(\lambda_i))=1$ for all $i$'s, and $e_i$ is a unitary vector in $E(\lambda_i)$, map $U:e_i\mapsto (f_i:\lambda_j\mapsto \delta_{ij})$, where $\delta_{ij}$ is the Kroenecker symbol, 
\[
Ue_i(\lambda_j)=\delta_{i,j},\ U(\sum_ia_ie_i)(\lambda_j)=a_j,
\]
defines a unitary map from $H$ to $L^2(\mathbb{R},\mu)$, where $\mu=\sum_i\delta_i$. 
We have then:
\begin{equation}\label{eqPVMspectralResolutionFour}
   (UAU^{-1}\varphi)(\lambda_j)=\lambda_j\varphi(\lambda_j),
\end{equation}
i.e. the operator $A$ is unitarily equivalent to the operator "multiplication times $\lambda$" on $L^2(\mu)$. Let's verify \eqref{eqPVMspectralResolutionFour}.
\begin{eqnarray*}
    \varphi&=&\sum_i\varphi(\lambda_i)f_i,\crcr 
    U^{-1}\varphi&=&\sum_i\varphi(\lambda_i)e_i,\crcr 
    AU^{-1}\varphi&=&\sum_i\lambda_i\varphi(\lambda_i)e_i,\crcr
    (UAU^{-1}\varphi)(\lambda_j)&=&\sum_i\lambda_i\varphi(\lambda_i)(f_i(\lambda_j))\crcr 
    &=&\lambda_j\varphi(\lambda_j).
\end{eqnarray*}
In the general case, eigenvalues $\lambda_i$ can have multiplicity $m_i\ge2$, and we have to split the above procedure accordingly. Let $m=\max(m_1,\dots.m_M)$ be the maximum of multiplicities. Rearrange so that $m=m_1\ge m_2\ge \dots\ge m_M$. For each eigenspace $E(\lambda_k)$ consider a orthonormal basis 
$\{e_j^k\}_{j=1}^{m_k}$ (the choice of which is not canonical); pick $m$ copies of $\mathbb{R}$, and on the $l^{th}$ copy consider the measure
\[
\mu_l=\sum_{k:m_k\ge l}\delta_{\lambda_k}.
\]
This way, $\delta_{\lambda_k}$ appears in the measures $\mu_1,\dots,\mu_{m_k}$, and not in $\mu_{m_k+1},\dots,\mu_m$. Let $k(l)$ to be the largest $k$ for which $m_k\ge l$.

Pick then a vector $x=\sum_{k}\sum_{i=1}^{m_k}a_i^k e_i^k$. Associate to it the family $Ux:=\{\varphi_l\}_{l=1}^{m}\in L^2(\mu_1)\oplus\dots\oplus L^2(\mu_m)$, where
\[
\varphi_l=\sum_{k:m_k\ge l}a_l^k f_l^k,
\]
and where, as before, $f_l^k(\lambda_j)=\delta_{k,j}$ for $1\le j\le k(l)$. This provides a constructive proof of the following.
\begin{theorem}[Spectral theorem finite dimensions: multiplicative form]\label{theoMultiplicativeSTfinite}
    Let $A$ be a self-adjoint operator on a finite dimensional space. Then, there are measures $\mu_1,\dots,\mu_M$ on $\mathbb{R}$ (sums of Dirac deltas at disjoint real points), and there is a unitary operator
    \[
    U:H\to L^2(\mu_1)\oplus\dots\oplus L^2(\mu_M),\ Ux=\varphi=(\varphi_1,\dots,\varphi_M),
    \]
    such that
    \[
    (UAU^{-1}\varphi)_l(\lambda)=\lambda\varphi_l(\lambda).
    \]
    Moreover, $\dim(H)$ is the total number of the Dirac's deltas involved.
\end{theorem}
We can rephrase the theorem in a language involving just one measure space. Let $X=\mathbb{R}\times\{1,\dots,M\}$, endowed with the topology generated by the open sets of each factors. Identify $\mu_j$ with the corresponding measure on $\mathbb{R}\times\{j\}$, and let $\mu=\mu_1+\dots+\mu_M$, which is a measure on $X$. Consider then the unitary map
\[
V:L^2(\mu_1)\oplus\dots\oplus L^2(\mu_M)\to L^2(\mu),\ V(\varphi_1,\dots,\varphi_M)(\lambda,j)=\varphi_j(\lambda).
\]
Then,
\[
VUAU^{-1}V^{-1}\psi(\lambda,j)=g(\lambda,j)\psi(\lambda,j), \text{ where }g(\lambda,j)=\lambda.
\]
This is a special case of the multiplicative form of the spectral theorem, stating that bounded, self-adjoint operators on Hilbert spaces can be always be represented as multiplication operators on (Borel) $L^2$ spaces. 

The presence of several summands is indication of how $\{A\}$ alone fails to be a complete system of  commuting observables. Multiplicity theory provides a more quantitative measure of the same.
\subsection{Uncertainty}\label{SSectUncertaintyFinite}
We saw that commuting self-adjoint operators $A$ and $B$ are compatible in the sense the measure of $A$ does not affect the measure of $B$. {\it Compatibility} should not be confused with {\it independence} of the random variables $X,Y$ associated with $A,B$ with respect to the probability measure $\mathbb{P}_h$ associated with a state $h\in H$, as the case $A=B$ clearly illustrates.

Lack of commutativity, however, has probabilistic consequences.
\begin{theorem}[Heisenberg inequality]\label{theoHeisenbergInequalityBaby}
Suppose that $A$ and $B$ are self-adjoint  and that they have zero mean with respect to a state $h$, $0=\langle h,Ah\rangle=\langle h,Bh\rangle$. Then,
    \begin{equation}\label{eqHeisenbergInequalityBaby}
       1/4 |\langle h,[A,B]h\rangle|^2\le \langle h,A^2h\rangle\langle h,B^2h\rangle=\text{Var}_h(X)\text{Var}_h(Y).
    \end{equation}
\end{theorem}
Actually, we will prove a slightly sharper inequality. The main consequence of \eqref{eqHeisenbergInequalityBaby} is that $[A,B]$ sets a lower bound to the simultaneous concentration of $X$ and $Y$ around their means.
\begin{proof}
    We set $\langle A\rangle_h=\langle h,Ah\rangle$ (in physics the subscript $h$ is implicit). If $A$ and $B$ are self-adjoint, $[A,B]$ is a measure of how much $AB$ fails to be self-adjoint: $AB-(AB)^\ast=AB-BA=[A,B]$. We can write $AB$ in "algebraic form",
    \[
    AB=\frac{\{A,B\}}{2}+i\frac{[A,B]}{2i}=S+iT, \text{ with }S=S^\ast \text{ and }T=T^\ast.
    \]
    where we use the notation $\{A,B\}=AB+BA$. 
    An easy calculation gives 
    \begin{equation*}
        \langle h,1/2\{A-\langle A\rangle_h,B-\langle B\rangle_h\}h \rangle =\langle h,[1/2\{A,B\}-\langle A\rangle_h\langle B\rangle_h ]h \rangle,
    \end{equation*}
    and clearly
    \[
    \langle h,[A-a,B-b]h\rangle=\langle h,[A,B]h\rangle
    \]
    for all constants $a,b$.
    We have then:
    \begin{eqnarray}\label{eqOggiNonSiVaAlCIRM}\left|1/(2i)\langle h,[A,B]h\rangle\right|^2&\le&\left|\langle h,[1/2\{A,B\}-\langle A\rangle_h\langle B\rangle_h ]h \rangle\right|^2+\left|1/(2i)\langle h,[A,B]h\rangle\right|^2\crcr 
        &=&|\langle h,(A-\langle A\rangle_h)(B-\langle B\rangle_h)h\rangle|^2\crcr
        &=&|\langle (A-\langle A\rangle_h)h,(B-\langle B\rangle_h)h\rangle|^2\crcr 
        &\le&\|(A-\langle A\rangle_h)h\|^2\|(B-\langle B\rangle_h)h\|^2\crcr 
        &=&\text{Var}_h(X)\text{Var}_h(Y),
    \end{eqnarray}
\end{proof}
The stronger inequality, starting from the right hand side of the first inequality in \eqref{eqOggiNonSiVaAlCIRM}, is known as the {\it Robertson-Schroedinger uncertainty relation}. One of its advantages is that it does not become trivial in the case $[A,B]=0$; e.g. when $A=B$. In a a sense it carries both classical and quantum probabilistic information.
\subsection{Evolution and the Schroedinger equation}\label{SSectSchroedingerEquation}
Quantum physical systems evolve, and, if they are isolated, they evolve according to the following rules. First a bit of notation. If $x$ is a state of the system, $U(t)x$ is the state of that system at time $t>0$. This way, $U(0)=I$ is the identity.
\begin{enumerate}
    \item[(a)] Distinguishable states remain distinguishable: 
    \[
    \langle x,y\rangle=0\implies\langle U(t)x,U(t)y\rangle=0.
    \]
    \item[(b)] If $\|x\|=1$, we normalize $U(t)x$ to have $\|U(x)\|=1$. 
    \item[(c)] $U(t+s)x=U(t)U(s)x$, which is self-explicatory: "the evolution of the evolution...".
    \item[(d)] $U(t)$ is linear. Together with (b), this implies that $\|U(t)x\|=\|x\|$. Linearity is a nontrivial requirement, whose roots are in experiments.
    \item[(e)] Some continuity is required...
\end{enumerate}
In purely mathematical terms, (a-e) mean:
\begin{enumerate}
    \item[(i)] each $U(t)$ is a unitary operator, $U(t)^\ast U(t)=U(t)U(t)^\ast=I$;
    \item[(ii)] $t\mapsto U(t)$ defines a semigroup of such operators;
    \item[(iii)] some continuity...
\end{enumerate}
We consider here the special case where
\begin{equation}\label{eqSemigroupsWithBoundedInfGen}
    U(t)=e^{itA},\ t\ge0,\text{ with $A$ bounded}.
\end{equation}
\begin{theorem}\label{theoSemigroupsWithBoundedInfGen}
    Let $A$ be a self-adjoint operator on a Hilbert space $H$.
    \begin{enumerate}
        \item[(i)] $U(t)$ is unitary.
        \item[(ii)] $t\mapsto U(t)$ is a semi-group for all bounded operators $A$.
        \item[(iii)] The Schroedinger equation\index{Schroedinger equation}
        \begin{equation}\label{eqSchroedingerEasy}
            \frac{1}{i}\frac{d U(t)}{dt}=A U(t)
        \end{equation}
    holds in the norm topology:
        \[
        \lim_{\Delta\to0}\frac{U(t+\Delta)-U(t)}{\Delta}=A U(t)
        \]
        in the operator norm.
        \item[(iv)] The semigroup is {\bf strongly continuous},
        \begin{equation}\label{eqStrongContSemigroup}
            \lim_{\Delta\to0}\|U(\Delta)x-x\|=0
        \end{equation}
        for all $x$ in $H$.
    \end{enumerate}
\end{theorem}
\begin{proof} Since the functions $z\mapsto e^{tz}$ are entire, we can use the holomorphic calculus. We can write
\[
U(t)=e^{itA}=\sum_{n=0}^\infty \frac{(itA)^n}{n!},
\]
where the series converges in the operator norm, and all usual manipulations of series hold, provided we take into account the eventual break of commutativity  (which here does not occur). We have $\|U(t)\|\le e^{t\|A\|}$. Item (ii) follows immediately. 

It is easily verified that
\[
U(t)^\ast=e^{-itA}.
\]
and this gives (i):
    \begin{eqnarray*}
        \langle U(t)x,U(t)y\rangle&=&\langle U(t)^\ast U(t)x,y\rangle\crcr 
        &=&\langle x,y\rangle.
    \end{eqnarray*}
Since $U(t+\Delta)=U(\Delta)U(t)$, it suffices to verify (iii) at $t=0$, which we can do using series:
\begin{eqnarray*}
    \frac{1}{i}\frac{U(\Delta)-I}{\Delta}-A&=&\sum_{n=1}^\infty\frac{A^n(i\Delta)^{n-1}}{n!}-A\crcr 
    &=&(i\Delta)\sum_{n=2}^\infty\frac{A^n(i\Delta)^{n-2}}{n!},
\end{eqnarray*}
which tends to zero in the operator norm as $\Delta\to0$. (iv) is an even weaker statement.
\end{proof}
All this is very nice, but also highly unsatisfactory. The Schroedinger equation for a free particle, in fact, is (but for the Plank constant) of the form \eqref{eqSchroedingerEasy}, but with $H=L^2(\mathbb{R}^3)$ and $A=\Delta_{\mathbb{R}^3}$, the Laplace operator, which is {\it not } bounded (and not even everywhere defined, in fact) on $H$. This provided the main motivation for von Neumann to develop the spectral theory for unbounded self-adjoint operators, which also has applications all over mathematics. One of main items in the unbounded case is the following theorem of Stone, which extends theorem \ref{theoSemigroupsWithBoundedInfGen}.
\begin{theorem}[Stone's theorem]\label{theStone}\index{theorem!Stone}
    Let $t\mapsto U(t)$ be a {\it strongly continuous semigroup} of unitary operators on $H$ (that is, (iv) in theorem \ref{theoSemigroupsWithBoundedInfGen} holds). Let $D(A)$ be the set of those $x$ in $H$ such that
    \[
    \lim_{\Delta\to0}\frac{U(\Delta)x-x}{i\Delta}=Ax
    \]
    exists in norm.
    
    Then, $D(A)$ is dense in $H$. The operator $A:D(A)\to H$ is self-adjoint, and $U(t)=e^{itA}$.

    Viceversa, if $A$ is self-adjoint (and possibly unbounded), then $U(t)=e^{itA}$ defines a strongly continuous unitary semigroup.
\end{theorem}
Of course, we should define what an {\it unbounded self-adjoint operator} is, and what we mean by $e^{itA}$, which requires a more general version of the spectral theorem. Such tools were developed by von Neumann, in fact, in the late 20's of the XX century, and they are the main topic we discuss in this chapter.

If the system starts in the state $x$ and the dynamics are given by $A$, at time $t$ the system is in the state $e^{itA}x$. Suppose we wish to perform a measure of the observable $B$ at time $t$. How does the expectation of the random variable $Y$ associate to $B$ change over time?
\[
\mathbb{E}_{U(t)x}Y=\langle U(t)x,BU(t)x\rangle=\rangle x,U(-t)BU(t)x\rangle=\langle x,e^{-itA}Be^{itA}x\rangle,
\]
hence (at least formally),
\[
\frac{d\ }{dt}\mathbb{E}_{U(t)x}Y=(-iu)\langle x,e^{-itA}(AB-BA)e^{itA}x\rangle=-i\langle U(t)x,[A,B]U(t)x\rangle.
\]
In particular, the mean remains constant if and only if $[A,B]=0$: the observables must be compatible. 
The case $A=B$ is interpreted in physics as {\it conservation of energy}.

Schroedinger's evolution operator can be written in terms of the pvm associated with $A$,
\[
A=\int_{-\infty}^{+\infty}\lambda \Pi_{d\lambda}.
\]
By the spectral theorem in pvm form,
\begin{equation}\label{eqFTandSchroedinger}
    e^{itA}=\int_{-\infty}^{+\infty}e^{i\lambda t} \Pi_{d\lambda},
\end{equation}
which makes sense for $t\in\mathbb{R}$, and it defines a notion of Fourier transform for the pvm $\Pi$, $\widehat{\Pi}(t)=e^{iAt}$. If you are familiar with Bochner's theorem, and you observe that $\Pi$ exhibits some positivity ($\Pi_E\ge0$ as operator), then you expect as well its "Fourier transform" to have a property similar to "positive definiteness":
\[
\sum_{j,k=1}^nc_j\overline{c_k}e^{i(t_j-t_k)A}\ge0
\]
as operator, for any choice of complex $\alpha_1,\dots,\alpha_n$ and real $t_1,\dots,t_n$. The guess is correct,
\[
\left\langle x,\sum_{j,k=1}^nc_j\overline{c_k}e^{i(t_j-t_k)A}x\right\rangle= \left\|\sum_{j=1}^nc_je^{it_jA}x\right\|^2\ge0.
\]

At some point one should talk about entenglement, but you will have to read it elsewhere.

\section{The spectral theorem for unitary operators}\label{SectUnitSpectral}
Here we follow one of the two proofs in Michael Taylor's lecture notes \href{https://mtaylor.web.unc.edu/wp-content/uploads/sites/16915/2018/04/specthm.pdf}{The Spectral Theorem for Self-Adjoint and Unitary Operators}. The proof only requires basic Fourier series.

Recall that a linear operator $H\xrightarrow{U}H$ is {\it unitary}\index{operator!unitary} iff and only if $U^\ast U=UU^\ast=I$ is he identity. This is equivalent to the property of preserving the inner product,
\begin{equation}\label{eqUnitaryIP}
    \langle Ux,Uy\rangle=\langle x,y\rangle,
\end{equation}
whenever $x,y\in H$. We have clearly that $\|U\|=1$, but there are plenty of operators having norm one which are not unitary.

An abundant source of unitary operators is provided by specific multiplication operators on $L^2$ spaces. If $(X,\nu)$ is a measure space and $\varphi:X\to\mathbb{R}$ is measurable, then 
\begin{equation}\label{eqUnitaryMultExe}
    M_{e^{i\varphi}}:f\mapsto e^{i\varphi}f
\end{equation}
is unitary on $L^2(\nu)$. Also observe that $M_{e^{i\varphi}}^\ast=M_{e^{-i\varphi}}$. At the end of the day, we learn that all unitary operators are unitarily equivalent to some $M_{e^{i\varphi}}$ for some measure space $(X,\nu)$.

\subsection{The spectrum of a unitary operator}\label{SSectUnitarySpec}\index{spectrum!unitary operator}
Recall that $\mathbb{T}=\{e^{it}:t\in (-\pi,\pi]\}$ is the unit circle in the complex plane (the {\it torus}).
\begin{proposition}\label{propSpecUO}
    If $H\xrightarrow{U}H$ is unitary, then $\sigma(U)\subseteq\mathbb{T}$.
\end{proposition}
\begin{proof}
    If $|\lambda|>1$, then $U-\lambda I=\lambda(U/\lambda-I)$ is invertible because $\|U/\lambda\|<1$, hence $\lambda\in \rho(U)$. If $|\lambda|<1$, then $U-\lambda I=U(I-\lambda U^\ast)$ is invertible because $U$ and $I-\lambda U^\ast$ are, so again $\lambda\in\rho(U)$.
\end{proof}
\subsection{The continuous functional calculus for a unitary operator}\label{SSectUnitaryFC}
If $p(e^{it})=p(e^{it})=\sum_{n=-N}^N a_n e^{int}$ is a trigonometric polynomial, it natural to define
\begin{equation}\label{eqPolCalc}
p(U)=\sum_{n=-N}^N a_n U^n.
\end{equation}
The obvious {\it polynomial calculus} defined by \eqref{eqPolCalc} generates a homomorphism of algebras, $(pq)(U)=p(U)q(U)$.

This definition immediately extends to functions $f:\mathbb{T}\to\mathbb{C}$ whose Fourier series converges absolutely,
\begin{equation}\label{eqWienerAlg}
    \frac{1}{2\pi}\sum_{n=-\infty}^\infty|\widehat{f}(n)|=\|f\|_{A(\mathbb{T})}<\infty.
\end{equation}
The space of such functions is called \href{https://en.wikipedia.org/wiki/Wiener_algebra}{Wiener algebra}. However, we just want to have a dense subspace of $L^2(\mathbb{T})$, so we consider $f\in C^2(\mathbb{T}):=C^2_{per}$. For such $f$'s, define
\begin{equation}\label{eqFCunitaryBis}
    f(U)=\frac{1}{2\pi}\sum_{n=-\infty}^\infty\widehat{f}(n)U^n,
\end{equation}
which converges in operator norm. For such "$C^2(\mathbb{T})$ calculus" the following properties are easily established.
\begin{enumerate}
    \item[(i)]  {\bf If $f\in C^2(\mathbb{T})$, then $f(U)^\ast=\overline{f}(U)$, where $\overline{f}(e^{it}):=\overline{f(e^{it})}$.} Using that
    \[
    \overline{\widehat{f}(n)}=\int_{-\pi}^\pi \overline{f(t)}e^{int}dt=\widehat{\overline{f}}(-n),
    \]
    in fact,
    \begin{eqnarray*}
        f(U)^\ast&=&\frac{1}{2\pi}\sum_{n=-\infty}^\infty\overline{\widehat{f}(n)}U^{\ast n}\\ 
        &=&\frac{1}{2\pi}\sum_{n=-\infty}^\infty\widehat{\overline{f}}(-n)U^{-n}\\ 
        &=&\overline{f}(U).
    \end{eqnarray*}
    \item[(ii)]  {\bf If $f,g\in C^2(\mathbb{T})$, then $(fg)(U)=f(U)g(U)$.} Absolute convergence of the corresponding Fourier series ensures that the following equalities hold:
    \begin{eqnarray*}
        (fg)(U)&=&\frac{1}{(2\pi)^2}\sum_{m,n\in\mathbb{Z}}\widehat{f}(n-m)\widehat{g}(m)U^n\\ 
        &=&\frac{1}{2\pi}\sum_{m\in\mathbb{Z}}\left(\frac{1}{2\pi}\sum_{n\in\mathbb{Z}}\widehat{f}(n-m)U^{n-m}\right)\widehat{g}(m)U^m\\ 
        &=&f(U)g(U).
        \end{eqnarray*}
        \item[(iii)] {\bf If $f\ge0$ and $f\in C^2(\mathbb{T})$, then $\langle x,f(U)x\rangle\ge0$ for all $x$ in $H$.} If $\epsilon>0$, then $f+\epsilon\in C^2(\mathbb{T})$ is strictly positive, hence $\sqrt{f+\epsilon}\in C^2(\mathbb{T})$. Since $\sqrt{f+\epsilon}$ is real valued, by (i) $(\sqrt{f+\epsilon})(U)$ is self-adjoint. Then,
        \begin{eqnarray*}
            \langle x,f(U)x\rangle+\epsilon\|x\|^2&=&\langle x,[f(U)+\epsilon I])x\rangle=\langle x,(\sqrt{f+\epsilon})(U)(\sqrt{f+\epsilon})(U)x\rangle\\ 
            &=&\langle (\sqrt{f+\epsilon})(U) x,(\sqrt{f+\epsilon})(U)x\rangle=\|(\sqrt{f+\epsilon})(U)x\|^2\\
            &\ge&0.
        \end{eqnarray*}
    Hence, $\langle x,f(U)x\rangle\ge0$.
    \item[(iv)]  {\bf If $f\in C^2(\mathbb{T})$, then $\|f(U)x\|\le\|f\|_{L^\infty_0}\|x\|$.} The function $g(e^{it})=\|f\|_{L^\infty_0}^2-|f|^2$ is positive, hence, by (iii),
    \begin{eqnarray*}
        0&\le&\langle x,(\|f\|_{L^\infty_0}^2-|f|^2)(U)x\rangle=\|f\|_{L^\infty_0}\|x\|^2-\langle x,(\overline{f}f)(U)x\rangle\\ 
        &=&\|f\|_{L^\infty_0}\|x\|^2-\langle x,f^\ast(U)f(U)x\rangle\\ 
        &=&\|f\|_{L^\infty_0}\|x\|^2-\|f(U)x\|^2.
    \end{eqnarray*}
\end{enumerate}
Item (iv) means that the regularity assumption $f\in C^2$ is somewhat redundant, since only the $L^\infty$ norm of $f$ is relevant in order to estimate $\|f(U)\|$.
Based on this, we can extend the $C^2(\mathbb{T})$ calculus to a {\it continuous calculus}. \index{continuous calculus!unitary operator}
\begin{theorem}\label{theoCCunitary} Let $U$ be a unitary operator on a Hilbert space $H$.
    There is a unique homomorphism of algebras $\Lambda:C(\mathbb{T})\to\mathcal{B}(H)$, $\Lambda(f)=f(U)$, such that
    \begin{enumerate}
        \item[(i)] if $e(e^{it})=e^{it}$ is the identity function on $\mathbb{T}$, then $\Lambda(e)=U$ (hence, $\Lambda(p)=p(U)$ when $p$ is a trigonometric polynomial);
        \item[(ii)] $\Lambda$ is continuous when $\mathcal{B}(H)$ is endowed with the uniform operator topology (the one induced by the operator norm), $\|\Lambda(f)\|\le\|f\|_{L^\infty_0}$.
    \end{enumerate}
    Moreover, if $f\in C(\mathbb{T})$,
    \begin{enumerate}
        \item[(iii)] $f(U)^\ast=\overline{f}(U)$;
        \item[(iv)] if $f$ is real valued, then $f(U)=f(U)^\ast$, and, if $f\ge0$, then $f(U)\ge 0$;
    \end{enumerate}
\end{theorem}
\begin{proof}
    The map $\Lambda$ was introduced above on the subspace of the trigonometric polynomial in accordance to (i). Given $f\in C^2(\mathbb{T})$,
    and 
    \[
    (S_Nf)(e^{it})=\frac{1}{2\pi}\sum_{|n|\le N}\widehat{f}(n)e^{int},   
    \]
    we have that $\|S_Nf-f\|_{L^\infty_0}\to0$ as $N\to\infty$, hence, by item (iv) above,
    \[
    \lim_{N\to\infty}\|\Lambda(S_Nf)-f(U)\|=\lim_{N\to\infty}\|(S_Nf)(U)-f(U)\|=0,
    \]
    where $f(U)$ was defined earlier; hence, $f(U)=\Lambda(f)$.

    If $g\in C(\mathbb{T})$ and $\{f_n\}$ is a sequence in $C^2(\mathbb{T})$ converging to $g$ uniformly, then $\{\Lambda(f_n)\}$ is a Cauchy sequence $\mathcal{B}(H)$ with respect to the norm topology, hence it converges to some operator $A\in\mathcal{B}(H)$. By requirement (ii), $A=\Lambda(g)$. It is readily verified that any other sequence converging to $g$ provides the same definition of $\Lambda(g)$. The fact that $\Lambda$ is a homomorphism of algebras is routine.
    
    This shows existence and uniqueness of a map $\Lambda$ with the properties (i), (ii). The proof of (iii) and (iv) is left as an exercise.
\end{proof}
The map $\Lambda$ needs not being injective. If $H$ is finite dimensional, for instance, $C(\mathbb{T})$ is infinite dimensional, whereas $\mathcal{L}(H)$ is finite dimensional.

\subsection{The spectral measures and the measurable calculus}\label{SSectSpecMeasUO}
Fix a unitary operator $U$ on a Hilbert space $H$.
By theorem \ref{theoCCunitary}, for each $x\in H$, the map
\begin{equation}
    \Lambda_{x,x}:f\mapsto \langle x,f(U)x\rangle,\ \Lambda_{x,x}: C(\mathbb{T})\to \mathbb{C},
\end{equation}
is linear, continuous, and positive. By the Riesz-Markov-Kakutani\index{theorem!Riesz-Markov-Kakutani for
positive measures} representation theorem for measures, there exists a (unique) finite Borel measure $\mu_{x,x}$ on $\mathbb{T}$ such that
\begin{equation}\label{eqSpecMeas}
\langle x,f(U)x\rangle=\int_{\mathbb{T}}f(e^{it})d\mu_{x,x}(e^{it}).
\end{equation}
By polarization, for each $x,y\in H$ we can find a complex valued, finite Borel measure $\mu_{x,y}$ such that
\begin{equation}\label{eqSpecMeasPol}
\langle x,f(U)y\rangle=\int_{\mathbb{T}}f(e^{it})d\mu_{x,y}(e^{it}).
\end{equation}
The measures $\mu_{x,y}$ are the {\it spectral measures} associated with $U$. Observe that the association is {\it canonical}: they are uniquely determined by $U$.\index{spectral measures!unitary operator}
We will need the following.
\begin{lemma}\label{lemmaMeasOUchange} If $U$ is unitary, $x,y\in H$, and $f\in C(\mathbb{T})$, then
    \begin{equation}\label{eqMeasOUchange}
        fd\mu_{x,y}=d\mu_{x,f(U)y}.
    \end{equation}
\end{lemma}
\begin{proof}
 For any other continuous $g$,
 \[
     \int_{\mathbb{T}}fgd\mu_{x,y}=\langle x,g(U)f(U)y\rangle= \int_{\mathbb{T}}gd\mu_{x,f(U)y}.
 \]
\end{proof}
The map $(x,y)\mapsto \mu_{x,y}$ is a sesquilinear map with values in the finite, complex measures on $\mathbb{T}$. 

We have, in particular,
\[
\int_{\mathbb{T}}e^{it}d\mu_{x,y}(e^{it})=\langle x,Uy\rangle,
\]
and $\mu_{x,y}(\mathbb{T})=\langle x,y\rangle$.
The advantage of the spectral measures is that we can define a {\it measurable calculus}.\index{measurable calculus!unitary operator}
\begin{theorem}\label{theoMeasCalcUO}[Measurable calculus for a unitary operator]
    There is a unique linear map $\Lambda:L^\infty_0(\mathbb{T})\to\mathcal{B}(H)$, $\Lambda(f)=f(U)$, such that:
    \begin{enumerate}
        \item[(i)] $\Lambda(U)=f(U)$ agrees with the previously defined one when $f\in C(\mathbb{T})$;
        \item[(ii)] $\Lambda$ is continuous in the sense that if $L^\infty(\mathbb{T})\ni f_n\to f\in L^\infty(\mathbb{T})$ pointwise and $\|f_n\|_{L^\infty_0}\le C$ is bounded, then $\Lambda(f_n)\to\Lambda(f)$ in the {\bf weak operator topology}.\index{operator topology!weak}
    \end{enumerate}
    Moreover, $\Lambda$ is a $\ast$-algebra homomorphism ($\Lambda(fg)=\Lambda(f)\Lambda(g)$, and $\Lambda(f)^\ast=\Lambda(\overline{f})$) and
    \begin{equation}\label{eqMeasCalcUO}
 \langle x,f(U)y\rangle=\int_{\mathbb{T}}fd\mu_{x,y}.     
    \end{equation}
\end{theorem}
\begin{proof}
We begin by showing that a linear map $\Lambda:L^\infty_0(\mathbb{T})\to\mathcal{B}(H)$ with the properties (i) and (ii) exists.  For $f\in L^\infty_0$ and $x,y\in H$, let
\[
[x,y]_f=\int_{\mathbb{T}}fd\mu_{x,y}.
\]
By uniqueness of the spectral measures, $[x,y]_f$ is sesquilinear, and 
\[
\left|[x,y]_f\right|\le C \|f\|_{L^\infty_0}
\]
if $\|x\|,\|y\|\le1$. By homogeneity, $\left|[x,y]_f\right|\le C \|f\|_{L^\infty_0}\|x\|\cdot \|y\|$ for all $x,y$, hence there exists a unique operator $\Lambda(U)=f(U)\in \mathcal{B}(H)$ such that
\[
[x,y]_f=\langle x,f(U)y\rangle.
\]
The map $\Lambda:f\mapsto f(U)$ just defined agrees with the previous definition when $f$ is continuous, so (i) holds. Moreover, it is clearly linear in the $f$ variable. 

Suppose $L^\infty(\mathbb{T})\ni f_n\to f\in L^\infty_0$ pointwise, and that they are uniformly bounded. For each $x,y$ in $H$, by dominated convergence,
\[
\lim_{n\to\infty}\langle x,f_n(U)y\rangle=\lim_{n\to\infty}\int_{\mathbb{T}}f_nd\mu_{x,y}=\int_{\mathbb{T}}fd\mu_{x,y}=\langle x,f(U)y\rangle,
\]
i.e. $f_n(U)\to f(U)$ is the weak operator topology, i.e. property (ii) holds.\index{operator topology!weak}

We show uniqueness. If (i) and (ii) hold for a linear map $\widetilde{\Lambda}$, the space $\mathcal{A}$ of functions $f:\mathbb{T}\to\mathbb{C}$ in $L^\infty_0(\mathbb{T})$ such that $\widetilde{\Lambda}(f)=\Lambda(f)$: (1) contains $C(\mathbb{T})$, and (2) is closed under pointwise limits of uniformly bounded functions. By the corollary \ref{corLebesgueHausdorff} to the Lebesgue-Hausdorff theorem\index{theorem!Lebesgue-Hausdorff}, $\mathcal{A}=L^\infty_0(\mathbb{T})$.

We are left with proving that $\Lambda$ is an $\ast$-algebra homomorphism. We start by extending lemma \ref{lemmaMeasOUchange} to $f\in L^\infty_0(\mathbb{T})$. For $n\in\mathbb{Z}$, let $e_n(t)=e^{int}$. We start by showing that
\[
\mu_{x,Uy}=\mu_{U^{-1}x,y}.
\]
In fact, for $n\in\mathbb{Z}$,
\[
\int_\mathbb{T}e^{imt}d\mu_{x,Uy}=\langle x,U^mUy\rangle=\langle U^{-1}x,U^my\rangle=\int_\mathbb{T}e^{imt}d\mu_{U^{-1}x,y},
\]
and a finite Radon measure $\mu$ on $\mathbb{T}$ is uniquely characterized by its Fourier coefficients $\widehat{\mu}(n)=\int_\mathbb{T}e^{-imt}d\mu$, $m\in\mathbb{Z}$.

Now, for $x,y\in H$, and $f\in L^\infty_0$,
\begin{eqnarray*}
    \int_\mathbb{T}e^{int}f(e^{it})d\mu_{x,y}(e^{it})&=&\int_\mathbb{T}f(e^{it})d\mu_{x,U^ny}(e^{it})=\int_\mathbb{T}f(e^{it})d\mu_{U^{-n}x,y}(e^{it})\\ 
    &=&\langle U^{-n}x,f(L)y\rangle=\langle x,U^nf(L)y\rangle\\ 
    &=&\int_\mathbb{T}e^{int}d\mu_{x,f(U)y}(e^{it}),
\end{eqnarray*}
which implies 
\begin{equation}\label{eqNoLH}
f(e^{it})d\mu_{x,y}(e^{it})=d\mu_{x,f(U)y}(e^{it}).
\end{equation}
With this at hands,
\begin{eqnarray*}
    \langle x,(fg)(U)y\rangle&=&\int_\mathbb{T}(fg)(e^{it})d\mu_{x,y}(e^{it})=\int_\mathbb{T}f(e^{it})g(e^{it})d\mu_{x,y}(e^{it})\\ 
    &=&\int_\mathbb{T}f(e^{it})d\mu_{x,g(U)y}(e^{it})=\int_\mathbb{T}f(e^{it})d\mu_{x,f(U)g(U)y}(e^{it})\\ 
    &=&\langle x,f(U)g(U)\rangle.
\end{eqnarray*}
\footnote{
{\bf A different proof of \eqref{eqNoLH}, using the Lebesgue-Hausdorff theorem.} For $x,y\in H$, let
\[
\mathcal{A}_{x,y}=\{f\in L^\infty_0(\mathbb{T}):\ fd\mu_{x,y}=d \mu_{x,f(U)y}\}.
\]
We know that $C(\mathbb{T})\subseteq\mathcal{A}_{x,y}$, and to apply the Lebesgue-Hausdorff theorem we only need to show that $f\in\mathcal{A}_{x,y}$ if $f$ is the pointwise limit of a uniformly bounded sequence $f_n\in\mathcal{A}_{x,y}$. For $g\in C(\mathbb{T})$, by dominated convergence and lemma \ref{lemmaMeasOUchange},
\begin{eqnarray*}
\int_{\mathbb{T}}fgd\mu_{x,y}&=&\lim_{n\to\infty}\int_{\mathbb{T}}f_ngd\mu_{x,y}=\lim_{n\to\infty}\int_{\mathbb{T}}gd\mu_{x,f_n(U)y}\\ 
&=&lim_{n\to\infty}\langle x,g(U)f_n(U)y\rangle=lim_{n\to\infty}\langle g(U)^\ast x,f_n(U)y\rangle\\ 
&=&\langle g(U)^\ast x,f(U)y\rangle=\langle x,g(U)f(U)y\rangle\\ 
&=&\int_{\mathbb{T}}gd\mu_{x,f(U)y},
\end{eqnarray*}
hence, $f\in \mathcal{A}_{x,y}$. Thus, lemma \ref{lemmaMeasOUchange} holds for $f\in L^\infty_0(\mathbb{T})$.
}

As a byproduct of the proof, we have that $(fg)(U)=f(U)g(U)$ if $f\in L^\infty(\mathbb{T})$ and $g\in C(\mathbb{T})$,
\[
\langle x,(fg)(U)\rangle=\int_{\mathbb{T}}fgd\mu_{x,y}=\langle x,g(U)f(U)y\rangle.
\]
We can now prove that $(fg)(U)=f(U)g(U)$ for $f,g\in L^\infty_0(\mathbb{T})$,
\begin{eqnarray*}
    \langle x,(fg)(U)y\rangle&=&\int_{\mathbb{T}}fgd\mu_{x,y}=\int_{\mathbb{T}}fd\mu_{x,g(U)y}=\int_{\mathbb{T}}d\mu_{x,f(U)g(U)y}\\ 
    &=&\langle x,f(U)g(U)y\rangle.
\end{eqnarray*}
To show that $\Lambda$ is a $\ast$-homomorphism,
\begin{eqnarray*}
\langle x,f(U)^\ast y\rangle&=&\langle f(U) x, y\rangle=\overline{\langle y,f(U) x\rangle}\\ 
&=&\overline{\int_{\mathbb{T}}fd\mu_{y,x} }=\int_{\mathbb{T}}\overline{f}d\overline{\mu_{y,x}}\\ 
&=&\int_{\mathbb{T}}\overline{f}d\mu_{x,y}=\langle x,\overline{f}(U) y\rangle,
\end{eqnarray*}
as wished.
\end{proof}
\subsection{From the measurable calculus to projection valued measures}\label{SSectUOpvm}
Let $E$ be a Borel subset of $\mathbb{T}$. We define $\mu(E)=\chi_E(U)$,
\begin{equation}\label{eqMissing}
    \langle x,\mu(E)y\rangle=\int_Ed\mu_{x,y}=\mu_{x,y}(E).
\end{equation}
The map $E\mapsto \mu(E)$ defines a function from the Borel $\sigma$-algebra into the set of the projections of $H$. In fact,
\begin{align*}\label{eqOUpvm}
&\mu(E)^\ast=\overline{\chi_E}(U)=\chi_E(U)=\mu(E),\\ 
&\mu(E)^2=\chi_E^2(U)=\chi_E(U)=\mu(E).
\end{align*}
\begin{theorem}\label{theoOUpvm}\index{projection valued measures!unitary operator}
    The map $E\mapsto \mu(E)$ is a {\bf projection valued measure} on $E$. That is,
    \begin{enumerate}
        \item[(i)] $\mu(\emptyset)=0$ and $\mu(\mathbb{T})=I$;
        \item[(ii)] if $\{E_n\}_{n=1}^\infty$ is a sequence of disjoint, Borel measurable sets in $\mathbb{T}$, then
    \[
    \mu\left(\bigcup_{n=1}^\infty E_n\right)=\sum_{n=1}^\infty\mu(E_n)
    \]
    converges in the {\bf strong operator topology}.\index{operator topology!strong}
    \end{enumerate} 
\end{theorem}
Having a measure, we expect to have an associated version of integral. In this case, we would like to integrate functions $f:\mathbb{T}\to \mathbb{C}$,
\begin{equation}\label{eqIntOUpvm}
\int_{\mathbb{T}}f(e^{it})d\mu(e^{it})\in\mathcal{L}(H).
\end{equation}
This can be done in a rigorous way. See e.g. \S 3.4 in \href{https://people.math.ethz.ch/~kowalski/spectral-theory.pdf}{ Spectral theory in Hilbert spaces (2009)},  ETH Zürich lecture notes, by Emmanuel Kowalski. We will  give the expression in \eqref{eqIntOUpvm} the weak meaning
\begin{equation}\label{eqIntOUpvmTwo}
\left\langle x,\left(\int_{\mathbb{T}}f(e^{it})d\mu(e^{it})\right)y\right\rangle:=\int_{\mathbb{T}}fd\mu_{x,y}.
\end{equation}
\begin{proof}
    Property (i) is clear. Let $\{E_n\}$ be like in (ii) and $E=\cup_{n=1}^\infty E_n$, and let $y\in H$. The vectors $\mu(E_n)y$, $n\ge1$, are orthogonal and 
    \[
    \sum_{n=1}^\infty\|\mu(E_n)y\|^2\le\|y\|^2,
    \]
    hence
    \[
    \sum_{n=1}^\infty \mu(E_n)y
    \]
    converges in $H$. Since for all $x\in H$
    \begin{eqnarray*}
        \langle x,\mu(E)y\rangle&=&\int_Ed\mu_{x,y}=\sum_{n=1}^\infty\int_{E_n}d\mu_{x,y}=\sum_{n=1}^\infty\langle x,\mu(E_n)y\rangle \\ 
        &=&\left\langle x,\sum_{n=1}^\infty\mu(E_n)y\right\rangle,
    \end{eqnarray*}
    we have that $\sum_{n=1}^\infty\mu(E_n)=\mu(E)$.
\end{proof}
The p.v.m. $\mu$ and the spectral measures $\mu_{x,y}$ are related by \eqref{eqMissing}. The advantage of $\mu$ is that we can use it to define null-sets and an $L^\infty$ norm. A Borel subset $E$ of $\mathbb{T}$ is a {\it null-set} for $U$ if $\mu(E)=0$. The countable union of null-sets is obviously a null-set. A property $\mathcal{P}(e^{it})$ ($e^{it}\in\mathbb{T}$) holds {\it $\mu$-almost everywhere} if the set where it fails is a null-set for $\mu$. A function $\varphi:\mathbb{T}\to\mathbb{C}$ belongs to $L^\infty(\mu)$ if there is $C\le 0$ such that $|\varphi(e^{it})|\le C$ $\mu-a.e.$. The minimum value of $C$ for which this holds is the $L^\infty(\mu)$-norm of $\varphi$, $\|\varphi\|_{L^\infty}(\mu)$.

Observe that other $L^p(\mu)$-norms can not be defined so easily. The advantage of the $L^\infty$-norm is that it is based on the primitive notion of null-set, not on the more sophisticated notion of the measure of a set ($\mu$ is a measure, but with values in the space of operators).

We now relate being null for $\mu$ with being null for $\mu_{x,x}$. In analogy with the case of finite positive Borel measures, we define the {\it support} of a p.v.m. $\mu$ on $\mathbb{T}$ to be the smallest closed subset $F$ of $\mathbb{T}$ such that $\mu(\mathbb{T}\setminus F)=0$.
\begin{theorem}\label{theoVariousNullities}[Nullsets for the p.v.m. associated with a unitary operator]\index{projection valued measures!unitary operator: nullsets}
    \begin{enumerate}
        \item[(i)] For a Borel subset $E$ of $\mathbb{T}$ the following are equivalent:
        \begin{itemize}
            \item[(a)] $\mu(E)=0$;
            \item[(b)] $\mu_{x,x}(E)=0$ for all $x\in H$;
            \item[(c)] $\mu_{x,y}(E)=0$ for all $x,y\in H$;
            \item[(d)] $\mu_{e_\alpha,e_\beta}(E)=0$ for all $e_\alpha,e_\beta$ in an orthonormal basis of $H$.
        \end{itemize}
        \item[(ii)] For a Borel function $\varphi:\mathbb{T}\to\mathbb{C}$,
        \[
        \|\varphi\|_{L^\infty(\mu)}=\sup_{x\in H}\|\varphi\|_{L^\infty(\mu_{x,x})}.
        \]
        \item[(iii)] If $f\in L^\infty_0(\mathbb{T})$, then $\|f(U)\|=\|f\|_{L^\infty(\mu)}$.
    \end{enumerate}
\end{theorem}
\begin{proof}
    (i) If $\mu(E)=0$, then $\mu_{x,y}(E)=\langle x,\mu(E)y\rangle=0$ for all $x,y$, which by polarization is the same as $\mu_{x,x}(E)=0$ for all $x$. If $\langle x,\mu(E)y\rangle=0$ for all $x,y$, on the other hand, $\mu(E)y=0$ for all $y$, then $\mu(E)=0$. This shows $(a)\iff(b)\iff(c)$. On the other hand, (d) states that all "matrix coefficients" $\langle e_\alpha,\mu(E)e_\beta\rangle=\mu_{e_\alpha,e_\beta}(E)$ of $\mu(E)$ vanish, hence that $\mu(E)=0$.

    (ii) For $\alpha>0$, using (i), we have
    \begin{eqnarray*}
        \|\varphi\|_{L^\infty(\mu)}\le\alpha&\iff&\mu\{t:|\varphi(e^{it})|>\alpha\}=0\\ 
        &\iff& \text{ for all }x\in H:\ \mu_{x,x}\{t:|\varphi(e^{it})|>\alpha\}=0\\ 
        &\iff&\sup_{x\in H}\|\varphi\|_{L^\infty(\mu)}\le\alpha,
    \end{eqnarray*}
    hence $\|\varphi\|_{L^\infty(\mu)}=\sup_{x\in H}\|\varphi\|_{L^\infty(\mu)}$.
    
    (iii) In one direction, for $x\in H$ we have
    \begin{eqnarray*}
        \|f(U)x\|^2&=&\langle x,|f|^2(U)x\rangle=\int_{\mathbb{T}}|f|^2d\mu_{x,x}\\ 
        &\le&\|f\|_{L^\infty(\mu_{x,x})}^2\|x\|^2\le \|f\|_{L^\infty(\mu)}^2\|x\|^2 \text{ by (ii)},
    \end{eqnarray*}
    hence, $\|f(U)\|\le \|f\|_{L^\infty(\mu)}$.

    In the other direction, let $0\le\lambda<\|f\|_{L^\infty(\mu)}$ (we can assume $\|f\|_{L^\infty(\mu)}>0$, since if $\|f\|_{L^\infty(\mu)}=0$ there is nothing to prove). By (ii), there is $x$ in $H$ such that $\lambda<\|f\|_{L^\infty(\mu_{x,x})}$. Hence, there exists a Borel set $E$ such that $\mu_{x,x}(E)>0$ and $|f(e^{it})|>\lambda$ for $e^{it}$ in $E$. Thus,
    \begin{eqnarray*}
        0&<&\int_E(|f|^2-\lambda^2)d\mu_{x,x}=\int_\mathbb{T}(|f|^2-\lambda^2)\chi_E^2d\mu_{x,x}\\ 
        &=&\int_\mathbb{T}(|f|^2-\lambda^2)d\mu_{\chi_E(U)x,\chi_E(U)x}\\ 
        &=&\langle \chi_E(U)x,|f|^2(U)\chi_E(U)x\rangle-\lambda^2  \langle \chi_E(U)x,\chi_E(U)x\rangle\\ 
        &=&\|f(U)[\chi_E(U)x]\|^2-\lambda^2 \|\chi_E(U)x\|^2.
    \end{eqnarray*}
    The strict inequality rules out the possibility that $\chi_E(U)x=0$, and the inequality itself implies that
    \[
    \lambda \|\chi_E(U)x\|<\|f(U)[\chi_E(U)x]\|,
    \]
    which implies that $\lambda< \|f(U)\|$. Hence, $\|f\|_{L^\infty(\mu)}\le \|f(U)\|$.
\end{proof}
\subsection{From the continuous calculus to the multiplicative form}\label{SSectOUmultform}
Let's go back to the continuous calculus. A vector $x\in H$ is {\it cyclic}\index{cyclic vector} for the unitary operator $U$ if the linear space
\[
\{f(U)x:\ f\in C(\mathbb{T})\}
\]
is dense in $H$. The space $V=\{f(U)x:\ f\in C(\mathbb{T})\}$ is {\it invariant for the action of $U$},\index{invariant subspace!unitary operator}
\[
UV\subseteq V.
\]
In fact, $V$ is invariant under the action of $f(U)$ if $f\in C(\mathbb{T})$, $[f(U)](V)\subseteq V$.
\begin{proposition}\label{propMultSpecUO}
    Suppose $U$ has a cyclic vector $x$. Then, the map $\Theta:C(\mathbb{T})\to H$, $\Theta(f)=f(U)x$, extends to a unitary map from $L^2(\mathbb{T},\mu_{x,x})$ to $H$, such that
    \begin{equation}\label{eqMultSpecUO}
        ((\Theta^{-1}U\Theta) f)(e^{it})=e^{it}f(e^{it}).
    \end{equation}
\end{proposition}
    That is, $U$ is unitarily equivalent to multiplication times $e^{it}$ on $L^2(\mathbb{T},\mu_{x,x})$,
    \[
    \begin{tikzcd}
L^2(\mathbb{T},\mu_{x,x}) \arrow{r}{\Theta} \arrow[swap]{d}{M_{e^{it}}} & H \arrow{d}{U} \\%
L^2(\mathbb{T},\mu_{x,x}) \arrow{r}{\Theta}& H
\end{tikzcd}
\]
\begin{proof}
    The map $\Theta$ is an $L^2(\mu_{x,x})$ isometry on $C(\mathbb{T})$,
    \[
    \|f(U)x\|^2=\langle x,f(U)^\ast f(U)x\rangle=\langle x,|f|^2(U)x\rangle=\int_{\mathbb{T}}|f|^2d\mu_{x,x}.
    \]
    As such, $\Theta$ extends to an isometry of $\overline{C(\mathbb{T})}^{L^2(\mu_{x,x})}$, the closure of $C(\mathbb{T})$ in the $L^2(\mu_{x,x})$ norm. The hypothesis that $x$ is cyclic means that $\Theta(C(\mathbb{T}))$ is dense in $H$, hence that $\Theta$ is a surjective isometry.
\end{proof}
We can then use some kind of finite or infinite induction, to have a similar result even when we do not have a cyclic vector at hands. We consider the case where $H$ is separable. The general case can be dealt with by Zorn's lemma.
\begin{theorem}\label{theoMultSpecUO}[Spectral theorem for unitary operators: multiplicative form]\index{theorem!spectral in multiplicative form for unitary operators}
    Let $U$ be a unitary operator on a separable Hilbert space $H$. Then, there exists a locally compact space $X$, a finite Borel measure $\nu$ on $X$, and a continuous function $\varphi:X\to\mathbb{R}$, and a unitary map $\Theta:L^2(\nu)\to H$, such that
    \begin{equation}\label{eqMultSpecUOter}
        (\Theta^{-1}U\Theta) f=e^{i\varphi}f,
    \end{equation}
    is multiplication times $e^{i\varphi}$.
\end{theorem}
\begin{proof}
    Let $\{x_n:\ n\ge1\}$ be a orthogonal basis for $H$, normalized to have $\sum_{n\ge1}\|x_n\|^2<\infty$. Let $P(\mathbb{T})$ be the space of the trigonometric polynomials. We construct a family of closed, mutually orthogonal subspaces $H_n$, whose sum exhausts $H$.
    \begin{itemize}
        \item Let $y_1=x_1$ and $H_1=\overline{\{p(U)y_1:\ p\in P(\mathbb{T})\}}^H$.
        \item Let $y_2=x_2-\pi_{H_1}x_2\in H_1^\perp$, and $H_2=\overline{\{p(U)y_2:\ p\in P(\mathbb{T})\}}^H$.
        \item If $H_1,\dots,H_{n-1}$ have been chosen, let $y_n=x_n-\pi_{H_1\oplus\dots\oplus H_{n-1}}x_n\in (H_1\oplus\dots\oplus H_{n-1})^\perp$, and $H_n=\overline{\{p(U)y_n:\ p\in P(\mathbb{T})\}}^H$.
    \end{itemize}
    Pick $p(U)y_n\in H_n$ and $q(U)y_m\in H_m$ with $m<n$. Then,
    \[
    \langle p(U)y_n,q(U)y_m\rangle=\langle y_n,(\overline{p}q)(U)y_m\rangle=0,
    \]
    because $\overline{p}q$ is a trigonometric polynomial, hence, $(\overline{p}q)(U)y_m\in H_m$, to which $y_n$ is orthogonal. Since the elements of the form $p(U)y_n$ are dense in $H_n$ and those of the form $q(U)y_m$ are dense in $H_m$, we have that $H_n\perp H_m$.

    After exhausting the list of the $x_n$'s, we have that $\bigoplus_{n\ge1}H_n$ is dense in $H$. Some of $H_n$'s might be trivial. In this case, we can remove them and renumber the sequence. The nontrivial $H_n$'s might be finitely many.

    For fixed $n$, we are in the situation described by proposition \ref{propMultSpecUO}. The map $\Theta_n:P(\mathbb{T})\to H_n$,
    \[
    \Theta_n(p)=p(U),
    \]
    extends to a unitary map from $\Theta_n:L^2(\mathbb{T},\mu_{y_n,y_n})\to H_n$, and
    \[
    [(\Theta_n^{-1}U\Theta_n)f](e^{it})=e^{it}f(e^it).
    \]
    Also, $\mu_{y_n,y_n}(\mathbb{T})=\|y_n\|^2\le\|x_n\|^2$.

    Let $A$ be the (discrete) set of the indices $n$ (which might be finite or infinite), and let $X=\mathbb{T}\times A$, with the product topology, which is locally compact. Let $\nu$ be the Borel measure which restricts to $\mu_{y_n,y_n}$ on the component $\mathbb{T}\times\{n\}$, and observe that $\nu$ is finite. We can naturally identify $L^2(X,\nu)=\bigoplus_{n\in A} L^2(\mathbb{T},\mu_{y_n,y_n})$. Finally, define $\Theta:L^2(X,\nu)\to H$ summand by summand,
    \[
    \Theta|_{L^2(\mathbb{T},\mu_{y_n,y_n})}=\Theta_n.
    \]
    If  $F:\mathbb{T}\times A\to\mathbb{C}$, then 
    \begin{equation}\label{eqSpMultUOb}
        [(\Theta^{-1}U\Theta)F](e^{it},n)=e^{it}F(e^{it},n).
    \end{equation}
    The operator $(\Theta^{-1}U\Theta)$, that is, acts like multiplication times $e^{it}$ on each "cl-open" component $\mathbb{T}\times\{n\}$ of $X$.
\end{proof}
The representation we have provided of $U$ as a multiplication operator is not canonical, even when $H=\mathbb{C}^2$. Consider for instance
\[
U=\begin{pmatrix}
    1&0\\ 
    0&-1
\end{pmatrix}.
\]
The vector $\begin{pmatrix}
    1\\-1
\end{pmatrix}$ is cyclic, and this provides a representation with just an "invariant subspace". The vectors $\begin{pmatrix}
    1\\0
\end{pmatrix}$, $\begin{pmatrix}
    0\\1
\end{pmatrix}$ are not, and they provide a representation with two "invariant subspaces".

The spaces $H_n$ in the proof of theorem \ref{theoMultSpecUO} are clearly invariant for $U$.

\subsection{The spectral theory of a multiplication operator}\label{SSectSTofMultO}
The multiplicative form of the spectral theorem tells us that any unitary operator is unitarily equivalent to a multiplication operator $M_\eta$ of the following form, for which we can explicitly write down the spectrum, the spectral measures, and so on.
\begin{itemize}
    \item[(a)] {\it We have a locally compact space $X$ and a finite positive Borel measure $\nu$, which we can assume to be finite, and a surjective local homeomorphism $\eta:X\to\mathbb{T}$, such that $M_\eta:L^2(\nu)\to L^2(\nu)$ is the operator $M_\eta f=\eta f$.}\index{spectrum!unitary operator}
    \item[(b)] {\it $\sigma(M_\eta)=\overline{\eta(\text{supp}(\nu))}\subseteq \mathbb{T}$ is the spectrum of $M_\eta$.}
    \begin{proof}
    Let $z\in\mathbb{T}$.
    \begin{eqnarray*}
        z\in\rho(M_\eta)&\iff&(z-M_\eta)^{-1}\in L^\infty(\nu)\\ 
        &\iff&|\eta(x)-z|\ge\epsilon\ \nu-a.e.\text{ for some }\epsilon>0\\ 
        &\iff&\eta^{-1}(D(z.\epsilon))\cap\text{supp}(\nu)=\emptyset \text{ for some }\epsilon>0\\ 
        &\iff&D(z,\epsilon)\cap \eta(\text{supp}(\nu))=\emptyset \text{ for some }\epsilon>0\\ 
        &\iff&z\notin \overline{\eta(\text{supp}(\nu))}.
    \end{eqnarray*}
    \end{proof}
    \item[(c)] {\it Let $f\in L^2(\nu)$. The spectral measure $\mu_{f,f}$ is such that, for any trigonometric polynomial $p:\mathbb{T}\to\mathbb{C}$,\index{spectral measures!unitary operator}
    \[
    \int_Xp(\eta)|f|^2d\nu=\langle f,p(M_\eta)f\rangle_{L^2(\nu)}=\int_{\mathbb{T}}pd\mu_{f,f}.
    \]
    That is, $d\mu_{f,f}=\eta_\ast(|f|^2d\nu)$.} (Recall that the {\it push-forward} of a measure $\nu$ on $X$ under a map $\eta:X\to Y$ is the measure defined by $\int_Y\varphi d(\eta_\ast\nu):=\int_X(\varphi\circ \eta)d\nu$. Of course, this needs some measurability assumptions, that is our case are easily verified). {\it Similarly, $d\mu_{f,g}=\eta_\ast(\overline{g}fd\nu)$.}
    \item[(d)] {\it The p.v.m. $\mu$ maps Borel measurable sets $E$ in $\mathbb{T}$ to projections of $L^2(\nu)$ in the following fashion:
    \[
    \int_X\overline{f}[\mu(E)g]d\nu=\langle f,\mu(E)g\rangle_{L^2(\nu)}=\mu_{f,g}(E)=\int_X\chi_E(\eta)\overline{f}gd\nu,
    \]
    i.e. ($\chi_E(\eta)=\chi_E\circ\eta$)
    \begin{equation}\label{eqMuForSPMF}
        \mu(E)g=\chi_E(\eta)g,\ \mu(E)=M_{\chi_E(\eta)}.
    \end{equation}}
    \item[(e)] {\it A Borel set $E\subseteq\mathbb{T}$ is a null-set for $\mu$ if and only if $0=\chi_E\circ\eta=\chi_{\eta^{-1}(E)}$ $\nu-a.e.$, i.e. if and only if $\nu(\eta^{-1}(E))=0$.}
    \item[(f)] {\it By (e) and the topological properties of $\eta$,  $\overline{\eta(\text{supp}(\nu))}=\text{supp}(\mu)$. }\index{projection valued measures!unitary operator}
    \begin{proof}
        \begin{eqnarray*}
            z\notin\text{supp}(\mu)&\iff&\mu(D(z,\epsilon))=0\text{ for some }\epsilon>0\\ 
            &\iff&\nu(\eta^{-1}(D(z,\epsilon)))=0\text{ for some }\epsilon>0,\text{ by }(e)\\ 
            &\iff&\eta^{-1}(D(z,\epsilon))\cap\text{supp}(\nu)=\emptyset \text{ for some }\epsilon>0\\ 
            &\iff&D(z,\epsilon)\cap\eta(\text{supp}(\nu))=\emptyset\text{ for some }\epsilon>0\\ 
            &\iff&z\notin\overline{\eta(\text{supp}(\nu))}.
        \end{eqnarray*}
    \end{proof}
    \item[(g)] {\it By (b) and (f),
    \begin{equation}\label{esSpectMultOprtr}
        \text{supp}(\mu)=\sigma(M_\eta).
    \end{equation}}
\end{itemize}

\subsection{Some applications of the spectral theorem in its multiplicative form}\label{SSectSTMFappied}
A notable advantage of the multiplicative form of the spectral theory is that we can now use measure theory to study functional analytic properties of unitary operators. We start by highlighting the role of the spectrum.
\begin{theorem}\label{theoSpectrumUO}[The spectrum and the support of the projection valued measure]\index{spectrum!unitary operator}\index{projection valued measures!unitary operator}
    Let $U:H\to H$ be unitary, $\{\mu_{x,x}:x\in H\}$ be the family of the spectral measures for $U$, and $\mu$ be the projection valued measure associated with $U$. We have
        \begin{equation}\label{eqSuppOfSMforUObis}
            \text{supp}(\mu)=\sigma(U).
        \end{equation}
\end{theorem}
\begin{proof}[Proof of theorem \ref{theoSpectrumUO}]
    By the multiplicative form of the spectral theorem, it suffices to show it for multiplication operators. This is done in \eqref{esSpectMultOprtr}.
\end{proof}
\begin{corollary}\label{corSpectrumUO}
With the same notation of theorem \ref{theoSpectrumUO}, we also have the following.
 \begin{enumerate}
        \item[(i)] For a function $f\in L^\infty(\mu)$, $f(U)=0$ if and only if $f$ vanishes $\mu$-a.e on $\sigma(U)$. Thus, the map $f\mapsto f\chi_{\sigma(U)}$ is the identity of $L^\infty(\mu)$. To highlight this fact, we also write $L^\infty(\mu)=L^\infty(\sigma(U))$, and for $\varphi:\sigma(U)\to\mathbb{C}$ measurable and bounded, we write $\varphi(U)=f(U)$, where $f$ is any measurable extension of $\varphi$ to $\mathbb{T}$. 
        \item[(ii)] If $\{f_n:n\ge1\}$ is a sequence in $L^\infty(\sigma(U))$, if $\|f_n\|_{L^\infty(\sigma)}\le C$ is uniformly bounded, and $f_n\to f$ $\mu-a.e.$, then $f_n(U)\to f(U)$ in the strong operator topology.\index{operator topology!strong}
        \item[(iii)] By Tietze extension theorem, the operator $f\mapsto f|_{\sigma(U)}=\Psi f$ maps $C(\mathbb{T})$ onto $C(\sigma(U))$. By (i), $f(U)=0$ if and only if $\Psi f=0$. For $g\in C(\sigma(U))$, then, we can write $g(U):=f(U)$. Now, if If $\{g_n:n\ge1\}$ is a sequence in $C(\sigma(U))$ converging to $g$ uniformly, then $g_n(U)$ converges to $g(U)$ uniformly.

        Moreover, $f\mapsto f(U)$ is an isometry of $C(\sigma(U))$ into $\mathcal{L}(H)$.
 \end{enumerate}   
\end{corollary}

Here is another sample application. We denote by $\mathcal{U}(H)$ the set of the unitary operators on $H$, which is also a group under composition.
\begin{proposition}\label{propUnitaryIsConnected}
    The unitary group $\mathcal{U}(H)$, endowed with the strong operator topology, is arc-connected.\index{operator topology!strong}
\end{proposition}
\begin{proof}
    We can assume that $H=L^2(\nu)$ and that $Uf=e^{i\varphi}f$, where $\varphi:X\to\mathbb{R}$ is measurable. For $t\in\mathbb{R}$, let
    \[
    U_tf=e^{it\varphi}f,
    \]
    $U=U_1$. Each $U_t$ is unitary, $U_0=I$, and it is easy to prove that $t\mapsto U_t$ is strongly continuous from $\mathbb{R}$ to $\mathcal{U}(L^2(\nu))$. (We have $U_{s+t}=U_sU_t$: $t\mapsto U_t$ is a {\it strongly continuous one parameter group in $\mathcal{U}(L^2(\nu))$}).

    Since any $U$ can be connected to $I$, two such $U$'s can be connected to each other.
\end{proof}
In general, there are many arcs of this group type connecting $U$ to the identity. We have in fact that
\[
e^{i\varphi}=e^{i(\varphi+2\pi n)},
\]
where $n:X\to\mathbb{Z}$ is any measurable function. The maps $t\mapsto e^{it(\varphi+2\pi n)}$ provide one parameter groups joining $U$ and $I$, and in general they are different from each other.
\subsection{Some final remarks}\label{SSectUOfinalremark}
The various forms we have seen of the spectral theorem for unitary operators, and their proofs, suggest a general framework, strategy, and objectives.
\begin{enumerate}
    \item[(i)] Given the operator $L$ on the Hilbert space $H$, we see if we can perform some simple algebraic operations on it. In this case, we can consider linear combinations of $U^n$ with complex coefficients, where $n$ can be both positive or negative. This provides a "restricted calculus", that in our case involves trigonometric polynomials $p(e^{it})$, which become the operators $p(L)$. 
    
    We will see when working with unbounded operators that finding a "restricted calculus" is not easy, since unbounded operators can not be in general composed without restricting more and more the domain. We will use instead the {\it Cauchy kernel} $u\mapsto\frac{1}{u-z}$ ($u\in\mathbb{R}$, $z\in\mathbb{C}\setminus\mathbb{R}$), and the operators $(L-z)^{-1}$ (the "resolvent map").
    \item[(ii)] If we have some information about convergence, we can extend the "restricted calculus" to a "continuous calculus". Here we had, for each $f\in C(\mathbb{T})$, a corresponding operator $f(L)$. The good news is that bounded, linear functionals on nice spaces of continuous functions are in a bijection with finite Borel measures by Riesz theorem. For instance, this is the case with the functionals
    \[
    f\mapsto\langle x,f(L)y\rangle.
    \]
    The "spectral" measures we have found are defined by the relation
    \[
    \langle x,f(L)y\rangle=\int_{\mathbb{T}}fd\mu_{x,y}.
    \]
    With unbounded operators things will be trickier. Lacking a polynomial, even less a continuous calculus, we have to resort to different results from real and complex analysis. Following Tao and one stream of literature, we will use {\it Herglotz representation theorem}, associating to each positive harmonic function defined on a half-plane (and having suitable decay) a positive measure on the half-plane's boundary. This will produce "spectral measures", satisfying
    \[
    \langle x,(L-z)^{-1}y\rangle=\int_{\mathbb{R}}\frac{d\mu_{x,y}(u)}{u-z}.
    \]
    In both cases,
    \begin{equation}\label{eqspecmultintuition}
        \langle x,Ly\rangle=\int_{\mathbb{X}}ud\mu_{x,y}(u):
    \end{equation}
    the operator $L$ becomes, in the world of each spectral measure, multiplication times the coordinate of the space.
    \item[(iii)] One way to exploit the state of things described in (ii) is by extending the integrals to Borel measurable functions,
    \[
    \langle x,m(L)y\rangle:=\int_{X}m(u)d\mu_{x,y}(u),
    \]
    for $m$ Borel measurable and bounded. Here $X$ is the "spectral set": the natural subset of the complex plane containing the spectra of the operators $L$ belonging to the class we are studying. The most interesting instance is that of the characteristic functions, $m=\chi_E$, where $E$ is measurable in $X$. Since $\chi_E^2=\chi_E$ and $\chi_E$ is real valued, if we have an effective calculus, then $\chi_E(L)^2=\chi_E(L)$ and $\chi_E(L)^\ast=\chi_E(L)$; i.e. $\mu(E):=\chi_E(L)$ is a projection onto a closed subspace of $H$! We have, that is, {\it projection valued measures}. By the way, these are of the uttmost importance in quantum mechanics. 
    \item[(iv)] A different strategy consists in starting again with (ii). We have realized that for elements of $H$ having the form $f(L)x$ (with $x\in H$ fixed), the map $\Theta:f\mapsto f(L)x$ is an isometry from $L^2(\mu_{x,x})$ onto a "cyclic" subspace of $H$, and that $\Theta$ intertwines $L$ with a multiplication operator on $L^2(\mu_{x,x})$. This is the main building block of the multiplicative form of the spectral theorem. 
\end{enumerate}
\section{Closed operators and their adjoints}\label{SectOpClosed}
Let $D$ be a subspace of $H$ and $D\xrightarrow{L} H$ be a linear operator. We say that $L$ is\index{densely defined!operator} {\it densely defined} (d.d.) id $D$ is dense in $H$. The {\it graph} $\Gamma(L)$ of $L$ is
\[
\Gamma(L)=\{(f,Lf):f\in D\}\subseteq H\times H.
\]
The Cartesian product $H\times H=H\oplus H$ is itself a Hilbert space with the inner product
\[
\langle(f,g),(h,k)\rangle_{H\oplus H}=\langle f,h\rangle_{H}+\langle g,k\rangle_{H}.
\]
The operator $L$ is {\it closed}\index{operator!closed} if $\Gamma(L)$ is closed in $H\times H$, i.e. if
\[
D\ni f_n\to f\text{ and }Lf_n\to g\text{ in }H, \text{ implies that }f\in D\text{ and }g=Lf.
\]
\begin{remark}\label{remInvBndd}
    Let $D\xrightarrow{L}H$ be closed, injective, and having closed range (this holds if, for instance, $L$ is a bijection). By the closed graph theorem, $\text{Ran}(L)\xrightarrow{L^{-1}}D$ is bounded.
\end{remark}
\begin{example}\label{exaClosedOp} The operator $(Lf)(x)=xf(x)$, having domain $D=\{f\in L^2(\mathbb{R}):\ xf\in L^2(\mathbb{R})\}$ is closed.
\end{example}
An {\it extension} $D_1\xrightarrow{L_1}H$ of an operator $D\xrightarrow{L}H$ is a linear operator $L_1$ with domain $D_1\supseteq D$ such that $D_1|_D=L$. 

The operator $D\xrightarrow{L}H$ is {\it symmetric}\index{operator!symmetric} if for all $f,g\in D$
\[
\langle f,Lg\rangle=\langle Lf,g\rangle.
\]
\begin{lemma}\label{lemmaClosureOp}
    If $D\xrightarrow{L}H$ is d.d. and symmetric, then
    \[
    \overline{\Gamma(L)}=\Gamma(\overline{L})
    \]
    is the graph of a closed, symmetric linear operator $L$, the {\bf closure} of $L$. The domain of $\overline{L}$ is denoted by $\overline{D}$. Beware! $\overline{D}$ is {\it not} the closure of $D$ in $H$ (which is $H$ itself), but rather $\pi_1(\overline{\Gamma(L)})$. That is,
    \[
    \overline{D}=\{f\in H:\exists D\ni f_n\to f\text{ such that }\exists\lim_{n\to\infty}f_n\in H\}.
    \]
\end{lemma}
If $L_1$ is a closed operator extending $L$, then $\Gamma(L_1)\supseteq \Gamma(\overline{L})\supseteq\Gamma(L)$. In this sense $\overline{L}$ is the minimal closed extension of $L$.
\begin{proof}
    Suppose $\Gamma(L)\ni(f_n,Lf_n)\to(f,g)$, and set $\overline{L}f=g$. The definition is well posed if it does not depend on the sequence $\{f_n\}$ in $D$, but just on $f$, and this happens if and only if $(h_n,Lh_n)\to(0,p)$ with $h_n\in D$ implies that $p=0$. For $k\in D$ we have that
    \[
    \langle k,p\rangle=\lim_n\langle k,Lh_n\rangle=\lim_n\langle Lk,Lh_h\rangle=0,
    \]
    hence $p\in D^\perp=H^\perp=0$, since $D$ is dense in $H$. Let $\pi_1(f,g)=f$ be the projection onto the first coordinate in $H\times H$. We have defined $\overline{D}\xrightarrow{\overline{L}}H$ on 
    \[
    \overline{D}:=\{f\in H:\exists D\ni f_n\to f \text{ such that }\lim_nLf_n \text{ exists in }H\}=\pi_1(\overline{\Gamma(L)}),
    \]
    and $\Gamma(\overline{L})=\overline{\Gamma(L)}$.

    The operator $\overline{L}$ is linear. If $D\ni f_n\to f$, $Lf_n\to g=\overline{L}f$, $D\ni h_n\to h$, and $Lh_n\to k=\overline{h}$, then $L(af_n+bh_n)\to ag+bk=\overline{af+bh}$, where $a,b$ are scalars.    
    Also, $\overline{L}$ is symmetric:
    \[
    \langle \overline{L}f,h\rangle=\lim_n\langle Lf_n,h_n\rangle=\lim_n\langle f_n,Lh_n\rangle=\langle f,\overline{L}h\rangle.
    \]
\end{proof}

The {\it adjoint}\index{operator!adjoint of a densely defined operator} of a d.d. linear operator $D\xrightarrow{L}H$ is defined as in the bounded operator case, although we have to be careful with the domains. Let
\[
D^\ast=\{g\in H:f\mapsto\langle g,Lf\rangle \text{ is bounded on }D\}.
\]
For each $g\in D^\ast$, by density of $D$ in $H$ and Riesz theorem, there exists a unique $L^\ast g\in H$ such that
\[
\langle L^\ast g,L\rangle=\langle g,Lf\rangle.
\]
It is an easy exercise showing that $D^\ast$ is a subspace of $H$ and $L^\ast$ is a linear operator.
\begin{remark}\label{remAdjSymm}
    If $L$ is d.d. and symmetric, then $D^\ast\supseteq D$ and $L^\ast$ extends $L$. In fact,
    for $g\in D$, $f\mapsto\langle Lf,g\rangle=\langle f,Lg\rangle$ is bounded, and $L^\ast g=Lg$.
\end{remark}
\begin{lemma}\label{lemmaLastClosed} If $D\xrightarrow{L}H$ is d.d., then $L^\ast$ is closed.
\end{lemma}
\begin{proof}
    Let $\Tilde{\Gamma}(T)=\{(Tf,f):f\in E\}$ be the co-graph of an operator $E\xrightarrow{T}H$. We prove that:
    \begin{equation}\label{eqLastClosed}
    \Tilde{\Gamma}(-L^\ast):=\{(-L^\ast g,g):g\in D^\ast\}=\Gamma(L)^\perp \text{ in }H\oplus H.
    \end{equation}
    On the one hand, for $f\in D$, we have
    \[
    \langle(-L^\ast g,g),(f,Lf)\rangle=-\langle L^\ast g,f\rangle+\langle g,Lf\rangle=0,
    \]
    which shows that $\Tilde{\Gamma}(-L^\ast)\subseteq \Gamma(L)^\perp$. Conversely, $(h,g)\in \Gamma(L)^\perp$ if and only if $0=\langle h,f\rangle+\langle g,Lf\rangle$ for all $f$ in $D$, but this implies that $f\mapsto \langle g,Lf\rangle$ is bounded on $D$ (hence, $g\in D^\ast$), and that $h=-L^\ast g$. Hence, $(h,g)\in \Tilde{\Gamma}(-L^\ast)$. By \eqref{eqLastClosed}, $\Tilde{\Gamma}(-L^\ast)$ is closed, which implies that $\Gamma(L^\ast)$ is closed, as wished.
\end{proof}
\begin{lemma}\label{lemmaAstAstSymm}
    Let $L$ be d.d. and symmetric. Then, $L^{\ast\ast}=\overline{L}$ and $L^\ast=\overline{L}^\ast$.
\end{lemma}
\begin{proof} By lemma \ref{lemmaClosureOp} we have
\begin{eqnarray*}
    \Tilde{\Gamma}(-L^\ast)&=&\Gamma(L)^\perp=\overline{\Gamma(L)}^\perp=\Gamma(\overline{L})^\perp\\ 
    &=&\Tilde{\Gamma}(-\overline{L}^\ast),
\end{eqnarray*}
    hence, $L^\ast=\overline{L}^\ast$, which is our second thesis. As a consequence, using lemma \ref{lemmaClosureOp} twice,
    \[
    \Gamma(L^{\ast\ast})=\Tilde{\Gamma}(-L^\ast)^\perp=\Tilde{\Gamma}(-\overline{L}^\ast)^\perp=\Gamma(\overline{L}),
    \]
    hence, $L^{\ast\ast}=\overline{L}$.
\end{proof}
\begin{example}\label{exaClosedOpProj}
    Unbounded operators provide examples of closed subspaces $M,N\subset \mathcal{H}$ of a Hilbert space $\mathcal{H}$ such that the orthogonal projection $\pi_M$ onto $M$ maps $N$ onto a dense subspace $\pi_M(N)\ne M$. This fact is not evident from pictures. A general machine to produce counterexamples is by considering $\mathcal{H}=H\times H$, where $H$ is infinite dimensional, $M=H$ (say, the first component), and $N=\Gamma(L)$, where $D\xrightarrow{L}H$ is a d.d., closed, unbounded operator (so that $D$, by the Hellinger-Toeplitz theorem, is strictly contained in $H$). Then, $\pi_1(\Gamma(L))=D$ is not the whole range of $\pi_1$, although it is dense in it.
\end{example}
A consequence of lemma \ref{lemmaLastClosed} is that, if $L$ is d.d. and symmetric, then $L^\ast$ is a closed extension of $L$, hence of $\overline{L}$. We might hope that, under these assumptions, $L^\ast=\overline{L}$. Unfortunately this is not the case.
\begin{remark}\label{remLsymmLastNot}
    There are d.d., symmetric operators such that $L^\ast$ is not symmetric. Since for all such operators $L^\ast=\overline{L}^\ast$ by lemma \ref{lemmaAstAstSymm}, such examples arise even for some d.d., symmetric, closed operators. See the section on the momentum operator for a concrete instance.
\end{remark}
A d.d. defined operator $D\xrightarrow{L}H$ is {\it self-adjoint}\index{self-adjoint!densely defined operator} (s.a.) if and only if $L=L^\ast$: $D=D^\ast$ and $L=L^\ast$. In particular, self-adjoint operators are closed (because $L^\ast$ is closed) and symmetric: if $f\mapsto\langle g, Lf\rangle$ is bounded on $D$, then $g\in D^\ast=D$, and $\langle g, Lf\rangle=\langle L^\ast g,f\rangle=\langle L g,f\rangle$. 

Self-adjointness is much stronger than it looks.
\begin{remark}\label{remSAmaximalDomain}
    Let $Dx\xrightarrow{L}H$ be self-adjoint, and let $D\subseteq E\xrightarrow{M}H$ be a symmetric extension of $L$. Then, $E=D$. That is, $L$ can not be extended to a larger domain while preserving symmetry.

    The proof is easy. If $y\in E$, then
    \[
    x\mapsto\langle y,Lx\rangle=\langle y,Mx\rangle=\langle My,x\rangle
    \]
    defines a bounded functional on $D$, hence, $y\in D=D^\ast$.
\end{remark}
Self-adjoint operators, that is, are maximal in the lattice of the symmetric operators, ordered by extension. In particular, the domain of a self-adjoint operator must be "best possible".

A d.d. symmetric operator $D\xrightarrow{L}H$ is just {\it essentially self-adjoint}\index{essentially self-adjoint operator} if its closure $\overline{L}$ is self-adjoint.
\begin{remark}\label{remSAclosed} A d.d. symmetric $D\xrightarrow{L}H$ is essentially self-adjoint if and only if $\overline{D}=D^\ast$. 
\end{remark}
The next proposition gives a first, abstract criterion for self-adjointness.
\begin{proposition}\label{propSArevealed}
    Let $D\xrightarrow{L}H$ be d.d., closed, symmetric. Then, $L$ is self-adjoint if and only if $L^\ast$ is symmetric.
\end{proposition}
\begin{proof}
    The only if direction is obvious. In the other direction, suppose $L^\ast$ is symmetric. We claim that, then, $D\subseteq D^\ast$ and $\Gamma(L)\subseteq\Gamma(L^\ast)$.
    If $g\in D$, in fact, then $|\langle g,Lf\rangle|=|\langle Lg,f\rangle|\le\|Lg\|\cdot\|f\|$, hence $g\in D^\ast$. For such $g$'s, $\langle Lg,f\rangle=\langle g,Lf\rangle=\langle L^\ast g,f\rangle$ for all $f$ in $D$, hence $L^\ast g=Lg$. This shows that $D\subseteq D^\ast$, hence that $\Gamma(L)\subseteq \Gamma(L^\ast)$. We want to show the opposite inclusion.

    Both $\Gamma(L)$ and $\Gamma(L^\ast)$ are closed. Let $(g,L^\ast g)\in \Gamma(L^\ast)\ominus\Gamma(L)$. For $f\in D$, 
    \[
    0=\langle(g,L^\ast g),(f,Lf)\rangle=\langle g,f\rangle+\langle L^\ast g,Lf\rangle,
    \]
    then, $f\mapsto \langle L^\ast g,Lf\rangle=-\langle g,f \rangle$ is bounded, hence $L^\ast g\in D^\ast$ and $L^{\ast\ast}g=-g$. Then,
    \[
    \|L^\ast g\|^2=\langle L^{\ast\ast}g,g\rangle=-\|g\|^2,
    \]
    which implies that $g=0$. Thus, $\Gamma(L)=\Gamma(L^\ast)$, and $L$ is self-adjoint.
\end{proof}
\section{Resolvent and spectrum}\label{sectRandS}

\subsection{Resolvent and spectrum for a general d.d. symmetric operator}\label{SSectResSpecGeneral}
\index{spectrum!symmetric, densely defined operator}
Consider $D\xrightarrow{L}H$ symmetric, d.d., and closed, and $z\in \mathbb{C}$. We say that $z\in \sigma(L)$ belongs to the {\it spectrum} of $L$ if $L-zI=L-z:D\to H$ is not a bijection. This can happen because $L-zI$ is not injective, or because it is not surjective. If it is injective, we can still define $R(z)=(L-z)^{-1}:\text{Ran}(L-z)\to D$. 
\begin{lemma}\label{lemmaInvBounded}
  Let $D\xrightarrow{L}H$ be closed, and $z\in\mathbb{C}$.
  If the inverse $(L-z)^{-1}:H\to D$ exists, then it is bounded.  
\end{lemma}
\begin{proof}
    It sufficed to prove that $\Gamma((L-z)^{-1})=\Tilde{\Gamma}(L-z)$ is closed, and apply the closed graph theorem. Suppose $f_n\to f$, and $(L-z)f_n\to g$ in $H$, so that $Lf_n\to zf+g$. Since $L$ is closed, $Lf=zf+g$, i.e. $g=(L-z)f$, showing that  $\Tilde{\Gamma}(L-z)$ is closed.
\end{proof}
As in the case of bounded operators, or of Banach algebras, we deduce that the {\it resolvent set}\index{resolvent!set} $\rho(L)=\mathbb{C}\setminus\sigma(L)$ of $L$ is open.
\begin{exercise}\label{exeResIsOpenForAllL}
    Let $D\xrightarrow{L}H$ be closed. Use a Neumann series to show that $\rho(L)$ is open. Deduce that the map $z\mapsto R(z)$, $R:\rho(L)\to\mathcal{B}(H)$, is holomorphic.
\end{exercise}
\begin{exercise}\label{exeSpctrumPositionOperator}
    Let $M_x$ be the {\it position operator}\index{position operator} $M_xf(x)=xf(x)$ on $L^2(\mathbb{R})$, with domain $D=\{f\in L^2:\ \int_{\mathbb{R}}(1+x^2)|f(x)|^2dx<\infty\}$. Show that $M_x$ is self-adjoint and that its spectrum is $\sigma(M_x)=\mathbb{R}$ (which is not compact).
\end{exercise}
The map $z\mapsto R(z)=(L-z)^{-1}$ is the {\it resolvent function}\index{resolvent!function} of $L$, $R:\rho(L)\to\mathcal{B}(H)$.

In order to define $R(z)=(L-z)^{-1}:\text{Ran}(L-z)\to D$ we just need $L-z$ to be injective, and not necessarily surjective. In the case of symmetric operators we have some quantitative information.
\begin{lemma}\label{lemmaResolventGood}
    Let $D\xrightarrow{L}H$ be d.d., symmetric and closed, and let $z\in\mathbb{C}\setminus\mathbb{R}$. Then, $\text{Ran}(L-z)$ is closed, $R(z):\text{Ran}(L-z)\to D$, and
    \begin{equation}\label{eqNormResolvent}
        \|R(z)\|\le\frac{1}{|\text{Im}(z)|}.
    \end{equation}
\end{lemma}
\begin{proof}
    By symmetry of $L$, for $f\in D$ we have $\langle f,Lf\rangle\in\mathbb{R}$, hence,
    \begin{eqnarray*}
        |-\text{Im}z|\cdot\|f\|^2&=&|\text{Im}\langle f,(L-z)f\rangle|\\ 
        &\le&\|(L-z)f\|\cdot\|f\|,
    \end{eqnarray*}
    showing that 
    \begin{equation}\label{eqCoercive}
    \|(L-z)f\|\ge |\text{Im}z|\cdot\|f\|.
    \end{equation}
    This implies that $L-z$ is injective.

    If $D\ni f_n$ and $(L-z)f_n\to g$, so that $\{(L-z)f_n\}$ is Cauchy, then  $\{f_n\}$ is Cauchy as well by \eqref{eqCoercive}. Thus $Lf_n\to g+zf=Lf$, because $L$ is closed, so $g=(L-z)f$. This shows that $\text{Ran}(L-z)$ is closed.  

    Finally, equation \eqref{eqCoercive} is equivalent to \eqref{eqNormResolvent}.
\end{proof}
\begin{proposition}\label{propResolventUnbndd}
    Let $D\xrightarrow{L}H$ be d.d., symmetric, and closed.
    \begin{enumerate}
        \item[(i)] If $z,w\in\mathbb{C}$ and $R(z),R(w)$ are everywhere defined (i.e. $z,w\in\rho(L)$), then
        \begin{equation}\label{eqResoEqUnbndd}
            R(z)-R(w)=(z-w)R(z)R(w).
        \end{equation}
        In particular $R(z)R(w)=R(w)R(z)$ and $R(z)R(w)$
        \item[(ii)] If $z\in\mathbb{C}\setminus\mathbb{R}$, and $R(z),R(\overline{z})\in\rho(L)$, then $R(z)^\ast=R(z)$.
    \end{enumerate}
\end{proposition}
\begin{proof}
   (i) Formally,
   \begin{eqnarray*}
       R(z)-R(w)&=&(L-z)^{-1}-(L-w)^{-1}\\ 
       &=&(L-z)^{-1}[(L-w)-(L-z)](L-w)^{-1}\\ 
       &=&(L-z)^{-1}(z-w)(L-w)^{-1}\\
       &=&(z-w)R(z)R(w).
   \end{eqnarray*}
   The expression after the first equality is well defined, since 
   \[
   H\xrightarrow{(L-w)^{-1}}D\xrightarrow{L}H\xrightarrow{(L-z)^{-1}}D.
   \]
   (ii) We have $D\ni u=R(z)x,v=R(\overline{z})y$ if and only if $x=(L-z)u$, $y=(L-\overline{z})v$, hence
   \begin{eqnarray*}
       \langle x,R(\overline{z})y\rangle&=&\langle (L-z)u,v\rangle=\langle u,(L-\overline{z})v\rangle\\ 
       &=&\langle R(z)x,y\rangle.
   \end{eqnarray*}
   The equalities hold for all $x,y$ in $H$ because $R(z)$ and $R(\overline{z})$ are defined on $H$. Hence $R(z)^\ast=R(\overline{z})$.
\end{proof}
The operator $D\xrightarrow{L}H$ is {\it positive}\index{operator!positive} if $\langle Lf,f\rangle\ge0$ for all $f\in H$. 
\begin{remark}\label{remPosSymm}
    If $H$ is, as in our case, a Hilbert space on $\mathbb{C}$, then positive operators are symmetric. In fact,
    \[
    \langle f,Lf\rangle=\overline{\langle f,Lf\rangle}=\langle Lf,f\rangle,
    \]
    and by polarization it follows that $\langle f,Lg\rangle=\langle f,Lg\rangle$. It suffices to show that the polarization identity holds even if an operator is sandwiched between the vectors,
    \[
    \langle x,Ly\rangle=\frac{1}{4}[\langle x+y,L(x+y)\rangle-\langle x-y,L(x-y)\rangle-i\langle x+iy,L(x+iy)\rangle+i\langle x-iy,L(x-iy)\rangle].
    \]
    The same is not true for Hilbert spaces over $\mathbb{R}$. For instance, the matrix
    $A=\begin{pmatrix}
       1&1\\
       0&1
    \end{pmatrix}$ defines a positive operator from $\mathbb{R}^2$ to itself, which is not symmetric.
\end{remark}
\begin{lemma}\label{lemmaResolventPosGood}
Let $D\xrightarrow{L}H$ be a d.d. positive, closed operator, and let $z\in\mathbb{C}$ with $\text{Re}z<0$. Then, $R(z)$ is defined on $\text{Ran}(L-z)$, which is closed, and
\[
\|R(z)\|\le\frac{1}{|\text{Re}z|}.
\]
\end{lemma}
The proof is similar to that of lemma \ref{lemmaResolventGood}.
\begin{proof} For $f\in D$,
   \begin{equation}\label{eqFPMF}
       -\text{Re}(z)\|f\|^2\le\langle(L-z)f,f \rangle\le\|(L-z)f\|\cdot\|f\|,\text{ hence, }|\text{Re}z|\cdot\|f\|\le\|(L-z)f\|.
   \end{equation} 
   As in the previous lemma this implies that $L-z$ is injective, that $\text{Ran}(L-z)$ is closed, and the estimate on $\|R(z)\|$.
\end{proof}
We have a first, useful criterion for self-adjointness.
\begin{proposition}\label{propSAcriterionOne}
    Let $D\xrightarrow{L}H$ be d.d., symmetric, and closed, let $D^\ast\xrightarrow{L^\ast}H$ be its adjoint, and let $z\in\mathbb{C}\setminus\mathbb{R}$.
    \begin{enumerate}
        \item[(i)] $z\in\rho(L)$ if and only if $L^\ast-\overline{z}$ is injective;
        \item[(ii)] if $L=L^\ast$ is self-adjoint, then $z\in\rho(L)$;
        \item[(iii)] if $z,\overline{z}\in \rho(L)$, then $L$ is self-adjoint.
    \end{enumerate}
\end{proposition}
\begin{corollary}\label{corSAcriterionOne}
    Let $D\xrightarrow{L}H$ be d.d., symmetric, and closed. Then, $L$ is self-adjoint if and only if $\pm i\in\rho(L)$, i.e. if and only if $R(i),R(-i)$ are defined on $H$.
\end{corollary}
    \begin{proof}
        (i) $R(z)$ is not defined on $H$ if and only if $\text{Ran}(L-z)\ne H$, i.e. if and only if there exists $0\ne v\perp \text{Ran}(L-z)$, since we proved in lemma \ref{lemmaResolventPosGood} that $\text{Ran}(L-z)$ is closed. For all $x\in D$, that is, $\langle v,(L-z)x\rangle=0$. Thus,
        \[
        |\langle v,Lx\rangle|=|z|\cdot|\langle v,x\rangle|\le |z|\cdot\|v\|\cdot\|x\|,
        \]
        which shows that $v\in D^\ast$, the domain of $L^\ast$. Moreover,
        \[
        \langle (L^\ast-\overline{z})v,x\rangle=\langle v,(L-x)x\rangle=0
        \]
        for all $x\in D$, hence $(L^\ast-\overline{z})v=0$ because $D$ is dense, showing that $L^\ast-\overline{z}$ is not injective. 
        
        Viceversa, if $0\ne v\in \ker(L^\ast-z)$ for some $v\in D^\ast$, then for $x\in D$
        \[
        0=\langle (L^\ast-\overline{z})v,x\rangle=\langle v,(L-x)x\rangle,
        \]
        i.e. $v\in D^\perp$, contradicting the fact that $L$ is d.d.

        (ii) By (i), if $z\notin \rho(L)$ then there is $0\ne v\in D^\ast=D$ such that $(L-\overline{z})v=(L^\ast-\overline{z})v=0$, which contradicts the inequality $\|(L-\overline{z})v\|\ge|\text{Im}z|\cdot\|v\|\ne0$ in lemma \ref{lemmaResolventGood}.

        (iii) Let $g\in D^\ast$. For all $f\in H$, then, using proposition \ref{propResolventUnbndd} in the first equality,
        \begin{eqnarray*}
            \langle f,R(\overline{z}(L^\ast-\overline{z})g)\rangle&=&\langle R(z)f,(L^\ast-\overline{z})g\rangle=\langle (L-z)R(z)f,g\rangle\\ 
            &=&\langle f,g\rangle.
        \end{eqnarray*}
        The second equality holds because $R(z)f\in D(L-z)$, the domain of $L-z$. Thus $R(\overline{z})(L^\ast-\overline{z})g=g$, hence $g\in \text{Ran}(R(\overline{z}))=D$. This shows that $D^\ast\subseteq D$, and the opposite inclusion holds for all closed, d.d., symmetric operators.
    \end{proof}
    The basic criterion of self-adjointness in corollary \ref{corSAcriterionOne} is less innocent that it seems. It requires, in fact, to show that for all $a$ in $H$ the equations
    \[
    (L\pm i)x=a
    \]
    have a solution $x\in D$. If $L$ is a differential operator, we need to solve differential equations, and the domain $D$ of the operator has to be chosen in such a way that the solution $x$ exists for all data $a$. Since $L\pm i$ must be invertible, the domain for which $L$ is self-adjoint, if it exists, it is uniquely determined. In concrete situations, this and other criteria require tools from "hard analysis".
    \begin{proposition}\label{propDomainOfR}
        Let $D\xrightarrow{L}H$ be d.d., closed and symmetric, and $z\in\mathbb{C}\setminus\mathbb{R}$.
        \begin{enumerate}
            \item[(i)] If $z\in\rho(L)$, then $w\in\rho(L)$ whenever $|z-w|<|\text{Im}z|$.
            \item[(ii)] If $z\in\rho(L)$ and $w\in\mathbb{C}\setminus\mathbb{R}$ is on the same side of $z$ with respect to the real line, then $w\in\rho(L)$.
        \end{enumerate}
    \end{proposition}
\begin{proof}
    Item (i) can be easily proved by a Neumann series argument, and item (ii) is a consequence of (i) and of the fact that half-planes are path-connected.
\end{proof}
For positive operators we have an analogous proposition.
\begin{proposition}\label{propDomainOfRpositive}
    Let $D\xrightarrow{L}H$ be d.d., closed, symmetric, and positive, and let $w\in\mathbb{C}\setminus[0,+\infty)$. Then, $w\in\rho(L)$ if and only if $L$ is self-adjoint.
\end{proposition}
\begin{proof} [If] By corollary \ref{corSAcriterionOne}, $L$ is self-adjoint if and only if $\pm i\in\rho(L)$. Draw a curve from $i$ to $w$ in the complex plane, which does not meet $[0,\infty)$. We can use Neumann series on overlapping discs centered on the curve to show that $w\in \rho(L)$, using the fact that $\|R(z)\|\le \frac{1}{|\text{Im}z|}$ if $\notin\mathbb{R}$, and that $\|R(z)\|\le\frac{1}{|\text{Re}z|}$ if $\text{Re}z<0$.

In the opposite direction, we use the same reasoning on curves starting at $w$ and ending at $\pm i$, respectively, showing that $\pm i\in\rho(L)$, hence that $L$ is self-adjoint. 
\end{proof}
Sometimes it is easier working with essentially self-adjoint operators, which do not require optimal information on the domain.
\begin{proposition}\label{propEssSA}
    Let $D\xrightarrow{L}H$ be d.d., symmetric. 
    \begin{enumerate}
        \item[(i)] $L$ is essentially self-adjoint if and only if $\text{Ran}(L\pm i)$ are dense in $H$.\index{essentially self-adjoint operator}
        \item[(ii)] If $L\ge0$, then $L$ is essentially self-adjoint if and only if $\text{Ran}(L+1)$ is dense in $H$.
        \item[(iii)] $L$ is essentially self-adjoint if and only if $\text{ker}(L^\ast\pm i)=0$.
    \end{enumerate}
\end{proposition}
\begin{proof}
    (i) $L$ is essentially self-adjoint $\iff$ $\overline{L}$ is self-adjoint $\iff$ $\overline{L}\pm i$ are onto $\iff$ for all $y$ in $H$ there is $\overline{x}$ in $\overline{D}$ such that $(\overline{L}\pm i)\overline{x}=y$ $\iff$ for all $y\in H$ there is a sequence $D\ni x_n\to\overline{x}$ such that $Lx_n\to z=L\overline{x}$, so that $y=z\pm i\overline{x}$. The latter holds if and only if for all $y$ in $H$ we have $y=\lim(Lx_n\pm i x_n)$, i.e. if and only if $y\in\overline{\text{Ran}(L\pm i)}$. This means that $\text{Ran}(L\pm i)$ are dense in $H$.

    (ii) The proof is similar, using the criterion of self-adjointness for positive operators.

    (iii) $\text{ker}(L^\ast\pm i)\ni h$ if and only if $f\mapsto\langle h,Lf\rangle$ is bounded on $D$ and
    \[
    \langle h,Lf\rangle=\langle L^\ast h,f\rangle=\mp i\langle h,f\rangle=\langle h,\pm if\rangle,
    \]
    i.e. $\langle h,(L\mp i)f\rangle=0$ for all $f\in D$: $\text{Ran}(L\pm i)$ is dense, then, if and only if $\text{ker}(L^\ast\pm i)=0$. 
\end{proof}
\subsection{The bounded operator case}\label{SSectResSpecBounded}
If $H\xrightarrow{L}H$ is bounded and self-adjoint (the Hellinger-Toeplitz theorem says that if $L$ is symmetric, then it is bounded), the results obtained so far are easily summarized. For $z\in\rho(L)$, let $R(z)=(L-z)^{-1}$.
\begin{enumerate}
    \item[(i)] $\sigma(L)\subset\mathbb{R}$ is compact (we already proved compactness in the context of Banach algebras; the proof consisted in showing that the Neumann series for $R(z)=(L-z)^{-1}=z^{-1}(L/z-I)^{-1}$ converges for large $z$). If $L\ge0$, then $\sigma(L)\subset[0,\infty)$. The map $z\mapsto R(z)$ is holomorphic from $\rho(L)$ to $\mathcal{H}$.
    \item[(ii)] We have the estimates:
    \begin{align*}
        &\|R(z)\|\le\frac{1}{|\text{Im}z|}\text{if $L$ is self-adjoint},\\ 
        &\|R(x)\|\le\frac{1}{|x|}\text{ if $x<0$ and $L\ge0$}.
    \end{align*}
    \item[(iii)] We still have the identities for the resolvent:
    \begin{align*}
        &R(z)-R(w)=(z-w)R(z)R(w),\\
        &R(z)^\ast=R(\overline{z}).
    \end{align*}
\end{enumerate}
There are two reasons for a real $\lambda$ to belong to $\sigma(L)$. 
\begin{enumerate}
    \item[(i)] $L-\lambda$ is not injective: there exists $v\ne0$ in $H$ such that $Lv=\lambda v$. In this case we say that $\lambda$ is an {\it eigenvalue} of $L$, and $v$ is an {\it eigenvector } corresponding to $\lambda$ (or that $v\in E_\lambda$ belongs to the {\it eigenspace} relative to $\lambda$).
    \item[(ii)] $L-\lambda$ is not surjective.
\end{enumerate}
\section{Intermezzo: the momentum operator $\mathbf{D}=i^{-1}\frac{d\ }{dx}$ on various domains}\label{SectMomOpOnDom}\index{momentum operator}
In this section we discuss the different operator theoretic properties of the {\it momentum operator} $\mathbf{D}=i^{-1}\frac{d\ }{dx}$ when it is defined on this or that domain. The right theoretical (as well as practical) framework to place the corresponding phenomenology is that of the {von Neumann index theorem}. An in depth discussion of the topic for physicists, with plenty of examples, is in 
\href{https://arxiv.org/pdf/quant-ph/0103153}{Self-adjoint extensions of operators and the teaching of quantum mechanics (2001)}, Guy BONNEAU Jacques FARAUT, Galliano VALENT. The proof of the index theorem can be found, e.g., in M. Reed, B. Simon \href{https://books.google.it/books?id=zHzNCgAAQBAJ&printsec=copyright&redir_esc=y#v=onepage&q&f=false}{Methods of modern mathematical physics: II Fourier Analysis, Self-Adjointness (1975)}, Academic Press, p. 135-141.

We will use the same symbol $\mathbf{D}$ for the operator on different domains. We will, however, always specify its domain, not to cause confusion.

\subsection{The real line}\label{SSectMOrealLine}
The operator $(\mathbf{D},C_c^1(\mathbb{R}))$, $C_c^1(\mathbb{R})\xrightarrow{\mathbf{D}=i^{-1}d/dx}L^2(\mathbb{R})$, is symmetric and has graph
\[
\Gamma(\mathbf{D},C_c^1(\mathbb{R}))=\{(\varphi,\mathbf{D}\varphi):\varphi\in C_c^1(\mathbb{R})\}.
\]
Since the Fourier transform is unitary (by Plancherel), we can consider the multiplication operator $M_\omega=\widehat{\mathbf{D}}$ in place of $\mathbf{D}$, where the starting domain is $\mathcal{F}(C_c^1(\mathbb{R}))$.

Let
\begin{equation}\label{eqHoneSob}
    H^1(\mathbb{R})=\{\varphi\in L^2(\mathbb{R}):M_\omega\widehat{\varphi}\in L^2(\mathbb{R})\}
\end{equation}
which is normed by
\[
\|\varphi\|_{H^1}^2=\int_\mathbb{R}(1+\omega^2)|\widehat{\varphi}(\omega)|^2d\omega.
\]
This is the {\it Sobolev space} $H^1(\mathbb{R})$ in disguise, but we do not need here being aware of it. For $\varphi\in C_c^1(\mathbb{R})$,
\[
\|\varphi\|_{H^1}^2=\|\varphi\|_{L^2}^2+\|\varphi'\|_{L^2}^2,
\]
which is the norm induced by the $L^2\times L^2$ norm on the graph $\Gamma(\mathbf{D},C_c^1(\mathbb{R}))$.
\begin{proposition}\label{propSAextOfDonR}
    The closure of $\mathcal{F}\left(C_c^1(\mathbb{R})\right)\xrightarrow{M_\omega}L^2(\mathbb{R})$ is 
    \[
    \mathcal{F}\left(H^1(\mathbb{R})\right)\xrightarrow{M_\omega}L^2(\mathbb{R}).
    \]
    Moreover, this extension is self-adjoint.
\end{proposition}
\begin{proof} We show first that $\overline{\Gamma\left(M_\omega,\mathcal{F}\left(C_c^1(\mathbb{R})\right)\right)}\subseteq \Gamma(M_\omega,\mathcal{F}(H^1))$.
    Suppose $(\psi,\eta)$ is in the closure of $\Gamma\left(M_\omega,\mathcal{F}\left(C_c^1(\mathbb{R})\right)\right)$: 
    \[
    \mathcal{F}\left(C_c^1(\mathbb{R}\right)\ni \psi_n\to\psi\text{ and }M_\omega\psi_n\to \eta\text{ in }L^2.
    \]
    Then, there is a subsequence $\psi_{n_j}(\omega)\to\psi(\omega)$ $a.e.\ \omega$, and a smaller subsequence $\omega\psi_{n_{j_l}}(\omega)\to \eta(\omega)$. Hence, $\omega\psi(\omega)=\eta(\omega)$ $a.e.\ \omega$. Hence, $(\psi,\eta)\in \Gamma(M_\omega,\mathcal{F}(H^1))$.

    In the opposite direction, consider $\psi=\widehat{\varphi}\in \mathcal{F}(H^1(\mathbb{R}))$. We use an approximation of the identity. Let $h\in C_c^\infty(\mathbb{R})$ with support in $[-1,1]$, $h\ge0$ and such that $\int_\mathbb{R}h(x)dx=1$. Let $h_n(x)=nh(nx)$, so that $\widehat{h_n}(0)=1$, $|\widehat{h_n}(\omega)|\le1$, and $\widehat{h_n}(\omega)=\widehat{h}(\omega/n)\to 1$ uniformly on compact sets. The function $h_n\ast\varphi$ lies in $L^2(\mathbb{R})\cap C^\infty(\mathbb{R})$. Also, $(h_n\ast\varphi)'\in L^2(\mathbb{R})$ because
    \[
    \|(h_n\ast\varphi)'\|_{L^2}^2=\frac{1}{2\pi}\|\omega\widehat{h_n}\widehat{\varphi}\|_{L^2}^2\le \frac{1}{2\pi}\|\omega\widehat{\varphi}\|_{L^2}^2<\infty.
    \]
    We have
    \[
    \|M_\omega\widehat{\varphi}-M_\omega\widehat{h_n\ast \varphi}\|_{L^2}=\|\widehat{\varphi}\omega (\widehat{h_n}-1)\|_{L^2}\to0\text{ as }n\to\infty
    \]
    by dominated convergence, since $\omega\widehat{\varphi}(\omega)\in L^2$. Then,
    \[
    \|h_n\ast \varphi-\varphi\|_{H^1}^2=\|h_n\ast \varphi-\varphi\|_{L^2}^2+\frac{1}{2\pi}\|M_\omega\widehat{\varphi}-M_\omega\widehat{h_n\ast\varphi}\|_{L^2}^2\to0\text{ as }n\to\infty.
    \]
    The function $\psi=h_n\ast\varphi$ might not lie in $C_c(\mathbb{R})$. Fix a function $\eta\in C^\infty(\mathbb{R})$ such that $0\le \eta\le 1$, $\eta(x)=1$ for $|x|\le1$ and $\eta(x)=0$ for $|x|\ge2$, and let $\psi_n(x)=\psi(x)\eta(x/n)=\psi(x)\eta_n(x)$. Then, $\psi_n\in C_c^1(\mathbb{R})$,
    \[
    \|\psi-\psi_n\|_{L^2}\to0\text{ as }n\to\infty
    \]
    by dominated convergence, and
    \begin{eqnarray*}
        \|(\psi-\psi_n)'\|_{L^2}&=&\|\psi'(1-\eta_n)-\psi\eta_n'\|_{L^2}\\ 
        &\le&\|\psi'(1-\eta_n)\|_{L^2}+\|\psi\eta_n'\|_{L^2}.
    \end{eqnarray*}
    Since $\psi'\in L^2$, the first term tends to zero as $n\to\infty$ by dominated convergence; while the second term goes to zero as $n\to\infty$ because 
    \[
    \|\psi\eta_n'\|_{L^2}^2=\frac{1}{n^2}\int_{n\le |x|\le 2n}|\psi(x)|^2\eta'(x/n)dx\to0\text{ as }n\to\infty.
    \]
    We have finished the proof that $\overline{\Gamma\left(M_\omega,\mathcal{F}\left(C_c^1(\mathbb{R})\right)\right)}= \Gamma(M_\omega,\mathcal{F}(H^1))$: $\Gamma(M_\omega,\mathcal{F}(H^1))$ is (modulo a unitary map) the closure of $(\mathbf{D},C_c^1(\mathbb{R}))$.

    To finish the proof, we have to show that $\mathcal{F}(H^1)\xrightarrow{\mathbf{D}}L^2(\mathbb{R})$ is self-adjoint. Let $g\in L^2(\mathbb{R})$ be such that the following map is bounded on $\mathcal{F}(H^1)$:
    \begin{eqnarray*}
        f&\mapsto&\frac{1}{2\pi}\langle g,M_\omega f\rangle_{L^2}=\frac{1}{2\pi}\int_{\mathbb{R}}\overline{g(\omega)}\omega f(\omega)d\omega\\ 
        &=&\frac{1}{2\pi}\int_{\mathbb{R}}\overline{\omega g(\omega)} f(\omega)d\omega.
    \end{eqnarray*}
    By density of $H^1(\mathbb{R})$ in $L^2(\mathbb{R})$, it must be $M_\omega g\in L^2(\mathbb{R})$, hence, $g\in \mathcal{F}(H^1)$, which proves self-adjointness.
\end{proof}
\subsection{$\mathbf{D}$ on some domains of functions on $[-\pi,\pi]$}\label{SSectMomemntumInterval}
\subsubsection{$\mathbf{D}$ on $C_c^1(-\pi,\pi)$}\label{SSSectMomemntumIntervalOne}
The operator $C_c^1(-\pi,\pi)\xrightarrow{\mathbf{D}}L^2[-\pi,\pi]$ is densely defined and it is symmetric,
\begin{equation}\label{eqMomemntumIntervalOne}
\int_{-\pi}^\pi \overline{\varphi}(x)i^{-1}\psi'(x)dx=i^{-1}[\overline{\varphi}(x)\psi(x)]_{-\pi}^\pi+\int_{-pi}^\pi\overline{i^{-1}\varphi'(x)}\psi(x)dx=\int_{-pi}^\pi\overline{i^{-1}\varphi'(x)}\psi(x)dx.
\end{equation}
Its closure is symmetric as well. We show that it is not self-ajoint.

For a complex number $k$, let $e_k(x)=e^{kx}$. Then, $e_k\in L^2[-\pi,\pi]$ and $\mathbf{D}e_k=i^{-1}ke_k$. The function $e_k$ belongs to the domain of $\mathbf{D}^\ast$. In fact,
\[
\varphi\mapsto\langle e_k,\mathbf{D}\varphi\rangle_{L^2}=\int_{-\pi}^\pi e^{\overline{k}x}i^{-1}\varphi'(x)dx= 
i\overline{k}\int_{-\pi}^\pi e^{\overline{k}x}i^{-1}\varphi(x)dx
\]
is bounded in $L^2$-norm on $C_0^1(-\pi,\pi)$. If $D^\ast$ where self-adjoint, it would be symmetric, but this would imply that
\[
i\overline{j}\langle e_j,e_k\rangle_{L^2}=\langle D^\ast e_j,e_k\rangle_{L^2}=\langle e_j,D^\ast e_k\rangle_{L^2}=i^{-1}k\langle e_j,e_k\rangle_{L^2},
\]
which implies that $j=-k$ for all $j,k$ real, which is absurd.

The closure of $C_c^1(-\pi,\pi)\xrightarrow{\mathbf{D}}L^2[-\pi,\pi]$ provides an example of closed, densely defined, symmetric operator which is not self-adjoint. We will see below, however, that this operator has infinitely many self-adjoint extensions.

\subsubsection{$\mathbf{D}$ on $C_{per}^1[-\pi,\pi]$}\label{SSectMomemntumIntervalTwo}
We consider here $C_{per}^1[-\pi,\pi]\xrightarrow{\mathbf{D}}L^2[-\pi,\pi]$, which is densely defined and symmetric in virtue of the periodicity: the boundary term in \eqref{eqMomemntumIntervalOne} vanish for different reasons. Clearly $C_{per}^1[-\pi,\pi]\supset C_c^1(-\pi,\pi)$, hence, the closure of $C_{per}^1[-\pi,\pi]\xrightarrow{\mathbf{D}}L^2[-\pi,\pi]$ is a closed, symmetric extension of the closure of $C_c^1(-\pi,\pi)\xrightarrow{\mathbf{D}}L^2[-\pi,\pi]$.

We can reason like we did with $C_c(\mathbb{R})$, with some simplifications. Let $H^1_{per}$ be the space of those $\varphi\in L^2[-\pi,\pi]$ such that $M_n\widehat{\varphi}\in \ell^2(\mathbb{Z})$, where $(M_nh)(n)=n\cdot h(n)$. We show first that $C_{per}^1$ is dense in $H^1$ with respect to the norm
\[
\|\varphi\|_{H^1}^2=\frac{1}{2\pi}\|\widehat{\varphi}\|_{\ell^2}^2+\frac{1}{2\pi}\|M_n\widehat{\varphi}\|_{\ell^2}^2=\frac{1}{2pi}\sum_n.(1+n^2)|\widehat{\varphi}(n)|^2.
\]
Let $\varphi_N(x)=\frac{1}{2\pi}\sum_{|n|\le N}\widehat{\varphi}(n)e^{inx}$, a trigonometric polynomial which surely sits in $C_{per}^1$. 
\begin{equation*}
    \|M_n\widehat{\varphi}-M_n\widehat{\varphi_N}\|_{\ell^2}^2=\sum_{|n|>N}n^2|\widehat{\varphi}(n)|^2\to0\text{ as }n\to\infty,
\end{equation*}
and {\it a fortiori} $\|\varphi-\varphi_N\|_{L^2}\to0$ as $N\to\infty$.

As we did before, one first shows that, modulo a passage to the frequency side, the closure of $\Gamma(\mathbf{D},C^1_{per})$ can be identified with $\Gamma(M_n,\mathcal{F}(H^1))$; then that $H^1\xrightarrow{\mathbf{D}}L^2$ is self-adjoint. Details are left to the non-indolent reader.

\subsubsection{{$\mathbf{D}$ on $C_\lambda^1[0,2\pi]$}\label{SSSectMomemntumIntervalThree}}
The operator $\mathbf{D}$ is symmetric also on the space
\[
C_\lambda^1[0,2\pi]=\{\varphi\in C^1[0,2\pi]:\varphi(2\pi)=e^{i\lambda}\varphi(0)\},
\]
where $\lambda\in\mathbb{R}$ (and we might consider $\lambda\in[0,2\pi)$). The map $\varphi\mapsto e^{-o\lambda x}\varphi=h$ maps $C_\lambda^1$ onto $C_{per}^1$ and we have
\[
(\mathbf{D}\varphi)(x)=(\lambda h(x)+(\mathbf{D}h)(x))e^{i\lambda x},
\]
hence, 
\begin{align*}
    \|\varphi\|_{L^2}^2&=\frac{1}{2\pi}\sum_{n}|\widehat{h}(n)|^2,\\ 
    \|\mathbf{D}\varphi\|_{L^2}^2&=\frac{1}{2\pi}\sum_{n}(\lambda+n)^2|\widehat{h}(n)|^2.
\end{align*}
For $\varphi\in C_\lambda^1$ we have then that the norm induced by the graph of $\mathbf{D}$ is
\[
\|\varphi\|_{H^1_\lambda}^2=\|\varphi\|_{L^2}^2+\|\mathbf{D}\varphi\|_{L^2}^2=\frac{1}{2\pi}\sum_n[1+(\lambda+n)^2]|\widehat{h}(n)|^2.
\]
We call $H^1_\lambda$ the completion of $C_\lambda^1$ with respect to this norm. The operator $\mathbf{D}$ acts on $\widehat{h}$ as $M_{\lambda+n}$, multiplication times $\lambda+n$. It is now routine verifying that $M_{\lambda+n}$ is unitarily equivalent to the closure of $\mathbf{D}$, and that it is self-adjoint.

\subsubsection{Final remarks on the finite interval case} The closure of the operator $\mathbf{D}_0=C_0^1\xrightarrow{\mathbf{D}} L^2$ is not self-adjoint, but we have found a one-parametr family of self-adjoint extensions of it. The adjoint $\mathbf{D}_0^\ast$ of $\mathbf{D}_0$, which is not symmetric, is a common extension of all these self-adjoint extensions we have found.

\subsubsection{The case of $[0,\infty)$}\label{SsectHalfLineMomentum}
Finally, we mention the case of $C_c^1(0,\infty)\xrightarrow{\mathbf{D}}L^2(0,\infty)$, which is symmetric and densely defined. 
To show that its closure is not self-adjoint, use a reasoning similar to that we used with $C_c(0,2\pi)$, using the functions $u_\lambda(x)=e^{-\lambda x}$ ($\lambda>0$). von Neumann's index theory shows that, in fact, this operator does not have any self-adjoint extension.

\section{The spectral theorem: measurable calculus form}\label{ScetSTFCF} In this section we discuss the {\it measurable calculus form} of the spectral theorem: if $D\xrightarrow{L}H$ is self-adjoint and $m\in L^\infty_0(\mathbb{R})$, we can make sense of the expression $H\xrightarrow{m(L)}H$ as a bounded operator. The map $m\mapsto m(L)$ is a $\ast$-homomorphism, and in fact a $\ast$-isomorphism after a notion of "almost everywhere" is set up. This is a far reaching extension of the holomorphic calculus in two respects: we go beyond holomorphic functions; we can deal with unbounded operators (in particular, with linear differential operators). The advantage of the measurable calculus form is that it is canonical: the operator $m(L)$ is uniquely determined by $m$ and $L$.

The statement, however, is long and somehow intricate. In the next section we will translate it into a more transparent form, showing that $L$ is unitarily equivalent to a multiplication operator on a nice measure space. The price to pay is that we will not have universality. We will show how to find a unitary map $\Theta:H\to L^2(\nu)$,  where $\nu$ is a suitable measure on a suitable measurable space, but neither $\nu$, nor $\Theta$, are uniquely determined: a number of arbitrary choices will be made.

\subsection{The measurable calculus for a general self-adjoint operator}\label{SSectMCSAOgeneral}

Let $D\xrightarrow{L}H$ be self-adjoint and $f\in H$. Define:
\begin{equation}\label{eqRisolTransf}
    F_{f,f}(z)=\langle f,R(z)f\rangle,\ z\in\rho(L).
\end{equation}
\begin{theorem}\label{theoRisolTransf}
    \begin{enumerate}
        \item[(i)] $F_{f,f}$ is holomorphic on $\rho(L)\supseteq\mathbb{C}\setminus\mathbb{R}$. Moreover, $F_{f,f}(\overline{z})=\overline{F_{f,f}(z)}$.
        \item[(ii)] We have the estimate
        \[
        \left|F_{f,f}(z)\right|\le\frac{\|f\|}{|\text{Im}(z)|}.
        \]
        \item[(iii)] If $z\in\mathbb{C}_+$, then $F_{f,f}(z)\in\mathbb{C}_+\cup\mathbb{R}$.
        \item[(iv)] We have
        \[
        \lim_{y\to+\infty}i^{-1}yF(iy)=\|f\|^2.
        \]
        \item[(v)] There exists $\mu_{f,f}\ge0$, a finite Borel measure on $\mathbb{R}$, with 
        \begin{equation}\label{eqSpecMeasArrive}
        \langle f,R(z)f\rangle=F_{f,f}(z)=\frac{1}{\pi}\int_{\mathbb{R}}\frac{d\mu_{f,f}(u)}{u-z}
        \end{equation}
        for all $z\in\mathbb{C}\setminus\mathbb{R}$. Moreover, $\frac{\mu_{f,f}(\mathbb{R})}{\pi}=\|f\|^2$.
        \item[(vi)] If $L\ge0$, then $\text{supp}(\mu_{f,f})\subseteq[0,+\infty)$.
    \end{enumerate}
\end{theorem}
The measures $\{\mu_{f,f}:f\in H\}$ are the {\it spectral measures}\index{spectral measures!self-adjoint operator} for the operator $L$. By polarization, we can extend the family to $\{\mu_{f,g}:f,g\in H\}$, where 
\begin{equation}\label{eqSpacMeas}
    F_{f,g}(z):=\langle f,R(z)g\rangle=\frac{1}{\pi}\int_{\mathbb{R}}\frac{d\mu_{f,g}(u)}{u-z},
\end{equation}
for $z\in \mathbb{C}\setminus\mathbb{R}$.
\begin{proof}
    (i) and (ii). We proved that the function $z\mapsto R(z)$ is holomorphic from $\mathbb{C}\setminus\mathbb{R}$ to $\mathcal(H)$, hence $F_{f,f}$ is holomorphic, and
    \[
    |F_{f,g}|\le\|R(z)\|\cdot\|f\|\cdot\|g\|\le\frac{\|f\|\cdot\|g\|}{|\text{Im}z|}.
    \]
    
    (iii).
    If $z=x+iy$ with $y>0$, then
    \[
    \text{Im}F_{f,f}(z)=\text{Im}\langle f,R(z)f\rangle=\text{Im}\langle (L-z)R(z)f,R(z)f\rangle=y\langle R(z)f,R(z)f\rangle\ge0.
    \]

    (iv) and (v). By (ii), (iii), and Herglotz theorem, there exists a positive, finite, Borel measure $\mu_{f,f}$ on the real line such that \eqref{eqSpecMeasArrive} holds for $z\in\mathbb{C}_+$. Since $R(z)^\ast=R(\overline{z})$ and the measure is real valued, the relation extends to $\mathbb{C}_-$. After splitting the Cauchy kernel in its real and imaginary parts, we have
    \[
    i^{-1}yF_{f,f}(iy)=\frac{i^{-1}}{\pi}\int_{\mathbb{R}}\frac{uy}{u^2+y^2}d\mu_{f,f}(u)+\frac{1}{\pi}\int_{\mathbb{R}}\frac{y^2}{u^2+y^2}d\mu_{f,f}(u)=i^{-1}A(y)+B(y).
    \]
    For the imaginary part we have
    \[
    A(y)=\frac{1}{\pi}\int_{\mathbb{R}}\frac{t}{t^2+1}d\mu_{f,f}(ty)\to0\text{ as }y\to\infty
    \]
    because $\mu_{f,f}$ is finite and $h(t)=\frac{1}{1+t^2}\in C_0(\mathbb{R})$.  We are left with the task of proving that $B(y)\to\|f\|^2$ as $y\to\infty$. We start by computing
    \begin{eqnarray*}
        2i\text{Im}F_{f,f}(iy)&=&\langle f,R(iy)f\rangle-\overline{\langle f,R(iy)f\rangle}=\langle f,[R(iy)-R(-iy)]f\rangle\\ 
        &=&\langle f,2iy R(-iy)R(iy)f\rangle\text{ by the resolvent identity}\\ 
        &=&2iy\|R(iy)f\|^2.
    \end{eqnarray*}
    Pick $f=R(i)g\in D$ (recall that $[R(i)](H)=D$):
    \begin{eqnarray*}
        y\text{Im}F_{f,f}(iy)&=&y^2\|R(iy)R(i)g\|^2\\ 
        &=&y^2\left\|\frac{[R(iy)-R(i)]g}{iy-i}\right\|^2\text{ by the resolvent identity}\\ 
        &=&\frac{y^2}{(y-1)^2}\|[R(iy)-R(i)]g\|^2\\ 
        &\to&\|R(i)g\|^2
    \end{eqnarray*}
    since $\|R(iy)\|\le \frac{1}{y}$. More generally, if $f,g\in H$, then
    \begin{eqnarray*}
        |y(\|R(iy)f\|-\|R(iy)R(i)g\|)|^2&\le&y^2\|R(iy)f-R(iy)R(i)g\|^2\\ 
        &\le&y^2\|R(iy)\|^2\|f-R(i)g\|^2\\ 
        &\le&\|f-R(i)g\|^2,
    \end{eqnarray*}
    because $\|R(iy)\|\le 1/y$. Thus,
    \begin{eqnarray*}
        \left|y\|R(iy)f\|-\|f\|\right|&\le&\left|y(\|R(iy)f\|-\|R(iy)R(i)g\|)\right|\\ 
        &\ &+\left|y\|R(iy)R(i)g\|-\|R(i)g\|\right|+\|R(i)g-f\|\\ 
        &\le&2\|f-R(i)g\|+\left|y\|R(iy)R(i)g\|-\|R(i)g\|\right|.
    \end{eqnarray*}
    For fixed $\epsilon>0$ choose first $g$ such that $\|f-R(i)g\|\le\epsilon$, then $y_\epsilon$ such that $\left|y\|R(iy)R(i)g\|-\|R(i)g\|\right|\le\epsilon$ for $y>y_\epsilon$. This shows that $y\|R(iy)f\|\to\|f\|$, as wished.

     (vi). This property will be more easily proved after we have the spectral theory in its multiplication form.
\end{proof}
\begin{lemma}\label{lemmaUniqSpecMeas}
    If $\mu,\nu$ are finite Borel measures on $\mathbb{R}$ such that $\int_{\mathbb{R}}\frac{d\mu(u)}{u-z}=\int_{\mathbb{R}}\frac{d\nu(u)}{u-z}$ for all $z\in\mathbb{C}\setminus\mathbb{R}$, then $\mu=\nu$.
\end{lemma}
\begin{proof}
    For $z=x+iy$, $y>0$, let $h_y(x)=\frac{1}{\pi}\int_{\mathbb{R}}\frac{d\mu(u)}{u-z}$, so that 
    \[
    k_y(x):=\frac{h_y(x)-h_{-y}(x)}{2i}=\frac{1}{\pi}\int_{\mathbb{R}}\frac{yd\mu(u)}{(u-x)^2+y^2}=(P_y\ast \mu)(x),
    \]
    hence, $\widehat{k_y}(\omega)=\widehat{\mu}(\omega)e^{-y|\omega|}$. If $h_y=h_{-y}=0$, then $\widehat{\mu}=0$, hence $\mu=0$.
\end{proof}
As usual, uniqueness leads to formulas.
\begin{proposition}\label{propFormSpecMeasUse}
    The map $(f,g)\mapsto\mu_{f,g}$ from $H\times H$ to the space of the finite Borel measures on $\mathbb{R}$ is sesquilinear and
    \[
    \overline{\mu_{f,g}}=\mu_{g,f}.
    \]
\end{proposition}
\begin{proof}
    We just prove the last item, the other proofs being similar.
    \begin{eqnarray*}
        \frac{1}{\pi}\int_{\mathbb{R}}\frac{d\mu_{f,g}(u)}{u-z}&=&\langle f,R(z)g\rangle=\overline{\langle g,R(\overline{z})f\rangle}\\ 
        &=&\overline{\frac{1}{\pi}\int_{\mathbb{R}}\frac{d\mu_{g,f}(u)}{u-\overline{z}}}\\ 
        &=&\frac{1}{\pi}\int_{\mathbb{R}}\frac{d\overline{\mu_{g,f}}(u)}{u-z}.
    \end{eqnarray*}
\end{proof}
Since $\frac{\mu_{f,f}(\mathbb{R})}{\pi}=\|f\|^2$, by polarization we have $\frac{\mu_{f,g}(\mathbb{R})}{\pi}=\langle f,g\rangle$. We denote by $L^\infty_0(\mathbb{R})$ the space of the Borel measurable, bounded functions on $\mathbb{R}$, normed with the sup norm. The subscript
means that we do not have a notion of almost everywhere, and neither of essential supremum. For $m\in L^\infty_0$, we define $\Gamma(f,g)$ by
\begin{equation}\label{eqMeasCalc}
    \Gamma(f,g):=\frac{1}{\pi}\int_{\mathbb{R}}m(u)d\mu_{f,g}(u),\ \Gamma:H\times H\to\mathbb{C}.
\end{equation}
By proposition \ref{propFormSpecMeasUse}, $\Gamma$ is sesquilinear, and 
\[
|\Gamma(f,g)|\le\|m\|_{L^\infty_0}\frac{|\mu|_{f,g}(\mathbb{R})}{\pi}\le C\|m\|_{L^\infty_0}\|f\|\cdot|g\|. 
\]
By a corollary of the Banach-Steinhaus theorem. there is a bounded operator $m(L)$ on $H$ such that $\Gamma(f,g)=\langle f,m(L)g\rangle$; i.e.
\begin{equation}\label{eqMeasCalcTwo}
    \langle f,m(L)g\rangle=\frac{1}{\pi}\int_{\mathbb{R}}m(u)d\mu_{f,g}(u).
\end{equation}
The map $m\mapsto m(L)$ is the {\it measurable calculus} for $L$. We already know how it forks for some maps $m$:
\begin{align}
    &m(x)=1\ \implies\ m(L)=I\text{ is the identity operator}:\ \langle f,g\rangle=\frac{1}{\pi}\int_{\mathbb{R}}d\mu_{f,g}(u);\\
    &m(x)=(u-z)^{-1}\ \implies\ m(L)=R(z)\text{ if }z\notin\mathbb{R}:\ \langle f,R(z)g\rangle=\frac{1}{\pi}\int_{\mathbb{R}}\frac{1}{u-z}d\mu_{f,g}(u).
\end{align}
We will split the measurable calculus theorem into two statements.
\begin{theorem}\label{theoMeasCalcI}[Spectral theorem, measurable calculus II]\index{measurable calculus!self-adjoint operator}
    There exists a unique linear map $\Lambda:L^\infty_0(\mathbb{R})\to\mathcal{L}(H)$ such that
    \begin{enumerate}
        \item[(i)] for $z\in\mathbb{C}\setminus\mathbb{R}$,
        \[
        \langle f,\Lambda(R(z))g\rangle=\frac{1}{\pi}\int_{\mathbb{R}}\frac{d\mu_{f,g}(u)}{u-z};
        \]
        \item[(ii)] if $m_n$ is a uniformly bounded sequence in $L^\infty_0(\mathbb{R})$, and $m_n\to m$ pointwise, then $\Lambda(m_n)\xrightarrow[WOT]{}\Lambda(m)$.  
    \end{enumerate}
  The map is the one in \eqref{eqMeasCalcTwo},
  \[
  \langle f,\Lambda(m)g\rangle=\frac{1}{\pi}\int_{\mathbb{R}}m(u)d\mu_{f,g}(u).
  \]
\end{theorem}
\begin{proof}
    The map $m\mapsto m(L)$ has properties (i,ii). To prove uniqueness, we check that the hypothesis of the Lebesgue-Hausdorff theorem are verified. By (i) and linearity,
        \begin{eqnarray*}
            \left\langle f,\Lambda\left(\frac{R(z)-R(\overline{z})}{2i}\right)g\right\rangle&=&\frac{1}{2\pi i}\int_{\mathbb{R}}\left(\frac{1}{u-z}-\frac{1}{u-\overline{z}}\right)d\mu_{f,g}(u)\\ 
            &=&\frac{1}{\pi}\int_{\mathbb{R}}\frac{y}{(x-u)^2+y^2}d\mu_{f,g}(u)\\ 
            &=&P_y[\mu_{f,g}](x).
        \end{eqnarray*}
        Let $\mathbf{P}(z)=\frac{R(z)-R(\overline{z})}{2i}=\Lambda(\tau_xP_y)$, where $\tau_x\psi(u)=\psi(u-x)$. Then,
        \[
        \langle f,\mathbf{P}(z)g\rangle=P_y[\mu_{f,g}](x).
        \]
        For $\varphi\in C_c(\mathbb{R})$, 
        \[
        P_y[\varphi](u)=\lim_{N\to\infty}\frac{1}{\pi}\sum_{n=-N}^N\frac{y}{(u-n/N)^2+y^2}\varphi(n/N)\frac{1}{N},
        \]
        where the convergence holds at least pointwise. Using (ii) in the last equality,
        \begin{eqnarray*}
            \frac{1}{\pi}\int_\mathbb{R}P_y[\varphi](u)d\mu_{f,g}(u)&=&\frac{1}{\pi}\int_\mathbb{R}P_y[\varphi](u)d\mu_{f,g}(u)\\ 
            &=&\frac{1}{\pi}\int_\mathbb{R}\left(\lim_{N\to\infty}\frac{1}{\pi}\sum_{n=-N}^N\frac{y}{(u-n/N)^2+y^2}\right.\\ 
            &\ &\left.\varphi(n/N)\frac{1}{N}\right)d\mu_{f,g}(u)\\ 
            &=&\left\langle f,\lim_{N\to\infty}\frac{1}{\pi}\sum_{n=-N}^N\mathbf{P}(-n/N+iy)\varphi(n/N)\frac{1}{N} g\right\rangle\\ 
            &=&\left\langle f,\lim_{N\to\infty}\frac{1}{\pi}\sum_{n=-N}^N\Lambda(\tau_{n/N}P_y\varphi)\varphi(n/N)\frac{1}{N} g\right\rangle\\ 
            &=&\left\langle f,\Lambda(P_y[\varphi]) g\right\rangle.
        \end{eqnarray*}
        Using again (ii),
        \begin{eqnarray*}
          \frac{1}{\pi}\int_\mathbb{R}\varphi(u)d\mu_{f,g}(u)&=&\frac{1}{\pi}\int_\mathbb{R}\lim_{y\to0^+}P_y[\varphi](u)d\mu_{f,g}(u)\\
          &=&\left\langle f,\Lambda(\lim_{y\to0^+}P_y[\varphi]) g\right\rangle\\ 
          &=&\left\langle f,\Lambda(\varphi) g\right\rangle.
        \end{eqnarray*}
        Thus, 
        \[
        \left\langle f,\Lambda(\varphi) g\right\rangle=\frac{1}{\pi}\int_\mathbb{R}\varphi(u)d\mu_{f,g}(u)
        \]
        for $\varphi\in C_c(\mathbb{R})$. By (ii) again, this extends to $C_0(\mathbb{R})$, and we can apply the Lebesgue-Hausdorff theorem.
\end{proof}
\begin{theorem}\label{theoMeasCalc}[Spectral theorem, measurable calculus II]\index{measurable calculus!self-adjoint operator}
 Let $D\xrightarrow{L}H$ be self-adjoint.
\begin{enumerate}
    \item[(i)] The map $m\mapsto m(L)$ satisfies $m(L)^\ast=\overline{m}(L)$.
    \item[(ii)] $d\mu_{f,R(z)g}(u)=\frac{1}{u-z}d\mu_{f,g}(u)$ for all $z\in\mathbb{C}\setminus\mathbb{R}$.
    \item[(iii)] $d\mu_{f,m(L)g}(u)=m(u)d\mu_{f,g}(u)$ for all $m\in L^\infty_0$.
    \item[(iv)] $m\mapsto m(L)$ is multiplicative (it is a $\ast$-morphism of $\ast$-algebras).
    \item[(v)] $\|m(L)\|\le\|m\|_{L^\infty_0}$.
    \item[(vi)] For $E\subseteq\mathbb{R}$, a Borel set, define $\mu(E)=\chi_E(L)$. Then, $\mu$ is a {\bf projection valued measure} (p.v.m.):\index{projection valued measure!self-adjoint operator}
    \begin{itemize}
        \item[(a)] $\mu(E):H\to H$ is an orthogonal projection;
        \item[(b)] if $E=\cup_{n=1}^\infty E_n$ with the $E_n$'s Borel, disjoint, then $\mu(E)=\sum_{n=1}^\infty\mu(E_n)$, where the series converges in the strong operator topology.\index{operator topology!strong}
        \item[(c)] $\mu(\emptyset)=0$ and $\mu(\mathbb{R})=I$. 
    \end{itemize}
    \item[(vii)] If $E$ is bounded, then $\text{Ran}(\mu(E))\subseteq D$.
    \item[(viii)] For $f\in H$ and $g\in D$, $d\mu_{f,Lg}(u)=ud\mu_{f,g}(u)$.
    \item[(ix)] Let $D_c=\cup_{N=1}^\infty\text{Ran}(\mu[-N,N])$. Then, $D_c$ is dense in $D$ and $D_c\xrightarrow{L}D_c$.
    \item[(x)] We have $D=\{f\in H:\frac{1}{\pi}\int_{\mathbb{R}}|u|^2d\mu_{f,f}(u)<\infty\}$. In particular, $D\xrightarrow{m(L)}D$ if $m\in L^\infty_0$.
    \item[(xi)] If $m,m'\in L^\infty_0$, $m'(u)=um(u)$, then
    \begin{itemize}
        \item[(a)] if $f\in H$, then $m(L)f\in D$ and $m'(L)f=Lm(L)f$;
        \item[(b)] if $f\in D$, then $m'(L)f=m(L)Lf$.
    \end{itemize}
\end{enumerate}
\end{theorem}
\begin{proof}
    (i). Linearity is clear. Also,
\begin{eqnarray*}
    \langle f,m(L)^\ast g\rangle&=&\overline{\langle g,m(L)f\rangle}=\overline{\frac{1}{\pi}\int_{\mathbb{R}}m(u)d\mu_{g,f}(u)}\\
    &=&\frac{1}{\pi}\int_{\mathbb{R}}\overline{m(u)}d\overline{\mu_{g,f}(u)}=\frac{1}{\pi}\int_{\mathbb{R}}\overline{m}(u)d\mu_{f,g}(u)\\
    &=&\langle f,\overline{m}(L)^\ast g\rangle.
\end{eqnarray*}

    (ii) We have
    \begin{eqnarray*}
        \frac{1}{\pi}\int_{\mathbb{R}}\frac{d\mu_{f,R(z)g}(u)}{u-w}&=&\langle f,R(w)R(z)f\rangle=\frac{\langle f,[R(w)-R(z)]f\rangle}{w-z}\\ 
        &=&\frac{1}{\pi}\int_{\mathbb{R}}\frac{1}{w-z}\left(\frac{1}{u-w}-\frac{1}{u-z}\right)d\mu_{f,g}(u)\\ 
        &=&\frac{1}{\pi}\int_{\mathbb{R}}\frac{1}{u-w}\frac{1}{u-z}d\mu_{f,g}(u),
    \end{eqnarray*}
    hence,
    \[
    d\mu_{f,R(z)g}(u)=\frac{1}{u-z}d\mu_{f,g}(u).
    \]
    
    (iii). From (ii),
    \[
    d\mu_{f,R(z)g}=\frac{d\mu_{f,g}(u)}{u-z}=\overline{\frac{d\mu_{g,f}(u)}{u-\overline{z}}}=\overline{d\mu_{g,R(\overline{z})f}(u)}=d\mu_{R(\overline{z})f,g},
    \]
    hence,
    \begin{eqnarray*}
        \frac{1}{\pi}\int_{\mathbb{R}}\frac{d\mu_{f,m(L)g}(u)}{u-z}&=&\langle f,R(z)m(L)g\rangle=\langle R(\overline{z})f,m(L)g\rangle\\ 
        &=&\frac{1}{\pi}\int_{\mathbb{R}}m(u)d\mu_{R(\overline{z})f,g}(u)=\frac{1}{\pi}\int_{\mathbb{R}}m(u)d\mu_{f,R(z)g}(u)\\ 
        &=&\frac{1}{\pi}\int_{\mathbb{R}}\frac{m(u)}{u-z}d\mu_{f,g}(u).
    \end{eqnarray*}
    (iv). By (iii),
    \[
    d\mu_{f,m'(L)m(L)g}(u)=m'(u)d\mu_{f,m(L)g}(u)=m'(u)m(u)d\mu_{f,g}(u)=d\mu_{f,[m'm](L)g}(u).
    \]
    (v). 
    \begin{eqnarray*}
        \|m(L)f\|^2&=&\langle f,m(L)^\ast m(L)f\rangle=\langle f,\overline{m}(L)m(L)f\rangle=\langle f,|m|^2(L)f\rangle\\ 
        &=&\frac{1}{\pi}\int_{\mathbb{R}}|m(u)|^2d\mu_{f,g}(u)\\
        &\le&\|m\|_{L^\infty_0}^2\frac{1}{\pi}\int_{\mathbb{R}}d\mu_{f,g}(u)=\|m\|_{L^\infty_0}^2\|f\|^2.
    \end{eqnarray*}
    (vi). Since
    \[
    \mu(E)^\ast=\overline{\chi_e}(L)=\chi_E(L)=\mu(E),
    \]
    the operator $\mu(E)$ is self-adjoint, and 
    \[
    \mu(E)^2=\chi_E^2(L)=\chi_E(L)=\mu(E),
    \]
    hence, $\mu(E)$ is an orthogonal projection. Also, if $E\cap F=\emptyset$, then 
    \[
    \mu(E\cup F)=\chi_{e\cup F}(L)=\chi_E(L)+\chi_F(L)=\mu(E)+\mu(F),
    \]
    thus
    \begin{eqnarray*}
        \|(\mu(E)-\sum_{n=1}^N\mu(E_n))f\|^2&=&\left\|\left(\chi_E(L)-\chi_{\cup_{n=1}^NE_n(L)}\right)f\right\|^2\\
        &=&\left\|\chi_{\cup_{n=N+1}^\infty E_n}(L)f\right\|^2\\    
        &=&\frac{1}{\pi}\int_{\cup_{n=N+1}^\infty E_n}d\mu_{f,f}=\sum_{n=N+1}^\infty \frac{1}{\pi}\int_{E_n}d\mu_{f,f}\\ &=&\sum_{n=N+1}^\infty\|\mu(E_n)f\|^2
        \to0\text{ as }N\to\infty 
\end{eqnarray*}
because 
\[
\sum_{n=1}^\infty\|\mu(E_n)f\|^2=\dots=\frac{1}{\pi}\int_Ed\mu_{f,f}=\langle f,\mu(E)f\rangle\le\|\mu(E)\|^2\cdot\|f\|^2\le\|f\|^2,
\]
hence the tail of the series tends to zero.

(vii). Recall that if $z\notin\mathbb{R}$, then $H\xrightarrow{L}D$ is a bijection. It suffices to show, then, that for each $f\in D$ there is $h\in H$ such that $\mu(E)f=R(z)h$. For now only formally, such $h$ is given by $h=R(z)^{-1}\mu(E)f$, so that
\[
\langle g,h\rangle=\langle g,R(z)^{-1}\mu(E)f\rangle=\langle g,(L-z)\mu(E)f\rangle=\frac{1}{\pi}\int_{\mathbb{R}}(u-z)\chi_E(u)d\mu_{g,f}(u).
\]
With this euristics at hands, we set $m(u)=(u-z)\chi_E(u)$, $m\in L^\infty_0$, and compute:
\begin{eqnarray*}
   \langle g,R(z)m(L)f\rangle&=&\frac{1}{\pi}\int_{\mathbb{R}}\frac{m(u)}{u-z}d\mu_{g,f}(u)=\frac{1}{\pi}\int_{\mathbb{R}}\chi_E(u)d\mu_{g,f}(u)\\ 
   &=&\langle g,\mu(E)f\rangle,
\end{eqnarray*}
which gives $h=m(L)f$.

(viii). If $g\in D$ and $z\notin\mathbb{R}$, then $g=R(z)h$ for some $h\in H$, and $Lg=LR(z)h$. Set $m(u)=\frac{u}{u-w}$, with $w\notin \mathbb{R}$. On the one hand,
\begin{eqnarray*}
    \langle f,R(w)h\rangle&=&\langle f,R(w)(L-z)g\rangle=\langle f,R(w)Lg\rangle-z\langle f,R(w)g\rangle\\
    &=&\frac{1}{\pi}\int_{\mathbb{R}}\frac{d\mu_{f,Lg}(u)}{u-w}-\frac{z}{\pi}\int_{\mathbb{R}}\frac{d\mu_{f,g}(u)}{u-w}.
\end{eqnarray*}
On the other hand,
\begin{eqnarray*}
\langle f,R(w)h\rangle&=&\frac{1}{\pi}\int_{\mathbb{R}}\frac{d\mu_{f,h}(u)}{u-w}=\frac{1}{\pi}\int_{\mathbb{R}}\left(\frac{u}{u-z}-\frac{z}{u-z}\right)\frac{d\mu_{f,h}(u)}{u-w}\\
&=&\frac{1}{\pi}\int_{\mathbb{R}}m(u)d\mu_{f,R(z)h}(u)-\frac{z}{\pi}\int_{\mathbb{R}}\frac{d\mu_{f,R(z)h}}{u-w}\\ 
&=&\frac{1}{\pi}\int_{\mathbb{R}}m(u)d\mu_{f,g}(u)-\frac{z}{\pi}\int_{\mathbb{R}}\frac{d\mu_{f,g}}{u-w}.
\end{eqnarray*}
Comparing the two expressions we have
\[
\frac{1}{\pi}\int_{\mathbb{R}}\frac{d\mu_{f,Lg}(u)}{u-w}=\frac{1}{\pi}\int_{\mathbb{R}}m(u)d\mu_{f,g}(u)=\frac{1}{\pi}\int_{\mathbb{R}}\frac{ud\mu_{f,g}(u)}{u-w}(u),
\]
which shows that $d\mu_{f,Lg}(u)=ud\mu_{f,g}(u)$, as wished.

As  byproduct of the proof, we have that
\begin{equation}\label{eqDOccasione}
\frac{1}{\pi}\int_{\mathbb{R}}d\mu_{f,R(w)Lg}(u)=\frac{1}{\pi}\int_{\mathbb{R}}\frac{d\mu_{f,Lg}(u)}{u-w}=\frac{1}{\pi}\int_{\mathbb{R}}m(u)d\mu_{f,g}(u)=\frac{1}{\pi}\int_{\mathbb{R}}d\mu_{f,m(L)g}(u),
\end{equation}
hence, $R(w)Lg=m(L)g$ for all $g\in D$. We have then:
\begin{equation}\label{eqMlMapstoInD}
    R(w)L=m(L)|_D.
\end{equation}
In particular, $m(L)(D)\subseteq D$, since the range of $R(w)$ is $D$.

(ix). Let $E_n=[-n,n]$, so that by (vii) we have that $\mu(E_n)(H)=\text{Ran}(\mu(E_n))\subseteq D$. If $f=\mu(E_n)h$, by (viii) 
\begin{eqnarray*}
    d\mu_{g,Lf}(u)&=&ud\mu_{g,f}(u)=ud_{g,\chi_{E_n}h}(u)=u\chi_{E_n}(u)d\mu_{g,h}(u)\\
    &=&m(u)d\mu_{g,h}(u) \text{ with $m(u)=u\chi_{E_n}(u)$, bounded,}\\ 
    &=&m(u)\chi_{E_n}(u)d\mu_{g,h}(u)=m(u)d\mu_{g,f}(u)\\ 
    &=&d\mu_{g,m(L)f},
\end{eqnarray*}
so that $Lf=m(L)f=\mu(E_n)m(L)f\in\text{Ran}(\mu(E_n))$. This shows that $\text{Ran}(\mu(E_n))\xrightarrow{L}\text{Ran}(\mu(E_n))$, hence that $D_c$ is mapped by $L$ into itself.

We have to show that $D_c$ is dense in $H$. For $f\in H$, in fact,
\[
\|\mu(E_n)f-f\|^2=\int_{\mathbb{R}\setminus E_n}d\mu_{f,f}\to0
\]
as $n\to\infty$ by dominated convergence.

(x). If $f\in D$, then, by (viii).
\begin{eqnarray*}
    \infty&>&\|Lf\|^2=\frac{1}{\pi}\int_{\mathbb{R}}d\mu_{Lf,Lf}(u)=\frac{1}{\pi}\int_{\mathbb{R}}ud\mu_{Lf,f}(u)\\
    &=&\frac{1}{\pi}\int_{\mathbb{R}}u^2d\mu_{f,f}(u).
\end{eqnarray*}
In the opposite direction, we could try the following:
\begin{eqnarray*}
    \infty&>&\frac{1}{\pi}\int_{\mathbb{R}}u^2d\mu_{f,f}(u)=\lim_{n\to\infty}\frac{1}{\pi}\int_{E_n}u^2d\mu_{f,f}(u)\\ 
    &=&\lim_{n\to\infty}\frac{1}{\pi}\int_{\mathbb{R}}u^2d\mu_{\mu(E_n)f,\mu(e_n)f}(u)\\ 
    &=&\lim_{n\to\infty}\frac{1}{\pi}\int_{\mathbb{R}}d\mu_{L\mu(E_n)f,L\mu(e_n)f}(u)\text{ by (viii), since }\mu(E_n)f\in D_c\subseteq D,\\ 
    &=&\lim_{n\to\infty}\|\mu(E_n)f\|^2.
\end{eqnarray*}
The idea is good, but short of showing that $f\in D$. We use the same calculation in estimating a Cauchy tail:
\begin{eqnarray*}
    \|L\mu(E_{n+m})f-L\mu(E_n)f\|^2&=&\|L\mu(E_{n+m}\setminus E_n)f\|^2\\ 
    &=&\frac{1}{\pi}\int_{\mathbb{R}}d\mu_{L\mu(E_{n+m}\setminus E_n)f,L\mu(E_{n-m}\setminus E_n)}(u)\\ 
    &=&\frac{1}{\pi}\int_{E_{n+m}\setminus E_n}u^2d\mu_{f,f}(u)\\ 
    &\to&0\text{ as $n\to\infty$ by dominated convergence}.
\end{eqnarray*}
We have shown the existence of $h=\lim_{n\to\infty}L\mu(E_n)f$. Since $L$ is closed and $\mu(E_n)f\to f$, we have that $f\in D$ and $h=Lf$.

(xi). We compute:
\[
\frac{1}{\pi}\int_{\mathbb{R}}u^2d\mu_{m(L)f,m(L)f}(u)=\frac{1}{\pi}\int_{\mathbb{R}}u^2|m(u)|^2d\mu_{f,f}(u)=\int_{\mathbb{R}}|m'(u)|^2d\mu_{f,f}(u)<\infty,
\]
hence, by (x), $m(L)f\in D$ if $f\in H$. We can now compute
\begin{eqnarray*}
    \langle f,m'(L)g\rangle&=&\frac{1}{\pi}\int_{\mathbb{R}}um(u)d\mu_{f,g}(u)=\frac{1}{\pi}\int_{\mathbb{R}}ud\mu_{f,m(L)g}(u)\\ 
    &=&\frac{1}{\pi}\int_{\mathbb{R}}d\mu_{f,Lm(L)g}(u)\text{ since }m(L)g\in D\\ 
    &=&\langle f,Lm(L)g\rangle.
\end{eqnarray*}
We have shown that $m'(L)=Lm(L)$.

In the other direction,
\begin{eqnarray*}
    \langle f,m'(L)g\rangle&=&\frac{1}{\pi}\int_{\mathbb{R}}um(u)d\mu_{f,g}(u)\\ 
    &=&\frac{1}{\pi}\int_{\mathbb{R}}m(u)d\mu_{f,Lg}(u)\text{ since }g\in D\\ 
    &=&\frac{1}{\pi}\int_{\mathbb{R}}d\mu_{f,m(L)Lg}(u)=\langle f,m(L)Lg\rangle,
\end{eqnarray*}
hence, $m'(L)=m(L)L$ on $D$.
\end{proof}

\section{The spectral theorem: multiplication form}\label{SectSTMultForm}
\subsection{The multiplicative form of the spectral theorem}\label{SSectSTMultForm}
We proceed like in the case of the unitary operators. Let $D\xrightarrow{L}H$ be self-adjoint. A subset $V$ of $H$ is {\it invariant for $L$ in the sense of the measurable claculus}\index{invariant subspace} if $m(L)(V)\subseteq V$ for all $m\in L^\infty_0(\mathbb{R})$.
\begin{lemma}\label{lemmaSTMultForm} Let $D\xrightarrow{L}H$ be a self-adjoint operator.
    \begin{enumerate}
        \item[(i)] If $V$ is invariant for $L$, then $V$ is linear and closed.
        \item[(ii)] If $V$ is invariant for $L$, then $V^\perp$ is invariant, and $\pi|_Vm(L)=m(L)\pi|_V$.
        \item[(iii)] If $V$ is invariant for $L$, then $D\cap V$ is dense in $V$ and $L(D\cap V)\subseteq V$ is closed and symmetric.
        \item[(iv)] If $V$ is invariant for $L$, then $D\cap V\xrightarrow{L}V$ is self-adjoint.
        \item[(v)] There is an index set $A$ and invariant subspaces $V_\alpha$, $\alpha\in A$, such that
        \[
        H=\bigoplus_{\alpha\in A}V_\alpha.
        \]
        If $H$ is separable, the index set is at most countable.
    \end{enumerate}
\end{lemma}
\begin{proof}
    (i) is clear. (ii) If $g\in V^\perp$ and $f\in V$, then
    \[
    \langle m(L)g,f\rangle=\langle g,m(L)^\ast f\rangle=\langle g,\overline{m}(L)f\rangle=0,
    \]
    thus $m(L)g\in V^\perp$. For $h\in H$,
    \[
    m(L)h=m(L)\pi_Vh+m(L)\pi_{V^\perp}\in V\oplus V^\perp,
    \]
    hence, $\pi_V(m(L)h)=m(L)\pi_Vh$, as wished.

    (iii) At first sight it may not seem obvious that $L|_V$ is densely defined in $V$. Let $\mu$ be the projection valued measure associated with $L$, and recall that $D_c=\cup_N>0\mu[-N,N]$ is dense in $D$. If $h\in V$, then $\mu[-N,N]h\in V\cap D_c$ because $V$ is invariant. Knowing this, the proof that $D_c\cap V$ is dense in $V$ is wholly similar to that that $D_c$ is dense in $H$.

    For $h=\mu[-N,N]f\in V\cap D_c$ ($f\in V$), we have that $Lh=m(L)f$ with $m(t)=t\chi_{[-N,N]}(t)$, hence $Lh\in V$ because $V$ is invariant. This shows that $L|_V$ maps $D_c\cap V$ into $V$. Since $L$ maps $D_c$ into $D_c$, $L$ maps $D_c\cap V$ into $D_c\cap V$.

    We now verify that $L(D\cap V)\subseteq V$. Suppose $g\in D\cap V\subseteq D$. Since $L$ is self-adjoint, there is $f$ in $H$, unique, such that $g=R(i)f$. Then,
    \[
    D\cap V\ni g=R(i)f=R(i)\pi_Vf+R(i)\pi_{V^\perp}f,
    \]
    with $R(i)\pi_Vf\in V$ and $R(i)\pi_{V^\perp}f\in V^\perp$ because $V$ and $V^\perp$ are both invariant. Thus $R(i)\pi_{V^\perp}f=0$, which in turns implies that $\pi_{V^\perp}f=0$, since $R(i)$ is a bijection. Hence, $f\in V$. We have shown that $g=R(i)f$ with $f\in V$. Then,
    \[
    Lg=(L-i)R(i)f+iR(i)f=f+R(i)f\in V,
    \]
    as we had to prove, because $V$ is invariant.

    The symmetry of $L|_V:D\cap V\to V$ follows immediately from what we have proved and the fact that $L$ is symmetric. Finally,
    \[
    \Gamma(L|_V)=\Gamma(L)\cap V\times V,
    \]
    which is closed because $\Gamma(L)$ and $V$ are closed.

    (iv) To show that $L|_D:D\cap V\to H$ is self-adjoint, it suffices to verify that $\pm i\in \rho(L|_V)$, i.e. that $L|_V\pm i:D\cap V\to V$ is invertible; something we did while proving (iii). 

    (v) The proof is step-by-step identical to the one for the unitary operator case.
\end{proof}
. 
We now proceed exactly like in \S \ref{SSectOUmultform}. A vector $f\in H$ is {\it cyclic for the self-adjoint operator $D\xrightarrow{L}H$} if $\{m(L)f:m\in L^\infty_0(\mathbb{R})\}$ \index{cyclic vector} is dense in $H$.
\begin{proposition}\label{propMultSpecSA}
    Suppose $L$ has a cyclic vector $f$. Then, the map $\Theta:L^\infty(\mathbb{R})\to H$, $\Theta(f)=f(L)x$, extends to a unitary map from $L^2(\mathbb{R},\mu_{f,f})$ to $H$, such that
    \begin{equation}\label{eqMultSpecSA}
        ((\Theta^{-1}L\Theta) m)(t)=t m(t).
    \end{equation}
\end{proposition}
    That is, $L$ is unitarily equivalent to multiplication times $t$ on $L^2(\mathbb{R},\mu_{f,f})$,
    \[
    \begin{tikzcd}
L^2(\mathbb{R},\mu_{f,f}) \arrow{r}{\Theta} \arrow[swap]{d}{M_{t}} & H \arrow{d}{L} \\%
L^2(\mathbb{R},\mu_{f,f}) \arrow{r}{\Theta}& H
\end{tikzcd}
\]
\begin{proof}
    The map $\Theta$ is an $L^2(\mu_{f,f})$ isometry on $L^\infty_0(\mathbb{R})$,
    \[
    \|m(L)f\|^2=\langle f,m(L)^\ast m(L)f\rangle=\langle f,|m|^2(L)f\rangle=\int_{\mathbb{R}}|m|^2d\mu_{f,f}.
    \]
    As such, $\Theta$ extends to an isometry of $\overline{L^\infty_0(\mathbb{R})}^{L^2(\mu_{f,f})}$, the closure of $L^\infty(\mathbb{R})$ in the $L^2(\mu_{f,f})$ norm. The hypothesis that $f$ is cyclic means that $\Theta(L^\infty_0(\mathbb{R}))$ is dense in $H$, hence that $\Theta$ is a surjective isometry.
    
    Suppose now that $m\in L^\infty_0(\mathbb{R})$ is such that $xm(x)=m'(x)$ defines a bounded function. Such $m$'s are dense in $L^2(\mu_{f,f})$, and, by theorem \ref{theoMeasCalc} (ix), 
    \[
    L\Theta(m)=Lm(L)f=m'(L)f=\Theta(M_xm)f.
    \]
    The relation extends by the unitary equivalence to the domain om $M_x$.
\end{proof}
\begin{theorem}\label{theoMultSpecSA}[Spectral theorem for self-adjoint operators: multiplicative form]\index{theorem!spectral in multiplicative form for self-adjoint operators}
    Let $D\xrightarrow{L}H$ be a self-adjoint operator on a separable Hilbert space $H$. Then, there exists a locally compact space $X$, a finite Borel measure $\nu$ on $X$, and a continuous function $\varphi:X\to\mathbb{R}$, and a unitary map $\Theta:L^2(\nu)\to H$, such that
    \begin{equation}\label{eqMultSpecUObis}
        (\Theta^{-1}U\Theta)m =\varphi m,
    \end{equation}
    is multiplication times $\varphi$.
\end{theorem}
The proof follows from proposition \ref{propMultSpecSA} and lemma \ref{lemmaSTMultForm} (v) (with an index set $A$ which is at most countable, in the separable case), line by line as in the corresponding statement for the unitary operators. In fact, the proof produces a rather concrete model for $X$ and $\varphi$.
\begin{enumerate}
    \item[(i)] $X=\mathbb{R}\times A$, with the product topology ($A$ is considered discrete);
    \item[(ii)] $d\nu(x,n)=d\mu_{f_n,f_n}(x)$, where $f_n$ is the cyclic vector for the $n^{th}$ summand in the decomposition of $H$ into invariant vectors;
    \item[(iii)] $\varphi(x,n)=x$: $\varphi$ multiplies times the coordinate $x$ on each copy $\mathbb{R}\times\{n\}$ of the real line.
\end{enumerate}
A picture of $\varphi$ can be obtained by changing coordinates in such a way $\mathbb{R}\times\{n\}$ is mapped onto $(-1/2+n,1/2+n)$.
\subsection{Some consequences}\label{SSectSTMultFormCons}

\subsubsection{The spectral theory of a multiplication operator}\label{SSectSTofMulSA}
We translate here in the world of self-adjoint operators what we saw in \S \ref{SSectSTofMultO} for the unitary operators. We omit the proofs, which are essentially the same.

The multiplicative form of the spectral theorem tells us that any self-adjoint operator is unitarily equivalent to a multiplication operator $M_\eta$ of the following form, for which we can explicitly write down the spectrum, the spectral measures, and so on.
\begin{itemize}
    \item[(a)] {\it We have a locally compact space $X$ and a finite, positive Borel measure $\nu$, which we can assume to be finite, and a surjective local homeomorphism $\eta:X\to\mathbb{R}$, such that $M_\eta:L^2(\nu)\to L^2(\nu)$ is the operator $M_\eta f=\eta f$.} Actually, we should specify the domain of $M_\eta$, 
    \[
    D(M_\eta)=\{f:X\to\mathbb{C}\text{ such that }\int_X(1+\eta^2)|f|^2d\nu<\infty\}.
    \]
    \item[(b)] {\it $\sigma(M_\eta)=\overline{\eta(\text{supp}(\nu))}\subseteq \mathbb{T}$ is the spectrum of $M_\eta$.}
    \item[(c)] {\it Let $f\in L^2(\nu)$. The spectral measure $\mu_{f,f}$ is such that, for any bounded function $m:\mathbb{R}\to\mathbb{C}$,
    \[
    \int_Xm(\eta)|f|^2d\nu=\langle f,m(M_\eta)f\rangle_{L^2(\nu)}=\int_{\mathbb{R}}md\mu_{f,f}.
    \]
    That is, $d\mu_{f,f}=\eta_\ast(|f|^2d\nu)$. Similarly, $d\mu_{f,g}=\eta_\ast(\overline{g}fd\nu)$.}
    \item[(d)] {\it The p.v.m. $\mu$ maps Borel measurable sets $E$ in $\mathbb{R}$ to projections of $L^2(\nu)$ in the following fashion:
    \[
    \int_X\overline{f}[\mu(E)g]d\nu=\langle f,\mu(E)g\rangle_{L^2(\nu)}=\mu_{f,g}(E)=\int_X\chi_E(\eta)\overline{f}gd\nu,
    \]
    i.e. ($\chi_E(\eta)=\chi_E\circ\eta$)
    \begin{equation}\label{eqMuForSPMFbis}
        \mu(E)g=\chi_E(\eta)g,\ \mu(E)=M_{\chi_E(\eta)}.
    \end{equation}}
    \item[(e)] {\it A Borel set $E\subseteq\mathbb{R}$ is a null-set for $\mu$ if and only if $0=\chi_E\circ\eta=\chi_{\eta^{-1}(E)}$ $\nu-a.e.$, i.e. if and only if $\nu(\eta^{-1}(E))=0$.}
    \item[(f)] {\it By (e) and the topological properties of $\eta$,  $\overline{\eta(\text{supp}(\nu))}=\text{supp}(\mu)$. }
    \item[(g)] {\it By (b) and (f),\index{spectrum!self-adjoint operator}\index{projection valued measure!self-adjoint operator}
    \begin{equation}\label{esSpectMultOprtrSA}
        \text{supp}(\mu)=\sigma(M_\eta).
    \end{equation}}
\end{itemize}

\subsubsection{Some applications of the spectral theorem in its multiplicative form}\label{SSectSTMFapplied}
\begin{theorem}\label{theoSpectrumSA}
    Let $L:D\to H$ be self-adjoint, $\{\mu_{x,x}:x\in H\}$ be the family of the spectral measures for $L$, and $\mu$ be the projection valued measure associated with $L$. We have
        \begin{equation}\label{eqSuppOfSMforUO}
            \text{supp}(\mu)=\sigma(L).
        \end{equation}
\end{theorem}
\begin{proof}[Proof of theorem \ref{theoSpectrumSA}]
    By the multiplicative form of the spectral theorem, it suffices to show it for multiplication operators. This is done in \eqref{esSpectMultOprtrSA}.
\end{proof}
\begin{corollary}\label{corSpectrumUObis}
With the same notation of theorem \ref{theoSpectrumSA}, we also have the following.
 \begin{enumerate}
        \item[(i)] For a function $f\in L^\infty(\mu)$, $f(L)=0$ if and only if $f$ vanishes $\mu$-a.e on $\sigma(L)$. Thus, the map $f\mapsto f\chi_{\sigma(L)}$ is the identity of $L^\infty(\mu)$. To highlight this fact, we also write $L^\infty(\mu)=L^\infty(\sigma(L))$, and for $\varphi:\sigma(L)\to\mathbb{C}$ measurable and bounded, we write $\varphi(L)=f(L)$, where $f$ is any measurable extension of $\varphi$ to $\mathbb{R}$. 
        \item[(ii)] If $\{f_n:n\ge1\}$ is a sequence in $L^\infty(\sigma(L))$, if $\|f_n\|_{L^\infty(\sigma)}\le C$ is uniformly bounded, and $f_n\to f$ $\mu-a.e.$, then $f_n(L)\to f(L)$ in the strong operator topology.\index{operator topology!strong}
 \end{enumerate}   
\end{corollary}
Another application is a relationship between unitary and self-adjoint operators. We will see some deeper result when we discuss Stone's theorem of strongly continuous, one parameter groups of unitary operators.
\begin{theorem}\label{theoUOSA}
    \begin{enumerate}
        \item[(i)] Let $H\xrightarrow{U}H$ be a unitary operator. Then there exist self-adjoint operators $L$, some of which are bounded, such that $U=e^{iL}$.
        \item[(ii)] Viceversa, if $D\xrightarrow{L}H$ is self-adjoint, then $e^{itL}$ is a unitary operator for all $t\in\mathbb{R}$.
    \end{enumerate}
    \begin{proof}
        (i) See proposition \ref{propUnitaryIsConnected}. (ii) By the measurable calculus, $L$ is unitarily equivalent to $M_\eta$ on some measure space, with $\eta$ measurable and real valued. It is easy to see, by the uniqueness of the measurable calculus, that $e^{itM_\eta}=M_{e^{it\eta}}$, and the latter is a unitary operator on $L^2(\nu)$ for all $t\in\mathbb{R}$.
    \end{proof}
\end{theorem}
\section{One parameter groups of unitary operators and Stone's theorem}\label{SectStone}\index{operator!unitary}
A {\it one parameter group of unitary operators}\index{strongly continuous one parameter group of unitary operators} on a Hilbert space $H$ is a map $t\mapsto U_t$ from $\mathbb{R}$ to the group $\mathcal{U}(H)$ of the unitary operators on $H$ such that
\begin{equation}\label{eqGoUO}
    U_{s+t}=U_sU_t, \text{ and }U_0=I.
\end{equation}
The group is {\it strongly continuous} if the map is continuous with respect to the strong operator topology,\index{operator topology!strong}
\begin{equation}\label{eqSCGoUO}
    \lim_{t\to0}U_th=h \text{ in }H
\end{equation}
for all $h\in H$.
\begin{proposition}\label{propStoneInv}
    Let $D\xrightarrow{L}H$ be self-adjoint, and for $t\in\mathbb{R}$ let $U_t:=e^{itL}$.\index{operator!self-adjoint} Then,
    \begin{enumerate}
        \item[(i)] $t\mapsto U_t$  is a strongly continuous group of unitary operators;
        \item[(ii)] for all $f\in D$,
    \begin{equation}\label{eqGenInfUOG}
    \lim_{t\to0}\frac{U_tf-f}{t}=iLf
     \end{equation}
     in $H$;
        \item[(iii)] if $f\in H\setminus D$, then the limit in \ref{eqGenInfUOG} does not exist.
    \end{enumerate}
\end{proposition}
The proposition might be read as providing a solution for the abstract, general version of the {\it Schroedinger equation}\index{Schroedinger equation}. The equation
\begin{equation}\label{eqSchroed}
    i^{-1}\dfrac{d\ }{dt}u(t)=Lu(t),
\end{equation}
where $u:\mathbb{R}\to H$, has a solution $u$ satisfying $u(0)=f\in D$, which is provided by
\begin{equation}\label{eqSchroedSol}
    u(t)=e^{itL}f.
\end{equation}
It is an important fact that the function $u$ in \eqref{eqSchroedSol} makes perfect sense, and it is strongly continuous, even when $f\in H\setminus D$.
\begin{proof}
    By the spectral theorem in multiplicative form, it suffices to prove the proposition for an operator of the form $L=M_g$ on $H=L^2(\nu)$, where $\nu$ is a finite Radon measure on a locally compact metric space $X$, and $g:X\to\mathbb{R}$ is continuous. In this framework, $U_t=M_{e^{itg}}$ is the operator of multiplication times $e^{itg}$. Since $|e^{itg}|=1$, the operator is unitary, as we remarked earlier, and $U_t^\ast$ is multiplication times $e^{-itg}$.

    (i) We only have to verify continuity. 
    \[
    \|U_{t+s}f-U_tf\|_{L^2(\nu)}^2=\int_X|e^{itg}-1|^2|f|^2d\nu\to0\text{ as }s\to0,
    \]
    by dominated convergence.

    (ii) We have
    \[
    \frac{U_tf-f}{t}-iLf=\left(\frac{e^{itg}-1}{t}-ig\right)f\to0\text{ as }t\to0\ \nu-a.e.\text{ on }X,
    \]
    and
    \[
    \left|\frac{e^{itg}-1}{t}\right|=|g|\left|\frac{e^{itg}-1}{tg}\right|\le C|g|
    \]
     because $|(e^{it}-1)/t|$ is bounded on $\mathbb{R}\setminus\{0\}$. We have then that
     \[
     \left|\left(\frac{e^{itg}-1}{t}-ig\right)f\right|\le (C+1)|fg|\in L^2(\nu),
     \]
     if $f\in D$, the domain of $M_g$. Then
     \[
     \lim_{t\to0}\left\|\left(\frac{e^{itg}-1}{t}-ig\right)f\right\|_{L^2(\nu)}=0,
     \]
     by dominated convergence.

     (iii) Viceversa, if there is $h\in L^2(\nu)$ such that
     \[
     \lim_{t\to0}\left\|\frac{e^{itg}-1}{t}f-h\right\|_{L^2(\nu)}=0,
     \]
     then, for some sequence $t_n\to0$,
     \[
     ifg-h=\lim_{n\to\infty}\frac{e^{it_ng}-1}{t_n}f-h=0\ \nu-a.e.,
     \]
     hence, $L^2(\nu)\ni fg,f$, which is the condition for membership of $f$ in $D$. 
\end{proof}
This simple proposition has a converse.
\begin{theorem}\label{theoStone}[Stone's theorem]\index{theorem!Stone}
    Let $\{U_t\}_{t\in\mathbb{R}}$ be a strongly continuous one parameter group of unitary operators on $H$. Then, there exists a unique self-adjoint operator $D\xrightarrow{L}H$ such that $U_t=e^{itL}$.
\end{theorem}
I present here two proofs. The first one, from Reed and Simon, does not require the Lebesgue-Hausdorff theorem. The second one, from  \href{https://terrytao.wordpress.com/2011/12/20/the-spectral-theorem-and-its-converses-for-unbounded-symmetric-operators/}{Tao's notes},  makes use of it.
\begin{proof}[Proof from \cite{RS1975}]
    For $\varphi\in C_c^\infty(\mathbb{R})$ and $f\in H$, let 
    \[
    f_\varphi=\int_{\mathbb{R}}\varphi(t)U_tfdt,
    \]
    which you might think of as a "smoothing" of the vector $f$, and let
    \[
    D=\text{span}\{f_\varphi:\varphi\in C_c^\infty,f\in H\}.
    \]
    Let $\eta_{\epsilon}(t)=1/\epsilon\eta(t/\epsilon)$ and approximation of the identity,
    \[
    \eta\in C^\infty_c(\mathbb{R}),\ \eta(t)=0\text{ if }|t|>1,\ \eta\ge0,\ \int_{\mathbb{R}}\eta(t)dx=1.
    \]
    We have then:
    \begin{eqnarray*}
       \|f_{\eta_\epsilon}-f\|&=&\left\|\int_{\mathbb{R}}\eta_\epsilon(t)[U_tf-f]dt\right\|\le\int_{\mathbb{R}}\eta_\epsilon(t)\sup_{|t|\le\epsilon}\|U_tf-f\|dt\\ 
       &\to&0\text{ as }\epsilon\to0.
    \end{eqnarray*}
    This shows that $D$ is dense in $H$.

    For $f_\varphi\in D$,
    \begin{eqnarray*}
        \frac{U_s-I}{s}f_\varphi&=&\int_{\mathbb{R}}\varphi(t)\frac{U_{s+t}-U_t}{s}fdt\\
        &=&\int_{\mathbb{R}}\frac{\varphi(\tau-s)-\varphi(\tau)}{s}U_\tau fd\tau\\ 
        &\to&-\int_{\mathbb{R}}\varphi'(\tau)U_\tau fd\tau\text{ in $H$ as }s\to0,
    \end{eqnarray*}
    since
    \[
    \lim_{s\to0}\frac{\varphi(\tau-s)-\varphi(\tau)}{s}=-\varphi'(\tau),
    \]
    uniformly with respect to $\tau$.

    For $f_\varphi\in D$, it is now natural defining
    \[
    Lf_\varphi:=i^{-1}f_{-\varphi'},
    \]
    so that $L:D\to D$. Also, it is easily verified that $U_s:D\to D$, and in fact that $U_s[f_\varphi]=[U_sf]_\varphi$. We now verify that
    \begin{equation}\label{eqSTcommute}
        U_sLf_\varphi=LU_sf_\varphi.
    \end{equation}
    In fact,
    \begin{eqnarray*}
       U_sLf_\varphi&=&U_si^{-1}f_{-\varphi'}=i^{-1}[U_sf]_{-\varphi'}=L[U_sf]_\varphi\\ 
       &=&LU_sf_\varphi.
    \end{eqnarray*}
    Also, $L$ is symmetric on $D$,
    \[
    \langle f_\varphi Lg_\psi\rangle=\lim_{s\to0}\left\langle f_\varphi,\frac{u_s-I}{is}g_\psi \right\rangle
    =\lim_{s\to0}\left\langle \frac{u_{-s}-I}{-is} f_\varphi,g_\psi \right\rangle=\langle L f_\varphi g_\psi\rangle.
    \]
    We have to prove that $D\xrightarrow{L}H$ is essentially self-adjoint\index{essentially self-adjoint operator}, i.e. that $\text{ker}(L^\ast- i)=0$. Suppose that $u\in\text{ker}(L^\ast- i)$. Then, for all $f\in D$,
    \begin{eqnarray*}
        \frac{d\ }{dt}\langle u, U_tf\rangle&=&\langle u,iL U_tf\rangle=i\langle L^\ast u,U_tf \rangle\\ 
        &=&i\langle iu,U_tf \rangle=\langle u,U_tf \rangle.
    \end{eqnarray*}
    i.e. $\alpha(t)=\langle u, U_tf\rangle$ satisfies $\alpha'(t)=\alpha(t),\ \alpha(0)=\langle u, f\rangle$. It must be
    \[
    \alpha(t)=e^t\langle u,f\rangle.
    \]
    On the other hand, it must be $|\alpha(t)|=|\langle u, U_tf\rangle|\le\|u\|\cdot\|f\|$, hence $0=\alpha(0)=\langle u, f\rangle$. Since $D$ is dense in $H$, $u=0$.
    \end{proof}

    The second proof is more spectral. We can not use directly the spectral theorem in its multiplicative form, because we do not know at this stage if the operators $U_t$ can be simultaneously "diagonalized". The proof has to go back to spectral measures, which have to be reconstructed from scratch to fit the present case, where we have a continuum of operators to work with simultaneously.
\begin{proof}[Proof from \cite{Tao2011}]
    \noindent{\bf Step I.} By proposition \ref{propStoneInv}, the domain of $L$ and its action are determined by the group,
    \[
    D=\left\{f\in H:\lim_{t\to0}\frac{U_tf-f}{t}f\text{ exists in }H\right\},\ Lf=i^{-1}\lim_{t\to0}\frac{U_tf-f}{t}f.
    \]
    \noindent{\bf Step II.} {\it We show that (simultaneous diagonalization) for each $f$ in $H$ there is a finite, positive Radon measure $\mu_{f,f}$ on $\mathbb{R}$ such that:
    \[
    \mu_{f,f}(\mathbb{R})=\|f\|^2, \text{ and }\langle f,U_tf\rangle=\int_{\mathbb{R}}e^{itx}d\mu_{f,f}(x).
    \]}
    The function $\varphi(t)=\langle f,U_tf\rangle$ is continuous on $\mathbb{R}$,
    \[
    |\varphi(t)-\varphi(s)|\le\|f\|\cdot\|[U_t-U_s]f\|=\|f\|\cdot\|[U_{t-s}f-f\|\to0\text{ as }t\to s,
    \]
    uniformly, by the hypothesis of continuity. Also, $\varphi$ is positive definite,
    \begin{eqnarray*}
        \sum_{l,m=1}^n\varphi(t_m-t_l)c_m\overline{c_l}&=&\sum_{l,m=1}^n\langle c_l U_{t_l}f,c_m U_{t_m}f\rangle\\ 
        &=&\sum_{m=1}^n\left\|c_mU_{t_m}f\right\|^2\ge0.
    \end{eqnarray*}
    By Bochner's theorem a measure with the desired properties exists and it is unique. Also,
    \[
    \mu_{f,f}(\mathbb{R})=\int_{\mathbb{R}}e^{i0x}d\mu_{f,f}(x)=\langle f,If\rangle=\|f\|^2.
    \]
    \noindent{\bf Step III.} By polarization, for all $f,g\in H$ we can construct Radon, finite, complex measures $\mu_{f,g}$ such that
    \[
    \langle f,U_tg\rangle=\int_{\mathbb{R}}e^{itx}d\mu_{f,g}(x),
    \]
    and the map $(f,g)\mapsto\mu_{f,g}$ is sesquilinear. Moreover,
    \begin{eqnarray*}
        \int_{\mathbb{R}}e^{itx}d\mu_{f,g}(x)&=&\langle f,U_tg\rangle=\langle U_{-t}f,g\rangle=\overline{\langle g,U_{-t}f\rangle}\\ 
        &=&\int_{\mathbb{R}}e^{itx}d\overline{\mu_{g,f}}(x),
    \end{eqnarray*}
    i.e. $\mu_{f,g}=\overline{\mu_{g,f}}$, by the injectivity of the Fourier transform on the space of the complex measures.

    \noindent{\bf Step IV.} {\it For each $m\in L^\infty_0(\mathbb{R})$  there is a unique operator $m(L)$ in $\mathcal{L}(H)$ such that 
    \[
    \langle f,m(L)g\rangle=\int_{\mathbb{R}}m(x)d\mu_{f,g}(x).
    \]
    Clearly, the map $m\mapsto m(L)$ is linear. Also, 
    \begin{equation}\label{eqEstStoneInt}
        \|m(L)\|\le \|m\|_{L^\infty_0}.
    \end{equation}
    Moreover, If $m_n\to m$ pointwise, $m_n,m\in L^\infty_0(\mathbb{R})$, and the $m_n$'s are uniformly bounded, then 
    \begin{equation}\label{eqPtwsStone}
        \lim_{n\to\infty}m_n(L)=m(L)
    \end{equation}
    in the weak operator topology.\index{operator topology!weak}
    }
    Elementary estimates, the definition of $\mu_{f,g}$, and the fact that $\mu_{f,f}(\mathbb{R})=\|f\|^2$, imply that
    \[
    \left|\int_{\mathbb{R}}m(x)d\mu_{f,g}(x)\right|\le C\|m\|_{L^\infty_0}(\|f\|+\|g\|).
    \]
    Passing to sup for $\|f\|,\|g\|\le1$, we see that the sesquilinear form
    \[
    (f,g)\mapsto \int_{\mathbb{R}}m(x)d\mu_{f,g}(x)
    \]
    is bounded, hence that a bounded operator $m(L)\in\mathcal{L}$ with the desired properties exists, and the estimate on its norm improves to $\|m(L)\|\le \|m\|_{L^\infty_0}$, by estimating the integral.

    The limit in \eqref{eqPtwsStone} follows from dominated convergence.

    \noindent{\bf Step V.} {\it We have
    \[
    d\mu_{f,m(L)g}(x)=m(x)d\mu_{f,g}(x).
    \]}
    We start with $m(x)=e^{isx}$, which corresponds to, we will see, $m(L)=U_s$. For all $t$,
    \begin{eqnarray}
        \int_{\mathbb{R}}e^{itx}e^{isx}d\mu_{f,g}(x)&=&\int_{\mathbb{R}}e^{i(t+s)x}d\mu_{f,g}(x)=\langle f,U_{t+s}g\rangle\\ 
        &=&\langle f,U_t(U_sg)\rangle=\int_{\mathbb{R}}e^{itx}d\mu_{U_sf,g}(x),
    \end{eqnarray}
    hence, $e^{isx}d\mu_{f,g}(x)=d\mu_{U_sf,g}(x)$. As a consequence, $\mu_{f,U_t}g=\mu_{U_{-t}f,g}$, since
    \[
    \int_\mathbb{R}e^{isx}d\mu_{f,U_tg}(x)=\langle f,U_sU_tg\rangle=\langle U_{-t}f,U_sg\rangle=\int_\mathbb{R}e^{isx}d\mu_{U_{-t}f,g}(x)
    \]
    and the Fourier transform is injective on the class of the finite Borel measures.
    Then, for $m\in L^\infty_0$ and $t\in\mathbb{R}$,
    \begin{eqnarray*}
        \int_{\mathbb{R}}e^{itx}m(x)d\mu_{f,g}(x)&=&\int_{\mathbb{R}}m(x)d\mu_{f,U_tg}(x)=\int_{\mathbb{R}}m(x)d\mu_{u_{-t}f,g}(x)\\ 
        &=&\langle U_{-t}f,m(L)g\rangle=\langle f,U_tm(L)g\rangle\\ 
        &=&\int_{\mathbb{R}}e^{itx}d\mu_{f,m(L)g}(x),
    \end{eqnarray*}
    which, again by injectivity of the Fourier transform, implies $m(x)d\mu_{f,g}(x)=d\mu_{f,m(L)g}(x)$.\footnote{I thank Terence Tao for providing me this argument, which shortened a more intricate one in a previous version of the notes.}
       
   \noindent{\bf Step VI.} {\it The map $m\mapsto m(L)$ is a $\ast$-homomorphism from $L^\infty_0(\mathbb{R})$ to $\mathcal{L}(H)$.} We proceed exactly as in the case of the unitary operators. Linearity is clear from Step IV, and the multiplicative property follows from Step V,
   \begin{eqnarray*}
       \langle f,(mn)(L)g\rangle&=&\int_{\mathbb{R}}mnd\mu_{x,y}=\int_{\mathbb{R}}md\mu_{x,n(L)y}\\
       &=&\int_{\mathbb{R}}d\mu_{x,m(L)n(L)y}=\langle f,m(L)n(L)g\rangle,
   \end{eqnarray*}
    if $f,g\in H$ and $m.n\in L^\infty_0(\mathbb{R})$. The fact that $m(L)^\ast=\overline{m}(L)$ is left as an exercise.

    \noindent{\bf Step VII.} {\it For each Borel measurable $E$ in $\mathbb{R}$, define $\mu(E)\in\mathcal{L}(H)$ by
    \begin{equation}\label{eqpvmSAO}
        \langle f,\mu(E)g\rangle=\int_Ed\mu_{f,g}.
    \end{equation}
    Then, $\mu:\mathbf{B}(\mathbb{R})\to\mathcal{L}(H)$ is a projection valued measure. 
    } The proof of this fact is wholly similar to the one in section \ref{SSectUOpvm}.

    \noindent{\bf Step VIII.} {\it There exists a self-adjoint operator $L$ having $\mu$ as its spectral measure.} We first define $L$ as a densely defined, symmetric operator, and show that it is essentially self-adjoint. \index{essentially self-adjoint operator}Let
    \[
    D_c=\cup_{N\ge1}\mu[-N,N],
    \]
    which is dense in $H$ (the proof is identical to that we have seen when proving the spectral theorem). For $g\in D_c$ define $Lg$ by
    \[
    \langle f,Lg\rangle:=\int_\mathbb{R}ud\mu_{f,g}(u).
    \]
    The integral converges absolutely because, for some $N\ge1$, $d\mu_{f,g}(u)=d\mu_{f,\mu[-N,N]g}(u)=\chi_{[-N,N]}(u)d\mu_{f,g}(u)$.

    Since the function $u\mapsto u$ is real valued, $L$ is symmetric. Moreover, $L(\mu[-N,N](H))\subseteq \mu[-N,N](H)$.  From this, it is easy to deduce that $(L\pm i)(H)$ is dense in $H$. In fact, on $\mu[-N,N](H)$,
    \[
    \langle f,(L\pm i)g\rangle=\int_{-N}^N(u\pm i)d\mu_{f,g}(u),
    \]
    and we have $1\le|u\pm i|^2\le N^2+1$, hence $L\pm i:\mu[-N,N](H)\to\mu[-N,N](H)$ is invertible with bounded inverse. Hence, $L$ is essentially self-adjoint.

    Next, we have to show that $\mu$ is the p.v.m. associated to $L$, or, equivalently, that $\mu_{f,g}$ are the spectral measures associated with $L$, let's call them $\nu_{f,g}$. If $z\in\mathbb{C}\setminus\mathbb{R}$, then 
    \[
    \langle f,(L-z)^{-1}g\rangle=\int_{\mathbb{R}}\frac{d\nu_{f,g}(u)}{u-z}.
    \]
    If we verify that
    \begin{equation}\label{eqSMareTheSame}
    \langle f,(L-z)^{-1}g\rangle=\int_{\mathbb{R}}\frac{d\mu_{f,g}(u)}{u-z},
    \end{equation}
    we have $\mu_{f,g}=\nu_{f,g}$ by the uniqueness of the Herglotz representation, and we are done.

    Let $m(u)=\frac{1}{u-z}$, so that $\langle f,m(L)g\rangle=\int_{\mathbb{R}}\frac{d\mu_{f,g}(u)}{u-z}$. Then,
    \[
    \langle f,(L-z)m(L)g\rangle=\int_{\mathbb{R}}(u-z)d\mu_{f,m(L)}(u)=\langle f,(u-z)m(z)d\mu_{f,g}(u)=\langle f,g\rangle.
    \]
    Hence $m(L)=(L-z)^{-1}$, and \eqref{eqSMareTheSame} holds.

    \noindent{\bf Step IX.} {\it Conclusion.} By the multiplicative form of the spectral theorem, after a unitary map our operator $L$ has the form $L\varphi=\eta\varphi$ on $L^2(X,\lambda)$ for some measure space $(X,\lambda)$ (and $\eta$ is measurable and real valued), and the operators $U_t$ have the form $U_t\varphi=e^{i\eta}\varphi$. Stone's theorem now follows from its converse, proposition \ref{propStoneInv}.
\end{proof}

    \section{Some applications to bounded and compact self-adjoint operators, and to the spectrum}\label{SectAppSPtoBOCO}
    Having at your disposal various versions of the spectral theorem for self-adjoint operators, especially the multiplicative one, by means of pictures, intuition, and measure theory, a number of conclusions can be drawn. We make a list of some of them, without striving for the greatest generality. The hints should be enough for you to come up with a proof.
    \subsection{Bounded operators, continuous calculus, spectrum}\label{SSectAppSPBO}
    \begin{enumerate}
        \item[(i)] {\it A self-adjoint operator $L$ is bounded if and only if the corresponding p.v.m. measure has bounded support in $\mathbb{R}$.} If $\mu$ has support in $[-N.N]$, then $L=L\mu[-N.N]$ is bounded. Viceversa, is the support of $\mu$ is unbounded you can find elements $f_n\in H$ with $\|f_n\|=1$ and $\|Lf_n\|\ge n$.
        \item[(ii)] {\bf Continuous calculus.} {\it Let $L$ be a self-adjoint operator, and for $\varphi\in C(\sigma(L))$, consider the operators $\varphi(L)$, as we did in the measurable calculus. Then, $\varphi\mapsto \varphi(L)$ is an homeomorphic isometry from $C(\sigma)$ endowed with the sup norm and a subalgebra of $\mathcal{L}(H)$ with respect to the operator norm.} 
        \item[(iii)] {\bf The spectral mapping theorem.} {\it Let $D\xrightarrow{L}H$ be a self-adjoint operator (possibly unbounded), and let $m\in L^\infty(\mu)$, where $\mu$ is the associated p.v.m. Then, $\sigma(m(L))=\text{ess-ran}(m)$, where the essential range is with respect to $\mu$. If $m\in C(\sigma(L))$, then $\sigma(m(L))=m(\sigma(L))$.}
        \item[(iv)] {\it An {\bf eigenvalue} $\lambda$ of $D\xrightarrow{L}H$ is a scalar such that $Lf=\lambda f$ for some $0\ne f\in D$. Any such $f$ is an {\bf eigenvector} relative to $\lambda$, and the subspace of the eigenvectors relative to $\lambda$ is the {\bf eigenspace} relative to $\lambda$. It is clear that $\lambda\in \sigma(L)$, hence $\lambda\in\mathbb{R}$ if $L$ is self-adjoint. In this case, if $\mu\ne\lambda$ are two eigenvectors, the corresponding eigenspaces are orthogonal. The set of the eigenvalues, $\sigma_d(L)$, is the {\bf discrete spectrum} of $L$.} To see this compute $\langle f,Lg\rangle=\langle Lf,g\rangle$ when $f,g$ are eigenvectors relative to $\lambda,\mu$, respectively.
        \item[(v)] {\it  If $f$ is an eigenvector relative to $\lambda$ and $m\in C_c(\sigma(L))$, then $m(\lambda)$ is an eigenvalue, having $f$ as eigenvector. If $\|f\|=1$, then
        \[
        m(\lambda)=\langle f,m(L)f\rangle=\int_{\mathbb{R}}m(x)d\mu_{f,f}(x)
        \]
        for all $m\in C_c(\mathbb{R})$; which correspond to the fact that $\mu_{f,f}$ is a Dirac delta at $\lambda$.}
        \item[(vi)] {\it If $\lambda\in\sigma(L)$ is not an eigenvalue, then $L-\lambda$ has dense range in $H$.} A simple way to prove this is by verifying the statement when $L=M_\eta$ is a multiplication operator on $L^2(X,\nu)$, where $X$ is locally compact and $\nu$ is a finite Radon measure. There, you can show that $\lambda$ is an eigenvalue if and only if $\nu$ "contains" a positive multiple of Dirac delta $a\delta_\zeta$ (i.e. $\nu\ge a\delta_\zeta$ and $\nu-a\delta_\zeta\perp\delta_\zeta$) and $\lambda=\eta(\zeta)$. If $\lambda$ is not an eigenvalue, then, $\eta^{-1}(\lambda)\subseteq X$ is a set which does not support point masses of $\nu$.
        \item[(vii)] {\it Let $\mu$ be the p.v.m. associated to $D\xrightarrow{L}H$, self-adjoint. Then, $\lambda\in\mathbb{R}\setminus\{0\}$ is an eigenvalue of $L$ is and only if $\mu(\{\lambda\})\ne0$.} Let $E=\mu(\{\lambda\})(H)\ni f$. Then
        \[
        \langle g,Lf\rangle=\langle g,L\mu(\{\lambda\})f\rangle=\int_{\{\lambda\}}xd\mu_{g,f}(x)=\lambda \int_{\{\lambda\}}d\mu_{g,f}(x)=\lambda\langle g,f\rangle.
        \]
        Thus, if $E\ne0$, then $\lambda$ is an eigenvalue. 
        
        Viceversa, if $\lambda$ is an eigenvalue with eigenvector $f$, then for all $m\in C_c(\mathbb{R})$
        \[
        \lambda\int_{\mathbb{R}}m(x)d\mu_{f,f}(x)=\lambda\langle f,m(L)f\rangle=\langle f,Lm(L)f\rangle=\int_{\mathbb{R}}xm(x)d\mu_{f,f}(x),
        \]
        This implies that $\mu_{f,f}$ is supported on $\{\lambda\}$, and 
        \[
        \|\mu(\{\lambda\})f\|^2=\int_{\{\lambda\}}d\mu_{f,f}(x)=\int_{\mathbb{R}}d\mu_{f,f}(x)=\|f\|^2,
        \]
        which shows that $f=\mu(\{\lambda\})f\in E$.
        \item[(viii)] {\it The {\bf Cayley map} $\Gamma:\mathbb{R}\to\mathbb{T}$,
        \[
        \Gamma(x)=\frac{x-i}{x+i},
        \]
        induces a bijection 
        \[
        L\mapsto U=\Gamma(L)=(L-i)(L+i)^{-1}
        \]
        from the family of the self-adjoint operators onto the family of the unitary operators. If $\mu_{f,g}$ are the spectral measures for $U$ and $\nu_{f,g}$ are those for $L$, and $m\in L^\infty_0(\mathbb{T})$, then
        \begin{eqnarray*}
            \langle f,m(U)g\rangle&=&\int_{\mathbb{T}}m(e^{it})d\mu_{f,g}(e^{it})\\ 
            &=&\int_{\mathbb{R}}m(\Gamma(x))d\mu_{f,g}(\Gamma(x))\\ 
            &=&\int_{\mathbb{R}}(\Gamma^\ast m)(x)d([\Gamma^{-1}]_\ast\mu_{f,g})(x)\\ 
            &=&\langle f,m(\Gamma(L))g\rangle.
        \end{eqnarray*}
        } 
        Recall that the {\it pull-back} of $\mathbb{T}\xrightarrow{m}\mathbb{C}$ by $\Gamma$ is $\mathbb{R}\xrightarrow{\Gamma^\ast m=m\circ\Gamma}\mathbb{C}$, while the push-forward of a measure $\mu$ on $\mathbb{T}$ by $\Phi=\Gamma^{-1}$ is defined by
        \[
        (\Phi_\ast\mu)(E)=\mu(\Phi^{-1}(E)).
        \]
        What we said can be rephrased as $\nu_{f,g}=[\Gamma^{-1}]_\ast\mu_{f,g}$.
    \end{enumerate}
    \subsection{Compact self-adjoint operators and their inverses}\label{SSectSPcompct}
    A bounded operator $L:H\to H$ is {\it compact} if the image $L(B_1)$ of the unit ball of $H$ has compact closure in $H$. Here we always assume $L$ to be self-adjoint.
    \begin{enumerate}
        \item[(i)] {\it Let $L$ be self-adjoint, bounded, and let $\mu$ be the corresponding p.v.m. Then, $L$ is compact if and only if $\mu|_{\mathbb{R}\setminus(-\epsilon,\epsilon)}$ reduces to a finite number of point masses.} The if part is easy. For the only if, suppose there is $\epsilon>0$ for which the requirement fails. Then we can find infinitely many, disjoint, measurable subsets $E_n$ in $\mathbb{R}\setminus(-\epsilon,\epsilon)$ with $\mu(E_n)\ne0$. Pick a unit vector $f_n$ in each $[\mu(E_n)](H)$: $\{f_n\}$ is bounded, but $\{Lf_n\}$ does not have convergent subsequences.
        \item[(ii)] {\it From (i) it is easy proving that self-adjoint $L$ is compact if and only if there are $L_n$ self-adjoint with finite rank such that $\|L-L_n\|\to0$ as $n\to\infty$.}
        \item[(iii)] {\it {\bf Hilbert-Schmidt Theorem.} If $L$ is self-adjoint, then $L$ is compact if and only if (when $H$ is separable and infinite dimensional) (a) $\sigma(L)=\{\lambda_n\}\cup\{0\}$, with finitely many $\lambda_n$'s, or $\lambda_n\to0$; (b) each $\lambda_n$ is an eigenvalue, and its eigenspace is finite dimensional.}
        \item[(iv)] {\it An self-adjoint operator $D\xrightarrow{\Delta}H$ has a compact inverse (separable, infinite dimensional case) if and only if: (a) $\sigma(M)=\{\mu_n\}\subset\mathbb{R}\setminus\{0\}$ with $|\mu_n|\to$; (b) each $\mu_n$ is an eigenvalue, and its eigenspace is finite dimensional.}
        \item[(v)] {\it With $M$ as in (iv), and $M$ {\bf positive}, the eigenvalues $0<\mu_1\le\dots\le\mu_n\le\dots$ can be computed using the {\bf Rayleigh quotients}.
        \begin{align*}
            \mu_1&=\min_{0\ne x\in H}\frac{\langle x,Mx\rangle}{\|x\|^2};\\
            &\text{pick $x_1$ such that the minimum is achieved}\\
            \mu_2&=\min_{0\ne x\perp x_1}\frac{\langle x,Mx\rangle}{\|x\|^2}\ge\mu_1;\\ 
            &\text{pick $x_2$ such that the minimum is achieved}\\
            &\dots\\
            \mu_n&=\min_{0\ne x\perp x_1,\dots,x_{n-1}}\frac{\langle x,Mx\rangle}{\|x\|^2}\ge\mu_{n-1};\\ 
            &\text{pick $x_n$ such that the minimum is achieved}\\
            &\dots\\
        \end{align*}
        } 
    \end{enumerate}

    \printindex

\end{document}